\documentclass[11pt,reqno]{amsart}
\usepackage[a4paper,margin=1.05in]{geometry}
\usepackage{amsmath,amssymb,amsthm,mathtools}
\usepackage{lmodern}
\usepackage{enumitem}
\usepackage{booktabs,longtable}
\usepackage{array}
\usepackage{needspace}
\usepackage{xcolor}
\usepackage[colorlinks=true,linkcolor=blue!50!black,citecolor=blue!50!black,urlcolor=blue!50!black]{hyperref}

\newtheorem{theorem}{Theorem}[section]
\newtheorem{proposition}[theorem]{Proposition}
\newtheorem{lemma}[theorem]{Lemma}
\newtheorem{corollary}[theorem]{Corollary}
\newtheorem*{conjS}{Conjecture S}
\newtheorem*{conjP}{Conjecture P (coefficient positivity)}
\theoremstyle{definition}
\newtheorem{definition}[theorem]{Definition}
\newtheorem{example}[theorem]{Example}
\theoremstyle{remark}
\newtheorem{remark}[theorem]{Remark}
\numberwithin{equation}{section}
\newcommand{\F}{\mathbb{F}}
\newcommand{\Z}{\mathbb{Z}}
\newcommand{\Q}{\mathbb{Q}}
\newcommand{\R}{\mathbb{R}}
\newcommand{\N}{\mathbb{N}}
\newcommand{\C}{\mathbb{C}}
\newcommand{\K}{\mathbb{K}}
\newcommand{\cK}{\mathcal{K}}

\newcommand{\cH}{\mathcal{H}}
\newcommand{\cL}{\mathcal{L}}
\newcommand{\cN}{\mathcal{N}}
\newcommand{\cP}{\mathcal{P}}
\newcommand{\cC}{\mathcal{C}}
\newcommand{\fo}{\mathfrak{o}}
\newcommand{\fp}{\mathfrak{p}}
\newcommand{\fn}{\mathfrak{n}}
\newcommand{\fm}{\mathfrak{m}}
\newcommand{\fc}{\mathfrak{c}}
\newcommand{\fH}{\mathfrak{H}}
\newcommand{\bO}{\mathbb{O}}
\newcommand{\tW}{\widetilde{W}}
\newcommand{\sH}{\mathsf{Hall}}
\newcommand{\sP}{\mathsf{P}}
\newcommand{\sN}{\mathsf{N}}
\newcommand{\sM}{\mathsf{H}}
\DeclareMathOperator{\GL}{GL}
\DeclareMathOperator{\Mat}{M}
\DeclareMathOperator{\val}{val}
\DeclareMathOperator{\vol}{vol}
\DeclareMathOperator{\Hom}{Hom}
\DeclareMathOperator{\End}{End}
\DeclareMathOperator{\Aut}{Aut}
\DeclareMathOperator{\Ext}{Ext}
\DeclareMathOperator{\Sur}{Sur}
\DeclareMathOperator{\Frac}{Frac}
\DeclareMathOperator{\Res}{Res}
\DeclareMathOperator{\Ch}{Ch}
\DeclareMathOperator{\ms}{ms}
\DeclareMathOperator{\eu}{eu}
\DeclareMathOperator{\mult}{mult}
\DeclareMathOperator{\mass}{mass}

\DeclareMathOperator{\dinv}{dinv}
\DeclareMathOperator{\diag}{diag}
\DeclareMathOperator{\rad}{rad}

\DeclareMathOperator{\pr}{pr}

\DeclareMathOperator{\Gr}{Gr}

\DeclareMathOperator{\im}{im}

\newcommand{\cE}{\mathcal E}

\newcommand{\qbinom}[2]{\genfrac{[}{]}{0pt}{}{#1}{#2}}
\newcommand{\ctd}{\trianglerighteq}
\newcolumntype{L}[1]{>{\raggedright\arraybackslash}p{#1}}

\author{Alexei Oblomkov}
\address{
A.~Oblomkov\\
Department of Mathematics and Statistics\\
University of Massachusetts at Amherst\\
Lederle Graduate Research Tower\\
710 N. Pleasant Street\\
Amherst, MA 01003 USA
}

\title[Cabling and equisingularity for lattice counts]{Cabling and equisingularity for lattice counts on\\ plane curve singularities}
\date{November 28, 2026}

\begin{document}

\begin{abstract}
Let $A\subset\F_q[[t]]$ be a plane branch with $\delta$-invariant
$\delta$, and let $E_A(n)$ count the $A$-submodules of
$\F_q[[t]]^n$ which generate $\F_q[[t]]^n$ as an $\F_q[[t]]$-module.
In every characteristic, we prove that $E_A(n)$ is a polynomial
in $q$ with integer coefficients, determined by $n$ and the
multiplicity sequence of $A$. We give a finite algorithm for
computing this polynomial and conjecture that its coefficients
are nonnegative; this holds in rank one and for torus branches.

The limit $\cN_A=\lim_{n\to\infty}q^{-n^2\delta}E_A(n)$ is
finite and equals the normalized mass of finite $A$-modules.
A cabling formula expresses it as a product of torus-branch
invariants determined by the multiplicity sequence. The torus
results of Huang--Jiang--Oblomkov and Huang--Lau--Ono give an
explicit infinite product of Rogers--Ramanujan type. The normalized
counting polynomials agree with this product through degree $n$
in $q^{-1}$; after multiplication by an explicit rational power,
the product is a modular unit of weight zero.

The proof realizes lattice and flag counts as orbital integrals.
Relations among invariant distributions on the Iwahori--Hecke
algebra give a universal transfer matrix for a quadratic
transform. Hall algebra identities identify its scalar
coefficients, and a uniform bound for an inverse triangular
matrix controls the passage to the limit.
\end{abstract}

\maketitle
\setcounter{tocdepth}{1}
\tableofcontents
\clearpage

\section{Introduction}

We count lattices over the complete local ring of a unibranch plane
curve singularity over a finite field. At finite rank the counts
depend only on the multiplicity sequence and are polynomials in the
size of the ground field (Theorem~\ref{thm:introB}). Their normalized
limit satisfies a cabling formula, factors into torus-branch
invariants, and has a modular interpretation
(Theorems~\ref{thm:introC}--\ref{thm:introM}).
We also formulate Conjectures~P and~S on coefficient positivity and
a fibrewise form of the cabling identity.

Our starting point is the rigid-cocenter framework of Kivinen and
Tsai \cite[\S2.3]{KT}. Its test functions count lattice chains for
framed plane branches in arbitrary module rank and characteristic.
The resulting transfer across quadratic transforms yields
an effective proof of equisingular invariance. A uniform estimate
for the inverse of a normalized Hall matrix then gives the stable
cabling formula. Combining this formula with the torus products of
\cite{HJO,HLO}, we obtain coefficientwise stabilization and modularity.
The method is described in Section~\ref{sub:method}.

\subsection{Lattice counts}
Let $F=\F_q$ be a finite field, $O=F[[t]]$ and $\cK=F((t))$. A \emph{branch} is an $F$-subalgebra $A\subseteq O$ with $\delta_A:=\dim_FO/A<\infty$, and such a ring is a complete local domain with residue field $F$, fraction field $\cK$ and normalization $O$. We call a branch $A$ a \emph{plane branch} if $A=F[[x,y]]$ is the closure of $F[x,y]$ for some $x,y\in tO$; equivalently, $A\cong F[[X,Y]]/(f)$ is a domain whose normalization is $F[[t]]$. Thus a plane branch is the complete local ring of a unibranch plane curve singularity with a rational branch.

An \emph{$A$-lattice} in $\cK^n$ is a finitely generated $A$-submodule that spans $\cK^n$, and the main objects of this paper are the finite sets:
\[
E_A(O^n)=\{L\subseteq O^n\mid L\text{ an }A\text{-lattice},\ OL=O^n\}
\]
together with the normalized counts:
\[
e_A(n)=q^{-n^2\delta_A}\,\#E_A(O^n),\qquad \cN_A=\lim_{n\to\infty}e_A(n).
\]
For $n=1$, normalizing the least valuation of a rank-one lattice
identifies $E_A(O)$ with the $F$-points of the local compactified
Jacobian, in the fixed-degree convention of Remark~\ref{rem:poly}(2).
This relates the count to knot invariants, orbital integrals and
combinatorics \cite{OS,ORS,Maulik,GM,Yun,KT}. For general $n$, it is
the set of points of the ``extension fibre'' of the punctual Quot
scheme of $O^n$ over $A$, studied by Huang and Jiang \cite{HJ} and
by Huang, Jiang and Oblomkov \cite{HJO}. The connection with
colored knot invariants is discussed in Section~\ref{sub:colored-slopes}.
Our first result gives the normalization and identifies the limit.

\begin{theorem}[Theorem~\ref{prop:cohnorm}, Corollaries~\ref{cor:limit-exists} and~\ref{cor:planarity-finiteness}]\label{thm:introA}
Let $A$ be a plane branch over $\F_q$, or more generally a branch that satisfies the Gorenstein condition \textup{(G)} of Section~\ref{sec:norm}. Then for every $n\ge1$ we have:
\[
\sum_{L\subseteq A^n}\big|A^n/L\big|^{-n}=e_A(n)\,\prod_{i=1}^n(1-q^{-i})^{-1},
\]
where the sum is over all $A$-submodules of finite colength. Hence $\cN_A$ exists in $(0,\infty]$, and we have:
\[
\cN_A=(q^{-1};q^{-1})_\infty\sum_{[M]}\frac1{|\Aut_A(M)|},
\]
where the sum is over the isomorphism classes of $A$-modules of finite cardinality. Among branches satisfying \textup{(G)}, one has $\cN_A<\infty$ if and only if $A$ is a plane branch.
\end{theorem}
\noindent Thus $\cN_A$ is the ratio of the Cohen--Lenstra sums \cite{CL} of $A$ and of its normalization, in other words the value at $1$ of the normalized Coh zeta function of $A$ in the sense of \cite{Huang,HJ}. When $A=F[[X,Y]]/(f)$, the sum $\sum_{[M]}|\Aut_A(M)|^{-1}$ is a weighted count of pairs of commuting nilpotent matrices satisfying $f(X,Y)=0$ (Remark~\ref{rem:coh}). For torus knot singularities the identity of the theorem is due to Huang, Jiang and Oblomkov \cite{HJO}, who deduce it from the functional equation of \cite{HJ}. Our proof is a direct computation with a Haar measure and residue duality, and it needs no separability hypothesis.

The finiteness of $\cN_A$ follows from a second formula, in which $e_A(n)$ is a sum over the finitely many $A$-lattices between $A^n$ and $O^n$ (Proposition~\ref{prop:finitesum}). By the duality (G), this sum is a count of the modules over the finite ring $A/t^{2\delta}O$ with the weights $|\Aut|^{-1}$, up to factors which tend to $1$ as $n\to\infty$ (Proposition~\ref{prop:conductor}). If $A=F[[x,y]]$, then a fixed power of $y$ vanishes on these modules, and such modules have finite total mass (Lemma~\ref{lem:ypower}). For torus knot singularities the finiteness was obtained in \cite{HJO} from the positive-definiteness theorem of Huang \cite{Huang2}. The converse is proved in Corollary~\ref{cor:planarity-finiteness}: a nonplane branch satisfying (G) has $\cN_A=\infty$.

\subsection{Main results}
For a plane branch $A$ we denote by $\mult(A)$ its multiplicity and by $A=A^{(0)}\subseteq A^{(1)}\subseteq\cdots$ the sequence of its quadratic transforms, and we call $\ms(A)=(\mult A^{(j)})_{j\ge0}$ the \emph{multiplicity sequence} of $A$. It is the characteristic-free form of the equisingularity class: in characteristic zero it determines, and is determined by, the Puiseux characteristic, the semigroup and the topological type of the link.

\begin{theorem}[Theorem~\ref{thm:E}, Corollary~\ref{cor:poly}]\label{thm:introB}
Let $A$ be a plane branch over $\F_q$, in any characteristic. For every $n\ge1$ the number $\#E_A(O^n)$ depends only on $q$, $n$ and $\ms(A)$. More precisely, for every multiplicity sequence $\underline\mu$ and every $n$ there is a polynomial $P_{\underline\mu,n}\in\Z[X]$ such that $\#E_A(O^n)=P_{\underline\mu,n}(q)$ for every plane branch $A$ over a finite field $\F_q$ with $\ms(A)=\underline\mu$.
\end{theorem}
\noindent The same holds for the finer counts of lattices and of chains of lattices with prescribed residual types (Theorem~\ref{thm:E}). For $n=1$ and tamely ramified branches the polynomiality, together with the nonnegativity of the coefficients, is a result of Kivinen and Tsai \cite[Cor.~1.11 and Cor.~8.3]{KT}, and for $n>1$ it was known for the singularities $y^2=x^{2k+1}$ \cite{HJ} and for torus knots \cite{HJO}.

The transfer theorem of Sections~\ref{sec:hecke}--\ref{sec:chains}
expresses chain counts using universal rational coefficients
(Theorem~\ref{thm:transfer}). Iterating the Euclidean transform of
Section~\ref{sec:qt} reduces them to finite flag counts on a smooth
branch and proves Theorem~\ref{thm:E}. Extension of constants then
gives integer values at an unbounded sequence of field sizes, so
Lemma~\ref{lem:ratpoly} shows that the rational functions are
polynomials over $\mathbb Q$. Their coefficients are integral by the
finite-type models of Appendix~\ref{app:integrality} and the
zeta-function argument of Lemma~\ref{lem:completion-integral}
(Corollary~\ref{cor:poly}).
Section~\ref{subsec:explicit-counting-algorithm} makes this recursion
explicit.

\begin{conjP}
For every multiplicity sequence $\underline\mu$ of a plane branch
over a finite field and every $n\ge1$, the polynomial of
Theorem~\ref{thm:introB} satisfies
\[
 P_{\underline\mu,n}\in\Z_{\ge0}[X].
\]
\end{conjP}
Conjecture~P concerns the coefficients of the total extension-fibre
polynomial. In rank one it holds for every multiplicity sequence
(Proposition~\ref{prop:rankone-positive}): the tame case is
\cite[Cor.~1.11 and Cor.~8.3]{KT}, with the fixed-degree
comparison in Remark~\ref{rem:poly}(2), and
Theorem~\ref{thm:introB} removes the tameness hypothesis, because
every multiplicity sequence is realized in large characteristic and
the polynomial does not depend on the realization. For the torus
sequences $\eu(a,b)$ it holds in every rank
(Proposition~\ref{prop:torus-positive}), by the finite-rank torus
formula.

The general case remains open; further nontorus families in every
rank are described in Section~\ref{sub:colored-slopes}.
The analogous assertion for the refined counts of
Theorem~\ref{thm:E} is false: Example~\ref{ex:refined-negative} gives a
residual-type count with a negative coefficient.

For coprime $a,b\ge1$, write $E_{a,b}(n)=\#E_{F[[t^a,t^b]]}(O^n)$ and $\cN_{a,b}=\cN_{F[[t^a,t^b]]}$, the invariant of the torus knot $T(a,b)$. Note that $\cN_{a,b}=1$ if $a=1$ or $b=1$.

\begin{theorem}[Theorems~\ref{thm:cable},~\ref{thm:limit} and~\ref{thm:product}]\label{thm:introC}
Let $A=F[[x,y]]$ be a plane branch over $F=\F_q$, in any characteristic, with $\val_tx=Da$, $\val_ty=Db$ and $\gcd(a,b)=1$. Choose integers $A_0,B_0$ with $A_0a+B_0b=1$, let $c\in F^\times$ be the constant term of $y^a/x^b$, and let us put:
\[
\varpi=x^{A_0}y^{B_0},\qquad\omega=c^{-1}y^a/x^b-1,\qquad T=F[[\varpi,\omega]],\qquad S=F[[\varpi,\varpi^{a+b-1}\omega]].
\]
Then $T$ and $S$ are plane branches,
and the following holds:
\begin{enumerate}[label=\textup{(\roman*)}]
\item For every $n$ we have:
\[
\#E_A(O^n)=\sum_{\mu\vdash Dn}B^T_n(\mu)\,\Phi^{a,b}_{q,Dn}(\mu).
\]
Here $B^T_n(\mu)$ counts the $\Lambda\in E_T(O^n)$ for which
$\omega$ acts on $\Lambda/\varpi\Lambda$ with Jordan type $\mu$.
The functions $\Phi^{a,b}_{q,N}$ depend only on $(a,b,q,N)$ and are
determined by the torus counts as follows: the assignments
\[
 u_\mu\longmapsto
 \frac{q^{-ab\binom{|\mu|}{2}}\Phi^{a,b}_{q,|\mu|}(\mu)}
      {|\Aut(M_\mu)|}
\]
form the Hall-algebra character whose value on $u_{(1^j)}$ is
$q^{-ab\binom j2}E_{a,b}(j)/|\GL_j(\F_q)|$.
\item The numbers $\cN_A$, $\cN_{a,b}$ and $\cN_S$ are finite, and we have:
\[
\cN_A=\cN_{a,b}^{\,D}\cdot\cN_S .
\]
\end{enumerate}
\end{theorem}
At finite rank, part~(i) retains the residual-type counts
$B_n^T(\mu)$. The product in part~(ii) arises after normalization
and passage to the limit in rank; the finite-rank counts themselves
are not multiplicative (see Section~\ref{sub:rankone}). 
If the generators are chosen with $\val x=\mult(A)$ and $\val x\nmid\val y$, then $\mult(S)=D<\mult(A)$, so the formula can be iterated, and we obtain the following.

\begin{theorem}[Theorem~\ref{thm:product}, Corollaries~\ref{cor:puiseux} and~\ref{cor:finalproduct}]\label{thm:introD}
Let $A$ be a plane branch over $\F_q$, in any characteristic. If $A=O$, take the empty cabling list and interpret both products below as $1$. If $A\ne O$, there is a list $(a_i,b_i,D_i)_{1\le i\le g}$ of triples of positive integers with $2\le a_i<b_i$, $\gcd(a_i,b_i)=1$ and $D_1>\dots>D_g=1$, computed from $\ms(A)$ alone by the algorithm of Section~\ref{sub:iterate}, such that
\[
\cN_A=\prod_{i=1}^g\cN_{a_i,b_i}^{\,D_i}=\prod_{i=1}^g\ \prod_{j\ge1}\big(1-q^{-j}\big)^{-D_i\,\rho_{a_i,b_i}(j)},
\]
where $\rho_{a,b}(j)=\min\{a,\,b,\,\operatorname{dist}(aj,(a+b)\Z)\}-1+\mathbf 1_{(a+b)\mid j}$.
If $\ms(A)$ is the multiplicity sequence of a branch with Puiseux characteristic $(\beta_0;\beta_1,\dots,\beta_g)$ and $e_i=\gcd(\beta_0,\dots,\beta_i)$, then $(a_i,b_i,D_i)=\big(e_{i-1}/e_i,\ (\beta_0+\beta_i-e_{i-1})/e_i,\ e_i\big)$.
\end{theorem}
We prove the first equality here. The second equality is the product formula $\cN_{a,b}=\prod_j(1-q^{-j})^{-\rho_{a,b}(j)}$, conjectured by Huang, Jiang and Oblomkov \cite{HJO} and proved by Huang, Lau and Ono \cite{HLO}; the case $a=3$ was proved earlier by Lau and Ono \cite{LO}. For example, for $A=\F_q[[t^4,t^6+t^7]]$ with $q$ odd we get $\cN_A=\cN_{2,3}^{\,2}\cN_{2,9}$, the product of a squared Rogers--Ramanujan product and an Andrews--Gordon product. If the characteristic $p$ divides the multiplicity, the list must be computed from the multiplicity sequence and not from the exponents of a parametrization, as the examples of Section~\ref{sec:examples} show.

The infinite product in Theorem~\ref{thm:introD} is the value at
$Q=q^{-1}$ of a power series determined by the multiplicity sequence.
This series is the coefficientwise limit of the normalized counting
polynomials of Theorem~\ref{thm:introB}; multiplication by an explicit
rational power of $Q$ gives a modular unit.

\Needspace{6\baselineskip}
\begin{theorem}[Theorems~\ref{thm:formal} and~\ref{thm:modular}]\label{thm:introM}
Let $A\ne\F_q[[t]]$ be a plane branch over $\F_q$, in any characteristic, with multiplicity sequence $\underline\mu$, $\delta$-invariant $\delta$ and cabling list $((a_i,b_i,D_i))_{1\le i\le g}$. Let us put:
\[
\cN_A(Q)=\prod_{i=1}^g\ \prod_{j\ge1}\big(1-Q^j\big)^{-D_i\,\rho_{a_i,b_i}(j)},\qquad c_A=\frac1{24}\sum_{i=1}^gD_i\,\frac{(a_i-1)(b_i-1)}{a_i+b_i},
\]
and let $L$ be the least common multiple of the numbers $a_i+b_i$.
\begin{enumerate}[label=\textup{(\roman*)}]
\item For every $n\ge1$ the function $Q^{n^2\delta}P_{\underline\mu,n}(Q^{-1})$ is a polynomial in $Q$ with constant term $1$, and it is congruent to $\cN_A(Q)$ modulo $Q^{\,n+1}$. In particular, the normalized counting polynomials converge coefficientwise to $\cN_A(Q)$.
\item The function $\Theta_A(\tau)=Q^{-c_A}\,\cN_A(Q)$, $Q=e^{2\pi i\tau}$, is holomorphic and nowhere zero on the upper half-plane, and it is invariant under the principal congruence subgroup $\Gamma(12L)$. Its power $\Theta_A^{12L}$ is a modular function on $\Gamma_1(L)$, and at the cusps $\infty$ and $0$ the function $\Theta_A$ has poles of the same order $c_A$, measured in $Q$.
\end{enumerate}
\end{theorem}
\noindent Thus
$\cN_A=q^{-c_A}\Theta_A\bigl(\tfrac{i\log q}{2\pi}\bigr)$ is,
up to its displayed factor, a value of a modular unit for
$\Gamma(12L)$. For the cusp,
$\Theta_A=Q^{-1/60}\prod_{n\ge1}(1-Q^{5n-1})^{-1}(1-Q^{5n-4})^{-1}$
is the Rogers--Ramanujan function, a modular unit for $\Gamma(60)$.
The torus products also have the character interpretation recalled
in Remark~\ref{rem:modular}(1). Part~(i) is proved $Q$-adically,
using the polynomiality of the refined counts
(Theorem~\ref{thm:introB}) since convergence at the values $Q=q^{-1}$
alone would not give it.

\subsection{Saturation}\label{sub:saturation}
The ring $T$ of Theorem~\ref{thm:introC} contains $A$. Indeed,
$\varpi^a=x\,\big(c(1+\omega)\big)^{B_0}$ implies that $x\in T$,
and the same argument gives $y\in T$. Thus every
$L\in E_A(O^n)$ has a \emph{saturation} $TL\in E_T(O^n)$, the
$T$-module generated by $L$. This defines the map:
\[
\sigma_T\colon E_A(O^n)\to E_T(O^n),\qquad L\mapsto TL,
\]
and the number of lattices decomposes over its fibres:
\[
\#E_A(O^n)=\sum_{\Lambda\in E_T(O^n)}\#\sigma_T^{-1}(\Lambda).
\]
The chains $(\Lambda_j)_{j\in\Z}$ of lattices with $x\Lambda_j=\Lambda_{j+a}$ and $y\Lambda_j=\Lambda_{j+b}$ are exactly the chains $\Lambda_j=\varpi^j\Lambda$ with $\Lambda$ a $T$-lattice (Lemma~\ref{lem:formal}).  This characterizes $T$ through the simultaneous shift relations for $x$ and $y$.

Write $\mu(\Lambda)$ for the Jordan type of $\omega$ on
$\Lambda/\varpi\Lambda$. Combining the definition of $B_n^T(\mu)$
with Theorem~\ref{thm:introC}(i) gives
\[
 \sum_{\Lambda\in E_T(O^n)}\#\sigma_T^{-1}(\Lambda)
 =\sum_{\Lambda\in E_T(O^n)}
        \Phi^{a,b}_{q,Dn}\bigl(\mu(\Lambda)\bigr).
\]
We expect that this identity holds fibre by fibre:

\begin{conjS}
For every $\Lambda\in E_T(O^n)$ we have $\#\sigma_T^{-1}(\Lambda)=\Phi^{a,b}_{q,Dn}(\mu(\Lambda))$.
\end{conjS}
\noindent Summing Conjecture~S over $\Lambda\in E_T(O^n)$ gives
Theorem~\ref{thm:introC}(i); the proof of the theorem uses only the
sum. Conjecture~S was numerically tested for several non-trivial examples over \(\F_2\) and \(\F_3\).



\subsection{Method and geometric interpretation}\label{sub:method}
\paragraph{\textbf{The finite moduli problem.}}
The set $E_A(O^n)$ is the set of $F$-points of a scheme of finite
type over $F$. If $t^hO\subseteq A$, then $L\mapsto L/t^hO^n$
identifies the lattices in $E_A(O^n)$ with the $A$-stable subspaces
of $(O/t^hO)^n$ that map onto $(O/tO)^n$, a locally closed subset of
a finite union of Grassmannians. Appendix~\ref{app:integrality}
constructs this model and its refinements for residual types and
chains of lattices; they give the integral coefficients of
Corollary~\ref{cor:poly}.

\paragraph{\textbf{The higher-rank affine Springer problem.}}
Choose nonzero generators $u,w\in tO$ with $A=F[[u,w]]$.
Put $K_u=F((u))$ and $D=[\cK:K_u]=\val(u)$, and regard
$V=\cK^n$ as a $K_u$-vector space of dimension $m=Dn$.
Let $G=\GL_{K_u}(V)$ and let $\gamma$ be multiplication by $w$.
Here $n$ is the module rank and $m$ the ambient matrix size.
Since $K_u[w]=\cK$, the minimal polynomial
$p_w(Y)\in K_u[Y]$ of $w$ has degree $D$, and
\[
 \begin{gathered}
 \operatorname{minpol}(\gamma)=p_w,\qquad
 \operatorname{charpol}(\gamma)=p_w^n,\\
 Z_G(\gamma)(K_u)=\operatorname{Aut}_{K_u[w]}(V)
                  =\GL_n(\cK).
 \end{gathered}
\]
If $\cK/K_u$ is separable, the $D$ distinct conjugates of $w$
are the eigenvalues of $\gamma$, each with multiplicity $n$.
Its algebraic centralizer over $K_u$ is
\[
 Z_G(\gamma)\cong
 \operatorname{Res}_{\cK/K_u}\GL_n,\qquad
 \dim Z_G(\gamma)=Dn^2.
\]
Thus $\gamma$ is regular semisimple when $n=1$, and semisimple
but nonregular when $n>1$. If the extension is inseparable,
the minimal-polynomial and centralizer calculations still hold,
but $\gamma$ is not semisimple. The counting argument below
also covers this case.

At the level of $F$-points, the affine Springer condition is
\[
 X_\gamma=\{\Lambda\subseteq V:
 \Lambda\text{ is an }F[[u]]\text{-lattice},\
 w\Lambda\subseteq\Lambda\}.
\]
These are precisely the $A$-lattices. The finite set counted here is
\[
 E_A(O^n)=\{\Lambda\in X_\gamma:O\Lambda=O^n\}.
\]
Fixing the $O$-saturation bounds the lattices between
$t^hO^n$ and $O^n$, as in the preceding Grassmannian model.
The full $X_\gamma$ can be much larger: for $A=O$ and $u=w=t$,
it is the set of $F$-points of the entire affine Grassmannian of
$\GL_n$, whereas $E_O(O^n)=\{O^n\}$.


\paragraph{\textbf{From moduli points to orbital integrals.}}
In the orbital-integral formulas, write $G$ for its group of
$K_u$-points and $Z=\GL_n(\cK)$ for the centralizer.
We normalize its Haar measure by $\vol(Z_0)=1$, where
$Z_0=\GL_n(O)$, and put
$k=n\val(w)=\val_{K_u}(\det w)$.
The orbital distribution is
\[
 D_w(f)=\int_{Z\backslash G}f(g^{-1}wg)\,d\dot g,
\]
with the quotient measure induced from the Haar measures on $G$
and $Z$; its convergence and invariance are proved in
Lemma~\ref{lem:orbital}.

The Haar measure on $G$ may be chosen arbitrarily: rescaling it
multiplies the quotient measure and the parahoric volumes below by
the same factor. For a composition $\alpha$ of $m$, the quotient
$G/P_\alpha$ is the set of $F$-points of a partial affine flag
variety, parametrizing indexed $u$-periodic lattice chains of type
$\alpha$. The conditions $w\Lambda_j\subseteq\Lambda_j$ define
$X_\alpha(w)$, the set of such chains whose lattices are $A$-modules.

Let $\fp_\alpha$ be the parahoric order with unit group $P_\alpha$.
The test function $F_{\alpha,k}$ is the characteristic function of
its elements of determinant valuation $k$, divided by
$\vol(P_\alpha)$. A conjugate $g^{-1}wg$ belongs to $\fp_\alpha$
exactly when $w$ preserves every lattice of the chain indexed by
$gP_\alpha$. Lemmas~\ref{lem:dictionary} and~\ref{lem:orbital} therefore give
\[
 D_w(F_{\alpha,k})
 =\sum_{[\Lambda_\bullet]\in Z\backslash X_\alpha(w)}
       \frac1{\vol(Z_{\Lambda_\bullet})}
 =\Ch_\alpha(A,u;n).
\]

The division by $\vol(P_\alpha)$ in the test function cancels
the volume of each right $P_\alpha$-coset. The remaining
centralizer mass becomes an ordinary finite count: every
$Z$-orbit meets the locus $O\Lambda_0=O^n$, and two normalized
chains in the same $Z$-orbit are in the same $Z_0$-orbit.
Their stabilizers lie in $Z_0$, so
$\vol(Z_{\Lambda_\bullet})^{-1}
 =[Z_0:Z_{\Lambda_\bullet}]$ is exactly the size of that orbit
in the normalized locus. This is Lemma~\ref{lem:massE} and
\eqref{eq:massfibre} applied to chains.

For a fixed first lattice $\Lambda$, the remaining chain is a
$\bar w$-stable partial flag in $\Lambda/u\Lambda$.
These are the $F$-points of a closed incidence subscheme of an
ordinary flag variety. Thus Lemma~\ref{lem:chainflag} says
\[
 \Ch_\alpha(A,u;n)
 =\sum_{\Lambda\in E_A(O^n)}H_\alpha\bigl(\mu(\Lambda)\bigr),
 \qquad
 D_w(F_{(m),k})=\#E_A(O^n).
\]
Here $\mu(\Lambda)$ is the Jordan type of $w$ on
$\Lambda/u\Lambda$. The next theorem gives universal relations
among these orbital integrals.

\begin{theorem}[{Theorem~\ref{thm:H}, Corollary~\ref{cor:basis}; cf.\ \cite[\S2.3]{KT}}]\label{thm:introH}
Let $\K=\F_q((\pi))$, $G=\GL_m(\K)$, and let $\cH$ be the Iwahori--Hecke algebra of $G$. For $k\in\Z$ let $J_k$ be the space of invariant distributions on $G$ supported on compact-mod-center elements, restricted to the elements of $\cH$ supported on $\{\val\det=k\}$. Then $\dim J_k=p(\gcd(m,k))$, the number of partitions of $\gcd(m,k)$. Put $e=m/\gcd(m,k)$. Moreover, let $f^*_\beta$, $\beta\vdash\gcd(m,k)$, be the normalized characteristic function of the coset $\tau^kP_{\beta^{[e]}}$, where $P_{\beta^{[e]}}$ is a certain standard parahoric subgroup normalized by $\tau^k$ (Section~\ref{sec:chains}). Then the evaluations at the functions $f^*_\beta$ form a basis of $J_k^*$, and the coordinates of the evaluation $D\mapsto D(T_w)$ in this basis
are rational functions of $q$ depending only on the combinatorial
data and not on the field, where
$T_w=\vol(I)^{-1}\mathbf 1_{IwI}$.
\end{theorem}
\noindent The upper bound $\dim J_k\le p(\gcd(m,k))$ and the
universality statement are contained in the work of Kivinen and Tsai
\cite[\S2.3]{KT}, who derive them from the theorems of He and Nie on
minimal length elements \cite{HN} and of He on the Newton decomposition
\cite{He}. The comparison with their results is discussed below.

We use the reduction to standard triples to express $D(T_w)$, for
$D\in J_k$, as a linear combination of the $D(T_{w_\rho})$ with
coefficients in $\mathbb Z[q]$. Here $\rho\vdash\gcd(m,k)$, and the
representatives $w_\rho$ are defined before Theorem~\ref{thm:H}.
Modulo commutators, every double coset reduces
to terminal ones. The gradient-flow argument of \cite{HN}, proved in
our conventions in Section~\ref{app:HN3}, reduces a terminal element,
without changing its length, to $ux$ with $x$ straight.

A terminal double coset with noncentral Newton point contains no
compact-mod-center element (Lemma~\ref{lem:cert}), which gives the
required vanishing. For central Newton points, a computation in the
Hecke algebra of the symmetric group and a conjugation interchanging
two blocks determine the cocenter classes
(Lemmas~\ref{lem:twist}--\ref{lem:swap}). Lattice counts for explicit
orbital distributions prove independence
(Proposition~\ref{prop:indep}); in particular, the $f^*_\beta$ are
independent in $\cH/[\cH,\cH]$ in every characteristic.

\paragraph{\textbf{Returning to the transformed branch.}}
The functions $f^*_\beta$ count chains for the transformed branch.
Put $n'=\gcd(m,k)$, $e=m/n'$, and $d_*=k/n'$. For a
partition $\beta\vdash n'$ with $l$ parts, the test function
$f^*_\beta$ selects chains of type $\beta^{[e]}$ satisfying
\[
 u\Lambda_j=\Lambda_{j+el},\qquad
 w\Lambda_j=\Lambda_{j+d_*l}.
\]
Choose $A_0e+B_0d_*=1$, and put
$\varpi=u^{A_0}w^{B_0}$ and $\theta=w^eu^{-d_*}$.
Lemma~\ref{lem:formal} identifies the same indexed chains of lattices
by the equalities
\[
 \varpi\Lambda_j=\Lambda_{j+l},\qquad
 \theta\Lambda_j=\Lambda_j.
\]
Writing $\omega=c^{-1}\theta-1$, with $c$ the leading
coefficient of $\theta$, these are exactly the
$\varpi$-periodic chains of $T=F[[\varpi,\omega]]$-lattices.
The identification is $Z$-equivariant and uses the same Haar
measure, so Lemma~\ref{lem:rigidchains} gives
$D_w(f^*_\beta)=\Ch_\beta(T,\varpi;n)$.

Evaluating the universal distribution identity of
Corollary~\ref{cor:basis} on $D_w$ therefore gives
\[
 \Ch_\alpha(A,u;n)=
 \sum_{\beta\vdash n'}c_\beta(F_{\alpha,k})
                 \Ch_\beta(T,\varpi;n),
\]
which is Theorem~\ref{thm:transfer}. The operator identities give a
bijection between the two sets of shifted chains, and the transfer
theorem relates their counts linearly to the ordinary stable-chain
counts. For the one-part composition, forgetting the transformed
flags gives Corollary~\ref{cor:orbitsum} and
Theorem~\ref{thm:introC}(i). The fibres in Lemma~\ref{lem:chainflag} parametrize stable flags
in a fixed lattice quotient. The fibres of the saturation map
$\sigma_T$ are different; their counts are the subject of Conjecture~S. Induction through the quadratic transforms of
Proposition~\ref{prop:transform} proves Theorem~\ref{thm:introB}.

\paragraph{\textbf{Relation with the work of Tsai and of Kivinen and Tsai.}}
For $\GL_m$ the framework of Theorem~\ref{thm:introH} is that of
\cite[\S2.3]{KT}. The cocenter $\cH/[\cH,\cH]$ is graded by the
valuation $k$ of the determinant and is the direct sum of Newton
strata, by the theorems of He and Nie \cite{HN} and of He \cite{He}.
The restriction to $\cH$ of an invariant distribution supported on
compact-mod-center elements factors through the dual of the
\emph{rigid} part, the sum of the strata with central Newton point
\cite[Thm.~2.25, Def.~2.26, Cor.~2.28, Prop.~2.29]{KT}.

In valuation
$k$ the rigid part has dimension $p(\gcd(m,k))$, with basis the classes
of the $T_{w_\rho}$ \cite[Prop.~2.31]{KT}. The classes of the functions
$f^*_\beta$, the translated parahoric test functions of Waldspurger
\cite[IV.1]{Wal91}, form another basis \cite[Def.~2.30, Cor.~2.33]{KT},
and the change of basis is given by rational functions of $q$ which do
not depend on the field \cite[the paragraph following Cor.~2.33]{KT};
they are not computed there, and
Section~\ref{subsec:explicit-counting-algorithm} gives an algorithm for
them.

The results of \cite[\S2.3]{KT} give the bound
$\dim J_k\le p(\gcd(m,k))$ and the universality statement for every
local field. They involve no regularity hypothesis on the orbit, so
they apply to our $D_w$ once Lemma~\ref{lem:orbital}
establishes convergence and the required support.

The equality
$\dim J_k=p(\gcd(m,k))$ is obtained in \cite{KT} for $m\mid k$ in every
characteristic, from the unipotent orbital integrals, and for all $k$
in characteristic zero, from Waldspurger's truncated Steinberg
characters \cite[Prop.~2.4]{Wal89}, \cite[Prop.~III.4]{Wal91},
\cite[Thm.~2.19, Rem.~2.21, Thm.~2.34]{KT}. In every characteristic it
can also be derived from the Newton decomposition of the cocenter of
$C_c^\infty(G)$ \cite[Thms.~A and~B]{He} and the injectivity of the
map from the cocenter of $\cH$ to that of $C_c^\infty(G)$, which follows
from the Morita equivalence between $\cH$ and the Iwahori block of the
category of smooth representations \cite{Borel},
\cite[Prop.~2.10]{BeDe}.

The counting theorems of \cite{KT} concern tamely ramified regular
semisimple elements: an orbital symmetric function is constructed from
the Newton pairs, and its Hall pairings give the parahoric orbital
integrals \cite[Thms.~1.6 and~1.10]{KT}. Their geometric application
has the same regularity hypothesis, and their application to
compactified Jacobians concerns torsion-free sheaves of rank one
\cite[\S\S8.1--8.2]{KT}; these theorems do not give the higher-rank
statement of Theorem~\ref{thm:introB}.

Tsai's algorithm for Shalika
germs assumes sufficiently large residue characteristic. It separates
nilpotent orbital integrals by finitely many test functions and
descends from a good semisimple leading term to its twisted Levi
centralizer \cite[\S\S3--4, Thm.~4.1]{Tsai}; in equal characteristic it
has a geometric form, a stratification of auxiliary affine Springer
spaces \cite[Thm.~6.4, Prop.~6.6, Lem.~6.8, Rem.~6.9]{Tsai}.

 Our recursion follows the Euclidean transform of a framed branch,
 retaining the module rank, the full chain vector and the framing
 valuation at each stage. It terminates at a finite stable-flag count
 on a smooth branch and applies also to inseparable framings.

\paragraph{\textbf{From finite-rank transfer to stable cabling.}}
To identify the universal coefficients we evaluate both sides on split elements, where lattices can be counted directly. This step uses the multiplicativity of lattice counts over distinct branches, and the result is the Hall-algebra character of Theorem~\ref{thm:introC}(i) (Theorem~\ref{thm:hallchar}).

To prove Theorem~\ref{thm:introC}(ii), we compare the identities of
part~(i) for $A$ and $S$. After normalization, they give two triangular
linear systems with the same matrix, indexed by partitions with at
most $D$ parts. A uniform Gaussian bound for the inverse matrix
(Lemma~\ref{lem:inverse}) controls the contribution of the difference
of their right-hand sides and permits passage to the limit over a
fixed finite field.

For Theorem~\ref{thm:introM}(i) we compare the same two systems in the ring of power series in $Q$. Since all their entries are rational functions of $q$ that do not depend on the field (Theorem~\ref{thm:introB}), the comparison can be made $Q$-adically, and there the estimate for the inverse matrix is immediate. Part (ii) is a computation with generalized Dedekind eta functions: the vanishing of the exponents $\rho_{a,b}$ on multiples of $a+b$ gives weight zero, a Fourier computation on $\Z/(a+b)$ gives the power of $Q$ (Lemma~\ref{lem:rho}), and the transformation formula of Yang \cite{Yang1} gives the invariance under $\Gamma(12L)$ (Section~\ref{app:siegel}).

The geometric comparisons in this paper concern point counts over
finite fields. They do not provide isomorphisms, pavings or motivic
identities between the moduli spaces, nor do they determine the
individual fixed-colength strata of the punctual Quot scheme.
The fibrewise cabling identity remains Conjecture~S.

\subsection{Colored knot invariants and slope operators}\label{sub:colored-slopes}
Colored refinements of the relation between plane curve singularities
and knot invariants have several geometric forms. Diaconescu, Hua
and Soibelman proposed a description of colored HOMFLY polynomials
by framed stable pairs on curve thickenings \cite[\S1.7.2]{DHS};
its Euler-characteristic version was proved by Maulik
\cite[Thm.~1.2]{Maulik}. More directly relevant here are the
higher-rank standard modules of Cherednik--Philipp
\cite[Conj.~2.12]{CP} and Cherednik
\cite[Thm.~6.2 and Conj.~9.1]{Ch25}.
In Cherednik's latter normalization, setting the colength variable
equal to $1$ and $a=0$ gives $\#E_A(O^n)$, and the conjectural
DAHA counterpart has row color $(n)$.

Cabling has an explicit realization through double affine Hecke
algebras. Samuelson constructs cabling maps in type $A_1$ and
proves their colored Jones specialization
\cite[Def.~5.7 and Thm.~5.8]{Sa}. Cherednik--Danilenko construct
colored DAHA invariants of iterated torus knots by successive
automorphisms and projections to the polynomial representation
\cite[Thm.~2.1]{CD}; Cherednik gives the corresponding construction
for algebraic links in \cite[\S3.3]{Ch25}.
The identification of these invariants with colored knot homology
is conjectural in general.

We use Negu\c{t}'s slope operators, constructed through the shuffle
realization of the elliptic Hall algebra
\cite[Thm.~1.1, \S6.11 and Prop.~6.13]{NegShuffle} and described
by explicit contour formulas \cite[Thm.~2.15 and \S2.9]{Neg}.
Negu\c{t} also proved integrality of their action in the slope
stable basis \cite[Thm.~1.2 and \S4.4]{NegPieri} and constructed
an integral form of the positive shuffle algebra
\cite[Thm.~1.2 and Lem.~1.3]{NegIntegral}.

In the geometric normalization, specializing Negu\c{t}'s contour
formulas at $t=1$ gives multiplication operators, and the resulting
slope maps are ring homomorphisms. Gorsky--Negu\c{t} studied this
specialization for uncolored torus-knot invariants and proved its
Dyck-path formula \cite[\S6.6 and Thm.~6.7]{GN}.
Kivinen--Tsai apply the specialized slope operators to their
Shalika-germ recursion. We use their normalization and the
rational-Dyck-path description of the images of elementary
symmetric functions \cite[Props.~4.9, 4.12 and Def.~5.13]{KT},
which we recall below.

Let $\Lambda=\Q(q)[e_1,e_2,\ldots]$ be the ring of symmetric
functions, with Hall pairing
$\langle m_\lambda,h_\mu\rangle=\delta_{\lambda\mu}$ and
involution $\omega(e_r)=h_r$. For coprime positive integers $d,e$,
the specialized slope map
$\phi_{d/e}^{(q)}\colon\Lambda\to\Lambda$ is the
$\Q(q)$-algebra homomorphism characterized by
\begin{equation}\label{eq:intro-slope}
 \phi_{d/e}^{(q)}(e_r)
 =\sum_{\pi\in\mathcal D_{er,dr}}
       q^{\operatorname{area}(\pi)}e_{\operatorname{runs}(\pi)}.
\end{equation}
Here $\mathcal D_{er,dr}$ consists of lattice paths from $(0,0)$
to $(er,dr)$ lying weakly below the diagonal,
$\operatorname{runs}(\pi)$ records the lengths of the horizontal
runs, and $\operatorname{area}(\pi)$ counts the whole unit squares
between the path and the diagonal. Thus $\phi_{d/e}^{(q)}$ sends
degree $N$ to degree $eN$.

In forthcoming work we will give the following consequence of
Kivinen--Tsai's slope recursion and the universal transfer theorem.
It identifies the matrices $\mathsf C_{m,k}(q)$ of
Section~\ref{subsec:explicit-counting-algorithm} with the normalized
operator of \cite[Thm.~6.21, Cor.~6.22 and Thm.~6.35]{KT}.
The additional comparison expresses this operator in our chain
coordinates and proves the identity on the full distribution
space of Theorem~\ref{thm:introH}.

\begin{theorem}[Slope formula in chain coordinates]\label{thm:intro-slope-matrix}
Let $e,d,N\ge1$ with $\gcd(e,d)=1$, and put
\[
 a_{e,d}(N)=\frac{edN^2-(e+d-1)N}{2},\qquad
 \mathcal T_{e,d,N}
   =q^{a_{e,d}(N)}\omega\phi_{d/e}^{(q^{-1})}\omega.
\]
The universal transfer matrix, with rows $\alpha\vdash eN$ and
columns $\beta\vdash N$, is
\begin{equation}\label{eq:intro-slope-matrix}
 \bigl(\mathsf C_{eN,dN}(q)\bigr)_{\alpha\beta}
   =\big\langle\mathcal T_{e,d,N}(m_\beta),h_\alpha\big\rangle.
\end{equation}
Its entries lie in $\Z[q]$, and the identity specializes over every
finite field. Equivalently, for every $D\in J_{dN}$ for
$\GL_{eN}$ and every composition $\alpha$ of $eN$,
\begin{equation}\label{eq:intro-slope-distribution}
 D(F_{\alpha,dN})
 =\sum_{\beta\vdash N}
    \big\langle\mathcal T_{e,d,N}(m_\beta),h_\alpha\big\rangle
    D(f^*_\beta).
\end{equation}
\end{theorem}
In this formula $\phi_{d/e}^{(q^{-1})}$ is $\Q(q)$-linear:
only the coefficients of the operator are specialized at $q^{-1}$.
A chain vector $\mathbf v$ is encoded as
$\sum_{\beta\vdash N}v_\beta m_\beta$, so
\eqref{eq:intro-slope-matrix} uses monomial coordinates.
The comparison follows by evaluating on a spanning family of
tame regular semisimple orbital distributions, identifying the
transformed chains, and applying the universality of
Corollary~\ref{cor:basis}. This gives
\eqref{eq:intro-slope-distribution} also for the nonregular
distributions that occur in higher module rank.

The Newton tower below ends at the smooth frame $(t,0)$, with
period of valuation one. Its chain vector is
\begin{equation}\label{eq:intro-slope-seed}
 \widetilde h_n[X;q]
   =(q;q)_n h_n\!\left[\frac{X}{1-q}\right]
   =\sum_{\alpha\vdash n}
       \frac{(q;q)_n}{\prod_i(q;q)_{\alpha_i}}\,m_\alpha.
\end{equation}
The coefficient of $m_\alpha$ is the Gaussian multinomial
counting flags of type $\alpha$ in $\F_q^n$. The brackets denote
plethystic substitution, with
$p_r[X/(1-q)]=p_r(X)/(1-q^r)$.
Let $\epsilon_e\colon\Lambda\to\Q(q)$ be the algebra homomorphism
given by $\epsilon_e(e_r)=1$ for all $r\ge1$; on degree $s$ it is
the Hall pairing with $e_s=s_{(1^s)}$.

To describe the iteration, choose a tame representative
$A=F[[t^m,y(t)]]$ of a multiplicity sequence
$\underline\mu$, with $m=\mult(A)>1$ and
$\operatorname{char}F\nmid m$. Write its Puiseux characteristic
as $(m;\beta_1,\ldots,\beta_g)$, set
$g_0=m$ and $g_i=\gcd(m,\beta_1,\ldots,\beta_i)$, and put
\[
 e_i=\frac{g_{i-1}}{g_i},\qquad
 d_1=\frac{\beta_1}{g_1},\qquad
 d_i=\frac{\beta_i-\beta_{i-1}}{g_i}\quad(i\ge2),\qquad
 \Phi^{(q)}=
   \phi_{d_1/e_1}^{(q)}\circ\cdots\circ\phi_{d_g/e_g}^{(q)}.
\]
The rightmost slope acts first. Iterating
\eqref{eq:intro-slope-matrix} and using
$\omega\widetilde h_n[X;q]
 =q^{\binom n2}\widetilde h_n[X;q^{-1}]$ gives
\begin{equation}\label{eq:intro-colored-count}
 \begin{split}
 P_{\underline\mu,n}(q)
   &=q^{\mathcal B_n}\,
       \epsilon_e\!\left(
          \Phi^{(q^{-1})}(\widetilde h_n[X;q^{-1}])\right),\\
 \mathcal B_n
   &=n^2\delta+\binom n2(m+\beta_g).
 \end{split}
\end{equation}
Here $\delta$ is the delta invariant and $\beta_g$ is the last
Puiseux characteristic exponent. Every multiplicity sequence has
a tame representative in sufficiently large characteristic, so
Theorem~\ref{thm:introB} makes \eqref{eq:intro-colored-count} an
identity for its universal counting polynomial, valid for all
realizations of the sequence. The smooth sequence gives
$P_{\varnothing,n}=1$.

The normalized limit of \eqref{eq:intro-colored-count} extends the
$t=1$ Rogers--Ramanujan identities of Huang, Jiang and Oblomkov
\cite{HJO} from torus branches to arbitrary plane branches.
Setting $Q=q^{-1}$, Theorems~\ref{thm:introD} and~\ref{thm:introM}
give
\[
 \begin{split}
 \lim_{n\to\infty}
 Q^{-(m+\beta_g)\binom n2}\,
 \epsilon_e\!\left(\Phi^{(Q)}(\widetilde h_n[X;Q])\right)
 &=\cN_A(Q)\\
 &=\prod_{i=1}^g\ \prod_{j\ge1}
       (1-Q^j)^{-D_i\rho_{a_i,b_i}(j)},
 \end{split}
\]
where $(a_i,b_i,D_i)$ is the cabling list of
Theorem~\ref{thm:introD}. The expression under the limit equals
$Q^{n^2\delta}P_{\underline\mu,n}(Q^{-1})$ and agrees with the
product modulo $Q^{n+1}$. Thus the limit is coefficientwise in
$\Z[[Q]]$.

For a torus branch, the product identity was conjectured in
\cite{HJO} and proved by Huang, Lau and Ono
\cite[Cor.~1.2]{HLO}; for the cusp $(a,b)=(2,3)$ it is the first
classical Rogers--Ramanujan identity. The stable cabling formula
expresses the limit for a general plane branch as a product of
these torus identities. These stable identities hold independently
of Conjecture~P.

The finite-rank identity \eqref{eq:intro-colored-count} also shows
that Conjecture~P is equivalent to the combinatorial assertion
\begin{equation}\label{eq:intro-slope-positivity}
 \epsilon_e\!\left(\Phi^{(q)}(\widetilde h_n[X;q])\right)
       \in\Z_{\ge0}[q]
\end{equation}
for every plane-branch slope sequence and every $n\ge1$.
Indeed, \eqref{eq:intro-colored-count} reverses the coefficients
and shifts the degrees, while Theorem~\ref{thm:introB} excludes
negative powers in the resulting count. Although each slope map
preserves nonnegative elementary coefficients, the seed need not
have that property: $\widetilde h_2=e_1^2+(q-1)e_2$.
Thus positivity of the Dyck-path weights alone does not settle
\eqref{eq:intro-slope-positivity}.

In the same forthcoming work we prove
\eqref{eq:intro-slope-positivity} in several families, in every
module rank. In fact the symmetric functions
\[
 \begin{gathered}
 \phi_{d/e}^{(q)}\phi_{1/D}^{(q)}(\widetilde h_n[X;q]),\qquad
 \phi_{d/e}^{(q)}\phi_{(2k+1)/2}^{(q)}(\widetilde h_n[X;q]),\\
 \phi_{d/e}^{(q)}\phi_{1/D}^{(q)}
       \phi_{(2k+1)/2}^{(q)}(\widetilde h_n[X;q])
 \end{gathered}
\]
are Schur-positive over $\Z_{\ge0}[q]$ for $n,D\ge1$, $k\ge0$
and coprime $d,e\ge1$. Since $\epsilon_e$ extracts the coefficient
of a column Schur function, this proves the required scalar
positivity. The first family includes the multiplicity sequences
of $F[[t^{aD},t^{bD}+t^{bD+1}]]$ with $1<a<b$ coprime and
$D\ge2$, represented in characteristic prime to $aD$.
The general case of \eqref{eq:intro-slope-positivity} remains open.

\subsection{External inputs and organization}
The general chain-transfer algorithm is independent of the torus
formulas. We use \cite{HJO} for the finite-rank torus formula and
\cite{HLO} for its infinite product, which enter coefficientwise
stabilization and modularity. Integral coefficients use rationality
of finite-field zeta functions \cite{Dwork} through
Lemma~\ref{lem:completion-integral}. Other inputs are Hall
polynomials, the Auslander--Buchsbaum formula, and standard facts
about local fields and Haar measures.

Section~\ref{sec:norm} proves the mass normalization and finiteness
criterion. Sections~\ref{sec:qt}--\ref{sec:chains} establish the
framed transform and universal transfer, then prove
Theorem~\ref{thm:introB}, including integral coefficients.
Section~\ref{subsec:explicit-counting-algorithm} gives the recursion. Section~\ref{sec:hall} identifies its scalar coefficients,
and Section~\ref{sec:limit} proves stable cabling and iteration.
Sections~\ref{sec:torus}--\ref{sec:modular} give the torus inputs,
coefficientwise limit, and modular interpretation.
Section~\ref{sec:examples} contains examples and computations.

Appendix~\ref{app:integrality} constructs the bounded moduli schemes
for the refined counts. Appendix~\ref{sec:algorithm-implementation}
contains the finite subroutines, their termination bounds, and a
worked symbolic example.
Appendix~\ref{app:conventions} proves the affine-Weyl reduction,
Bruhat decomposition and Hecke relations, and modular transformation
formulas in our conventions. It includes an exact regularity test
and a finite construction of standard triples, and concludes with
the character normalization of the cabling product.

\subsection{Conventions}
Throughout the paper $F=\F_q$, $Q=q^{-1}$, $O=F[[t]]$, $\cK=F((t))$ and $\val=\val_t$. From Section~\ref{sec:qt} on we also use the products $\mathbf O=O^r$ and $\mathbf K=\cK^r$, and for the groups $\GL_m(\K)$ an auxiliary local field $\K$, equal to $F((\pi))$ from Section~\ref{sec:chains} on.

Partitions are written $\mu=(\mu_1\ge\mu_2\ge\cdots)$, $\mu'$ is the conjugate partition, $\mu\vdash N$ means that $|\mu|=N$, $\ctd$ is the dominance order and $p(N)$ is the number of partitions of $N$. The \emph{Jordan type} of a nilpotent endomorphism is the partition formed by the sizes of its Jordan blocks, and we denote by $\bO_\mu$ the set of nilpotent matrices of type $\mu$ and by $M_\mu$ the corresponding $F[z]$-module.

We write $(x;y)_j=\prod_{i=0}^{j-1}(1-xy^i)$ for $j\le\infty$, thus $(x;x)_j=\prod_{i=1}^j(1-x^i)$, and $\qbinom nk_x$ stands for the Gaussian binomial coefficient.

Haar measures are normalized as indicated. For a locally compact totally disconnected unimodular group $Z$ acting on a set $S$ with compact open stabilizers we put $\mass_Z(S)=\sum_{s\in Z\backslash S}\vol(Z_s)^{-1}\in[0,\infty]$, and if $S\to S'$ is a $Z$-equivariant map with finite fibres, then
\begin{equation}\label{eq:massfibre}
\mass_Z(S)=\sum_{s'\in Z\backslash S'}\frac{\#\text{fibre over }s'}{\vol(Z_{s'})},
\end{equation}
because $\sum_{s\in Z_{s'}\backslash\text{fibre}}\vol(Z_s)^{-1}=\vol(Z_{s'})^{-1}\sum_s[Z_{s'}:Z_s]$.

\noindent{\bf Acknowledgements.} I would like to thank Yifeng Huang and Ruofan Jiang for many discussions and, most importantly, for introducing me to the circle of ideas related to
the Cohen--Lenstra series. I would like to thank Cheng-Chiang Tsai for explaining his paper \cite{Tsai} back in 2016. His paper ignited my long interest in the mathematics of local fields. The current paper is my attempt to understand better the paper \cite{Tsai} and the subsequent work \cite{KT}.
I would like to thank Eugene Gorsky, Oscar Kivinen, Lev Rozansky and Dimitri Wyss for discussions related to various aspects of the paper.
This material is based upon work supported by the National Science Foundation under Grant No. DMS-2424139 while the author was in residence at the Simons Laufer Mathematical Sciences Institute in Berkeley,
California, during the Fall 2026 semester.

\noindent{\bf Declaration of AI use.} During the preparation of this work, the author used Claude and Claude Code to assist with mathematical discussions, the development and exposition of proofs, literature searches, and the drafting and editing of the manuscript. The author independently verified all mathematical arguments and references, revised all AI-assisted material, and take full responsibility for the content of the paper.

\section{Lattices and the Cohen--Lenstra normalization}\label{sec:norm}

Residue duality gives the normalization of the lattice count.
Theorem~\ref{prop:cohnorm} identifies the normalized count with a
special value of the Quot series, and Corollary~\ref{cor:limit-exists}
expresses its limit as a mass of finite modules. To prove that this mass is finite, we
rewrite the normalized count using modules over the finite conductor
quotient (Proposition~\ref{prop:conductor}). For a plane branch, a
commuting-matrix estimate then gives a bound uniform in the rank
(Theorem~\ref{thm:finite}); this bound will also control the passage
to the limit in Section~\ref{sec:limit}.

\subsection{Branches and lattices}
Let $A\subseteq O$ be a branch and let $\delta=\delta_A$. The action map $A\to\End_F(O/A)$ has a nonzero kernel,
since $A$ is infinite-dimensional over $F$ whereas
$\End_F(O/A)$ is finite-dimensional. For a nonzero element $a$
of this kernel, $aO\subseteq A$. Thus $t^cO\subseteq A$ for
$c=\val(a)$, and in particular $A$ is closed in $O$. A \emph{lattice} in $\cK^n$ is an $F$-subspace $\Lambda$ with $t^NO^n\subseteq\Lambda\subseteq t^{-N}O^n$ for some $N$, and an $A$-lattice is a lattice which is an $A$-submodule, in agreement with the definition in the introduction. For lattices $\Lambda,\Lambda'$ the index is defined as follows:
\[
[\Lambda:\Lambda']=\dim_F\Lambda/(\Lambda\cap\Lambda')-\dim_F\Lambda'/(\Lambda\cap\Lambda').
\]
Let us recall that the index is additive, $[\Lambda:\Lambda'']=[\Lambda:\Lambda']+[\Lambda':\Lambda'']$, and that $[\Lambda:g\Lambda]=\val\det g$ for $g\in\GL_n(\cK)$. We normalize the Haar measure of $\cK^n$ by $\vol(O^n)=1$, and then $\vol(\Lambda)=q^{-[O^n:\Lambda]}$. Note that the set $E_A(O^n)$ is finite, because each of its elements contains $t^cO^n$.

For $\xi,\lambda\in\cK^n$ we write $\langle\xi,\lambda\rangle=\sum_j\xi_j\lambda_j\in\cK$, and for a lattice $\Lambda$, respectively an $A$-lattice $\Lambda$, we define $\Lambda^\perp$, respectively $\Lambda^\vee$, as follows:
\begin{align*}
\Lambda^\perp&=\{\xi\in\cK^n:\ \Res_{t=0}\big(\langle\xi,\lambda\rangle\,dt\big)=0\ \text{ for all }\lambda\in\Lambda\},\\
\Lambda^\vee&=\{\xi\in\cK^n:\ \langle\xi,\lambda\rangle\in A\ \text{ for all }\lambda\in\Lambda\}.
\end{align*}
We say that $A$ satisfies condition (G) if
\begin{itemize}[leftmargin=3em]
\item[(G)] there is $\kappa\in\cK^\times$ with $A^\perp=\kappa A$ (here $n=1$).
\end{itemize}

\begin{lemma}\label{lem:duality}
\begin{enumerate}[label=\textup{(\alph*)}]
\item We have $(t^jO^n)^\perp=t^{-j}O^n$, and if $\Lambda\subseteq\Lambda'$ are lattices, then $\Lambda'^\perp\subseteq\Lambda^\perp$ and $\dim\Lambda^\perp/\Lambda'^\perp$ equals $\dim\Lambda'/\Lambda$. Hence $[O^n:\Lambda^\perp]=-[O^n:\Lambda]$ for every lattice $\Lambda$.
\item If $A$ satisfies \textup{(G)}, then $\Lambda^\perp=\kappa\Lambda^\vee$ for every $A$-lattice $\Lambda$.
\item If $A$ satisfies \textup{(G)}, then $\val(\kappa)=-2\delta$, and $[O^n:\Lambda^\vee]=2\delta n-[O^n:\Lambda]$ for every $A$-lattice $\Lambda$.
\end{enumerate}
\end{lemma}
\begin{proof}
(a) The first statement holds because the residue $\Res(t^it^j\,dt)$ is nonzero exactly when $i+j=-1$. Next, if $N$ is such that $t^NO^n\subseteq\Lambda\subseteq\Lambda'\subseteq t^{-N}O^n$, then the residue pairing induces a perfect pairing of the finite-dimensional space $t^{-N}O^n/t^NO^n$ with itself, and the lattice $\Lambda^\perp$ is the preimage of the annihilator of $\Lambda/t^NO^n$. Since annihilators reverse inclusions and have complementary dimensions, the second statement follows. The last statement follows by comparing $\Lambda$ and $O^n$ with $\Lambda\cap O^n$.

(b) Since $\Lambda$ is an $A$-module, $\xi\in\Lambda^\perp$ if and only if $\Res(\langle\xi,\lambda\rangle r\,dt)=0$ for all $r\in A$ and $\lambda\in\Lambda$, that is, if and only if $\langle\xi,\lambda\rangle\in A^\perp=\kappa A$ for all $\lambda$, or equivalently $\kappa^{-1}\xi\in\Lambda^\vee$.

(c) By (a) we have $[O:A^\perp]=-[O:A]=-\delta$, and on the other hand $[O:\kappa A]=[O:\kappa O]+[\kappa O:\kappa A]=\val(\kappa)+\delta$. The equality $A^\perp=\kappa A$ gives $\val(\kappa)=-2\delta$, and then (b) and (a) imply:
\[
[O^n:\Lambda^\vee]=[O^n:\kappa^{-1}O^n]+[\kappa^{-1}O^n:\kappa^{-1}\Lambda^\perp]=2\delta n+[O^n:\Lambda^\perp]=2\delta n-[O^n:\Lambda].\qedhere
\]
\end{proof}

\begin{lemma}\label{lem:symmetric}
Let $\Gamma=\val(A\smallsetminus\{0\})$ be the value semigroup of $A$, and let $c$ be its conductor, the least integer with $c+\N\subseteq\Gamma$. Then the set of valuations of the nonzero elements of $A^\perp$ is $\{j\in\Z:\ -1-j\notin\Gamma\}$. If $\Gamma$ is symmetric, that is, $j\in\Gamma\Leftrightarrow c-1-j\notin\Gamma$ for all $j\in\Z$, then \textup{(G)} holds.
\end{lemma}
\begin{proof}
Let us observe first that $O=O^\perp\subseteq A^\perp$ and that $\dim A^\perp/O=\delta$ by Lemma~\ref{lem:duality}(a). If $\xi\in A^\perp$ has valuation $j$ and $-1-j=\val(r)$ with $r\in A$, then $\xi r$ has valuation $-1$ and nonzero residue, a contradiction. Thus the valuations of the nonzero elements of $A^\perp$ lie in $\Gamma'=\{j:-1-j\notin\Gamma\}$. Since the set $\N\smallsetminus\Gamma$ has $\delta$ elements, the set $\Gamma'$ contains every $j\ge0$ and exactly $\delta$ negative integers. On the other hand, $\dim A^\perp/O$ is the number of negative valuations occurring in $A^\perp$, hence every element of $\Gamma'$ occurs.

Now let $\Gamma$ be symmetric. Then $\Gamma'=\Gamma-c$, and for $\kappa\in A^\perp$ of valuation $-c$ the lattices $\kappa A\subseteq A^\perp$ have the same set of valuations, hence the same index in $t^{-c}O$, and thus $\kappa A=A^\perp$.
\end{proof}

\begin{remark}\label{rem:gorenstein}
Every plane branch satisfies (G), in every characteristic, and a direct proof is given in Lemma~\ref{lem:planeG}. More generally, for a branch $A$ the lattice $A^\perp$ is a canonical module, thus (G) says that $A$ is Gorenstein. By a theorem of Kunz \cite{Kunz}, this holds if and only if the value semigroup of $A$ is symmetric, and Lemma~\ref{lem:symmetric} is the easy half of this statement. We use (G) only through Lemma~\ref{lem:duality}.
\end{remark}

\subsection{The finite-rank identity}
\begin{theorem}\label{prop:cohnorm}
Let $A$ be a branch satisfying \textup{(G)}, for instance a plane branch. Then for every $n\ge1$ we have:
\begin{equation}\label{eq:cohnorm}
\sum_{L\subseteq A^n}\big|A^n/L\big|^{-n}\;=\;q^{-n^2\delta}\,\#E_A(O^n)\,\prod_{i=1}^n(1-q^{-i})^{-1},
\end{equation}
where the sum is over the $A$-submodules of $A^n$ of finite colength.
\end{theorem}
\begin{proof}
Let us denote by $dg$ the Haar measure of $\GL_n(\cK)$ with $\vol(\GL_n(O))=1$ and by $d^+g$ the additive Haar measure of $\Mat_n(\cK)$ with $\vol^+(\Mat_n(O))=1$, and we put $|\det g|=q^{-\val\det g}$. Since the measure $|\det g|^{-n}d^+g$ is a Haar measure of $\GL_n(\cK)$, we have:
\begin{equation}\label{eq:haar}
|\det g|^n\,dg=c_n\,d^+g\ \text{ on }\GL_n(\cK),\qquad c_n=\vol^+(\GL_n(O))^{-1}=\prod_{i=1}^n(1-q^{-i})^{-1}.
\end{equation}
Note also that the complement of $\GL_n(\cK)$ in $\Mat_n(\cK)$ has measure zero.

\emph{Step 1. Parametrization.} If $L\subseteq A^n$ is an $A$-submodule of finite colength, then $OL$ is an $O$-lattice. Thus $OL=gO^n$ for some $g\in\GL_n(\cK)$, unique up to right multiplication by $\GL_n(O)$, and the lattice $\Lambda=g^{-1}L$ belongs to $E_A(O^n)$. Conversely, if $g\in\GL_n(\cK)$ and $\Lambda\in E_A(O^n)$ satisfy $g\Lambda\subseteq A^n$, then $L=g\Lambda$ is such a submodule and $OL=gO^n$.

For a fixed $L$, the admissible matrices $g$ are exactly those with
$gO^n=OL$, and form one right coset of $\GL_n(O)$. For each such
$g$, the lattice $\Lambda=g^{-1}L$ is uniquely determined.
Our normalization gives this coset measure~$1$.
All summands are nonnegative, so Tonelli's theorem permits the
interchange of the sum and integral without a prior convergence
assumption. Consequently, for every nonnegative function $\phi$ on the finite-colength submodules
of $A^n$, we have:
\[
\sum_L\phi(L)=\int_{\GL_n(\cK)}\ \sum_{\Lambda\in E_A(O^n)}\mathbf 1[g\Lambda\subseteq A^n]\ \phi(g\Lambda)\ dg .
\]

\emph{Step 2. The index.} For $g$ and $\Lambda$ as above we have:
\[
\dim_FA^n/g\Lambda=[A^n:O^n]+[O^n:gO^n]+[gO^n:g\Lambda]=-n\delta+\val\det g+[O^n:\Lambda],
\]
and therefore $|A^n/g\Lambda|^{-n}=q^{n^2\delta}\,|\det g|^n\,q^{-n[O^n:\Lambda]}$.

\emph{Step 3. A volume.} For $g\in\Mat_n(\cK)$ the condition $g\Lambda\subseteq A^n$ means that every row of $g$ lies in $\Lambda^\vee$. Hence the set of these $g$ has $d^+$-volume $\vol(\Lambda^\vee)^n=q^{-n[O^n:\Lambda^\vee]}$, and by Lemma~\ref{lem:duality}(c) this volume is equal to $q^{-2\delta n^2+n[O^n:\Lambda]}$.

\emph{Step 4.} We apply Steps 1 and 2 with $\phi(L)=|A^n/L|^{-n}$, then \eqref{eq:haar} and Step 3, and obtain:
\[
\sum_L\big|A^n/L\big|^{-n}=\sum_{\Lambda\in E_A(O^n)}q^{n^2\delta-n[O^n:\Lambda]}\cdot c_n\cdot q^{-2\delta n^2+n[O^n:\Lambda]}=c_n\,q^{-n^2\delta}\,\#E_A(O^n).\qedhere
\]
\end{proof}

\begin{corollary}\label{cor:limit-exists}
Let $A$ be a branch satisfying \textup{(G)}, with maximal ideal $\fm$, and let $e_A(n)=q^{-n^2\delta}\#E_A(O^n)$.
\begin{enumerate}[label=\textup{(\roman*)}]
\item For all $n\ge1$ we have:
\[
e_A(n)\prod_{i=1}^n(1-q^{-i})^{-1}=\sum_{[M]:\ g(M)\le n}\frac1{|\Aut_A(M)|}\prod_{i=0}^{g(M)-1}(1-q^{i-n}),
\]
where the sum is over the isomorphism classes of $A$-modules $M$ of finite length with $g(M):=\dim_FM/\fm M\le n$.
\item The limit $\cN_A=\lim_{n\to\infty}e_A(n)$ exists in $(0,\infty]$, and it is given by the formula:
\[
\cN_A=(q^{-1};q^{-1})_\infty\sum_{[M]}\frac1{|\Aut_A(M)|}.
\]
\item $\cN_A<\infty$ if and only if the sequence $e_A(n)$ is bounded.
\end{enumerate}
\end{corollary}
\begin{proof}
(i) The number of submodules $L\subseteq A^n$ with $A^n/L\cong M$ is $|\Sur_A(A^n,M)|/|\Aut_A(M)|$. By Nakayama's lemma an $n$-tuple of elements generates $M$ if and only if its image spans $M/\fm M\cong F^{g(M)}$, and there are $\prod_{i=0}^{g(M)-1}(q^n-q^i)$ spanning $n$-tuples in $F^{g(M)}$, each of which has $|\fm M|^n$ lifts. Hence $|\Sur_A(A^n,M)|=|M|^n\prod_{i=0}^{g(M)-1}(1-q^{i-n})$, and (i) is \eqref{eq:cohnorm} with the submodules $L$ grouped according to the isomorphism class of $A^n/L$.

(ii) Each term on the right of (i) is nonnegative and nondecreasing in $n$, with limit $|\Aut_A(M)|^{-1}$, and the product on the left converges to $(q^{-1};q^{-1})_\infty^{-1}$, so (ii) follows from (i) by monotone convergence. The term $M=0$ shows that $\cN_A\ge(q^{-1};q^{-1})_\infty>0$.

(iii) This follows from (ii).
\end{proof}

\begin{remark}\label{rem:coh}
(1) In the notation of \cite{HJ,HJO}, the left side of \eqref{eq:cohnorm} is the value at $t=q^{-n}$ of the Quot zeta function $\sum_L t^{\dim A^n/L}$ of $A^n$, and the theorem says that its normalization by the corresponding function of $O^n$ takes the value $e_A(n)$ there. In \cite{HJO} this is proved for $A=F[[t^a,t^b]]$ from the functional equation of \cite{HJ} and the value at $t=1$, while the argument above gives the special value directly. It uses the residue pairing of $\cK$ rather than a trace form, thus no separability condition enters, but it does not give the functional equation.

(2) If $A=F[[X,Y]]/(f)$, then an $A$-module structure on $F^m$ is a pair of commuting nilpotent matrices $(X,Y)$ with $f(X,Y)=0$, and by the orbit--stabilizer formula we have:
\[
\sum_{[M]}|\Aut_A(M)|^{-1}=\sum_{m\ge0}\#\{(X,Y)\}/|\GL_m(F)|.
\]
Thus Corollary~\ref{cor:limit-exists}(ii) expresses $\cN_A$ as a normalized count of matrix points of the curve.

(3) By Corollary~\ref{cor:limit-exists}, the existence of the limit $\cN_A$ requires no further input, whereas its finiteness is not formal. For plane branches we prove it in Theorem~\ref{thm:finite} below, and it fails for all other branches satisfying (G) (Corollary~\ref{cor:planarity-finiteness}). For $A=O$ we have $E_O(O^n)=\{O^n\}$ and $\cN_O=1$, and in this case (ii) is the classical identity:
\[
\sum_\lambda|\Aut(M_\lambda)|^{-1}=(q^{-1};q^{-1})_\infty^{-1}.
\]
\end{remark}

\subsection{The lattices between \texorpdfstring{$A^n$}{A^n} and \texorpdfstring{$O^n$}{O^n}}\label{sub:finite}
A finite sum over the $A$-lattices between $A^n$ and $O^n$ gives a
second expression for $e_A(n)$, from which we will deduce the
finiteness of $\cN_A$ for plane branches.

\begin{proposition}\label{prop:finitesum}
Let $A$ be a branch. Then for every $n\ge1$ we have:
\begin{equation}\label{eq:finitesum}
\#E_A(O^n)=\sum_{A^n\subseteq L\subseteq O^n}\big|O^n/L\big|^n,\qquad\text{and hence}\qquad e_A(n)=\sum_{A^n\subseteq L\subseteq O^n}\big|L/A^n\big|^{-n},
\end{equation}
where the sums are over the $A$-submodules $L$ of $O^n$ which contain $A^n$.
\end{proposition}
\begin{proof}
Let us choose $c\ge1$ with $t^cO\subseteq A$ and put $R=O/t^cO$. A matrix in $\Mat_n(R)$ is invertible if and only if its reduction modulo $t$ is invertible, thus $|\GL_n(R)|=\gamma_n|R|^{n^2}$, where $\gamma_n=q^{-n^2}|\GL_n(F)|$. Every $L\in E_A(O^n)$ contains $t^cO^n$, because $t^cO^n=(t^cO)\,OL=(t^cO)L\subseteq AL=L$, and we denote by $\mathcal G(L)$ the set of matrices $g\in\GL_n(R)$ all of whose columns lie in $L/t^cO^n$. Since $O^n=OL=(F+tO)L=L+tO^n$, the reduction modulo $t$ is a surjective $F$-linear map $L/t^cO^n\to F^n$. Hence the $n$-tuples of elements of $L/t^cO^n$ whose reductions form a basis of $F^n$ make up the fraction $\gamma_n$ of all $n$-tuples, that is:
\[
|\mathcal G(L)|=\gamma_n\,\big|L/t^cO^n\big|^n .
\]

Let us now count the pairs $(g,L)$ with $L\in E_A(O^n)$ and $g\in\mathcal G(L)$, each of them with the weight $|\mathcal G(L)|^{-1}$. Summing over $g$ for a fixed $L$, we see that the weighted number of pairs is $\#E_A(O^n)$.

On the other hand, let $g\in\GL_n(R)$, and let $\tilde g\in\Mat_n(O)$ be a matrix with reduction $g$, which lies in $\GL_n(O)$ because its determinant is a unit. For an $A$-submodule $L$ of $O^n$ which contains $t^cO^n$, the columns of $g$ lie in $L/t^cO^n$ if and only if $\tilde gA^n\subseteq L$. The $A$-submodules $L$ of $O^n$ with $\tilde gA^n\subseteq L$ are the modules $\tilde gL'$, where $L'$ runs over the $A$-submodules of $O^n$ which contain $A^n$. Each of them contains $\tilde g(t^cO^n)=t^cO^n$, belongs to $E_A(O^n)$, because $O\tilde gA^n=O^n$, and satisfies $|\tilde gL'/t^cO^n|=|L'/t^cO^n|$. Therefore the weighted number of pairs is also equal to:
\[
\sum_{g\in\GL_n(R)}\ \sum_{A^n\subseteq L'\subseteq O^n}\frac1{\gamma_n\,|L'/t^cO^n|^n}=\sum_{A^n\subseteq L'\subseteq O^n}\Big(\frac{|R|^n}{|L'/t^cO^n|}\Big)^n=\sum_{A^n\subseteq L'\subseteq O^n}\big|O^n/L'\big|^n .
\]
This proves the first formula, and the second one follows, because $|O^n/L|\cdot|L/A^n|=|O^n/A^n|=q^{n\delta}$.
\end{proof}

For a branch $A$ satisfying (G) we put $\fc=t^{2\delta}O$, and for an $A$-module $M$ of finite length we denote by $s(M)$ the dimension over $F$ of its socle $\{m\in M:\fm m=0\}$, where $\fm$ is the maximal ideal of $A$.

The quotients of $A^n$ in Corollary~\ref{cor:limit-exists} are
controlled by their number of generators. For the submodules of
$(O/A)^n$ in Proposition~\ref{prop:finitesum}, the corresponding
invariant is the socle dimension. Indeed, residue duality identifies
$O/A$ with the dual of $A/\fc$, and injectivity into $(O/A)^n$ can
be tested on the socle. This gives the factors involving $s(M)$
in the next formula.

\begin{proposition}\label{prop:conductor}
Let $A$ be a branch satisfying \textup{(G)}.
\begin{enumerate}[label=\textup{(\roman*)}]
\item We have $\fc\subseteq A$, and the pairing $(\xi,r)\mapsto\Res_{t=0}(\kappa\xi r\,dt)$ induces an isomorphism of $A$-modules $O/A\cong\Hom_F(A/\fc,F)$. Consequently $|\Hom_A(M,O/A)|=|M/\fc M|$ for every $A$-module $M$ of finite length.
\item For all $n\ge1$ we have:
\[
e_A(n)=\sum_{[M]:\ \fc M=0}\frac1{|\Aut_A(M)|}\prod_{i=0}^{s(M)-1}(1-q^{i-n}),
\]
where the sum is over the isomorphism classes of $A$-modules $M$ of finite length with $\fc M=0$.
\item $\displaystyle\cN_A=\sum_{[M]:\ \fc M=0}\frac1{|\Aut_A(M)|}$, and $e_A(n)$ increases to $\cN_A$ as $n\to\infty$.
\end{enumerate}
\end{proposition}
\begin{proof}
(i) By Lemma~\ref{lem:duality}(c) we have $\val(\kappa)=-2\delta$, that is, $\kappa^{-1}O=\fc$, and since $A\subseteq O$, Lemma~\ref{lem:duality}(a) gives $O=O^\perp\subseteq A^\perp=\kappa A$, hence $\fc\subseteq A$. For $\xi\in O$ and $r\in A$ the residue $\Res(\kappa\xi r\,dt)$ vanishes if $r\in\fc$, because then $\kappa\xi r\in O$, and it vanishes if $\xi\in A$, because then $\kappa\xi r\in\kappa A=A^\perp$ and $1\in A$. Thus $(\xi,r)\mapsto\Res(\kappa\xi r\,dt)$ is a well-defined pairing between $O/A$ and $A/\fc$, and the multiplication by an element of $A$ in the first argument corresponds to the multiplication by the same element in the second one.

If $\xi\in O$ pairs to zero with all $r\in A$, then $\kappa\xi\in A^\perp=\kappa A$, that is, $\xi\in A$. Since $\dim O/A=\delta=\dim A/\fc$, the pairing is perfect, and it gives the isomorphism of $A$-modules $O/A\cong\Hom_F(A/\fc,F)$, where $a\in A$ acts on a linear form $\phi$ by $(a\phi)(r)=\phi(ar)$. Finally, for an $A$-module $M$ of finite length we have $\Hom_A(M,\Hom_F(A/\fc,F))=\Hom_F(M\otimes_AA/\fc,F)$, and $M\otimes_AA/\fc=M/\fc M$.

(ii) By Proposition~\ref{prop:finitesum}, $e_A(n)$ is the sum of the numbers $|M'|^{-n}$ over the $A$-submodules $M'$ of $(O/A)^n$, and all these modules are annihilated by $\fc$, because $\fc O\subseteq A$. Let $M$ be an $A$-module of finite length with $\fc M=0$. The submodules of $(O/A)^n$ isomorphic to $M$ are the images of the injective homomorphisms $M\to(O/A)^n$, and two of these have the same image if and only if they differ by an automorphism of $M$.

By (i), the homomorphisms $\phi=(\phi_1,\dots,\phi_n)\colon M\to(O/A)^n$ correspond to the $n$-tuples $(\lambda_1,\dots,\lambda_n)$ of linear forms on $M$, where $\phi_i(m)$ is the form $r\mapsto\lambda_i(rm)$. The kernel of $\phi$ is the largest $A$-submodule of $M$ on which all $\lambda_i$ vanish. A nonzero submodule of $M$ contains a nonzero element $m$ of the socle, and for such an element we have $Am=Fm$. Thus $\phi$ is injective if and only if the linear map $(\lambda_1,\dots,\lambda_n)$ from the socle of $M$ to $F^n$ is injective.

Since the restriction of linear forms from $M$ to its socle is surjective, and the number of injective linear maps $F^s\to F^n$ is $\prod_{i=0}^{s-1}(q^n-q^i)$, the number of injective homomorphisms $M\to(O/A)^n$ is equal to:
\[
|M|^n\prod_{i=0}^{s(M)-1}(1-q^{i-n}) .
\]
Dividing by $|\Aut_A(M)|$, multiplying by $|M|^{-n}$ and summing over the isomorphism classes, we obtain (ii).

(iii) Each term of the sum in (ii) is nonnegative and nondecreasing in $n$, with limit $|\Aut_A(M)|^{-1}$, so (iii) follows from (ii) by monotone convergence.
\end{proof}

The next lemma computes the mass of finite-length modules over
$F[[X,Y]]/(Y^k)$ by counting pairs of commuting nilpotent matrices.
In particular, this mass is finite.

\begin{lemma}\label{lem:ypower}
For integers $d\ge0$, $k\ge1$ and $\ell\ge0$ let $\mathcal X_{d,k,\ell}$ be the set of pairs $(X,Y)$ of commuting nilpotent matrices in $\Mat_d(F)$ with $Y^k=0$ and $\dim\ker Y=\ell$, and let $Q=q^{-1}$. Then we have:
\[
\sum_{d\ge0}\frac{|\mathcal X_{d,k,\ell}|}{|\GL_d(F)|}=Q^{\ell}\sum_{m_1+\dots+m_k=\ell}\ \prod_{i=1}^k\frac1{(Q;Q)_{m_i}},\qquad\sum_{\ell\ge0}\ \sum_{d\ge0}\frac{|\mathcal X_{d,k,\ell}|}{|\GL_d(F)|}=(Q;Q)_\infty^{-k},
\]
where the first sum on the right is over the vectors $(m_1,\dots,m_k)$ of nonnegative integers with sum $\ell$.
\end{lemma}
\begin{proof}
\emph{Step 1. Nilpotent matrices.} Let $\nu_d$ be the number of nilpotent matrices in $\Mat_d(F)$. By Fitting's lemma, every $Z\in\Mat_d(F)$ gives a decomposition $F^d=\ker Z^d\oplus Z^dF^d$, where $Z$ is nilpotent on the first summand and invertible on the second one, and conversely, a decomposition $F^d=V_0\oplus V_1$ together with a nilpotent endomorphism of $V_0$ and an automorphism of $V_1$ defines a matrix $Z$ with $\ker Z^d=V_0$ and $Z^dF^d=V_1$. The group $\GL_d(F)$ acts transitively on the decompositions with $\dim V_0=j$, with stabilizer $\GL_j(F)\times\GL_{d-j}(F)$. Counting all matrices in this way and using $|\GL_d(F)|=q^{d^2}(Q;Q)_d$, we get:
\begin{equation}\label{eq:fitting}
\sum_{j=0}^d\frac{\nu_j}{|\GL_j(F)|}=\frac{q^{d^2}}{|\GL_d(F)|}=\frac1{(Q;Q)_d}\qquad(d\ge0).
\end{equation}
Taking the difference of two consecutive identities, we obtain:
\[
\frac{\nu_d}{|\GL_d(F)|}=\frac{Q^d}{(Q;Q)_d}\quad(d\ge0),\qquad\text{and by \eqref{eq:fitting}}\qquad\sum_{m\ge0}\frac{Q^m}{(Q;Q)_m}=\frac1{(Q;Q)_\infty} .
\]

\emph{Step 2. Centralizers.} Let $Y\in\Mat_d(F)$ be a nilpotent matrix with $Y^k=0$ which has $m_i$ Jordan blocks of size $i$ for $1\le i\le k$, and let $\mathcal Z$ be its centralizer in $\Mat_d(F)$. The subspaces $K_j=\ker Y\cap Y^jF^d$, $0\le j\le k$, form a decreasing filtration of $K_0=\ker Y$ with $K_k=0$, and $\dim K_{i-1}/K_i=m_i$, because $\dim K_j$ is the number of blocks of size greater than $j$. Every $Z\in\mathcal Z$ preserves this filtration, and we denote by $\operatorname{gr}(Z)$ the family of the endomorphisms induced by $Z$ on the spaces $K_{i-1}/K_i$. The map
\[
\operatorname{gr}\colon\mathcal Z\to\prod_{i=1}^k\End_F(K_{i-1}/K_i)\cong\prod_{i=1}^k\Mat_{m_i}(F)
\]
is a homomorphism of algebras. It is surjective: in a Jordan basis we have $F^d=\bigoplus_i\big(F[y]/(y^i)\big)\otimes F^{m_i}$, where $Y$ acts as $y\otimes1$, and for $h_i\in\Mat_{m_i}(F)$ the endomorphism $\bigoplus_i1\otimes h_i$ commutes with $Y$ and induces $h_i$ on $K_{i-1}/K_i=y^{i-1}\otimes F^{m_i}$.

Moreover, every element $Z$ of the kernel of $\operatorname{gr}$ is nilpotent. Indeed, such a $Z$ maps $K_{i-1}$ into $K_i$ for all $i$, hence $Z^k$ vanishes on $\ker Y$. For $j\ge1$ the map $Y^j$ sends $\ker Y^{j+1}$ into $\ker Y$ and commutes with $Z^k$, thus $Z^k$ maps $\ker Y^{j+1}$ into $\ker Y^j$, and we get $Z^{k^2}=0$ on $\ker Y^k=F^d$.

Consequently, an element $Z\in\mathcal Z$ is nilpotent if and only if $\operatorname{gr}(Z)$ is nilpotent, and it is invertible if and only if $\operatorname{gr}(Z)$ is invertible. Indeed, if $\operatorname{gr}(Z)^N=0$, then $Z^N$ lies in the kernel of $\operatorname{gr}$, hence $Z^N$ is nilpotent, and so is $Z$. If $\operatorname{gr}(Z)$ is invertible and $Z'\in\mathcal Z$ is an element with $\operatorname{gr}(Z')=\operatorname{gr}(Z)^{-1}$, then $ZZ'-1$ lies in the kernel of $\operatorname{gr}$, thus $ZZ'$ is invertible, and so is $Z$. The converse statements hold because $\operatorname{gr}$ is a homomorphism of algebras and the inverse of an invertible element of $\mathcal Z$ lies in $\mathcal Z$. Since all fibres of $\operatorname{gr}$ have the same number of elements, Step~1 gives:
\[
\frac{\#\{Z\in\mathcal Z:\ Z\text{ nilpotent}\}}{\#\{Z\in\mathcal Z:\ Z\text{ invertible}\}}=\prod_{i=1}^k\frac{\nu_{m_i}}{|\GL_{m_i}(F)|}=\prod_{i=1}^k\frac{Q^{m_i}}{(Q;Q)_{m_i}} .
\]

\emph{Step 3. The count.} The number of matrices conjugate to $Y$ is $|\GL_d(F)|$ divided by the number of invertible elements of $\mathcal Z$, and the matrices $X$ with $(X,Y)\in\mathcal X_{d,k,\ell}$ are the nilpotent elements of $\mathcal Z$, where $\ell=\dim\ker Y=m_1+\dots+m_k$. Hence, by Step~2, the pairs $(X,Y')$ with $Y'$ conjugate to $Y$ contribute $\prod_iQ^{m_i}/(Q;Q)_{m_i}$ to the number $|\mathcal X_{d,k,\ell}|/|\GL_d(F)|$, where $d=\sum_iim_i$. The conjugacy classes of nilpotent matrices $Y$ with $Y^k=0$ and $\dim\ker Y=\ell$, of all sizes $d$, correspond to the vectors $(m_1,\dots,m_k)$ with sum $\ell$, and the first formula follows. The second formula follows from the first one and from the last identity of Step~1.
\end{proof}

\begin{theorem}\label{thm:finite}
Let $A=F[[x,y]]$ be a plane branch over $F=\F_q$, in any characteristic, let $\delta=\delta_A$ and $\fc=t^{2\delta}O$, and let $k\ge1$ be an integer with $y^k\in\fc$, for instance any $k\ge1$ with $k\val(y)\ge2\delta$. Then for all $n\ge1$ we have:
\[
e_A(n)\ \le\ \cN_A=\sum_{[M]:\ \fc M=0}\frac1{|\Aut_A(M)|}\ \le\ (Q;Q)_\infty^{-k},\qquad Q=q^{-1} .
\]
In particular $\cN_A<\infty$. Moreover, for every $s\ge0$ we have:
\begin{equation}\label{eq:soclebound}
\sum_{[M]:\ \fc M=0,\ s(M)\ge s}\frac1{|\Aut_A(M)|}\ \le\ Q^s\,(Q;Q)_\infty^{-k}\sum_{j\ge0}\binom{j+s+k-1}{k-1}Q^j .
\end{equation}
\end{theorem}
\begin{proof}
By Lemma~\ref{lem:planeG} the branch $A$ satisfies (G), so Proposition~\ref{prop:conductor} applies: its part (iii) is the formula for $\cN_A$, and the first inequality follows from parts (ii) and (iii), because the products in (ii) lie in $[0,1]$.

Let $M$ be an $A$-module of finite length with $\fc M=0$ and of dimension $d$ over $F$, and let $X$ and $Y$ be the endomorphisms of the $F$-vector space $M$ given by the multiplications by $x$ and $y$. They commute, $X$ is nilpotent, because $x^{2\delta}\in\fc$, and $Y^k=0$, because $y^k\in\fc$. Since $A$ is the closure of $F[x,y]$ and $\fc$ is open, we have $A=F[x,y]+\fc$. Hence the $A$-module structure of $M$ is determined by the pair $(X,Y)$, and an $F$-linear map between two such modules is $A$-linear if and only if it intertwines the corresponding pairs.

Thus, after the choice of a basis of $M$, the isomorphism class of $M$ becomes an orbit of $\GL_d(F)$ acting by simultaneous conjugation on the set of all pairs $(X,Y)$ of commuting nilpotent matrices in $\Mat_d(F)$ with $Y^k=0$. The stabilizer of a point of this orbit is isomorphic to $\Aut_A(M)$, and distinct isomorphism classes give distinct orbits. Counting the points of these orbits, we get for every $d\ge0$:
\[
\sum_{[M]:\ \fc M=0,\ \dim_FM=d}\frac1{|\Aut_A(M)|}\ \le\ \sum_{\ell\ge0}\frac{|\mathcal X_{d,k,\ell}|}{|\GL_d(F)|},
\]
and the sum over $d$ is at most $(Q;Q)_\infty^{-k}$ by Lemma~\ref{lem:ypower}.

Finally, the socle of $M$ is contained in $\ker Y$, because $y\in\fm$, so the modules with $s(M)\ge s$ give pairs in the sets $\mathcal X_{d,k,\ell}$ with $\ell\ge s$. By Lemma~\ref{lem:ypower} the left side of \eqref{eq:soclebound} is therefore at most:
\[
\sum_{\ell\ge s}Q^\ell\sum_{m_1+\dots+m_k=\ell}\ \prod_{i=1}^k\frac1{(Q;Q)_{m_i}}\ \le\ (Q;Q)_\infty^{-k}\sum_{\ell\ge s}\binom{\ell+k-1}{k-1}Q^\ell,
\]
because there are $\binom{\ell+k-1}{k-1}$ vectors $(m_1,\dots,m_k)$ with sum $\ell$, and \eqref{eq:soclebound} follows.
\end{proof}

\begin{corollary}[Finiteness characterizes planarity]\label{cor:planarity-finiteness}
Let $A\subseteq F[[t]]$ be a branch over $F=\F_q$ satisfying
the Gorenstein condition \textup{(G)}. Then
\[
 \cN_A<\infty\quad\Longleftrightarrow\quad A\text{ is a plane branch}.
\]
\end{corollary}
\begin{proof}
If $A$ is a plane branch, Theorem~\ref{thm:finite} gives
$\cN_A<\infty$. For the converse, first note that a
branch with $\dim_F\fm/\fm^2\le2$ is a plane branch. Indeed, choose
$x,y\in\fm$ whose classes span $\fm/\fm^2$. From
$\fm=Fx+Fy+\fm^2$ one obtains $A=F[x,y]+\fm^N$ for all $N$ by
induction. Since $\fm^N\subseteq t^NO$ and $A$ is closed in $O$,
the closure of $F[x,y]$ is $A$.

Suppose now that $A$ is not a plane branch. Then
$\dim_F\fm/\fm^2\ge3$, so the ring
$A/\fm^2=F\oplus\fm/\fm^2$ has
$B=F[u,v,w]/(u,v,w)^2$ as a quotient. A $B$-module structure on
$F^{2r}$ is a triple of matrices all of whose products vanish.
Consider the triples whose images span an $r$-dimensional subspace
$W$. Each matrix vanishes on $W$ and maps $F^{2r}/W$ into $W$;
the sum of the three maps is surjective from $(F^{2r}/W)^3$ onto
$W$. Conversely, each such surjection specifies a triple with these
properties, and its images determine $W$ uniquely. Thus the number
of these triples is
\[
 \qbinom{2r}r_q\prod_{i=0}^{r-1}(q^{3r}-q^i)
 =\frac{(Q;Q)_{3r}}{(Q;Q)_{2r}(Q;Q)_r^2}\,|\GL_{2r}(F)|
 \ge (Q;Q)_\infty\,|\GL_{2r}(F)|,
 \qquad Q=q^{-1}.
\]
Here we used
\[
 \qbinom{2r}r_q=q^{r^2}\frac{(Q;Q)_{2r}}{(Q;Q)_r^2},\qquad
 \prod_{i=0}^{r-1}(q^{3r}-q^i)
   =q^{3r^2}\frac{(Q;Q)_{3r}}{(Q;Q)_{2r}},\qquad
 |\GL_{2r}(F)|=q^{4r^2}(Q;Q)_{2r}.
\]
Restriction along $A\twoheadrightarrow B$ preserves the
automorphism groups of these modules. By the orbit--stabilizer
formula of Remark~\ref{rem:coh}(2), the $A$-modules of dimension
$2r$ therefore have total mass at least $(Q;Q)_\infty>0$ for every
$r\ge1$. The sum $\sum_{[M]}|\Aut_A(M)|^{-1}$ diverges.
Condition (G) and Corollary~\ref{cor:limit-exists}(ii) now imply
$\cN_A=\infty$.
\end{proof}

\begin{remark}\label{rem:finite}
(1) Together with Corollary~\ref{cor:limit-exists}(ii), Theorem~\ref{thm:finite} shows that for every plane branch $A$ the sum of the numbers $|\Aut_A(M)|^{-1}$ over all $A$-modules $M$ of finite length, which is the value at $1$ of the Coh zeta function of $A$ \cite{Huang,HJ}, is finite and equal to $(Q;Q)_\infty^{-1}$ times the sum over the modules annihilated by $\fc$. In particular the Coh zeta function $\sum_{[M]}|\Aut_A(M)|^{-1}t^{\dim_FM}$ of a plane branch has radius of convergence at least $1$. For $A=F[[t^a,t^b]]$ the finiteness was obtained in \cite{HJO} from the positive-definiteness theorem of Huang \cite{Huang2}, which is Theorem~\ref{thm:torus}(ii) below.

(2) For the cusp $A=F[[t^2,t^3]]$ we have $\fc=t^2O$ and $A/\fc=F$, and Proposition~\ref{prop:conductor}(iii) reads:
\[
\cN_A=\sum_{g\ge0}\frac1{|\GL_g(F)|}=\sum_{g\ge0}\frac{Q^{g^2}}{(Q;Q)_g},
\]
which is the sum side of the first Rogers--Ramanujan identity. For $A=F[[t^a,t^b]]$ with coprime $2\le a<b$ we can take $y=t^b$ and $k=a-1$ in Theorem~\ref{thm:finite}, because $(a-1)b\ge(a-1)(b-1)=2\delta$. The resulting bound $\cN_{a,b}\le(Q;Q)_\infty^{-(a-1)}$ agrees with the product formula of Theorem~\ref{thm:torus}(iii), in which all exponents $\rho_{a,b}(j)$ are at most $a-1$.

(3) An example with $\cN_A=\infty$ is the closure $A=F+Ft^4+Ft^5+Ft^6+t^8O$ of $F[t^4,t^5,t^6]$, which satisfies (G) by Lemma~\ref{lem:symmetric}. Here $\fc=t^8O=\fm^2$, the classes of $t,t^2,t^3,t^7$ form a basis of $O/A$, and we have $t^{7-j}\cdot t^j=t^7$ for $1\le j\le3$, while all other products of $t^4,t^5,t^6$ with these classes vanish. Hence a subspace $M'$ of $(O/A)^n=F^n\otimes(O/A)$ is a submodule if and only if its image $P$ in $F^n\otimes\langle t,t^2,t^3\rangle$ is contained in $S\otimes\langle t,t^2,t^3\rangle$, where $S\subseteq F^n$ is the subspace with $M'\cap(F^n\otimes t^7)=S\otimes t^7$. For given $S$ and $P$, of dimensions $s$ and $p$, these submodules correspond to the linear maps $P\to(F^n/S)\otimes t^7$, thus there are $q^{p(n-s)}$ of them, and Proposition~\ref{prop:finitesum} gives:
\[
e_A(n)=\sum_{s=0}^n\sum_{p=0}^{3s}\qbinom ns_q\qbinom{3s}p_q\,q^{p(n-s)-n(s+p)}=\sum_{s=0}^n\qbinom ns_Q\sum_{p=0}^{3s}\qbinom{3s}p_Q\,Q^{(p-s)^2}\ \ge\ n+1,
\]
where the last inequality comes from the terms with $p=s$.
\end{remark}

\section{Quadratic transforms and framed curves}\label{sec:qt}

The induction uses a framed pair consisting of a ring with two
specified generators: the first determines the period of a lattice
chain, and the second its stability condition. The Euclidean
construction of Proposition~\ref{prop:transform} performs a block of
quadratic transforms and determines the resulting framing and
multiplicity sequence. Corollary~\ref{cor:cabling-data} applies it to
the two branches that enter the cabling formula.

In Section~\ref{sec:hall} we need to treat several branches simultaneously, and for this reason we work in the semilocal ring $\mathbf O=O^r$, $r\ge1$, whose total fraction ring is $\mathbf K=\cK^r$ and whose radical is $\rad\mathbf O=(tO)^r$. For $z\in\mathbf K^\times$ we denote by $\val_i(z)$ the valuation of its $i$-th component, and we put $|z|=\sum_i\val_i(z)$. If $x\in\rad\mathbf O$ is a nonzerodivisor, then $F[[x]]\subseteq\mathbf O$ is a power series ring and $F((x))\subseteq\mathbf K$ is a field. For $x,y\in\rad\mathbf O$ the closure of $F[x,y]$ in $\mathbf O$ is the image of $F[[X,Y]]$, and we denote it by $F[[x,y]]$.

\subsection{Two generators}
\begin{lemma}\label{lem:gen}
Let $x,y\in\rad\mathbf O$, and assume that $x$ is a nonzerodivisor.
\begin{enumerate}[label=\textup{(\alph*)}]
\item The ring $\mathbf O$ is a free $F[[x]]$-module of rank $|x|$.
\item We have $F[[x,y]]=F[[x]][y]$, and this ring is a free $F[[x]]$-module with basis $1,y,\dots,y^{d-1}$, where $d=\dim_{F((x))}F((x))[y]$.
\item The following conditions are equivalent: $F((x))[y]=\mathbf K$, $\dim_F\mathbf O/F[[x,y]]<\infty$, and $d=|x|$.
\item If these conditions hold, then $y^{|x|}\in x\,F[[x,y]]$.
\end{enumerate}
\end{lemma}
\begin{proof}
(a) For every $i$ the ring $O$ is free over $F[[x_i]]$ with basis $1,t,\dots,t^{\val(x_i)-1}$, where $x_i$ is the $i$-th component of $x$.

(b) Let $P\in F[[x]][Y]$ be the characteristic polynomial of multiplication by $y$ on the free module $\mathbf O$. Since $P$ is monic and $P(y)=0$, the ring $F[[x]][y]$ is a finitely generated $F[[x]]$-module. Thus it is compact, hence closed, and it is equal to $F[[x,y]]$. The minimal polynomial $\chi$ of $y$ over the field $F((x))$ divides $P$, and since $F[[x]]$ is integrally closed, we have $\chi\in F[[x]][Y]$. Thus $F[[x]][y]=\sum_{i<d}F[[x]]y^i$ with $d=\deg\chi$, and the elements $1,\dots,y^{d-1}$ are linearly independent over $F((x))$.

(c) By (a) and (b), $F[[x,y]]$ is a free submodule of rank $d$ of a free module of rank $|x|$, and $F((x))[y]=F((x))\otimes_{F[[x]]}F[[x,y]]$.

(d) If $d=|x|$, then $\chi=P$, and modulo $x$ the polynomial $P$ is the characteristic polynomial of $y$ on $\mathbf O/x\mathbf O$, where $y$ is nilpotent, thus $P\equiv Y^{|x|}$ modulo $x$, and $P(y)=0$ gives $y^{|x|}\in x\,F[[x]][y]=x\,F[[x,y]]$.
\end{proof}

\begin{lemma}\label{lem:planeG}
Every plane branch satisfies condition \textup{(G)} of Section~\ref{sec:norm}.
\end{lemma}
\begin{proof}
Let $r=1$, and let $A=F[[x,y]]$ be a plane branch, where we may assume $x\ne0$ after exchanging $x$ and $y$, because $A\ne F$. Let us put $R=F[[x]]$, $K_x=F((x))$ and $\mu=\val x$. By Lemma~\ref{lem:gen} we have $A=\bigoplus_{i<\mu}Ry^i$, and $1,y,\dots,y^{\mu-1}$ is a $K_x$-basis of $\cK$. Let $\lambda:\cK\to K_x$ be the $K_x$-linear form with $\lambda(y^i)=0$ for $i<\mu-1$ and $\lambda(y^{\mu-1})=1$. For $0\le i,j<\mu$ the element $\lambda(y^iy^j)$ lies in $R$, because $y^{i+j}\in A$, and it is $0$ if $i+j<\mu-1$ and $1$ if $i+j=\mu-1$. Hence the matrix $(\lambda(y^iy^j))_{i,j}$ is invertible over $R$, and therefore:
\begin{equation}\label{eq:euler}
\{\xi\in\cK:\ \lambda(\xi a)\in R\ \text{for all }a\in A\}=A .
\end{equation}
Next we consider the $F$-linear form $\varrho(\xi)=\Res_{x=0}(\lambda(\xi)\,dx)$ on $\cK$. Since $A$ is an $R$-module and an element $g\in K_x$ with $\Res(gx^j\,dx)=0$ for all $j\ge0$ lies in $R$, the identity \eqref{eq:euler} says that $\{\xi:\varrho(\xi a)=0\text{ for all }a\in A\}=A$. In particular, $\varrho\ne0$ and $\varrho$ vanishes on $A\supseteq t^cO$ for some $c$, and for $\kappa=\sum_{j<c}\varrho(t^j)\,t^{-j-1}\in\cK$ we then have $\Res_{t=0}(\kappa t^i\,dt)=\varrho(t^i)$ for all $i\in\Z$. Both sides of $\varrho(\xi)=\Res_{t=0}(\kappa\xi\,dt)$ depend only on the coefficients of $\xi$ in degrees below $c$, hence this equality holds for all $\xi\in\cK$. Thus $\kappa\ne0$, and $\xi\in A$ if and only if $\kappa\xi\in A^\perp$, that is, $A^\perp=\kappa A$.
\end{proof}

\begin{lemma}\label{lem:qtlength}
Let $x,y$ be as in Lemma~\ref{lem:gen} with $F((x))[y]=\mathbf K$, assume that $y/x\in c+\rad\mathbf O$ for some $c\in F$, and let us put $B=F[[x,y]]$ and $B^\sharp=F[[x,\,y/x-c]]$. Then $B\subseteq B^\sharp$, $F((x))[y/x-c]=\mathbf K$, and
\[
\dim_FB^\sharp/B=\binom{|x|}2 .
\]
\end{lemma}
\begin{proof}
For $z=y/x$ and $\mu=|x|$ we have $y=xz\in B^\sharp$ and $F((x))[z-c]=F((x))[y]$, and Lemma~\ref{lem:gen} gives:
\[
B^\sharp=\bigoplus_{i<\mu}F[[x]]z^i,\qquad B=\bigoplus_{i<\mu}F[[x]]y^i=\bigoplus_{i<\mu}x^iF[[x]]z^i,
\]
hence $B^\sharp/B\cong\bigoplus_{i<\mu}F[[x]]/(x^i)$, and its dimension is $\sum_{i<\mu}i$.
\end{proof}

\subsection{Branches and multiplicity sequences}
Let $r=1$, and for a branch $A$ with maximal ideal $\fm_A=A\cap tO$ let us put $\mult(A)=\min\{\val(g):g\in\fm_A\smallsetminus0\}$. By the geometric series, an element of $A$ of valuation $0$ is a unit of $A$.

\begin{lemma}\label{lem:qt-intrinsic}
Let $A$ be a plane branch and let $g\in\fm_A$ with $\val(g)=\mult(A)$. Then the ring $A':=A[\fm_A/g]$ does not depend on $g$, and it is a plane branch containing $A$. Moreover, if $A=F[[x,y]]$ with $\val x\le\val y$ and $c$ is the constant term of $y/x$, then $\mult(A)=\val x$ and $A'=F[[x,\,y/x-c]]$.
For every finite field extension $E/F$, the branch
$A_E=A\otimes_F E\subseteq E[[t]]$ has the same multiplicity
as $A$, and $(A_E)'=A'\otimes_F E$. In particular, extension
of constants preserves the entire multiplicity sequence.
\end{lemma}
\begin{proof}
Every element of $\fm_A$ is a power series in $x,y$ without constant term, hence $\mult(A)=\val x$, and $A[\fm_A/x]=A[y/x]$ because $\fm_A=xA+yA$. The element $y/x\in O$ is integral over $A$, so $A[y/x]$ is a finitely generated $A$-module and therefore closed in $O$. It contains $x$ and $y/x-c$, thus also $F[[x,y/x-c]]$, and conversely, the latter ring contains $x$, $y/x$ and $y=x\cdot y/x$. Finally, if $g\in\fm_A$ is another element of valuation $\val x$, then $g/x\in B=A[\fm_A/x]$ has valuation $0$, so it is a unit of $B$ and $\fm_A/g=(x/g)\,\fm_A/x\subseteq B$, and by symmetry we get $A[\fm_A/g]=B$.

To check extension of constants, choose an $F$-basis
$b_1,\ldots,b_s$ of $E$. Every element of $A_E$ is uniquely
$\sum_i b_i a_i$ with $a_i\in A$. Independence of the leading
coefficients gives
\[
 \val_t\Bigl(\sum_i b_i a_i\Bigr)=\min_i\val_t(a_i)
 \qquad\text{for a nonzero sum}.
\]
Thus the valuation filtration, including the maximal ideal, is
obtained by scalar extension. A minimal-valuation element
$g\in\fm_A$ remains minimal over $E$, and adjoining the
generators of $\fm_A/g$ commutes with scalar extension:
\[
 A_E[\fm_{A_E}/g]
   =\bigl(A[\fm_A/g]\bigr)\otimes_F E.
\]
Ordinary tensor products suffice because $E/F$ is finite.
Iteration proves the assertion about multiplicity sequences.
\end{proof}
\noindent The ring $A'$ is the local ring of the strict transform after one blow-up, and we call it the \emph{quadratic transform} of $A$. The iterated quadratic transforms of $A$ are defined by $A^{(0)}=A$ and $A^{(j+1)}=(A^{(j)})'$.

\begin{lemma}\label{lem:noether}
Let $A$ be a plane branch. If $\mult(A)=1$, then $A=O$. In general we have $\delta_A-\delta_{A'}=\binom{\mult A}2$, hence $A^{(j)}=O$ for $j\gg0$ and $\delta_A=\sum_{j\ge0}\binom{\mult A^{(j)}}2$.
\end{lemma}
\begin{proof}
The first statement holds because $A\supseteq F[[g]]=O$ for every element $g\in A$ of valuation $1$, the second statement is Lemma~\ref{lem:qtlength} with $r=1$, and the third statement follows, because $\delta$ strictly decreases as long as the multiplicity is at least $2$.
\end{proof}
\noindent The \emph{multiplicity sequence} $\ms(A)=(\mult A^{(j)})_{j\ge0}$ of $A$ is eventually equal to $1$ by Lemma~\ref{lem:noether}, and we identify sequences which differ by trailing $1$'s.

For coprime integers $a,b\ge1$ the \emph{Euclidean sequence} $\eu(a,b)$ is defined recursively as follows:
\[
\eu(1,1)=(1),\qquad\eu(a,b)=(a,\eu(a,b-a))\ \text{ if }a<b,\qquad\eu(a,b)=(b,\eu(a-b,b))\ \text{ if }a>b,
\]
The final entry in $\eu(1,1)=(1)$ is retained when computing its
length. In Proposition~\ref{prop:transform}, this entry records the
last quadratic transform, when the two current generators have equal
valuation and one subtracts the common residue from their quotient.
The convention of omitting trailing $1$'s from a multiplicity
sequence therefore does not apply to the length of a Euclidean block.

We denote by $\ell(a,b)$ its length. By induction we have:
\begin{equation}\label{eq:eusums}
\eu(a,b)=\eu(b,a),\qquad\sum_{m\in\eu(a,b)}m=a+b-1,\qquad\sum_{m\in\eu(a,b)}m^2=ab,
\end{equation}
hence $\sum_{m\in\eu(a,b)}\binom m2=\frac{(a-1)(b-1)}2$. If $2\le a<b$, then $\eu(a,b)$ begins with $\lfloor b/a\rfloor$ entries equal to $a$, followed by the entry $b-\lfloor b/a\rfloor a<a$, and $\eu(a,b+ja)=(a^j,\eu(a,b))$ for $j\ge0$.

\begin{lemma}\label{lem:normalgen}
If $A$ is a plane branch with $m=\mult(A)\ge2$, then $A=F[[x,y]]$ for some $x,y$ with $\val x=m$ and $m\nmid\val y$.
\end{lemma}
\begin{proof}
Let us write $A=F[[x,y]]$ with $\val x=m\le\val y$. If $\val y=km$, then there is $c\in F^\times$ with $\val(y-cx^k)>\val y$ and $F[[x,y-cx^k]]=A$, and we repeat this step with $y-cx^k$ in place of $y$. If the process does not stop, then $y$ is a limit of elements of $F[[x]]$, thus $A=F[[x]]$ and $\Frac A=F((x))$, which is impossible, since $F((x))\ne\cK$ for $m\ge2$.
\end{proof}

\subsection{Framed pairs and their transforms}
\begin{definition}\label{def:framed}
A \emph{framed pair} in $\mathbf O=O^r$ is a pair $(u,w)$ of elements of $\rad\mathbf O$ such that $u$ is a nonzerodivisor and $F((u))[w]=\mathbf K$, and its \emph{ring} is $T=F[[u,w]]$. For coprime integers $e,d_*\ge1$ the pair is \emph{balanced of type} $(e,d_*)$ if there are integers $g_1,\dots,g_r\ge1$ and $c\in F^\times$ such that:
\[
\val_i(u)=g_ie,\qquad\val_i(w)=g_id_*\quad(1\le i\le r),\qquad \theta:=\frac{w^e}{u^{d_*}}\in c\,(1+\rad\mathbf O),
\]
and in this case its \emph{weight} is $g=\sum_ig_i$. For $r=1$ a framed pair with $w\ne0$ is called a \emph{framed branch}, written as $(T;u,w)$, and every framed branch is balanced, with $g=\gcd(\val u,\val w)$.
\end{definition}

\begin{proposition}\label{prop:transform}
Let $(u,w)$ be a balanced framed pair of type $(e,d_*)$ and weight $g$ in $\mathbf O$, with ring $T$. Let us choose integers $A_0,B_0$ with $A_0e+B_0d_*=1$ and put:
\[
\varpi=u^{A_0}w^{B_0},\qquad\omega=c^{-1}\theta-1,\qquad T_1=F[[\varpi,\omega]].
\]
\begin{enumerate}[label=\textup{(\roman*)}]
\item We have $\val_i(\varpi)=g_i$ and the following two identities:
\[
u=\varpi^e\theta^{-B_0},\qquad w=\varpi^{d_*}\theta^{A_0}.
\]
The pair $(\varpi,\omega)$ is a framed pair, and the ring $T_1$ does not depend on the choice of $(A_0,B_0)$, since a different choice multiplies $\varpi$ by a unit of $F[[\omega]]$, while $F[[\varpi]]$ changes accordingly.
\item We have $T\subseteq T_1$ and
\[
\dim_FT_1/T=\sum_{m\in\eu(e,d_*)}\binom{gm}2=g^2\,\frac{(e-1)(d_*-1)}2+\binom g2(e+d_*-1).
\]
\item If $r=1$, then $T_1=T^{(\ell(e,d_*))}$ and the multiplicities of $T^{(j)}$ for $0\le j<\ell(e,d_*)$ are the entries of $g\cdot\eu(e,d_*)$, hence $\ms(T)=(g\cdot\eu(e,d_*),\ \ms(T_1))$.
\end{enumerate}
\end{proposition}
\begin{proof}
(i) Since all components of $u$ and $w$ are nonzero, we have $u,w\in\mathbf K^\times$, and $\varpi,\theta\in\mathbf K^\times$ are defined. The valuations and the two identities follow directly from $A_0e+B_0d_*=1$, and in particular $\varpi\in\rad\mathbf O$ is a nonzerodivisor and $\omega\in\rad\mathbf O$. By Lemma~\ref{lem:gen}(b), the closed ring $F[[\varpi,\omega]]$ equals $F[[\varpi]][\omega]$, and it contains the units $\theta^{\pm1}$ of $F[[\omega]]$. By the two identities it contains $u$ and $w$, hence also $F[[u]]$, and since $u^{-1}=\varpi^{-e}\theta^{B_0}\in F((\varpi))[\omega]$, we get $F((\varpi))[\omega]\supseteq F((u))[w]=\mathbf K$. Finally, another choice of $(A_0,B_0)$ replaces $\varpi$ by $\varpi\theta^s$, $s\in\Z$.

(ii) and (iii) We run the subtractive Euclidean algorithm on Laurent monomials, where a \emph{state} is a pair $(x,y)=(u^{p}w^{q},u^{p'}w^{q'})$ with $\left(\begin{smallmatrix}p&q\\p'&q'\end{smallmatrix}\right)\in\GL_2(\Z)$ and $a:=pe+qd_*\ge1$, $b:=p'e+q'd_*\ge1$. For a state we have $\val_i(x)=g_ia$, $\val_i(y)=g_ib$, $|x|=ga$ and $\gcd(a,b)=1$, and the initial state is $(u,w)$, with $(a,b)=(e,d_*)$. Suppose inductively that $B=F[[x,y]]$ has finite colength in $\mathbf O$, in which case Lemma~\ref{lem:gen}(c) applies to $(x,y)$ and to $(y,x)$.

If $a<b$, then $y/x\in\rad\mathbf O$, and Lemma~\ref{lem:qtlength} with $c=0$ gives the ring $B^\sharp=F[[x,y/x]]$ of the state $(x,y/x)$, in which $(a,b)$ is replaced by $(a,b-a)$, with $\dim B^\sharp/B=\binom{ga}2$. If $a>b$, we exchange the roles of $x$ and $y$: the new state is $(x/y,y)$ and the index is $\binom{gb}2$. For $r=1$ this step is the quadratic transform, by Lemma~\ref{lem:qt-intrinsic}, and the multiplicity of $B$ is $g\min(a,b)$.

Finally, let $a=b=1$. Then $(p'-p,q'-q)$ is orthogonal to $(e,d_*)$, so it equals $k(-d_*,e)$ for some $k\in\Z$, and since $k$ is the determinant of the exponent matrix, we have $k=\pm1$ and $y/x=\theta^{\pm1}\in c^{\pm1}(1+\rad\mathbf O)$. Moreover, $(p,q)=(A_0,B_0)+s(-d_*,e)$ for some $s$, so $x=\varpi\theta^s$, and Lemma~\ref{lem:qtlength} gives $B^\sharp=F[[\varpi\theta^s,\theta^{\pm1}-c^{\pm1}]]$ with $\dim B^\sharp/B=\binom g2$. Now we have:
\[
\theta-c=c\,\omega,\qquad\theta^{-1}-c^{-1}=-c^{-1}\omega(1+\omega)^{-1},
\]
so each of them generates the same closed subalgebra as $\omega$. Since this subalgebra contains the unit $\theta^s$, we get $B^\sharp=F[[\varpi,\omega]]=T_1$.

The pairs $(a,b)$ that we visit are the states of the Euclidean algorithm for $(e,d_*)$, which ends at $(1,1)$ after $\ell(e,d_*)-1$ steps, and the recorded numbers $\min(a,b)$ are the entries of $\eu(e,d_*)$. Summing the indices and using $\binom{gm}2=g^2\binom m2+\binom g2m$ and \eqref{eq:eusums}, we obtain (ii), and for $r=1$ we obtain (iii).
\end{proof}

\begin{lemma}\label{lem:invtransform}
Let $(\varpi,\omega)$ be a framed pair in $\mathbf O$ with $\omega\in\rad\mathbf O$, let $e,d_*\ge1$ be coprime, let $A_0e+B_0d_*=1$, let $c\in F^\times$, and let $\theta=c(1+\omega)$. Then the following pair is a balanced framed pair of type $(e,d_*)$ with $g_i=\val_i(\varpi)$:
\[
u=\varpi^e\theta^{-B_0},\qquad w=\varpi^{d_*}\theta^{A_0},
\]
and its transform with respect to $(A_0,B_0)$ is $(\varpi,\omega)$.
\end{lemma}
\begin{proof}
The elements $u,w\in\rad\mathbf O$ are nonzerodivisors with the stated valuations, and we have $w^e/u^{d_*}=\theta$ and $u^{A_0}w^{B_0}=\varpi$. The algebra $F((u))[w]$ is a finite-dimensional $F((u))$-subspace of $\mathbf K$, thus it is closed, and multiplication by the nonzerodivisor $w$ is injective, hence bijective, on it. Therefore it contains $w^{-1}$, hence also $\varpi^{\pm1}$, $\theta$, $\omega$ and $F((\varpi))$, and thus $F((u))[w]\supseteq F((\varpi))[\omega]=\mathbf K$.
\end{proof}

\begin{example}\label{ex:torus}
Let $r=1$, $u=t^a$ and $w=t^b$ with $\gcd(a,b)=1$. Then $g=1$, $\varpi=t$, $\theta=1$, $\omega=0$ and $T_1=O$, and we obtain:
\[
\delta_{F[[t^a,t^b]]}=(a-1)(b-1)/2=:\delta_{a,b},\qquad\ms(F[[t^a,t^b]])=\eu(a,b).
\]
\end{example}

\subsection{The cabling data of a plane branch}\label{sub:cabling-data}
\begin{corollary}\label{cor:cabling-data}
Let $A=F[[x,y]]$ be a plane branch with $\val x=Da$, $\val y=Db$ and $\gcd(a,b)=1$, and let $c,\varpi,\omega,T,S$ be as in Theorem~\ref{thm:introC}, with $h=a+b-1$, namely:
\[
\varpi=x^{A_0}y^{B_0},\qquad 1+\omega=c^{-1}y^a/x^b,\qquad T=F[[\varpi,\omega]],\qquad S=F[[\varpi,\varpi^h\omega]].
\]
\begin{enumerate}[label=\textup{(\roman*)}]
\item $(T;\varpi,\omega)$ is a framed branch (or $T=O$ and $\omega=0$) with $\val\varpi=D$. It is the transform of the framed branch $(A;x,y)$, of type $(a,b)$ and weight $D$, and also the transform of the framed branch $(S;\varpi,\varpi^h(1+\omega))$, of type $(1,h)$ and weight $D$.
\item $T=A^{(\ell(a,b))}=S^{(h)}$, $\ms(A)=(D\cdot\eu(a,b),\ms(T))$ and $\ms(S)=(D^h,\ms(T))$, where $D^h$ stands for $h$ entries equal to $D$.
\item For the $\delta$-invariants we have:
\[
\delta_A-\delta_T=D^2\delta_{a,b}+\binom D2h,\qquad\delta_S-\delta_T=\binom D2h,\qquad\delta_A-\delta_S=D^2\delta_{a,b}.
\]
\item Every element of $E_T(O^n)$ is a free $F[[\varpi]]$-module of rank $Dn$, and $\omega^D\Lambda\subseteq\varpi\Lambda$ for every $T$-lattice $\Lambda$.
\end{enumerate}
\end{corollary}
\begin{proof}
(i)--(iii) We apply Proposition~\ref{prop:transform} to $(x,y)$, with $(e,d_*,g)=(a,b,D)$, and for $S$ we observe that $F[[\varpi,\varpi^h(1+\omega)]]=F[[\varpi,\varpi^h\omega]]$ and that $(\varpi,\varpi^h(1+\omega))$ is a framed pair, because $F((\varpi))[\varpi^h(1+\omega)]=F((\varpi))[\omega]$, and that with $(e,d_*)=(1,h)$ and $(A_0,B_0)=(1,0)$ its transform is again $(\varpi,\omega)$, since $\theta=1+\omega$. If $\omega=0$, then $D=1$, and all three rings are equal to $O$ after the indicated transforms.

(iv) This follows from Lemma~\ref{lem:gen}(a),(d), applied to the pair $(\varpi,\omega)$ in place of $(x,y)$.
\end{proof}
\noindent The branch $S$ depends on the choice of $(A_0,B_0)$ only through the framing $\varpi$, and we shall see that its lattice counts depend only on $\ms(S)$.


\section{The Iwahori--Hecke algebra and compact-mod-center elements}\label{sec:hecke}

For invariant distributions supported on compact-mod-center elements,
evaluations on a fixed-valuation part of the Iwahori--Hecke algebra
are determined by finitely many elements indexed by partitions
(Theorem~\ref{thm:H}). We prove this by reducing to terminal elements
and then to finite Hecke algebras. Section~\ref{sec:chains} proves
that the resulting dimension bound is sharp and gives a basis of
test functions that count lattice chains.

In this section $\K$ is a non-archimedean local field with ring of integers $\fo$, uniformizer $\pi$ and residue field $\F_q$, and we fix a Haar measure on the group $G=\GL_m(\K)$.

\subsection{Notation}
Let us denote by $b_1,\dots,b_m$ the standard basis of $\K^m$, extended to all $i\in\Z$ by $b_{i+m}=\pi b_i$, and set:
\[
\cL_i=\sum_{j>i}\fo\,b_j\qquad(i\in\Z).
\]
Then $(\cL_i)_{i\in\Z}$ is a complete periodic chain of lattices with $\cL_i\supsetneq\cL_{i+1}$ and $\cL_{i+m}=\pi\cL_i$, and its stabilizer $I$ is an Iwahori subgroup.

An \emph{affine permutation} is a bijection $w$ of $\Z$ with $w(i+m)=w(i)+m$, and it acts on $\K^m$ by $b_i\mapsto b_{w(i)}$. The affine permutations form the group $\tW\cong\Z^m\rtimes S_m$ of monomial matrices with entries powers of $\pi$. Let $\tau(i)=i+1$ and $\val(w)=\frac1m\sum_{i=1}^m(w(i)-i)=\val\det w$, and set $W_a=\ker\val$ and $\Omega=\langle\tau\rangle$.
For $m=1$, the set of simple reflections is empty and $W_a=\{1\}$.
For $m\ge2$ and $0\le i<m$, the simple reflection $s_i$ exchanges
$i+mj$ and $i+1+mj$ for all $j$ (indices modulo $m$), and
$W_a=\langle s_0,\dots,s_{m-1}\rangle$.
We define $\ell(w)$ as the number of inversions of $w$:
\[
\ell(w)=\#\{(i,j):1\le i\le m,\ i<j,\ w(i)>w(j)\}.
\]
Then $\tW=W_a\rtimes\Omega$, $\tau s_i\tau^{-1}=s_{i+1}$, $\ell(\tau w)=\ell(w\tau)=\ell(w)$, and $\ell(w)$ is the smallest number $r$ such that $w=s_{i_1}\cdots s_{i_r}\tau^k$ for some $i_1,\dots,i_r$ and $k$. Lemma~\ref{lem:IMlength} in Section~\ref{app:IM} proves these length facts in our conventions. We use the Bruhat decomposition and the Iwahori--Matsumoto relations \cite{IM}: for a simple reflection $s$ we have
\begin{gather}
G=\bigsqcup_{w\in\tW}IwI,\qquad I\tau=\tau I,\notag\\
IsI\cdot IwI=IswI\ \text{ if }\ell(sw)>\ell(w),\qquad IsI\cdot IwI=IwI\sqcup IswI\ \text{ if }\ell(sw)<\ell(w),\label{eq:IM}
\end{gather}
and the mirror statements hold for $IwI\cdot IsI$. They are proved in Section~\ref{app:IM}, in the present conventions (Lemma~\ref{lem:IMbruhat} and Proposition~\ref{prop:IM}).

Let us denote by $\cH=C_c(I\backslash G/I)$ the Iwahori--Hecke algebra of complex-valued functions, with convolution defined by the chosen Haar measure, and put $T_w=\vol(I)^{-1}\mathbf 1_{IwI}$. Thus $T_1=\vol(I)^{-1}\mathbf 1_I$ is the unit. The family $(T_w)_{w\in\tW}$ is a basis of $\cH$, and by Proposition~\ref{prop:IM} we have $T_wT_{w'}=T_{ww'}$ if $\ell(ww')=\ell(w)+\ell(w')$, and $T_s^2=(q-1)T_s+q$ for a simple reflection $s$. In particular, every $T_w$ is invertible (Proposition~\ref{prop:IM}(ii)). The algebra $\cH$ is graded by $\val$, and we denote by $\cH^{\val=k}$ the span of the $T_w$ with $\val(w)=k$. Finally, we write $\equiv$ for congruence modulo $[\cH,\cH]$.

An element $g\in G$ is \emph{compact-mod-center} if all its eigenvalues in an algebraic closure of $\K$ have the same valuation. By Lemma~\ref{lem:Gc} below, this holds if and only if the image of $g$ in $G/\K^\times$ lies in a compact subgroup. The set $G_c$ of these elements is conjugation-invariant and closed in $G$: indeed, the condition says that the Newton polygon of the characteristic polynomial $\sum_ic_iX^{m-i}$ is a segment, that is, $m\val(c_i)\ge i\val(c_m)$ for all $i$, while $\val(c_m)=\val\det g$ is locally constant. A distribution $D$ is a linear form on $C_c^\infty(G)$, and it is \emph{invariant} if $D(f^g)=D(f)$, where $f^g(y)=f(gyg^{-1})$.

\begin{lemma}\label{lem:trace}
An invariant distribution $D$ satisfies $D(f_1*f_2)=D(f_2*f_1)$ for all $f_1,f_2\in C_c^\infty(G)$, and in particular $D|_{\cH}$ vanishes on $[\cH,\cH]$.
\end{lemma}
\begin{proof}
Let us write both convolutions as integrals; the second formula below holds because $G$ is unimodular:
\[
(f_1*f_2)(y)=\int_Gf_1(h)f_2(h^{-1}y)\,dh,\qquad(f_2*f_1)(y)=\int_Gf_1(h)f_2(yh^{-1})\,dh.
\]
For fixed $h$ the function $R_hf_2\colon y\mapsto f_2(yh^{-1})$ is the conjugate by $h$ of $L_hf_2\colon y\mapsto f_2(h^{-1}y)$, hence $D(R_hf_2)=D(L_hf_2)$. Since $f_2$ is invariant under a compact open subgroup on both sides, the functions $L_hf_2$ and $R_hf_2$ are locally constant in $h$, and both convolutions are finite linear combinations $f_1*f_2=\sum_jc_jL_{h_j}f_2$ and $f_2*f_1=\sum_jc'_jR_{h'_j}f_2$. Thus we get:
\[
\sum_jc_jD(L_{h_j}f_2)=\int_Gf_1(h)D(L_hf_2)\,dh=\int_Gf_1(h)D(R_hf_2)\,dh=\sum_jc'_jD(R_{h'_j}f_2).
\]
\end{proof}

For $k\in\Z$ let us denote by $J_k$ the space of linear forms on $\cH^{\val=k}$ obtained by restricting the invariant distributions supported on $G_c$, and put $n'=\gcd(m,k)$. For a partition $\rho=(\rho_1\ge\dots\ge\rho_l)$ of $n'$ we split $\{1,\dots,n'\}$ into the consecutive blocks $\{p+1,\dots,p+c\}$ of sizes $\rho_1,\dots,\rho_l$, denote by $u_\rho\in S_{n'}:=\langle s_1,\dots,s_{n'-1}\rangle\subseteq W_a$ the product of the Coxeter elements $s_{p+c-1}\cdots s_{p+2}s_{p+1}$ of these blocks, and set:
\[
w_\rho=u_\rho\,\tau^k .
\]

\begin{theorem}\label{thm:H}
Let $k\in\Z$ and $n'=\gcd(m,k)$.
\begin{enumerate}[label=\textup{(\roman*)}]
\item $\dim J_k\le p(n')$.
\item For every $w\in\tW$ with $\val(w)=k$ and every $\rho\vdash n'$ there is a polynomial $a_{w,\rho}\in\Z[X]$, depending only on $(m,w,\rho)$ and not on the field $\K$, such that
\[
D(T_w)=\sum_{\rho\vdash n'}a_{w,\rho}(q)\,D(T_{w_\rho})\qquad\text{for all }D\in J_k .
\]
\end{enumerate}
\end{theorem}
\noindent For $m=1$ the statement is trivial, so we assume $m\ge2$. The proof occupies the rest of this section: it combines the reduction method of Geck--Pfeiffer and He--Nie with two elementary facts about the group $G$ (Lemmas~\ref{lem:Gc} and~\ref{lem:cert}) and three computations in Hecke algebras (Lemmas~\ref{lem:twist}--\ref{lem:swap}).

\subsection{Conjugacy classes and compact-mod-center elements}
Let us write $w\in\tW$ as $w=t^\lambda\sigma$ with $\sigma\in S_m$ and $\lambda\in\Z^m$, that is, $w(i)=\sigma(i)+m\lambda_{\sigma(i)}$ for $1\le i\le m$. Then $\val(w)=\sum_i\lambda_i$, and if $L$ is a multiple of the order of $\sigma$, then $w^L=t^{\lambda^{(L)}}$ is a translation. The \emph{Newton point} of $w$ is $\nu_w=\lambda^{(L)}/L\in\Q^m$, and it is \emph{central} if all its coordinates are equal, in which case they are necessarily equal to $\val(w)/m$.
In Section~\ref{app:HN3} we denote this same vector by
$\nu_w^{\mathrm{full}}$ and write $\nu_w^{\mathrm{red}}$ for its
image modulo the diagonal line, represented by a vector of coordinate
sum zero.

\begin{lemma}\label{lem:classes}
The conjugacy classes of $\tW$ correspond to the multisets $\{(\ell_c,s_c)\}_c$, where $c$ runs over the cycles of $\sigma$, $\ell_c$ is the length of $c$ and $s_c=\sum_{i\in c}\lambda_i$. The coordinate of $\nu_w$ at $i\in c$ is $s_c/\ell_c$, and the classes with central Newton point and valuation $k$ correspond to the partitions of $\gcd(m,k)$.
\end{lemma}
\begin{proof}
Conjugation by $t^\mu$ replaces $\lambda$ by $\lambda+\mu-\sigma\mu$, which does not change the $s_c$, and conjugation by an element of $S_m$ permutes the cycles. Conversely, if $\lambda,\lambda'$ have the same cycle sums for $\sigma$, then $\lambda-\lambda'$ has sum zero on each cycle, hence lies in $(1-\sigma)\Z^m$. Next, since $w^L=t^{\sum_{j<L}\sigma^j\lambda}$, the coordinate of $\lambda^{(L)}$ at $i\in c$ is $(L/\ell_c)s_c$. With $e=m/\gcd(m,k)$ and $d_*=k/\gcd(m,k)$, the Newton point is central if and only if $s_c/\ell_c=d_*/e$ for all $c$, that is, $e\mid\ell_c$ and $s_c=d_*\ell_c/e$. Such classes correspond to the partitions $(\ell_c/e)_c$ of $\gcd(m,k)$.
\end{proof}

\begin{lemma}\label{lem:Gc}
An element $g\in G$ is compact-mod-center if and only if there are a periodic lattice chain $(M_j)_{j\in\Z}$ and $c\in\Z$ with $gM_j=M_{j+c}$ for all $j$. Here periodic means that $M_j\supseteq M_{j+1}$ and $\pi M_j=M_{j+r}$ for some $r\ge1$. In that case $g^r\in\pi^c\,\cC$ for a compact open subgroup $\cC$ of $G$, and the image of $g$ in $G/\K^\times$ lies in a compact subgroup. Conversely, an element whose image in $G/\K^\times$ lies in a compact subgroup is compact-mod-center.
\end{lemma}
\begin{proof}
First, if $gM_j=M_{j+c}$ and $\pi M_j=M_{j+r}$, then $\pi^{-c}g^r$ stabilizes the lattice $M_0$. Thus it lies in a compact open subgroup and has unit eigenvalues, and the eigenvalues of $g$ have valuation $c/r$.

Conversely, let all eigenvalues of $g$ have valuation $d_*/e$ in lowest terms. Then $h=\pi^{-d_*}g^e$ has unit eigenvalues, so $h$ and $h^{-1}$ are integral over $\fo$. Thus for any lattice $\Lambda$ the module $\Lambda_0=\fo[h,h^{-1}]\Lambda$ is an $h$-stable lattice, and we put:
\[
M_j=\sum\pi^ag^i\Lambda_0,
\]
where the sum is over the pairs $(a,i)$ with $a\in\Z$, $0\le i<e$ and $ae+id_*\ge j$. These are lattices, decreasing in $j$, with $\pi M_j=M_{j+e}$. Multiplication by $g$ increases the weight of every generator by $d_*$: it sends the generator $\pi^ag^i\Lambda_0$ of weight $ae+id_*$ to $\pi^ag^{i+1}\Lambda_0$ if $i<e-1$, and to $\pi^{a+d_*}h\Lambda_0=\pi^{a+d_*}\Lambda_0$ if $i=e-1$. Every generator of weight at least $j+d_*$ arises in this way from one of weight at least $j$, thus $gM_j=M_{j+d_*}$.

To prove the remaining statements, let $\mathsf N$ be the group of all $g'\in G$ with $g'M_j=M_{j+c(g')}$ for all $j$ and some $c(g')\in\Z$. Then $c:\mathsf N\to\Z$ is a homomorphism whose kernel $\mathsf P$ is the stabilizer of the chain, and $\mathsf P$ is compact and open. Since $c(\pi)=r$, the subgroup $\K^\times\mathsf P$ has index at most $r$ in $\mathsf N$, and $\mathsf N/\K^\times$ is compact.

Conversely, let the image of $g$ lie in a compact subgroup of $G/\K^\times$, and denote by $\Delta(g')$ the difference of the largest and the smallest valuation of the eigenvalues of $g'$. The function $\Delta$ is locally constant on $G$, because the Newton polygon of the characteristic polynomial is locally constant, and it is invariant under $\K^\times$. Thus $\Delta$ is bounded on the powers of $g$, and $\Delta(g)=0$ because $\Delta(g^N)=N\Delta(g)$.
\end{proof}

\subsection{Reduction modulo commutators}\label{sub:reduction}
For $w\in\tW$ and a simple reflection $s$ we have the following congruences in $\cH$:
\begin{equation}\label{eq:GP}
T_w\equiv T_{sws}\ \ \text{if }\ell(sws)=\ell(w),\qquad T_w\equiv(q-1)\,T_{sw}+q\,T_{sws}\ \ \text{if }\ell(sws)=\ell(w)-2,
\end{equation}
and $T_{\tau w\tau^{-1}}\equiv T_w$, because $T_\tau T_wT_\tau^{-1}=T_{\tau w\tau^{-1}}$ (Proposition~\ref{prop:IM}(ii)). Indeed, if $\ell(sws)=\ell(w)$ and $sws\ne w$, then one of $sw$, $ws$ is shorter than $w$ (Lemma~\ref{lem:IMlength}(iv)), say $w=sw_1$ with $\ell(w)=\ell(w_1)+1$, and then we have:
\[
T_w=T_sT_{w_1}\equiv T_{w_1}T_s=T_{sws}.
\]
Next, if $w=sw_1s$ with $\ell(w)=\ell(w_1)+2$, we get:
\[
T_w=T_sT_{w_1}T_s\equiv T_{w_1}T_s^2=(q-1)T_{w_1s}+qT_{w_1},
\]
where $w_1s=sw$ and $w_1=sws$. We write $w\approx w'$ if $w'$ is obtained from $w$ by a sequence of conjugations by simple reflections which preserve the length and conjugations by powers of $\tau$, and in this case $T_w\equiv T_{w'}$. An element $w$ is \emph{terminal} if no $w'\approx w$ admits a simple reflection $s$ with $\ell(sw's)<\ell(w')$.

\begin{lemma}\label{lem:shift}
If $w\approx w'$, then every element of $IwI$ is $G$-conjugate to an element of $Iw'I$.
\end{lemma}
\begin{proof}
For conjugation by $\tau$ the statement holds because $\tau$ normalizes $I$. Next, let $\ell(sws)=\ell(w)$ and $sws\ne w$, and say $w=sw_1$ with lengths adding, the case $w=w_2s$ being symmetric. Then $IwI=IsI\cdot Iw_1I$ and $IswsI=Iw_1I\cdot IsI$ by \eqref{eq:IM}, hence for $g=ab$ with $a\in IsI$ and $b\in Iw_1I$ we have $a^{-1}ga=ba\in IswsI$.
\end{proof}

A \emph{standard triple} is a triple $(K,x,u)$ with the following properties: $K$ is a proper subset of the set of simple reflections, $x\in\tW$ satisfies $xKx^{-1}=K$ and $\ell(sx)>\ell(x)$ for $s\in K$, and $u\in W_K$. Here $W_K$ stands for the subgroup generated by $K$, which is finite by Lemma~\ref{lem:parabolic}(i). An element $x\in\tW$ is \emph{straight} if $\ell(x^n)=n\ell(x)$ for all $n\ge1$. The following form of the reduction adapts and combines
He--Nie \cite[Propositions~2.5 and~2.8]{HN}.
A proof in our conventions is given in Section~\ref{app:HN3},
including the gradient-flow argument and the treatment of the
infinite group of length-zero elements.

\begin{proposition}[He--Nie]\label{prop:HN3}
For every terminal $w\in\tW$ there is a standard triple $(K,x,u)$ with $x$ straight and $w\approx ux$. The Newton point of $x$ is central if and only if that of $w$ is central, and in this case $x=\tau^{\val(w)}$.
\end{proposition}
\begin{proof}
This is Corollary~\ref{cor:terminal-standard}, proved in
Section~\ref{app:HN3} from the all-element reduction of
Proposition~\ref{prop:HN3app}.
\end{proof}

\begin{lemma}\label{lem:cert}
Let $(K,x,u)$ be a standard triple with $x$ straight and $\nu_x$ not central, and let $w\approx ux$. Then $IwI\cap G_c=\emptyset$, and consequently $D(T_w)=0$ for every distribution $D$ supported on $G_c$.
\end{lemma}
\begin{proof}
By Lemma~\ref{lem:shift} it suffices to treat $w=ux$, and we put $P_K=\bigcup_{y\in W_K}IyI$, which is a group by Lemma~\ref{lem:parabolic}(iii). Since $x$ is the shortest element of $W_Kx$ (Lemma~\ref{lem:standard}(S5)), we have $\ell(ux)=\ell(u)+\ell(x)$ and $IuxI=IuI\cdot IxI\subseteq P_K\cdot IxI$. For $s\in K$ we have $\ell(xs)=\ell((xsx^{-1})x)>\ell(x)$, and (S5) applied to $x^{-1}$, together with Lemma~\ref{lem:IMlength}(i), shows that $x$ is also the shortest element of $xW_K$ and that $\ell(xy)=\ell(x)+\ell(y)$ for $y\in W_K$. Thus for $y\in W_K$ we have:
\[
IxI\cdot IyI=IxyI=I(xyx^{-1})xI=I(xyx^{-1})I\cdot IxI,\qquad xyx^{-1}\in W_K.
\]
Hence $IxI\cdot P_K\subseteq P_K\cdot IxI$, and we obtain:
\[
(IuxI)^n\subseteq P_K\cdot(IxI)^n=P_K\cdot Ix^nI\subseteq\bigcup_{y\in W_K}Iyx^nI ,
\]
where the equality uses straightness. For the last inclusion,
successive left multiplication by simple reflections in $K$ keeps
every Iwahori index in the left coset $W_Kx^n$, by \eqref{eq:IM};
apply this to the finite parabolic decomposition of $P_K$.

Let $x^{n_0}=t^\lambda$, where $\lambda$ is not a multiple of $(1,\dots,1)$, and suppose that $g\in IuxI$ is compact-mod-center. By Lemma~\ref{lem:Gc}, $g^r\in\pi^c\cC$ for a compact open subgroup $\cC$. Here $\cC\subseteq\bigcup_{z\in Z_0}IzI$ for a finite set $Z_0\subseteq\tW$, and for every $l\ge1$ we have:
\[
g^{rn_0l}\in\bigcup_{y\in W_K}Iyt^{rl\lambda}I\qquad\text{and}\qquad g^{rn_0l}\in\pi^{cn_0l}\cC\subseteq\bigcup_{z\in Z_0}It^{cn_0l(1,\dots,1)}zI.
\]
Since double cosets of distinct elements of $\tW$ are disjoint, we get $t^{l(r\lambda-cn_0(1,\dots,1))}\in W_KZ_0$ for all $l\ge1$. This is a contradiction, because the right side is a finite set independent of $l$, while the left side takes infinitely many values.
\end{proof}

\subsection{Blocks}\label{sub:blocks}
We index the simple reflections by $\Z/m$. Then $\tau s_i\tau^{-1}=s_{i+1}$, and the image of $s_i$ in $S_m$ is the transposition of the positions $i,i+1\in\Z/m$. We fix $k$ and set $n'=\gcd(m,k)$, $e=m/n'$ and $\delta=\mathrm{Ad}(\tau^k):s_i\mapsto s_{i+k}$. A \emph{block structure} is a nonempty set $\Sigma\subseteq\Z/m$ of \emph{separators} which is stable under $i\mapsto i+n'$, or equivalently under $i\mapsto i+k$. Its \emph{blocks} are the arcs of positions $B=\{p+1,\dots,p+c\}$ between cyclically consecutive separators $p$ and $p+c$, and the sizes of the blocks in one period form a cyclic composition of $n'$.

We put $K_\Sigma=\{s_i:i\notin\Sigma\}$ and note that every $\delta$-stable proper subset of the simple reflections is of this form. Then $W_{K_\Sigma}=\prod_BS_B$, where $S_B$ is generated by $s_{p+1},\dots,s_{p+c-1}$ (Lemma~\ref{lem:parabolic}(i)). The map $\delta$ sends $S_B$ onto $S_{B+k}$ and simple reflections to simple reflections, and the blocks $B+jk$, $0\le j<e$, are pairwise distinct. We put $\gamma_B=s_{p+c-1}\cdots s_{p+2}s_{p+1}$, with $\gamma_B=1$ if $c=1$, and for a set $R$ of representatives of the $\delta$-orbits of blocks we define:
\[
u(\Sigma,R)=\prod_{B\in R}\gamma_B\in W_{K_\Sigma},\qquad[\Sigma,R]=\text{the class of }T_{u(\Sigma,R)\tau^k}\ \text{in }\cH/[\cH,\cH].
\]
Let $\Sigma_\rho$ be the block structure with separators $0,\rho_1,\rho_1+\rho_2,\dots$ modulo $n'$. If $R$ is the set of its blocks contained in $\{1,\dots,n'\}$, then $u(\Sigma_\rho,R)\tau^k=w_\rho$.

For a set $K$ of simple reflections let us denote by $\cH_K\subseteq\cH$ the span of the $T_y$, $y\in W_K$. For $K=K_\Sigma$ it is the product of the commuting subalgebras $\cH_B$ spanned by the $T_y$, $y\in S_B$ (Lemma~\ref{lem:parabolic}(iv)). By Proposition~\ref{prop:genericHecke}(ii), $\cH_B$ is the specialization at $q$ of the Hecke algebra $H_{|B|}$ of the symmetric group (Lemma~\ref{lem:finiteA} below). We put $\vartheta=T_{\tau^k}$ and observe that $\vartheta T_y\vartheta^{-1}=T_{\delta(y)}$ by Proposition~\ref{prop:IM}(ii), thus for $\delta$-stable $K$ the formula $\delta(h)=\vartheta h\vartheta^{-1}$ extends $\delta$ to $\cH_K$. We also have $T_{y\tau^k}=T_y\vartheta$, again by Proposition~\ref{prop:IM}(ii), and we consider the map:
\[
\phi:\cH_K\longrightarrow\cH/[\cH,\cH],\qquad\phi(h)=h\,\vartheta .
\]

\begin{lemma}\label{lem:twist}
\begin{enumerate}[label=\textup{(\roman*)}]
\item $\phi(h_1h_2)=\phi(h_2\,\delta(h_1))$ for $h_1,h_2\in\cH_K$.
\item Let $K=K_\Sigma$ and let $R$ be a set of representatives of the $\delta$-orbits of blocks. Then for $y\in W_K$ the element $\phi(T_y)$ is a linear combination of the elements $\phi(\prod_{B\in R}T_{y_B})$, $y_B\in S_B$. The coefficients are the values at $q$ of polynomials in $\Z[X]$ which depend only on $m$, $k$, $\Sigma$, $R$ and $y$, and not on $\K$.
\item $\phi(hh'g)=\phi(h'hg)$ for $h,h'\in\cH_B$, $B\in R$, and $g\in\prod_{B'\in R\smallsetminus\{B\}}\cH_{B'}$.
\end{enumerate}
\end{lemma}
\begin{proof}
(i) Since $\delta(h_1)\vartheta=\vartheta h_1$, we have $h_1h_2\vartheta-h_2\delta(h_1)\vartheta=h_1(h_2\vartheta)-(h_2\vartheta)h_1$.

(ii) If we write $y=\prod_By_B$ over all blocks, then $T_y=\prod_BT_{y_B}$ and the factors commute (Lemma~\ref{lem:parabolic}(iv)). Let $B_0\in R$ and $B=B_0+jk$ with $1\le j\le e-1$. Applying (i) with $h_1=T_{y_B}$ and $h_2$ the product of the other factors, we see that $\phi$ does not change if the factor at $B$ is replaced by $1$ and the factor at $B+k$ is multiplied on the right by $T_{\delta(y_B)}\in\cH_{B+k}$. We do this first for $j=e-1$, where $B+k=B_0$ since $ek\equiv0$ modulo $m$, then for $j=e-2$ and again for $j=e-1$, and so on. In this way all factors are moved into $\cH_{B_0}$, and these moves do not depend on $\K$.

Finally, by Proposition~\ref{prop:genericHecke}, a product of basis elements of $\cH_{B_0}$ is a combination of basis elements whose coefficients are the values at $q$ of structure constants of $H_{|B_0|}$, and these structure constants are polynomials in $\Z[X]$.

(iii) By (i) with $h_1=h$ and $h_2=h'g$ we have $\phi(hh'g)=\phi(h'g\,\delta(h))$ with $\delta(h)\in\cH_{B+k}$. If $e=1$, then $\tau^k$ is central, $\delta$ is the identity, and $h'gh=h'hg$. If $e>1$, we move the factor $\delta(h)$ around the orbit as in (ii), each time on the right, and we obtain $\phi(h'g\,\delta^e(h))$. Here $\delta^e$ is the identity because $\tau^{ke}$ is central, and again $h'gh=h'hg$.
\end{proof}

\begin{lemma}\label{lem:finiteA}
Let $H_n$ be the Hecke algebra of $S_n$ over $\Z[X]$ of Proposition~\ref{prop:genericHecke}, in which $T_{s_i}^2=(X-1)T_{s_i}+X$, and let $H_0=H_1=\Z[X]$. For a composition $\lambda=(\lambda_1,\dots,\lambda_r)$ of $n$ we put $a_i=\lambda_1+\dots+\lambda_{i-1}$, denote by $w_\lambda$ the product over the blocks $\{a_i+1,\dots,a_i+\lambda_i\}$ of the elements $s_{a_i+\lambda_i-1}\cdots s_{a_i+1}$, and set $T_\lambda=T_{w_\lambda}$. Then we have:
\[
H_n=\sum_{\lambda}\Z[X]\,T_\lambda+[H_n,H_n],
\]
where the sum is over all compositions $\lambda$ of $n$.
\end{lemma}
\begin{proof}
We use induction on $n$, the cases $n\le1$ being clear. We write $\ell$ for the number of inversions, and we denote by $H_j\subseteq H_n$ the subalgebra generated by $T_{s_1},\dots,T_{s_{j-1}}$. By the rules of Proposition~\ref{prop:genericHecke}(i), the span of the $T_w$ with $w\in S_j$ is stable under right multiplication by these generators.

Conversely, if $w\in S_j$ with $\ell(w)>0$ and $i$ is such that $\ell(ws_i)<\ell(w)$, then $w(i)>w(i+1)$, hence $i<j$, because $w(c)=c$ for $c>j$ and $w(j)\le j$. By induction on $\ell(w)$, the element $T_w=T_{ws_i}T_{s_i}$ lies in $H_j$, so $H_j$ is the span of the $T_w$ with $w\in S_j$, and by the uniqueness in Proposition~\ref{prop:genericHecke}(i) it is the Hecke algebra of $S_j$.

\emph{Claim.} For $2\le j\le n$ we have:
\[
H_j=H_{j-1}+T_{s_{j-1}}H_{j-1}+[H_j,H_{j-1}].
\]
Indeed, let $w\in S_j$ and $i=w^{-1}(j)$, and put $d=1$ if $i=j$ and $d=s_{j-1}s_{j-2}\cdots s_i$ if $i<j$. Then $d(i)=j$ and $d$ maps $\{1,\dots,j\}\smallsetminus\{i\}$ increasingly onto $\{1,\dots,j-1\}$, hence $a=wd^{-1}$ fixes $j$, that is, $a\in S_{j-1}$. The inversions of $w=ad$ are the $j-i$ pairs $(i,c)$ with $c>i$, and the pairs $(c,c')$ with $c,c'\ne i$ for which $(d(c),d(c'))$ is an inversion of $a$. Thus $\ell(w)=\ell(a)+\ell(d)$. If $d\ne1$, then $d=s_{j-1}b$ with $b\in S_{j-1}$ and $\ell(d)=1+\ell(b)$, so $T_w=T_aT_{s_{j-1}}T_b$ by Proposition~\ref{prop:genericHecke}(i), and we have:
\[
T_aT_{s_{j-1}}T_b=T_{s_{j-1}}T_bT_a+[T_a,T_{s_{j-1}}T_b].
\]

Next we put $E_1=1$ and $E_{i+1}=T_{s_{n-1}}T_{s_{n-2}}\cdots T_{s_{n-i}}$ for $1\le i\le n-1$, and we observe that $E_{i+1}$ commutes with $H_{n-i-1}$. Indeed, $s_{n-1},\dots,s_{n-i}$ commute with $s_1,\dots,s_{n-i-2}$ and the lengths add, so the corresponding elements $T_{s_c}$ commute by Proposition~\ref{prop:genericHecke}(i). Multiplying the claim for $j=n-i$ by $E_{i+1}$ on the left, and using $E_{i+1}T_{s_{n-i-1}}=E_{i+2}$ and $E_{i+1}(hg-gh)=[E_{i+1}h,g]$ for $g\in H_{n-i-1}$, we get:
\begin{equation}\label{eq:finiteA-step}
E_{i+1}H_{n-i}\subseteq E_{i+1}H_{n-i-1}+E_{i+2}H_{n-i-1}+[H_n,H_n]\qquad(0\le i\le n-2).
\end{equation}
Starting from $H_n=E_1H_n$ and applying this inclusion for $i=0,1,\dots,n-2$, we obtain:
\[
H_n=\sum_{j=1}^nE_jH_{n-j}+[H_n,H_n].
\]
Since $E_j[h_1,h_2]=[E_jh_1,h_2]$ for $h_1,h_2\in H_{n-j}$, the induction hypothesis for $H_{n-j}$ gives:
\begin{equation}\label{eq:finiteA-ind}
E_jH_{n-j}\subseteq\sum_{\mu}\Z[X]E_jT_\mu+[H_n,H_n],
\end{equation}
where the sum is over the compositions $\mu$ of $n-j$. Finally, $E_jT_\mu=T_{(\mu,j)}$, because $E_j=T_{s_{n-1}\cdots s_{n-j+1}}$ and the lengths add.
\end{proof}

\begin{lemma}\label{lem:swap}
The class $[\Sigma,R]$ depends only on the partition of $n'$ formed by the sizes of the blocks of $\Sigma$ in one period.
\end{lemma}
\begin{proof}
\emph{Independence of $R$.} Suppose that $R'$ is obtained from $R$ by replacing $B$ by $B+k$. Applying Lemma~\ref{lem:twist}(i) with $h_1=T_{\gamma_B}$ and $h_2$ the product of the other factors, and noting that $\delta(\gamma_B)=\gamma_{B+k}$, we get $[\Sigma,R']=[\Sigma,R]$. Since the $\delta$-orbit of a block $B$ consists of the blocks $B+jk$, any two sets of representatives are related by a sequence of such replacements, and thus we may write $[\Sigma]$.

\emph{Rotation.} We have $\tau\,u(\Sigma,R)\tau^k\,\tau^{-1}=u(\Sigma+1,R+1)\tau^k$, thus $[\Sigma+1]=[\Sigma]$. Note that two block structures with the same cyclic sequence of sizes differ by a rotation.

\emph{Interchange.} Let $B=\{p+1,\dots,p+a\}$ and $B'=\{p+a+1,\dots,p+a+b\}$ be cyclically consecutive blocks in different $\delta$-orbits. Then $a+b\le n'$. Let $\widetilde\Sigma$ be obtained from $\Sigma$ by replacing the separators $p+a+jn'$ by $p+b+jn'$, $j\in\Z$, and note that its blocks $\widetilde B=\{p+1,\dots,p+b\}$ and $\widetilde B'=\{p+b+1,\dots,p+a+b\}$ have the interchanged sizes.

We put $\Sigma_J=\Sigma\smallsetminus(p+a+n'\Z)$ and $J=K_{\Sigma_J}=K_\Sigma\cup\{s_{p+a+jn'}\}_j$. Since $1\le a,b<n'$, the separators $p$ and $p+a+b$ of $\Sigma$ are not congruent to $p+a$ modulo $n'$, so they lie in $\Sigma_J$. Thus $J$ is a proper subset of the simple reflections, and the arcs $(B\cup B')+jn'$, which are pairwise disjoint because $a+b\le n'$, are blocks of $\Sigma_J$. By Lemma~\ref{lem:parabolic}(i) there is a unique $X\in W_J$ with the following two properties: on every window $\{\tilde p+1,\dots,\tilde p+a+b\}$ of $\Sigma_J$ over one of these arcs it acts by
\[
X(\tilde p+i)=\tilde p+b+i\ \ (1\le i\le a),\qquad X(\tilde p+a+i)=\tilde p+i\ \ (1\le i\le b),
\]
and it is the identity on the other windows of $\Sigma_J$.

By Lemma~\ref{lem:parabolic}(i),(ii), lengths and conjugations of elements of $W_J$ can be computed in the product of the symmetric groups of the windows, and by construction $\tau^kX\tau^{-k}=X$. We choose $R\ni B,B'$, put $\widetilde R=(R\smallsetminus\{B,B'\})\cup\{\widetilde B,\widetilde B'\}$, and observe that:
\[
Xs_{p+i}X^{-1}=s_{p+b+i}\ \ (1\le i\le a-1),\qquad Xs_{p+a+i}X^{-1}=s_{p+i}\ \ (1\le i\le b-1).
\]
Hence $X\gamma_BX^{-1}=\gamma_{\widetilde B'}$ and $X\gamma_{B'}X^{-1}=\gamma_{\widetilde B}$. The other factors of $u=u(\Sigma,R)$ have support disjoint from that of $X$, thus $XuX^{-1}=\widetilde u:=u(\widetilde\Sigma,\widetilde R)$ and $Xu\tau^k=\widetilde u\tau^kX$.

The element $X$ is increasing on every window of $\Sigma$, and $X^{-1}$ is increasing on every window of $\widetilde\Sigma$. By Lemma~\ref{lem:IMlength}(ii) we get $\ell(Xs)>\ell(X)$ for $s\in K_\Sigma$, and $\ell(sX)>\ell(X)$ for $s\in K_{\widetilde\Sigma}$. By (S5) applied to $X^{-1}$ and to $X$, the element $X$ is the shortest element of $XW_{K_\Sigma}$ and of $W_{K_{\widetilde\Sigma}}X$. Again by (S5) and Lemma~\ref{lem:IMlength}(i), this gives $\ell(Xu)=\ell(X)+\ell(u)$ and $\ell(\widetilde uX)=\ell(\widetilde u)+\ell(X)$, and therefore:
\[
T_X\,T_{u\tau^k}=T_{Xu\tau^k}=T_{\widetilde u\tau^kX}=T_{\widetilde u\tau^k}\,T_X .
\]
Since $T_X$ is invertible, we get $T_{\widetilde u\tau^k}=T_XT_{u\tau^k}T_X^{-1}\equiv T_{u\tau^k}$, and thus $[\Sigma]=[\widetilde\Sigma]$.

Finally, rotations and interchanges of consecutive blocks generate all rearrangements of the cyclic sequence of sizes.
\end{proof}

\subsection{Proof of Theorem~\ref{thm:H}}\label{sub:proofH}
Let us fix $D\in J_k$ and recall that, by Lemma~\ref{lem:trace}, $D$ vanishes on $[\cH,\cH]$.
The reduction identities in Steps~1 and~3 hold in the full
cocenter. Only Step~2 uses support on $G_c$ to discard
noncentral terminal terms after applying $D$.

\emph{Step 1. Reduction to terminal elements.} If $w$ is not terminal, then there are $w_1\approx w$ and $s$ with $\ell(sw_1s)=\ell(w_1)-2$, and by \eqref{eq:GP} we have:
\[
T_w\equiv(q-1)T_{sw_1}+qT_{sw_1s},
\]
which is a combination of shorter elements of the same valuation. Fixing such a choice for every $w$ and using induction on the length, we get:
\[
T_w\equiv\sum_{w'}b_{w,w'}(q)T_{w'},
\]
where $w'$ is terminal, $\val(w')=\val(w)$ and $b_{w,w'}\in\Z[X]$ is independent of $\K$.

\emph{Step 2. Terminal elements with non-central Newton point.} If $w'$ is terminal, then by Proposition~\ref{prop:HN3} we have $w'\approx ux$ for a standard triple with $x$ straight. If $\nu_{w'}$ is not central, then $\nu_x$ is not central, and $D(T_{w'})=0$ by Lemma~\ref{lem:cert}.

\emph{Step 3. Terminal elements with central Newton point.} In this case $x=\tau^k$, and $w'\approx u\tau^k$ with $u\in W_K$, where $K$ is a proper $\delta$-stable set of simple reflections. Thus $K=K_\Sigma$ for a block structure $\Sigma$, and $T_{w'}\equiv T_{u\tau^k}=\phi(T_u)$. By Lemma~\ref{lem:twist}(ii), this is a combination of elements $\phi(\prod_{B\in R}h_B)$ with $h_B$ basis elements of $\cH_B$, with coefficients as in that lemma, and by Lemma~\ref{lem:twist}(iii), $\phi(\prod_Bh_B)$ depends on each $h_B$ only modulo $[\cH_B,\cH_B]$.

Applying Lemma~\ref{lem:finiteA} in $H_{|B|}$ and the isomorphism of Proposition~\ref{prop:genericHecke}(ii) for each block $B$, we see that $\phi(\prod_Bh_B)$ is a combination of the elements $\phi(\prod_{B\in R}T_{\lambda^B})$, where $\lambda^B$ is a composition of $|B|$. The coefficients are the values at $q$ of polynomials in $\Z[X]$ given by Lemma~\ref{lem:finiteA}. Now $\prod_{B\in R}T_{\lambda^B}=T_{u(\Sigma',R')}$, where $\Sigma'\supseteq\Sigma$ is the block structure which subdivides every block $B+jk$ according to $\lambda^B$, and $R'$ consists of the sub-blocks of the blocks in $R$. By Lemma~\ref{lem:swap} we have:
\[
\phi(T_{u(\Sigma',R')})=[\Sigma']=[\Sigma_\rho],
\]
where $\rho$ is the partition formed by the sizes of $\Sigma'$, and finally $[\Sigma_\rho]$ is the class of $T_{w_\rho}$. All the polynomials which occur in Steps 1 and 3 are determined by $m$, $k$, $w$ and by choices made in $\tW$, and they do not depend on $\K$.

This proves (ii), and part (i) follows because $D$ is determined by its values on the $T_w$. \qed

\begin{remark}\label{rem:terminal}
(1) For $w_\rho=u_\rho\tau^k$ we have:
\[
w_\rho^e=u_\rho\,\delta(u_\rho)\cdots\delta^{e-1}(u_\rho)\,\tau^{ke},
\]
where $\tau^{ke}$ is central and the other factor is the product of the Coxeter elements of all blocks of $\Sigma_\rho$, hence $\nu_{w_\rho}$ is central. Since $w_\rho$ permutes the $e$ translates of each block cyclically, the cycles of its underlying permutation have lengths $e\rho_i$. Thus, by Lemma~\ref{lem:classes}, the $w_\rho$ represent the $p(n')$ conjugacy classes with central Newton point and valuation $k$.

Already for $m=k=6$, the elements
$(s_1s_3)\tau^6$ and $(s_1s_4)\tau^6$ have the same partition
$(2,2,1,1)$ but are not related by $\approx$.
A length-preserving simple conjugation fixes each product,
while $\tau$-conjugation rotates the indices modulo $6$ and
preserves their circular distance, respectively $2$ and $3$.
Their cocenter classes are equal by Lemma~\ref{lem:swap}.
Thus the block-sorting step supplies an identity beyond
the elementary reductions in \eqref{eq:GP}.

(2) Theorem~\ref{thm:H}(i) is contained in \cite[\S2.3]{KT}; the relation is described in Section~\ref{sub:method}. For $(m,k)=(4,2)$ the elements $w_{(1,1)}=\tau^2$ and $w_{(2)}=s_1\tau^2$ are the elements $w_1,w_2$ of \cite[Ex.~2.32]{KT}. The polynomials $a_{w,\rho}$ of Theorem~\ref{thm:H}(ii) are not considered there. From this theory the proof above uses only Proposition~\ref{prop:HN3}. Lemma~\ref{lem:finiteA} is a form of the classical description of the cocenter of the Hecke algebra of the symmetric group \cite{Ram}, \cite[\S8.2]{GP}, of which only the spanning property is needed here.

(3) Section~\ref{subsec:explicit-counting-algorithm} gives a finite algorithm for these three steps.
\end{remark}

\section{Orbital integrals and chains of lattices}\label{sec:chains}

The relations of Section~\ref{sec:hecke} become identities between
lattice counts: ordinary parahoric test functions count stable chains,
and rigid test functions count chains over the transformed ring.
The resulting transfer formula (Theorem~\ref{thm:transfer}) gives
equisingular invariance and polynomiality by induction through
quadratic transforms. Section~\ref{subsec:explicit-counting-algorithm}
gives a finite algorithm for the counting polynomial. Throughout this
section we take $\K=F((\pi))$, $\fo=F[[\pi]]$, and
$G=\GL_m(\K)$, with the notation of Section~\ref{sec:hecke}.

\subsection{Framed modules}\label{sub:framedmod}
Let $(u,w)$ be a framed pair in $\mathbf O=O^r$ with ring $T=F[[u,w]]$ (Definition~\ref{def:framed}), let $\mathbf n=(n_1,\dots,n_r)$ be a tuple of positive integers, and set:
\[
V=\prod_{i=1}^r\cK^{n_i},\qquad\mathbf O^{\mathbf n}=\prod_{i=1}^rO^{n_i},\qquad Z=\prod_{i=1}^r\GL_{n_i}(\cK),\qquad Z_0=\prod_{i=1}^r\GL_{n_i}(O),
\]
and $K_u=F((u))$, $\fo_u=F[[u]]$. By Lemma~\ref{lem:gen}(a), $V$ is a $K_u$-vector space of dimension $m=\sum_in_i\val_i(u)$. Since $K_u[w]=\mathbf K$, the group $Z=\Aut_{\mathbf K}(V)$ is the centralizer of $w$ in $G_u:=\GL_{K_u}(V)$. We equip $Z$ with the Haar measure such that $\vol(Z_0)=1$.

A \emph{lattice} in $V$ is a compact open $\fo_u$-submodule. A \emph{$T$-lattice} is a lattice stable under $w$, or equivalently a compact open $T$-submodule. For a $T$-lattice $\Lambda$ the module $\mathbf O\Lambda$ is an $\mathbf O$-lattice, so it is of the form $z\mathbf O^{\mathbf n}$ with $z\in Z$. Let us denote by $E_T(\mathbf O^{\mathbf n})$ the set of $T$-lattices $\Lambda$ with $\mathbf O\Lambda=\mathbf O^{\mathbf n}$. This set is finite: indeed, $T\supseteq t^c\mathbf O$ for some $c$, and then every $\Lambda\in E_T(\mathbf O^{\mathbf n})$ satisfies $t^c\mathbf O^{\mathbf n}\subseteq\Lambda\subseteq\mathbf O^{\mathbf n}$.

\begin{lemma}\label{lem:massE}
For every $Z$-stable set $\mathcal S$ of $T$-lattices in $V$ we have:
\[
\mass_Z(\mathcal S)=\#\big(\mathcal S\cap E_T(\mathbf O^{\mathbf n})\big).
\]
\end{lemma}
\begin{proof}
Every $Z$-orbit in $\mathcal S$ meets $E_T(\mathbf O^{\mathbf n})$, because $\mathbf O\Lambda$ is of the form $z\mathbf O^{\mathbf n}$ with $z\in Z$ for every $T$-lattice $\Lambda$. Next, two elements of $E_T(\mathbf O^{\mathbf n})$ in the same $Z$-orbit are in the same $Z_0$-orbit, and the stabilizer in $Z$ of $\Lambda\in E_T(\mathbf O^{\mathbf n})$ stabilizes $\mathbf O\Lambda$, so it lies in $Z_0$. Thus, using $\vol(Z_0)=1$, we get:
\[
\mass_Z(\mathcal S)=\sum_{\Lambda\in Z_0\backslash(\mathcal S\cap E_T)}[Z_0:Z_{0,\Lambda}],
\]
and the right-hand side is the number of elements of $\mathcal S\cap E_T(\mathbf O^{\mathbf n})$.
\end{proof}

Let us recall that a composition $\alpha=(\alpha_1,\dots,\alpha_l)$ of $m$ is a sequence of positive integers with sum $m$. A \emph{$u$-periodic chain of type $\alpha$} in $V$ is a family $(\Lambda_j)_{j\in\Z}$ of lattices such that $\Lambda_{j-1}\supseteq\Lambda_j$ and $\Lambda_{j+l}=u\Lambda_j$ for all $j$, and $\dim_F\Lambda_{j-1}/\Lambda_j=\alpha_j$ for $1\le j\le l$.

\begin{definition}\label{def:chainmass}
The \emph{chain mass} $\Ch_\alpha(u,w;\mathbf n)$ is the $Z$-mass of the set of $u$-periodic chains of type $\alpha$ of $T$-lattices in $V$. For $\mu\vdash m$ the \emph{refined count} $B_{(u,w)}(\mu;\mathbf n)$ is the number of $\Lambda\in E_T(\mathbf O^{\mathbf n})$ such that the nilpotent endomorphism induced by $w$ on $\Lambda/u\Lambda\cong F^m$ has Jordan type $\mu$. For $r=1$ we write $\Ch_\alpha(T,u;n)$ and $B_{(u,w)}(\mu;n)$.
\end{definition}
\noindent Note that the chain mass depends only on $T$, $u$ and $\mathbf n$, and the same is true for the refined count. Indeed, by Lemma~\ref{lem:gen} the ring $T/uT=F[\bar w]$ is a truncated polynomial ring, and the Jordan type of a generator of its maximal ideal on a module does not depend on the generator:
if $N$ is the action of one generator, the action of another one is $NU(N)$ with
$U(0)\ne0$. The operator $U(N)$ is invertible and commutes
with $N$, so $\ker((NU(N))^j)=\ker(N^j)$ for every $j$;
these dimensions determine the Jordan type. If the maximal
ideal is zero, both actions are zero.

For a composition $\alpha$ of $m$ and $\mu\vdash m$ we set:
\[
H_\alpha(\mu)=\#\{F^m=V_0\supseteq V_1\supseteq\dots\supseteq V_l=0:\ \dim V_{j-1}/V_j=\alpha_j,\ XV_j\subseteq V_j\}\qquad(X\in\bO_\mu).
\]

\begin{lemma}\label{lem:chainflag}
For every composition $\alpha$ of $m$ we have:
\[
\Ch_\alpha(u,w;\mathbf n)=\sum_{\mu\vdash m}H_\alpha(\mu)\,B_{(u,w)}(\mu;\mathbf n).
\]
In particular $\Ch_{(m)}(u,w;\mathbf n)=\#E_T(\mathbf O^{\mathbf n})$, and if $r=1$ and $T=O$, then $\Ch_\alpha(O,u;n)=H_\alpha((D^n))$, where $D=\val u$.
\end{lemma}
\begin{proof}
A chain is determined by $\Lambda_0$ and by the flag $(\Lambda_j/u\Lambda_0)_{0\le j\le l}$ in $\Lambda_0/u\Lambda_0$. Every subspace between $u\Lambda_0$ and $\Lambda_0$ is a lattice, and it is a $T$-lattice if and only if its image is stable under $\bar w$. Hence the fibre of the map $(\Lambda_j)\mapsto\Lambda_0$ over a $T$-lattice $\Lambda_0$ has $H_\alpha(\mu)$ elements, where $\mu$ is the Jordan type of $\bar w$ on $\Lambda_0/u\Lambda_0$, and the formula follows from \eqref{eq:massfibre} applied to this map and from Lemma~\ref{lem:massE}.

If $T=O$, then $E_O(O^n)=\{O^n\}$. For $D>1$, the
equality $O=F[[u,w]]$ forces $\val(w)=1$, because $u,w$ generate
the maximal ideal of $O$. Hence $O=F[[w]]$ and
$u=w^Dv(w)$ for a unit $v(w)$, so multiplication by $w$ on
$O/uO\cong F[w]/(w^D)$ is one Jordan block of size $D$.
For $D=1$, the quotient is $F$ and $w$ acts by zero. On $n$
copies the Jordan type is therefore $(D^n)$ in both cases.
\end{proof}

Thus the chain mass counts pairs consisting of a normalized lattice
$\Lambda_0\in E_T(\mathbf O^{\mathbf n})$ and a compatible flag in
$\Lambda_0/u\Lambda_0$. In particular, it is a nonnegative integer;
we will use this when passing from rational transfer formulas to
counting polynomials.

\subsection{Parahoric test functions}
Let $\alpha=(\alpha_1,\dots,\alpha_l)$ be a composition of $m$, and let $\sigma_\alpha(j)=\alpha_1+\dots+\alpha_j$ for $0\le j\le l$, extended to all $j\in\Z$ by $\sigma_\alpha(j+l)=\sigma_\alpha(j)+m$. The \emph{standard chain of type $\alpha$} is $\Lambda^\alpha_j=\cL_{\sigma_\alpha(j)}$, and it is $\pi$-periodic of type $\alpha$. We denote by $P_\alpha\supseteq I$ its stabilizer, which is a parahoric subgroup, and set:
\[
\fp_\alpha=\{g\in\Mat_m(\K):\ g\Lambda^\alpha_j\subseteq\Lambda^\alpha_j\ \text{for all }j\},\qquad F_{\alpha,k}=\frac{\mathbf 1_{\fp_\alpha\cap G^{\val=k}}}{\vol(P_\alpha)}\in\cH^{\val=k},
\]
where $G^{\val=k}=\{g:\val\det g=k\}$.

Now let $n'$ be a divisor of $m$ and $e=m/n'$, and for a composition $\beta=(\beta_1,\dots,\beta_l)$ of $n'$ let us write $\beta^{[e]}$ for the composition $(\beta,\beta,\dots,\beta)$ of $m$ with $e$ copies of $\beta$. Then $\sigma_{\beta^{[e]}}(j+l)=\sigma_{\beta^{[e]}}(j)+n'$, hence $\tau^{n'}\Lambda^{\beta^{[e]}}_j=\Lambda^{\beta^{[e]}}_{j+l}$ and $\tau^{n'}$ normalizes $P_{\beta^{[e]}}$. For $k=d_*n'$, $d_*\in\Z$, we set:
\[
f^*_\beta=f^{*,k}_\beta=\frac{\mathbf 1_{\tau^kP_{\beta^{[e]}}}}{\vol(P_{\beta^{[e]}})}\in\cH^{\val=k}.
\]
Membership in the support of $F_{\alpha,k}$ requires an inclusion
at each lattice of the chain; membership in the support of $f^*_\beta$
requires equality with a specified translate. Lemma~\ref{lem:formal}
identifies the latter condition with periodicity over the transformed
ring. Division by the parahoric volumes gives the same mass
normalization for both orbital integrals.

To expand these functions in the basis $(T_w)$, let $\Sigma_\alpha\subseteq\Z/m$ be the set of residues of the numbers $\sigma_\alpha(j)$ and $W_\alpha=W_{K_{\Sigma_\alpha}}$, and note that the blocks of $\Sigma_\alpha$ have the sizes $\alpha_1,\dots,\alpha_l$. An element $w\in\tW$ lies in $P_\alpha$ if and only if $w\cL_p=\cL_p$ for every $p$ whose residue lies in $\Sigma_\alpha$, that is, if and only if $w$ maps every window of $\Sigma_\alpha$ onto itself. Since $I\subseteq P_\alpha$, Lemmas~\ref{lem:IMbruhat}(iv) and~\ref{lem:parabolic}(i) give $P_\alpha=IW_\alpha I$, and we obtain:
\begin{equation}\label{eq:testfunctions}
F_{\alpha,k}=\frac1{[P_\alpha:I]}\sum_{w\in\tW,\ \val(w)=k,\ w\Lambda^\alpha_j\subseteq\Lambda^\alpha_j\ \forall j}T_w,\qquad f^*_\beta=\frac1{[P_{\beta^{[e]}}:I]}\sum_{y\in W_{\beta^{[e]}}}T_{\tau^ky},
\end{equation}
where the index sets are finite and do not depend on $q$, and $[P_\alpha:I]=\prod_{j=1}^l\prod_{i=1}^{\alpha_j}\frac{q^i-1}{q-1}$ by Lemma~\ref{lem:parabolic}(iii).

\begin{lemma}\label{lem:dictionary}
Let $\gamma\in G$, and assume that its centralizer $Z_\gamma$
is unimodular and equipped with a Haar measure. For a nonnegative
$f\in C_c^\infty(G)$ define
\[
 O_\gamma(f)=\int_{Z_\gamma\backslash G}
                  f(g^{-1}\gamma g)\,d\dot g\in[0,\infty].
\]
The following identities hold as equalities of nonnegative extended
real numbers; no convergence hypothesis is needed.
\begin{enumerate}[label=\textup{(\alph*)}]
\item If $\val\det\gamma=k$, then $O_\gamma(F_{\alpha,k})$ is the $Z_\gamma$-mass of the set of $\pi$-periodic chains $(\Lambda_j)$ of type $\alpha$ in $\K^m$ with $\gamma\Lambda_j\subseteq\Lambda_j$ for all $j$.
\item $O_\gamma(f^{*,k}_\beta)$ is the $Z_\gamma$-mass of the set of $\pi$-periodic chains $(\Lambda_j)$ of type $\beta^{[e]}$ with $\gamma\Lambda_j=\Lambda_{j+d_*l}$ for all $j$.
\end{enumerate}
\end{lemma}
\begin{proof}
The group $G$ acts transitively on the set of $\pi$-periodic chains of type $\alpha$: indeed, it suffices to lift a basis of $\Lambda_0/\pi\Lambda_0$ adapted to the flag $(\Lambda_j/\pi\Lambda_0)_j$. The stabilizer of the standard chain is $P_\alpha$, thus this set is $G/P_\alpha$, and $Z_\gamma$ acts on it with stabilizers $Z_\gamma\cap gP_\alpha g^{-1}$. Now let $C\subseteq G$ be a compact open subset with $P_\alpha CP_\alpha=C$. Since the function $g\mapsto\mathbf 1_C(g^{-1}\gamma g)$ is right $P_\alpha$-invariant and the double coset $Z_\gamma gP_\alpha$ has measure $\vol(P_\alpha)/\vol(Z_\gamma\cap gP_\alpha g^{-1})$ in $Z_\gamma\backslash G$, Tonelli's theorem gives the following
identity without first assuming finiteness:
\[
O_\gamma\Big(\frac{\mathbf 1_C}{\vol(P_\alpha)}\Big)=\mass_{Z_\gamma}\{gP_\alpha:\ g^{-1}\gamma g\in C\}.
\]
For $C=\fp_\alpha\cap G^{\val=k}$ the condition $g^{-1}\gamma g\in C$ means that $\gamma g\Lambda^\alpha_j\subseteq g\Lambda^\alpha_j$ for all $j$. For $C=\tau^kP_{\beta^{[e]}}$ it means that $\gamma g\Lambda^0_j=g\tau^k\Lambda^0_j=g\Lambda^0_{j+d_*l}$ for all $j$, where $\Lambda^0_\bullet$ is the standard chain of type $\beta^{[e]}$.
\end{proof}

\subsection{The orbital distribution of a balanced pair}
\begin{lemma}\label{lem:orbital}
Let $(u,w)$ be a balanced framed pair of type $(e,d_*)$ in $\mathbf O$ with $\val_i(u)=g_ie$, and let $\mathbf n$ be as above. We set $n'=\sum_in_ig_i$, so that $m=en'$, and $k=d_*n'$. Choosing a $K_u$-basis of $V$, we identify $G_u$ with $G=\GL_m(\K)$ and $u$ with $\pi$.
\begin{enumerate}[label=\textup{(\alph*)}]
\item The element $w\in G$ is compact-mod-center, $\val\det w=k$, $\gcd(m,k)=n'$, and the centralizer of $w$ is $Z$.
\item For every $f\in C_c^\infty(G)$ the integral $D_w(f)=\int_{Z\backslash G}f(g^{-1}wg)\,d\dot g$ converges absolutely. The functional $D_w$ is an invariant distribution supported on $G_c\cap G^{\val=k}$, and its restriction to $\cH^{\val=k}$ belongs to $J_k$.
\item For every composition $\alpha$ of $m$ we have:
\[
\Ch_\alpha(u,w;\mathbf n)=D_w(F_{\alpha,k}).
\]
\end{enumerate}
\end{lemma}
\begin{proof}
(a) The ring $O$ in the $i$-th factor is free over $F[[u_i]]$
by Lemma~\ref{lem:gen}(a). Smith normal form for multiplication
by $w_i$ gives
\[
 \val_{u_i}\det(w_i\mid\cK)
 =\operatorname{length}_{F[[u_i]]}(O/w_iO)
 =\dim_F(O/w_iO)=\val_t(w_i)=g_id_*.
\]
Taking $n_i$ copies and summing yields
$\val\det w=\sum_i n_i g_id_*=k$.
This calculation does not require separability. Next, we have $w^e=u^{d_*}\theta$, and the eigenvalues of the unit $\theta$ of $\mathbf O$ are units, because $\theta$ and its inverse preserve the lattice $\mathbf O^{\mathbf n}$. Since $w$ and $\theta$ commute, all eigenvalues of $w$ have valuation $d_*/e$. Finally, $\gcd(m,k)=n'\gcd(e,d_*)=n'$, and since $K_u[w]=\mathbf K$, the centralizer of $w$ is $Z=\Aut_{\mathbf K}(V)$.

(b) For $M\ge0$ the ring $R_M=\fo_u[u^Mw]\subseteq\mathbf O$ is a finitely generated $\fo_u$-module which spans $K_u[w]=\mathbf K$. Hence it has finite index in $\mathbf O$ and contains $t^c\mathbf O$ for some $c$, and as in Lemma~\ref{lem:massE}, the $Z$-mass of the set of lattices stable under $u^Mw$ is the number of those with saturation $\mathbf O^{\mathbf n}$, which is finite.

Now let $f\in C_c^\infty(G)$, put $L_0=\mathbf O^{\mathbf n}$,
and let $C_0$ be its stabilizer in $G$. Choose $M\ge0$ such that
$hL_0\subseteq u^{-M}L_0$ for all $h$ in the support of $f$, and set
\[
 C_M=\{h\in G:hL_0\subseteq u^{-M}L_0,\ \val\det h=k\}.
\]
This is a compact open $C_0$-biinvariant subset of $G$: it is
the intersection of a compact matrix lattice with the inverse
image, under the determinant, of $u^k\fo_u^\times$.
Every conjugate of $w$ has determinant valuation $k$.
Choose an $\fo_u$-basis of $L_0$ and apply the nonnegative
double-coset formula of Lemma~\ref{lem:dictionary} to $C_M$.
It gives
\[
 \int_{Z\backslash G}|f(g^{-1}wg)|\,d\dot g
 \le \|f\|_\infty\,\vol(C_0)\,
       \mass_Z\{\Lambda:u^Mw\Lambda\subseteq\Lambda\}<\infty.
\]
The same estimate is uniform for functions supported in a fixed
compact set, so $D_w$ is a distribution.
The invariance of $D_w$ follows by a change of variables. Finally, since the set $G_c\cap G^{\val=k}$ is closed, conjugation-invariant and contains $w$, it contains the support of $D_w$.

(c) This follows from Lemma~\ref{lem:dictionary}(a), since a lattice is a $T$-lattice if and only if it is stable under $w$.
\end{proof}

\subsection{Flags, nilpotent orbits and the Hall algebra}\label{sub:hall}
Let $\sH$ be the Hall algebra of finite nilpotent $F[z]$-modules over $\Q$ \cite[Ch.~II, III]{Mac}. Let us recall that it has the basis $(u_\mu)$ indexed by all partitions, with the multiplication:
\[
u_\alpha u_\beta=\sum_\mu g^\mu_{\alpha\beta}u_\mu,
\]
where $g^\mu_{\alpha\beta}$ stands for the number of submodules $N\subseteq M_\mu$ with $N\cong M_\beta$ and $M_\mu/N\cong M_\alpha$. The algebra $\sH$ is associative, and it is commutative because the duality $M\mapsto\Hom_F(M,F)$ preserves isomorphism types and exchanges submodules and quotients. The coefficient of $u_\mu$ in a product $u_{\alpha^1}\cdots u_{\alpha^l}$ is the number of filtrations of $M_\mu$ by submodules with successive quotients of types $\alpha^1,\dots,\alpha^l$. The numbers $g^\mu_{\alpha\beta}$ are the values at $q$ of polynomials with integer coefficients, the Hall polynomials \cite[Ch.~II, \S4]{Mac}. Finally, we write $a_\mu=|\Aut_{F[z]}(M_\mu)|$, which is the value at $q$ of a polynomial with integer coefficients as well \cite[Ch.~II, \S1]{Mac}.

For a composition $\mathbf n=(n_1,\dots,n_r)$ of $N$ let $e_2(\mathbf n)=\sum_{i<j}n_in_j$, and denote by $\fn_{\mathbf n}(F)\subseteq\Mat_N(F)$ the space of strictly lower block triangular matrices for the block sizes $\mathbf n$. It has $q^{e_2(\mathbf n)}$ elements, and it is the set of $X$ with $XV_i\subseteq V_{i-1}$ for a flag $0=V_0\subseteq\dots\subseteq V_r=F^N$ of coordinate subspaces with $\dim V_i/V_{i-1}=n_{r+1-i}$. We define the numbers $F_{\mathbf n}(\mu)$ by the expansion:
\[
u_{(1^{n_1})}\cdots u_{(1^{n_r})}=\sum_{\mu\vdash N}F_{\mathbf n}(\mu)\,u_\mu .
\]
Thus $F_{\mathbf n}(\mu)$ is the number of flags $0=V_0\subseteq V_1\subseteq\dots\subseteq V_r=F^N$ with $\dim V_i/V_{i-1}=n_i$ and $XV_i\subseteq V_{i-1}$, for a fixed $X\in\bO_\mu$. By commutativity it does not change if $\mathbf n$ is permuted.

\begin{lemma}\label{lem:hallbasics}
\begin{enumerate}[label=\textup{(\alph*)}]
\item $\#\big(\fn_{\mathbf n}(F)\cap\bO_\mu\big)=q^{e_2(\mathbf n)}\prod_{i}|\GL_{n_i}(F)|\cdot F_{\mathbf n}(\mu)/a_\mu$.
\item If $F_{\mathbf n}(\mu)\ne0$, then $\mu'\ctd\mathbf n^+$, where $\mathbf n^+$ is the partition obtained by sorting $\mathbf n$. Moreover, $F_\lambda(\lambda')=1$ for every partition $\lambda$, and if $0=V_0\subseteq\dots\subseteq V_l=F^N$ is a flag with $\dim V_i/V_{i-1}=\lambda_i$ and $X$ satisfies $XV_i\subseteq V_{i-1}$ for all $i$, then $X\in\bO_{\lambda'}$ if and only if all the induced maps $V_i/V_{i-1}\to V_{i-1}/V_{i-2}$, $2\le i\le l$, are injective.
\item The elements $u_{(1^j)}$, $j\ge1$, are algebraically independent and generate $\sH$.
\item The matrix $\sM_N=(H_\alpha(\mu))_{\alpha,\mu\vdash N}$, in which the partition $\alpha$ is regarded as a composition, is invertible. Its determinant is, up to sign, a power of $q$.
\end{enumerate}
\end{lemma}
\begin{proof}
(a) Let us count in two ways the pairs consisting of an element $X\in\bO_\mu$ and a flag with $\dim V_i/V_{i-1}=n_{r+1-i}$ and $XV_i\subseteq V_{i-1}$. If we fix $X$, we obtain $|\GL_N(F)|F_{\mathbf n}(\mu)/a_\mu$ pairs, and on the other hand, since the stabilizer $P_{\mathbf n}(F)$ of the flag has order $q^{e_2(\mathbf n)}\prod_i|\GL_{n_i}(F)|$, we obtain $|\GL_N(F)|\,\#(\fn_{\mathbf n}\cap\bO_\mu)/|P_{\mathbf n}(F)|$ pairs if we fix the flag.

(b) By commutativity we may assume that $\mathbf n=\lambda$ is a partition. If $XV_i\subseteq V_{i-1}$ for all $i$, then $V_i\subseteq\ker X^i$, and thus:
\[
\lambda_1+\dots+\lambda_i\le\dim\ker X^i=\mu'_1+\dots+\mu'_i,
\]
and if $\mu'=\lambda$, the inclusions $V_i\subseteq\ker X^i$ are equalities, so the flag is unique. The equalities $V_i=\ker X^i$ hold for all $i$ if and only if every map $V_i/V_{i-1}\to V_{i-1}/V_{i-2}$ is injective. Indeed, if these maps are injective, then by induction a vector outside $V_{i-1}$ is not killed by $X^{i-1}$, and conversely, a nonzero class in the kernel of such a map gives a vector outside $V_{i-1}$ killed by $X^{i-1}$.

(c) By (b), the products $u_{(1^{\lambda_1})}\cdots u_{(1^{\lambda_l})}$, $\lambda\vdash N$, are related to the basis $(u_{\lambda'})_{\lambda\vdash N}$ by a matrix which is unitriangular for the dominance order.

(d) Let us put $v_j=\sum_{\lambda\vdash j}u_\lambda$ and $v_0=1$, and note that $H_\alpha(\mu)$ is the coefficient of $u_\mu$ in $v_{\alpha_1}\cdots v_{\alpha_l}$. If $j,j'\ge0$ and $M$ is a module of dimension $j+j'$ whose socle $\ker z$ has dimension $s$, then also $\dim_FM/zM=s$, and the coefficient of $[M]$ in $u_{(1^j)}v_{j'}$ is the number of submodules $N\supseteq zM$ of codimension $j$, that is, the number $\qbinom sj_q$ of subspaces of codimension $j$ of $M/zM$. By the $q$-binomial theorem we have:
\[
\sum_{j=0}^s(-1)^jq^{\binom j2}\qbinom sj_q=\prod_{i=0}^{s-1}(1-q^i),
\]
and this product vanishes for $s\ge1$. Thus we obtain the following identity in $\sH[[x]]$:
\[
\Big(\sum_{j\ge0}(-1)^jq^{\binom j2}u_{(1^j)}x^j\Big)\Big(\sum_{j'\ge0}v_{j'}x^{j'}\Big)=1.
\]
Therefore $v_j=(-1)^{j+1}q^{\binom j2}u_{(1^j)}+(\text{a polynomial in }u_{(1^i)},\ i<j)$. Thus the products $v_{\alpha_1}\cdots v_{\alpha_l}$, $\alpha\vdash N$, are related to the products $u_{(1^{\alpha_1})}\cdots u_{(1^{\alpha_l})}$ by a matrix which is triangular for the refinement order, with diagonal entries $\pm q^{\sum_i\binom{\alpha_i}2}$, and combining this with (c) we obtain (d).
\end{proof}

\subsection{Split modules}\label{sub:split}
For a composition $\mathbf n=(n_1,\dots,n_r)$ of $N$ and pairwise distinct elements $\eta_1,\dots,\eta_r\in tO$ we set in $\mathbf O=O^r$:
\[
\varpi=(t,\dots,t),\qquad\omega_\eta=(\eta_1,\dots,\eta_r),\qquad k_{ij}=\val(\eta_i-\eta_j).
\]
Then $F((\varpi))$ is the field $\cK$ embedded diagonally, and since $F((\varpi))[\omega_\eta]=\mathbf K$ by Lagrange interpolation, $(\varpi,\omega_\eta)$ is a framed pair. In this case $V=\bigoplus_i\cK^{n_i}=\cK^N$ has dimension $m=N$ over $K_\varpi=\cK$, the lattices are the $O$-lattices in $\cK^N$, and $\omega_\eta$ acts as the diagonal matrix $\diag(\eta_i1_{n_i})$.

\begin{proposition}\label{prop:split}
For every $\mu\vdash N$ the refined count of the pair $(\varpi,\omega_\eta)$ is:
\[
B_{(\varpi,\omega_\eta)}(\mu;\mathbf n)=q^{\sum_{i<j}n_in_j(k_{ij}-1)}\ \#\big(\fn_{\mathbf n}(F)\cap\bO_\mu\big).
\]
\end{proposition}
\begin{proof}
Let $U\subseteq\GL_N(\cK)$ be the group of lower block unipotent matrices for the block sizes $\mathbf n$, and let $U(O)=U\cap\GL_N(O)$. By the Iwasawa decomposition $\GL_N(\cK)=ZU\GL_N(O)$, every $Z$-orbit of lattices meets $\{cO^N:c\in U\}$. If $zc_1O^N=c_2O^N$ with $z\in Z$ and $c_i\in U$, then $c_2^{-1}zc_1$ is an integral invertible lower block triangular matrix whose diagonal blocks are those of $z$, so these blocks are integral with unit determinants, and $z\in Z_0$. Thus the $Z$-orbits of lattices correspond to the orbits of $Z_0$ acting by conjugation on $U/U(O)$, with the same stabilizers. Hence, for a $Z$-stable set $\mathcal S$ of lattices, $\mass_Z(\mathcal S)$ is the number of cosets $cU(O)$ with $cO^N\in\mathcal S$.

Let us write $c=1+C$ and $c^{-1}\omega_\eta c=\omega_\eta+X$, where $X$ is strictly lower block triangular. For blocks $i>j$ the identity $(1+C)X=\omega_\eta C-C\omega_\eta$ then reads:
\[
X_{ij}=(\eta_i-\eta_j)\,C_{ij}-\sum_{i>\ell>j}C_{i\ell}X_{\ell j}.
\]
In the order of increasing $i-j$ this is a triangular change of variables between $C$ and $X$, whose linear diagonal part is multiplication by $\eta_i-\eta_j\ne0$. Hence $C\mapsto X$ is a bijection of the space of strictly lower block triangular matrices over $\cK$, and it multiplies the additive Haar measure by $q^{-\sum_{i<j}n_in_jk_{ij}}$. We normalize this measure by giving the integral matrices volume $1$, and then it is the Haar measure of $U$ with $\vol(U(O))=1$.

The lattice $cO^N$ is stable under $\omega_\eta$ if and only if $\omega_\eta+X$ is integral, or equivalently $X$ is integral, and in this case the type of $\omega_\eta$ on $cO^N/t\,cO^N$ is the type of $X$ modulo $t$, because $\omega_\eta\equiv0$ modulo $t$. The set of integral $X$ with reduction in $\bO_\mu$ has volume $q^{-e_2(\mathbf n)}\#(\fn_{\mathbf n}(F)\cap\bO_\mu)$, hence the corresponding set of $C$ has volume:
\[
q^{\sum_{i<j}n_in_jk_{ij}-e_2(\mathbf n)}\#(\fn_{\mathbf n}(F)\cap\bO_\mu).
\]
This set is stable under right multiplication by $U(O)$, so its volume is its number of cosets. By the first paragraph and Lemma~\ref{lem:massE} this number is $B_{(\varpi,\omega_\eta)}(\mu;\mathbf n)$.
\end{proof}

\subsection{The functions \texorpdfstring{$f^*_\beta$}{f*} count chains of the transform}
\begin{lemma}\label{lem:formal}
Let $u,w$ be commuting invertible operators on a set, let $e,d_*\ge1$ be coprime integers with $A_0e+B_0d_*=1$, and let $\varpi=u^{A_0}w^{B_0}$ and $\theta=w^eu^{-d_*}$. Then $u=\varpi^e\theta^{-B_0}$ and $w=\varpi^{d_*}\theta^{A_0}$, and for every family $(\Lambda_j)_{j\in\Z}$ of subsets and every $l\ge1$ we have:
\[
\big(u\Lambda_j=\Lambda_{j+el}\ \text{and}\ w\Lambda_j=\Lambda_{j+d_*l}\ \text{ for all }j\big)\iff\big(\varpi\Lambda_j=\Lambda_{j+l}\ \text{and}\ \theta\Lambda_j=\Lambda_j\ \text{ for all }j\big).
\]
\end{lemma}
\begin{proof}
The two identities follow from $A_0e+B_0d_*=1$. If the left side of the equivalence holds, then $\varpi\Lambda_j=\Lambda_{j+l(A_0e+B_0d_*)}$ and $\theta\Lambda_j=\Lambda_{j+ed_*l-d_*el}$. Conversely, if the right side holds, then $u\Lambda_j=\varpi^e\Lambda_j=\Lambda_{j+el}$ and $w\Lambda_j=\Lambda_{j+d_*l}$.
\end{proof}

\begin{lemma}\label{lem:rigidchains}
Let $(u,w)$ be a balanced framed pair of type $(e,d_*)$ in $\mathbf O$ with transform $(\varpi,\omega)$ (Proposition~\ref{prop:transform}), and let $\mathbf n$, $n'$, $m=en'$, $k=d_*n'$ and $D_w$ be as in Lemma~\ref{lem:orbital}. Then for every composition $\beta$ of $n'$ we have:
\[
D_w\big(f^{*,k}_\beta\big)=\Ch_\beta(\varpi,\omega;\mathbf n).
\]
\end{lemma}
\begin{proof}
If $l$ denotes the number of parts of $\beta$, then by Lemma~\ref{lem:dictionary}(b) the number $D_w(f^*_\beta)$ is the $Z$-mass of the set of $u$-periodic chains of type $\beta^{[e]}$ with $w\Lambda_j=\Lambda_{j+d_*l}$ for all $j$. Since such a chain satisfies $u\Lambda_j=\Lambda_{j+el}$, Lemma~\ref{lem:formal} shows that these chains are the descending families of compact open $F$-subspaces $(\Lambda_j)$ with the following properties:
\[
\varpi\Lambda_j=\Lambda_{j+l}\ \text{ and }\ \theta\Lambda_j=\Lambda_j\ \text{ for all }j,\qquad\dim_F(\Lambda_{j-1}/\Lambda_j)=\beta_j\ \text{ for }1\le j\le l.
\]
In such a family $\varpi\Lambda_j=\Lambda_{j+l}\subseteq\Lambda_j$.
Closure then gives stability under $F[[\varpi]]$, so each
$\Lambda_j$ is an $F[[\varpi]]$-lattice. The condition $\theta\Lambda_j=\Lambda_j$ is equivalent to $\omega\Lambda_j\subseteq\Lambda_j$, because $\omega=c^{-1}\theta-1$ and, conversely, $\theta^{\pm1}$ are convergent power series in the topologically nilpotent element $\omega$. Thus these families are exactly the $\varpi$-periodic chains of type $\beta$ of $F[[\varpi,\omega]]$-lattices, and since the group $Z$ and its measure are the same for both framed pairs, the lemma follows.
\end{proof}

\begin{proposition}\label{prop:indep}
Let $m\ge1$, $k\in\Z$ and $n'=\gcd(m,k)$. There are invariant distributions $D_\beta$, $\beta\vdash n'$, on $G=\GL_m(\K)$, supported on compact-mod-center elements of determinant valuation $k$, such that the matrix $\big(D_\beta(f^{*,k}_{\beta'})\big)_{\beta',\beta\vdash n'}$ is invertible.
\end{proposition}
\begin{proof}
Multiplication by the central element $\pi=\tau^m$ maps $f^{*,k}_\beta$ to $f^{*,k+m}_\beta$, and it transports invariant distributions on $G_c\cap G^{\val=k}$ to those on $G_c\cap G^{\val=k+m}$. Thus we may assume that $d_*=k/n'\ge1$, and we put $e=m/n'$ and fix $A_0,B_0$ with $A_0e+B_0d_*=1$.

For $\beta=(\beta_1\ge\dots\ge\beta_r)\vdash n'$ let us consider the split pair $(\varpi,\omega_\beta)$ of Section~\ref{sub:split} in $O^r$ with $\eta_i=t^i$ and ranks $\mathbf n=\beta$. Let $(U,W)=(\varpi^e\theta^{-B_0},\varpi^{d_*}\theta^{A_0})$, where $\theta=1+\omega_\beta$, be the balanced pair of type $(e,d_*)$ with transform $(\varpi,\omega_\beta)$ given by Lemma~\ref{lem:invtransform}. Here $g_i=1$, $\sum_in_ig_i=n'$ and $\dim_{K_U}V=m$, and we let $D_\beta=D_W$ be the distribution of Lemma~\ref{lem:orbital}. Since $k_{ij}=\val(t^i-t^j)=i$ for $i<j$, Lemmas~\ref{lem:rigidchains} and~\ref{lem:chainflag} and Proposition~\ref{prop:split} give:
\[
D_\beta(f^*_{\beta'})=\Ch_{\beta'}(\varpi,\omega_\beta;\beta)=\sum_{\mu\vdash n'}H_{\beta'}(\mu)\ q^{c_\beta}\,\#\big(\fn_\beta(F)\cap\bO_\mu\big),\qquad c_\beta=\sum_{i<j}\beta_i\beta_j(i-1).
\]
Thus the matrix of the proposition is the product:
\[
\sM_{n'}\cdot\sN\cdot\diag(q^{c_\beta}),
\]
where $\sN=\big(\#(\fn_\beta(F)\cap\bO_\mu)\big)_{\mu,\beta}$. The matrix $\sM_{n'}$ is invertible by Lemma~\ref{lem:hallbasics}(d), and by Lemma~\ref{lem:hallbasics}(a),(b) the matrix $\sN$ with its rows indexed by $\mu'$ is triangular for the dominance order with nonzero diagonal entries.
\end{proof}

\begin{corollary}\label{cor:basis}
Let $m\ge1$, $k\in\Z$ and $n'=\gcd(m,k)$.
\begin{enumerate}[label=\textup{(\roman*)}]
\item We have $\dim J_k=p(n')$, and the linear forms $D\mapsto D(f^{*,k}_\beta)$, $\beta\vdash n'$, are a basis of $J_k^*$. Hence for every $f\in\cH^{\val=k}$ there are unique numbers $c_\beta(f)$, $\beta\vdash n'$, with the following property:
\[
D(f)=\sum_{\beta\vdash n'}c_\beta(f)\,D(f^{*,k}_\beta)\qquad\text{for all }D\in J_k .
\]
\item Let $f$ be one of the following functions: $T_w$ with $\val(w)=k$, $F_{\alpha,k}$ with $\alpha$ a composition of $m$, or $f^{*,k}_{\beta'}$ with $\beta'$ a composition of $n'$. Then $c_\beta(f)=c_{\beta,f}(q)$ for a rational function $c_{\beta,f}\in\Q(X)$ which depends only on $(m,k,\beta)$ and on $w$, $\alpha$, $\beta'$ respectively, and not on the field $\F_q$.
\end{enumerate}
\end{corollary}
\begin{proof}
(i) By Theorem~\ref{thm:H}(i) we have $\dim J_k\le p(n')$, while the restrictions of the distributions $D_\beta$ of Proposition~\ref{prop:indep} belong to $J_k$, and on their span the $p(n')$ evaluations at the $f^*_\beta$ are linearly independent. Hence $\dim J_k=p(n')$, and these evaluations are a basis of $J_k^*$.

(ii) By Theorem~\ref{thm:H}(ii) we have $D(T_w)=\sum_\rho a_{w,\rho}(q)D(T_{w_\rho})$ for $D\in J_k$, with $a_{w,\rho}\in\Z[X]$, and by \eqref{eq:testfunctions} the same holds for $F_{\alpha,k}$ and $f^*_{\beta'}$ with coefficients in $\Q(X)$. For $\beta'\vdash n'$ let us write accordingly:
\[
D(f^*_{\beta'})=\sum_\rho\sP_{\beta'\rho}(q)D(T_{w_\rho}).
\]
For every prime power $q$ the $p(n')$ evaluations at the $T_{w_\rho}$ span $J_k^*$, hence they are a basis, and by (i) the evaluations at the $f^*_{\beta'}$ are another basis. Thus $\sP(q)$ is invertible for every prime power $q$, and $\det\sP\ne0$ in $\Q(X)$. Finally, the row vector $(c_\beta(f))_\beta$ is $(a_\rho(f))_\rho\,\sP^{-1}$.
\end{proof}

\subsection{The transfer theorem}
The chain masses of a balanced pair are related to those of its transform as follows.

\begin{theorem}\label{thm:transfer}
Let $(u,w)$ be a balanced framed pair of type $(e,d_*)$ in $\mathbf O=O^r$ with $\val_i(u)=g_ie$, and let $(\varpi,\omega)$ be its transform. Let $\mathbf n=(n_1,\dots,n_r)$ be a tuple of positive integers, and let $n'=\sum_in_ig_i$, $m=en'$ and $k=d_*n'$. Then for every composition $\alpha$ of $m$ we have:
\[
\Ch_\alpha(u,w;\mathbf n)=\sum_{\beta\vdash n'}c_{\alpha\beta}\ \Ch_\beta(\varpi,\omega;\mathbf n),\qquad c_{\alpha\beta}=c_\beta(F_{\alpha,k}),
\]
with the coefficients of Corollary~\ref{cor:basis} for $(m,k)$. The coefficients $c_{\alpha\beta}$ depend only on $(e,d_*,n',\alpha,\beta)$ and $q$, and they are the values at $q$ of rational functions which do not depend on $\F_q$.
\end{theorem}
\begin{proof}
By Lemma~\ref{lem:orbital} we have $\Ch_\alpha(u,w;\mathbf n)=D_w(F_{\alpha,k})$ with $D_w\in J_k$ and $\gcd(m,k)=n'$, and Corollary~\ref{cor:basis} gives:
\[
D_w(F_{\alpha,k})=\sum_\beta c_\beta(F_{\alpha,k})D_w(f^*_\beta).
\]
Lemma~\ref{lem:rigidchains} gives $D_w(f^*_\beta)=\Ch_\beta(\varpi,\omega;\mathbf n)$, proving the theorem.
\end{proof}

For $\alpha=(m)$ the transfer theorem can be stated in terms of the following functions.

\begin{definition}\label{def:Phi}
For coprime $a,b\ge1$, $N\ge1$ and $\mu\vdash N$ we set:
\[
\Phi^{a,b}_{q,N}(\mu)=\sum_{\beta\vdash N}c_\beta\big(F_{(aN),bN}\big)\,H_\beta(\mu),
\]
where the coefficients are those of Corollary~\ref{cor:basis} for $(m,k)=(aN,bN)$. Here $(aN)$ is the composition with one part, and thus $F_{(aN),bN}$ is the characteristic function of $\Mat_{aN}(\fo)\cap G^{\val=bN}$ divided by $\vol(\GL_{aN}(\fo))$.
\end{definition}

\begin{corollary}[Orbit-sum identity]\label{cor:orbitsum}
In the situation of Theorem~\ref{thm:transfer} we have:
\[
\#E_{F[[u,w]]}(\mathbf O^{\mathbf n})=\sum_{\mu\vdash n'}B_{(\varpi,\omega)}(\mu;\mathbf n)\ \Phi^{e,d_*}_{q,n'}(\mu).
\]
\end{corollary}
\begin{proof}
We apply Theorem~\ref{thm:transfer} for $\alpha=(m)$, and then Lemma~\ref{lem:chainflag} to both sides of the resulting identity.
\end{proof}

\subsection{Equisingular invariance and polynomiality}
As a first application of the transfer theorem we prove the equisingular invariance of the chain masses.

\begin{theorem}\label{thm:E}
Let $T\subseteq O$ be a plane branch over $\F_q$, in any characteristic. Let $u\in\fm_T\setminus\{0\}$ satisfy $T=F[[u,w]]$ for some $w\in\fm_T$, and put $D=\val u$. Then for every $n\ge1$ and every composition $\alpha$ of $Dn$ the chain mass $\Ch_\alpha(T,u;n)$ depends only on $q$, $n$, $\alpha$, $D$ and $\ms(T)$. More precisely, there is a rational function $R\in\Q(X)$, depending only on $(n,\alpha,D,\ms(T))$, such that $\Ch_\alpha(T,u;n)=R(q)$, and the same holds for the refined counts $B_{(u,w)}(\mu;n)$, $\mu\vdash Dn$. In particular, $\#E_T(O^n)$ depends only on $(q,n,\ms(T))$.
\end{theorem}
\begin{proof}
We use induction on the number $N(T)$ of quadratic transforms needed
to reach $O$. The rank $n$ stays fixed, but the framed valuation $D$
is part of the induction data and may change at each step.

\emph{Step 1. The smooth branch.}
If $T=O$, its only normalized lattice is $O^n$. By
Lemma~\ref{lem:chainflag},
\[
\Ch_\alpha(O,u;n)=H_\alpha((D^n)),
\]
which is a polynomial in $q$.

\emph{Step 2. Choose a generator of minimal valuation.}
Suppose $T\ne O$ and write $d=\mult(T)$. Since $u,w$ generate the
maximal ideal, $d=\min\{\val u,\val w\}$. If $\val w>D$, replace
$w$ by $w+u$. This preserves $T$ and the chain masses, and now
$\val w=d$. Set
\[
g=\gcd(D,d),\qquad e=D/g,\qquad d_*=d/g.
\]
Proposition~\ref{prop:transform} constructs the transformed frame
$(T_1;\varpi,\omega)$ with $\val\varpi=g$. Its multiplicity sequence
is obtained by deleting the initial block $g\,\eu(e,d_*)$ from
$\ms(T)$. Equivalently, $T_1=T^{(\ell)}$ for
$\ell=\ell(e,d_*)\ge1$. Thus $N(T_1)<N(T)$, regardless of whether
the framed valuation itself decreases.

\emph{Step 3. Transfer the full chain vector.}
Theorem~\ref{thm:transfer} gives
\[
\Ch_\alpha(T,u;n)=
\sum_{\beta\vdash gn}c_{\alpha\beta}\,
                 \Ch_\beta(T_1,\varpi;n).
\]
The coefficients are values at $q$ of rational functions determined
by $(D,d,n,\alpha,\beta)$. The induction hypothesis applies to the
frame $(T_1;\varpi,\omega)$ with framed valuation $g$ and the remaining
multiplicity sequence. This proves the asserted rationality and
invariance of every chain mass.

\emph{Step 4. Recover residual types and the total count.}
Restricting to partition indices $\alpha\vdash Dn$,
Lemma~\ref{lem:chainflag} expresses the resulting chain vector as
$\sM_{Dn}B$, where $B$ is the vector of refined counts. The flag
matrix $\sM_{Dn}$ has polynomial entries and is invertible by
Lemma~\ref{lem:hallbasics}(d), so the same assertion holds for $B$.
Finally, choose $u$ of valuation $\mult(T)$ and take the one-part
composition $\alpha=(Dn)$. Its chain mass is $\#E_T(O^n)$.

\end{proof}

To pass from rational functions to polynomials we need the following elementary lemma.

\begin{lemma}\label{lem:ratpoly}
If $R\in\Q(X)$ takes integer values on an unbounded set of positive integers, then $R\in\Q[X]$.
\end{lemma}
\begin{proof}
Let us write $R=P+S$ with $P\in\Q[X]$ and $S(x)\to0$ as $x\to\infty$, and choose $d\ge1$ with $dP\in\Z[X]$. At the given integers $x$ the number $dS(x)=dR(x)-dP(x)$ is an integer whose absolute value is less than one for large $x$, thus $S$ has infinitely many zeros and $S=0$.
\end{proof}

Integer values do not by themselves imply integer coefficients.
For the second step we use the fact that one fixed finite-type
scheme realizes all extension-field counts.

\begin{lemma}\label{lem:completion-integral}
Let $X$ be a separated scheme of finite type over $k=\F_{q_0}$.
Suppose $P\in\Q[T]$ satisfies
\[
 \#X(\F_{q_0^r})=P(q_0^r)\qquad\text{for every }r\ge1.
\]
Then $P\in\Z[T]$.
\end{lemma}
\begin{proof}
By Dwork's theorem \cite{Dwork} we have
\[
 Z_X(z)=\exp\left(\sum_{r\ge1}\#X(\F_{q_0^r})\frac{z^r}{r}\right)
 \in\Q(z):
\]
Dwork treats affine varieties, and the general case follows because
$Z_X$ is multiplicative over a finite stratification of $X$ by
locally closed affine subschemes. Write
$P(T)=\sum_{i=0}^d c_iT^i$. Formal differentiation gives
\[
 \frac{Z_X'(z)}{Z_X(z)}
   =\sum_{i=0}^d\frac{c_iq_0^i}{1-q_0^iz}.
\]
The points $z=q_0^{-i}$ are distinct. The residue of the right side at
$z=q_0^{-i}$ is $-c_i$. The logarithmic derivative of a nonzero rational
function has integral residues: locally write
$Z_X(z)=(z-a)^mU(z)$ with $m\in\Z$ and $U(a)\ne0$.
Thus $-c_i=m\in\Z$ for every $i$.
\end{proof}

\begin{corollary}\label{cor:poly}
For every multiplicity sequence $\underline\mu$ of a plane branch over a finite field and every $n\ge1$ there is a polynomial $P_{\underline\mu,n}\in\Z[X]$ such that $\#E_A(O^n)=P_{\underline\mu,n}(q)$ for every plane branch $A$ over every finite field $\F_q$ with $\ms(A)=\underline\mu$. The same holds for the refined counts and the chain masses of Theorem~\ref{thm:E}.
\end{corollary}
\begin{proof}
By Theorem~\ref{thm:E} there is $R\in\Q(X)$, depending only on $(\underline\mu,n)$, such that $\#E_A(O^n)=R(q)$ for all such $A$ and $q$. We may assume that a plane branch $A_0=\F_{q_0}[[x,y]]$ with $\ms(A_0)=\underline\mu$ exists, since otherwise there is nothing to prove. For every $r\ge1$ the ring $A_0\otimes_{\F_{q_0}}\F_{q_0^r}=\F_{q_0^r}[[x,y]]\subseteq\F_{q_0^r}[[t]]$ is then a plane branch, and since its quadratic transforms are obtained from those of $A_0$ by extension of scalars (Lemma~\ref{lem:qt-intrinsic}), it has the same multiplicity sequence. Hence $R(q_0^r)\in\Z$ for all $r\ge1$, and $R$ is a polynomial
with rational coefficients by Lemma~\ref{lem:ratpoly}. The same
argument applies to the refined and chain counts of Theorem~\ref{thm:E}.

Proposition~\ref{cor:completion-integral-main} constructs, independently
of polynomiality, a finite-type scheme over $\F_{q_0}$ for each of
these counts. Its number of $\F_{q_0^r}$-points is $R(q_0^r)$ for
every $r$, so Lemma~\ref{lem:completion-integral} gives $R\in\Z[X]$.
\end{proof}

\begin{remark}\label{rem:poly}
(1) The recursion of Theorem~\ref{thm:E} retains the full chain
vector, the remaining multiplicity sequence, the rank $n$, and the
framed valuation $D=\operatorname{val}(u)$, which need not equal the
multiplicity. Its terminal vector is
$\bigl(H_\alpha((D^n))\bigr)_{\alpha\vdash Dn}$.
Section~\ref{subsec:explicit-counting-algorithm} gives the algorithm
and matrix conventions. The scalar functions $\Phi^{a,b}_{q,N}$
computed in Section~\ref{sec:hall} give a shorter calculation when
the transformed residual counts are known.

(2) For $n=1$, normalizing the least valuation of a lattice
to zero gives a set of representatives for the $t^{\mathbb Z}$-action
on the affine Springer lattice points. This gives the local
compactified-Jacobian point count, with a fixed-degree convention:
if $r=\dim_F(O/L)$, then $L\mapsto t^{\delta_A-r}L$
bijects $E_A(O)$ with the $A$-lattices in $O$ of colength
$\delta_A$. Indeed, a unit $f\in L$ gives $fA\subseteq L$,
so $r\le\delta_A$, and normalization by the least valuation
gives the inverse.

For tamely ramified branches the polynomiality, with nonnegative
coefficients, is \cite[Cor.~1.11 and Cor.~8.3]{KT}; the global
statement for compactified Picard schemes is
\cite[\S8.2, Thm.~8.12]{KT}. Proposition~\ref{prop:rankone-positive}
below extends the nonnegativity to all plane branches. The
relation of Theorem~\ref{thm:E} with \cite{KT} and \cite{Tsai} is
described in Section~\ref{sub:method}; here
Lemmas~\ref{lem:formal}--\ref{lem:rigidchains} interpret the test
functions of \cite[\S2.3]{KT} for arbitrary framed generators, and
Lemma~\ref{lem:orbital} permits arbitrary module rank and
inseparable extensions.

(3) Corollary~\ref{cor:poly} proves integral coefficients.
Coefficient nonnegativity for the total count is Conjecture~P.
It is proved in rank one (Proposition~\ref{prop:rankone-positive})
and for torus sequences (Proposition~\ref{prop:torus-positive}), and
it fails for the refined counts (Example~\ref{ex:refined-negative});
this is why the conjecture concerns the total count only.
\end{remark}

\begin{proposition}[Conjecture~P in rank one]\label{prop:rankone-positive}
For every multiplicity sequence $\underline\mu$ of a plane branch over
a finite field we have $P_{\underline\mu,1}\in\Z_{\ge0}[X]$.
\end{proposition}
\begin{proof}
By Corollary~\ref{cor:poly}, $P_{\underline\mu,1}\in\Z[X]$ is a
single polynomial that counts every plane branch with multiplicity
sequence $\underline\mu$, over every finite field. We will determine
this polynomial from a realization in sufficiently large
characteristic and its extensions of constants.

\emph{Step 1: realization over another field.} We claim that for
every framed branch $(T;u,w)$ over a finite field $F$ and every finite field $F'$
there is a framed branch $(\widetilde T;\tilde u,\tilde w)$ over $F'$ with $\ms(\widetilde T)=\ms(T)$
and $\val \tilde u=\val u$. We use induction on $\delta_T$. If $T=O$, take
$(\tilde u,\tilde w)=(t^{\val u},t)$.

Let $T\ne O$, let $(e,d_*)$ be the type and
$g$ the weight of $(u,w)$, and let $(\varpi,\omega)$ be its transform,
with ring $T_1$. By Proposition~\ref{prop:transform}(ii) we have
$\dim_FT_1/T=\sum_{m\in\eu(e,d_*)}\binom{gm}2$. This sum vanishes only
if $g=1$ and $\min(e,d_*)=1$, that is, only if $u$ or $w$ has
valuation $1$, which is excluded by $T\ne O$. Hence
$\delta_{T_1}<\delta_T$.

If $\omega\ne0$, the induction hypothesis
applied to the framed branch $(T_1;\varpi,\omega)$ gives a framed
branch $(\widetilde T_1;\tilde\varpi,\tilde\omega)$ over $F'$ with $\ms(\widetilde T_1)=\ms(T_1)$
and $\val\tilde\varpi=g$. If $\omega=0$, then $\cK=F((\varpi))$, hence
$g=1$ and $T_1=O$, and we take the framed pair
$(\tilde\varpi,\tilde\omega)=(t,0)$, with ring $\widetilde T_1=F'[[t]]$. In both cases
Lemma~\ref{lem:invtransform}, with the same $(e,d_*)$ and with $c=1$,
gives a balanced framed pair $(\tilde u,\tilde w)$ over $F'$ of type $(e,d_*)$
and weight $g$ whose transform is $(\tilde\varpi,\tilde\omega)$. Then
$\tilde w\ne0$ and $\val \tilde u=ge=\val u$, and
Proposition~\ref{prop:transform}(iii), applied to both pairs, gives
$\ms(F'[[\tilde u,\tilde w]])=(g\cdot\eu(e,d_*),\ms(T_1))=\ms(T)$.

\emph{Step 2: a tame model.} Choose a plane branch
$A=\F_q[[x,y]]$ with $\ms(A)=\underline\mu$, and put
$\mu_0=\mult(A)$. For $A=O$ the polynomial is $1$, so let $A\ne O$. Then $x$ and $y$ are nonzero,
because $F[[x]]$ has finite colength in $O$ only if $\val x=1$, and
$(A;x,y)$ is a framed branch by Lemma~\ref{lem:gen}(c).

Fix a prime
$p>2\mu_0$. By Step~1 there is a plane branch $\widetilde A=\F_p[[\tilde x,\tilde y]]$ with
$\ms(\widetilde A)=\underline\mu$, and by Lemma~\ref{lem:qt-intrinsic} we may
exchange the generators so that $\val \tilde x=\mult(\widetilde A)=\mu_0$. For
$r\ge1$ put $F_r=\F_{p^r}$, $\widetilde A_r=F_r[[\tilde x,\tilde y]]$ and
$K_r=F_r((\tilde x))\subseteq\cK_r=F_r((t))$. By
Lemma~\ref{lem:qt-intrinsic}, $\ms(\widetilde A_r)=\underline\mu$.

The
extension $\cK_r/K_r$ is totally ramified of degree $\mu_0<p$, hence
tamely ramified and separable, and $\cK_r=K_r[\tilde y]$ by
Lemma~\ref{lem:gen}(c). Thus multiplication by $\tilde y$ is an element
$\gamma_r$ of $\Mat_{\mu_0}(F_r[[\tilde x]])$, for the lattice $F_r[[t]]$,
which is regular semisimple and elliptic, and $K_r(\gamma_r)=\cK_r$
is totally and tamely ramified over $K_r$: it is tamely ramified and
inertially elliptic in the sense of \cite[\S1.1]{KT}. The bound
$p>2\mu_0$ is the one required in \cite[Thm.~2.8]{KT}, on which the
Lie algebra statements of \cite{KT} rely; tameness alone needs
$p\nmid\mu_0$.

\emph{Step 3: the theorems of Kivinen and Tsai.} Let $k=\val \tilde y$. By
Lemmas~\ref{lem:orbital}(c) and~\ref{lem:chainflag}, applied to the
framed branch $(\widetilde A_r;\tilde x,\tilde y)$ with $n=1$, and by
Corollary~\ref{cor:poly},
\[
 D_{\tilde y}\bigl(F_{(\mu_0),k}\bigr)=\Ch_{(\mu_0)}(\widetilde A_r,\tilde x;1)
 =\#E_{\widetilde A_r}(F_r[[t]])=P_{\underline\mu,1}(p^r).
\]
Here $F_{(\mu_0),k}$ is the characteristic function of
$\Mat_{\mu_0}(F_r[[\tilde x]])\cap G^{\val=k}$ divided by the volume of
$\GL_{\mu_0}(F_r[[\tilde x]])$, and the centralizer $\cK_r^\times$ of
$\gamma_r$ carries the Haar measure for which $F_r[[t]]^\times$ has
volume one. Since every conjugate of $\gamma_r$ has determinant
valuation $k$, the left-hand side is the orbital integral
$I_{\gamma_r}(\mathbf 1_{(\mu_0)})$ of \cite[Thm.~1.10]{KT}, with the
measures of \cite[Def.~2.1]{KT}. By \cite[Thm.~1.10]{KT}, proved
there as Theorem~6.48, it equals the Hall pairing
$\langle f_{\gamma_r},h_{(\mu_0)}\rangle$, and by
\cite[Thm.~1.6]{KT}, proved there as Corollary~6.22,
\[
 f_{\gamma_r}=q^{\Xi}\,\omega\bigl(f_{\vec p,\vec q}\bigr)\big|_{q\mapsto q^{-1}}
 \qquad\text{at }q=p^r .
\]
The Newton pairs $(\vec p,\vec q)$ and the integer $\Xi$ are defined
by the valuations of the differences of the conjugates of $\tilde y$ over
$K_r$ \cite[Defs.~1.2 and~6.12]{KT}. The extension $K_r/K_1$ is
unramified, whereas $\cK_1/K_1$ is totally ramified, so these
extensions are linearly disjoint. The minimal polynomial of
$\tilde y$ therefore remains irreducible over $K_r$, with the same
roots and difference valuations in a common separable closure.
Consequently $(\vec p,\vec q)$ and $\Xi$ are independent of $r$.
By \cite[Prop.~4.12 and Thm.~5.14]{KT},
$f_{\vec p,\vec q}=\sum_{\lambda\vdash\mu_0}c_\lambda e_\lambda$ with
$c_\lambda\in\Z_{\ge0}[q]$: each $c_\lambda$ is a sum of powers
$q^{\operatorname{area}}$ over tuples of lattice paths. Since
$\omega(e_\lambda)=h_\lambda$ and
$\langle h_\lambda,h_{(\mu_0)}\rangle=1$, we obtain
\[
 P_{\underline\mu,1}(p^r)=\sum_{\lambda\vdash\mu_0}
     p^{r\Xi}\,c_\lambda(p^{-r})\qquad\text{for all }r\ge1 .
\]
Put
\[
 L(X)=X^\Xi\sum_{\lambda\vdash\mu_0}c_\lambda(X^{-1})
       \in\Z_{\ge0}[X,X^{-1}].
\]
We have proved that $P_{\underline\mu,1}(p^r)=L(p^r)$ for every
$r\ge1$. Multiplying their difference by a sufficiently large power
of $X$ gives a polynomial with infinitely many distinct zeros, so
$P_{\underline\mu,1}=L$. Since Corollary~\ref{cor:poly} gives
$P_{\underline\mu,1}\in\Z[X]$, this identity shows that its
coefficients are nonnegative.

The polynomial $P_{\underline\mu,1}$ is the one fixed at the start
of the proof. In particular, for the original branch $A$ over
$\F_q$, Corollary~\ref{cor:poly} gives
\[
 \#E_A(\F_q[[t]])=P_{\underline\mu,1}(q)=L(q).
\]
Thus the polynomial identified using the tame model counts the
original branch as well, completing the proof.
\end{proof}

\begin{example}[A refined count with a negative coefficient]\label{ex:refined-negative}
Conjecture~P cannot be extended to the refined counts of
Theorem~\ref{thm:E}. Let $A=F[[t^4,t^7]]$, $u=t^4$, $w=t^7$ and
$n=1$, so that $\ms(A)=\eu(4,7)=(4,3)$. In the notation of
Section~\ref{subsec:explicit-counting-algorithm} below, the chain
vector is $\mathbf V(4,(4,3);1)=\mathsf C_{4,4}\mathsf C_{4,3}\,(1)$,
and multiplication by $\sM_4^{-1}$ (Step~4 of the proof of
Theorem~\ref{thm:E}) gives the refined counts:
\[
\begin{array}{c|ccccc}
\mu&(4)&(3,1)&(2,2)&(2,1,1)&(1^4)\\\hline
B_{(u,w)}(\mu;1)\vphantom{\big|}&q^9&3q^8+3q^7&2q^7+2q^6-q^5&3q^6+6q^5+4q^4+2q^3&q^3+2q^2+q+1
\end{array}
\]
Their sum is
$P_{(4,3),1}(q)=q^9+3q^8+5q^7+5q^6+5q^5+4q^4+3q^3+2q^2+q+1$, with
nonnegative coefficients.


\end{example}

\subsection{An explicit algorithm for the counting polynomial}
\label{subsec:explicit-counting-algorithm}

The algorithm follows the induction in Theorem~\ref{thm:E}. Affine
length reduction and a finite Hecke recursion compute the
coefficients $a_{w,\rho}(q)$, from which we form the transfer matrices.
Appendix~\ref{sec:algorithm-implementation} supplies the finite
support enumeration, block refinements, and the Hall calculation
for the terminal flags.

The input is a realized multiplicity sequence $\underline\mu$ and a
positive rank $n$; the output is $P_{\underline\mu,n}(q)$.
Throughout the computation $q$ is an indeterminate and arithmetic is
exact in $\Q(q)$. Represent partitions by decreasing tuples of positive
integers, in one fixed order at each size. Trailing $1$'s in a
multiplicity sequence may be omitted; they are understood to continue
indefinitely. A smooth sequence has output $1$.

\paragraph{\textbf{State and recursion.}}
For a framed branch $(T;u,w)$, retain the column vector
\[
 \mathbf V(D,\underline\nu;n)
 =\bigl(\Ch_\alpha(T,u;n)\bigr)_{\alpha\vdash Dn},
 \qquad D=\val(u),\quad\underline\nu=\ms(T).
\]
Partitions suffice as indices because the stable-flag count is
unchanged by permuting the parts of a composition: $H_\alpha(\mu)$
is the coefficient of $u_\mu$ in the product
$v_{\alpha_1}\cdots v_{\alpha_l}$ in the commutative algebra $\sH$
(Section~\ref{sub:hall} and the first sentence of the proof of
Lemma~\ref{lem:hallbasics}(d)). The data at one singular step are:
\begin{center}
\begin{tabular}{@{}ll@{}}
\toprule
Symbol & Meaning\\
\midrule
$n$ & Module rank over the normalization, unchanged by transfer\\
$D$ & Current framing valuation $\val(u)$\\
$d=\nu_0$ & Current branch multiplicity\\
$m=Dn$ & Ambient dimension before transfer\\
$g=\gcd(D,d)$ & Next framing valuation\\
$N=gn$ & Ambient dimension after transfer\\
\bottomrule
\end{tabular}
\end{center}
The frame $D$ is part of the state: it need not equal the current
multiplicity. In particular, the smooth terminal value is
\begin{equation}
 \mathbf V(D,\varnothing;n)_\alpha=H_\alpha((D^n)),
 \label{eq:alg-smooth}
\end{equation}
where $(D^n)$ has $n$ parts equal to $D$
(Lemma~\ref{lem:chainflag}).

At a singular step arrange $\val(w)=d$ as in Theorem~\ref{thm:E}
and put
\begin{equation}
 g=\gcd(D,d),\qquad e=D/g,\qquad d_*=d/g,\qquad L=\ell(e,d_*).
 \label{eq:alg-stage}
\end{equation}
Proposition~\ref{prop:transform} deletes the prefix $g\,\eu(e,d_*)$
of length $L$ from $\underline\nu$, leaving
$\underline\nu_{\mathrm{next}}$, and replaces $D$ by $g$.
For $m,k\ge1$ let $N=\gcd(m,k)$ and
\[
 \mathsf C_{m,k}(q)
 =\bigl(c_\beta(F_{\alpha,k})\bigr)_{
       \alpha\vdash m,\ \beta\vdash N}.
\]
The recursion and the requested output are
\begin{align}
 \mathbf V(D,\underline\nu;n)
   &=\mathsf C_{Dn,dn}(q)
       \mathbf V(g,\underline\nu_{\mathrm{next}};n),
       \label{eq:alg-recursion}\\
 P_{\underline\mu,n}(q)
   &=\mathbf V(\mu_0,\underline\mu;n)_{(\mu_0n)}.
       \label{eq:alg-output}
\end{align}
The transfer matrix has $p(Dn)$ rows and $p(gn)$ columns.
Each step removes at least one singular entry, so the recursion
terminates. The frame never increases, and all partition sizes are
at most $n\mu_0$, but the frame need not strictly decrease.
If $D=d$, one must still apply $\mathsf C_{Dn,Dn}$ and delete
the single entry $D$; this matrix need not be the identity.

\paragraph{\textbf{The cusp.}}
Take $A=F[[t^2,t^3]]$, $n=1$, and the initial state
$(D,\underline\nu;n)=(2,(2);1)$. With $u=t^2$, replace the other
generator by $w=t^2+t^3$, so that the balanced transform is
$\varpi=u=t^2$, $\omega=w/u-1=t$. The transformed branch is
smooth, while its frame still has valuation two. With both rows
and columns ordered as $(2),(1,1)$, the matrix construction below
gives
\[
 \mathsf C_{2,2}=
 \begin{pmatrix}q+2&-1\\2q+2&-1\end{pmatrix},
 \qquad
 \mathbf V(2,\varnothing;1)=\begin{pmatrix}1\\1\end{pmatrix}.
\]
Their product is $(q+1,2q+1)^{\mathsf T}$, and the one-part entry
is $P_{(2),1}(q)=q+1$.
For comparison, the canonical recursion for $(4,2,2)$ has states
\[
 (4;(4,2,2))\longrightarrow(4;(2,2))
                  \longrightarrow(2;\varnothing)
\]
and contracts $\mathsf C_{4n,4n}\mathsf C_{4n,2n}$ with the
column $H_\alpha((2^n))$. Sections~\ref{sub:rankone}
and \ref{sub:hecke46}
give the resulting rank-one polynomial and
its Hecke comparison; the symbolic calculation is recorded in
Section~\ref{subsec:alg-symbolic}.

\paragraph{\textbf{Computing the coefficients $a_{w,\rho}(q)$.}}
Fix $m,k\ge1$ and put $N=\gcd(m,k)$ and $e=m/N$. The transfer
matrix $\mathsf C_{m,k}$ is built below, in \eqref{eq:alg-A} and
\eqref{eq:alg-M}, from the rows of polynomials
\[
 \mathbf a(w)=\bigl(a_{w,\rho}(q)\bigr)_{\rho\vdash N}
 \in\Z[q]^{p(N)},\qquad
 D(T_w)=\sum_{\rho\vdash N}a_{w,\rho}(q)D(T_{w_\rho})
 \quad(D\in J_k),
\]
of finitely many affine permutations $w$ with $\val(w)=k$. Here
$T_w=\vol(I)^{-1}\mathbf1_{IwI}$, the polynomials $a_{w,\rho}$ are
those of Theorem~\ref{thm:H}(ii), and $w_\rho=u_\rho\tau^k$ is the
representative defined before Theorem~\ref{thm:H}: $u_\rho$ is the
product of the descending Coxeter words on consecutive blocks of
sizes $\rho_i$ in $\{1,\ldots,N\}$. The row $\mathbf a(w)$ is
determined by $w$, because the evaluations $D\mapsto D(T_{w_\rho})$
form a basis of $J_k^*$ by Theorem~\ref{thm:H} and
Corollary~\ref{cor:basis}(i) (see Section~\ref{subsec:alg-finite-block}).

We compute $\mathbf a(w)$ by the three steps in the proof of
Theorem~\ref{thm:H} (Section~\ref{sub:proofH}). Retain previously
computed rows, starting with $\mathbf a(w_\rho)=\mathbf e_\rho$,
the unit row at $\rho$. For $m=1$ these are all the rows, since
$\tau^k=w_{(1)}$ is the only element of valuation $k$.
The bounds, block refinements and commutator identities used below
are proved in Sections~\ref{subsec:alg-discrete}
and~\ref{subsec:alg-finite-block}.

\emph{Terminal elements.}
Recall from Section~\ref{sub:reduction} that $v\approx w$ means that $v$ is
obtained from $w$ by a sequence of moves of the following two kinds:
\[
 x\longmapsto sxs\quad\text{for a simple reflection $s$ with }
 \ell(sxs)=\ell(x),\qquad
 x\longmapsto\tau^{\pm1}x\,\tau^{\mp1}.
\]
These moves preserve the length and the valuation, and
$T_v\equiv T_w$ modulo $[\cH,\cH]$, hence $\mathbf a(v)=\mathbf a(w)$.
We denote by $\mathcal C(w)$ the set of all $v\approx w$, and we say
that a simple reflection $s$ \emph{shortens} $x$ if
$\ell(sxs)<\ell(x)$. In this case the elements $sx$ and $sxs$ have
lengths $\ell(x)-1$ and $\ell(x)-2$ by Lemma~\ref{lem:IMlength}(ii),
and \eqref{eq:GP} gives:
\[
 \boxed{\ \mathbf a(x)=(q-1)\,\mathbf a(sx)+q\,\mathbf a(sxs).\ }
\]
The element $w$ is \emph{terminal} if no element of $\mathcal C(w)$
is shortened by a simple reflection.

For example, for $m=3$ no
simple reflection shortens $(6,1,2)$, but $s_1(6,1,2)s_1=(2,6,1)$ has
the same length and is shortened by $s_2$, which gives the step
$a(6,1,2)=a(2,6,1)$ of the coefficient examples below. The terminal
elements are described by Proposition~\ref{prop:HN3}, which
Stages~2 and~3 use.

\emph{Affine windows.}
Represent $w$ by the integer tuple $(w(1),\ldots,w(m))$, with distinct
residues modulo $m$, and extend by $w(i+m)=w(i)+m$. Write
\[
 w(i)=\sigma(i)+mz_i,\qquad \sigma\in S_m,\qquad
 \sum_i z_i=k,\qquad
 \ell(w)=\sum_{i<j}
 \left|\left\lfloor\frac{w(j)-w(i)}m\right\rfloor\right|.
\]
Products mean $(vw)(i)=v(w(i))$. Right multiplication by $s_i$,
$1\le i<m$, swaps entries $i,i+1$; for $s_0$ it replaces the first
and last entries by $w(m)-m,w(1)+m$. Left multiplication applies
the periodic permutation $s_i$ to each entry. Conjugation by $\tau$
has window
\[
 \tau w\tau^{-1}
 =\bigl(w(m)-m+1,w(1)+1,\ldots,w(m-1)+1\bigr).
\]

\emph{Stage 1: length reduction.}
Let us regard $\mathcal C(w)$ as the vertex set of a graph whose
edges are the moves above. This graph is finite: if $\ell(w)=L$,
then by Section~\ref{subsec:alg-discrete} the shifts of every vertex
lie in the box
\[
 \left\lceil\frac{k-(m-1)(L+1)}m\right\rceil
 \le z_i\le
 \left\lfloor\frac{k+(m-1)(L+1)}m\right\rfloor .
\]
Explore the graph, for instance breadth-first, starting at $w$.
At each vertex $x$, check whether its row is already known and
compute $sxs$ and its length for every simple reflection $s$.
Add the vertices $sxs$ with $\ell(sxs)=\ell(x)$, together with
$\tau x\tau^{-1}$ and $\tau^{-1}x\tau$, to the list of discovered
vertices if they have not appeared before. The search ends in one
of three ways:
\begin{enumerate}[label=\textup{(\alph*)}]
\item an examined vertex $x$ has a known row, and then
$\mathbf a(w)=\mathbf a(x)$;
\item a simple reflection $s$ shortens an examined vertex $x$, and
then $\mathbf a(w)=\mathbf a(x)$ is given by the boxed formula, where
the rows of $sx$ and $sxs$ are computed by calling the routine on
these elements, of lengths $\ell(w)-1$ and $\ell(w)-2$;
\item all discovered vertices have been examined without (a) or (b);
then they form the whole set $\mathcal C(w)$, no element of which is
shortened, so that $w$ is terminal, and the routine passes to
Stages~2 and~3.
\end{enumerate}

Assign the resulting row to every discovered vertex, since all lie
in $\mathcal C(w)$. Every recursive call in (b) has smaller length,
so Stage~1 ends after at most $\ell(w)$ levels.

\emph{Stage 2: terminal elements with non-central Newton point.}
For each cycle $c$ of $\sigma$ put $r_c=|c|$ and
$s_c=\sum_{i\in c}z_i$. By Lemma~\ref{lem:classes}, where
$z_i=\lambda_{\sigma(i)}$, the Newton point of $w$ has the
coordinate $s_c/r_c$ at $i\in c$, so centrality is the integer test
\[
 m s_c=k r_c\qquad\text{for every cycle }c.
\]
If this fails for a terminal element, return $\mathbf a(w)=0$: by
Proposition~\ref{prop:HN3} and Lemma~\ref{lem:cert} we have
$D(T_w)=0$ for every $D\in J_k$, as in Step~2 of the proof of
Theorem~\ref{thm:H}. This vanishing is used after evaluation on
$J_k$; it is not a zero identity in the full affine cocenter. In
particular, the Newton test must not discard a nonterminal input:
the element $(6,2,1)$ of the coefficient examples has a non-central
Newton point and a nonzero row.

\emph{Stage 3: terminal elements with central Newton point.}
In this case, search the vertices $v$ of $\mathcal C(w)$ and the
integers $0\le p<N$. Put $u=v\tau^{-k}$ and test:
\[
 \boxed{\quad
 \left\lfloor\frac{u(x)-p-1}{N}\right\rfloor
 =\left\lfloor\frac{x-p-1}{N}\right\rfloor
 \quad(1\le x\le m).\quad}
\]
This says that $u$ preserves every interval
$\{p+tN+1,\ldots,p+(t+1)N\}$, $t\in\Z$.

A pair $(v,p)$ which passes
the test exists. Indeed, by Proposition~\ref{prop:HN3}, as in Step~3
of the proof of Theorem~\ref{thm:H}, some vertex has the form
$v=u\tau^k$ with $u\in W_K$, where $K$ is a proper set of simple
reflections stable under $\delta=\operatorname{Ad}(\tau^k)$. By
Section~\ref{sub:blocks} we have $K=K_\Sigma$ for a block structure
$\Sigma$, which is a nonempty set of separators stable under
$i\mapsto i+N$. Hence $\Sigma$ contains some $p$ with $0\le p<N$,
together with all of $p+N\Z$ modulo $m$, so every block of $\Sigma$
lies in one of the above intervals, and $u\in\prod_BS_B$ preserves
them. Thus only $N$ separator tests are needed at each vertex.

\emph{Folding to a finite Hecke algebra.}
Let $(v,p)$ pass the interval test, let $u=v\tau^{-k}$, and put
$B_j=\{p+jk+1,\ldots,p+jk+N\}$ for $0\le j<e$. Since
$\gcd(k/N,e)=1$, the intervals $B_j$ represent the $e$ classes modulo
$m$ of the intervals of the test, and $\delta=\operatorname{Ad}(\tau^k)$
maps $S_{B_j}$ onto $S_{B_{j+1}}$ for $j<e-1$ and $S_{B_{e-1}}$ onto
$S_{B_0}$. Hence $u=u_0u_1\cdots u_{e-1}$, where $u_j$ acts on
$B_j+m\Z$ as $u$ does and fixes all other integers, and these factors
commute. After translating $B_j$ by $-(p+jk)$, the factor $u_j$
becomes the permutation
\[
 \pi_j(r)=u(p+jk+r)-(p+jk),\qquad 1\le r\le N,
\]
of $\{1,\ldots,N\}$.

Let $H_N$ be the Hecke algebra of
Lemma~\ref{lem:finiteA} with $X$ specialized to $q$, and let
$\iota\colon H_N\to\cH$ be the homomorphism which maps $T_y$ to
$T_{\tau^py\tau^{-p}}$, where $y\in S_N$ is regarded as an element of
$\tW$ permuting $\{1,\ldots,N\}$. Thus $\iota(T_{s_i})=T_{s_{p+i}}$,
and $\iota$ identifies $H_N$ with $\cH_{B_0}$
(Section~\ref{sub:blocks}). Then $T_{u_j}=\delta^j(\iota(T_{\pi_j}))$,
and with $\phi(h)=h\vartheta$ as in Section~\ref{sub:blocks} we have
$T_w\equiv T_v=\phi(T_{u_0}T_{u_1}\cdots T_{u_{e-1}})$.

We move the
factor at $B_{e-1}$ into $B_0$ by Lemma~\ref{lem:twist}(i), where it
is multiplied on the right, then the factor at $B_{e-2}$ through
$B_{e-1}$ into $B_0$, and so on, as in the proof of
Lemma~\ref{lem:twist}(ii). Since $\delta^e$ is the identity, this
gives:
\[
 T_w\equiv\phi(\iota(F)),\qquad
 \boxed{\quad F=T_{\pi_0}T_{\pi_{e-1}}\cdots T_{\pi_1}\in H_N,\quad}
\]
with $F=T_{\pi_0}$ if $e=1$. By Lemma~\ref{lem:twist}(iii) with
$g=1$, the class $\phi(\iota(h))$ depends only on $h$ modulo
$[H_N,H_N]$.

Suppose now that the finite routine below gives
\[
 F\equiv\sum_\eta c_\eta(q)\,T_\eta\pmod{[H_N,H_N]},
\]
where $\eta$ runs over the compositions of $N$ and $T_\eta$ is the
element of Lemma~\ref{lem:finiteA}. For each $\eta$ we have
$\iota(T_\eta)\vartheta=T_{u(\Sigma,R)\tau^k}$, where $\Sigma$ is the
block structure with separators $p,p+\eta_1,p+\eta_1+\eta_2,\ldots$
modulo $N$ and $R$ is the set of its blocks contained in $B_0$. By
Lemma~\ref{lem:swap} the class of this element is $[\Sigma_\rho]$,
the class of $T_{w_\rho}$, where $\rho=\bar\eta$ is the partition
formed by the parts of $\eta$. Therefore
$T_w\equiv\sum_\eta c_\eta(q)T_{w_{\bar\eta}}$ modulo $[\cH,\cH]$,
and since the row $\mathbf a(w)$ is determined by $w$, we obtain:
\[
 a_{w,\rho}(q)=\sum_{\eta:\ \bar\eta=\rho}c_\eta(q).
\]
This is Step~3 of the proof of Theorem~\ref{thm:H} for the block
structure $p+N\Z$ and $R=\{B_0\}$; the normalization introduces
no further power of $q$.

\emph{The finite routine.}
For $c\ge0$ let $H_c$ be the Hecke algebra of $S_c$ of
Lemma~\ref{lem:finiteA}, with $X=q$, and for a composition
$\eta=(\eta_1,\ldots,\eta_l)$ of $c$ let $T_\eta=T_{w_\eta}$ be the
element of that lemma: $w_\eta$ is the product of the descending
Coxeter words $s_{a+\eta_i-1}\cdots s_{a+1}$ of the consecutive
blocks $\{a+1,\ldots,a+\eta_i\}$, where $a=\eta_1+\cdots+\eta_{i-1}$.
For a partition $\rho$ of $N$, the element $T_\rho$ of $H_N$ is thus
$T_{u_\rho}$. The finite routine takes an element $h$ of $H_c$ with
coefficients in $\Z[q]$ and returns polynomials $c_\eta(q)$ such that:
\[
 h\equiv\sum_\eta c_\eta(q)\,T_\eta\pmod{[H_c,H_c]},
\]
where $\eta$ runs over the compositions of $c$. Following the proof
of Lemma~\ref{lem:finiteA}, we give the procedure for $h=T_v$,
$v\in S_c$, and extend it linearly. Applied to $F\in H_N$, it gives
the coefficients of the folding step.

Write each element of $H_c$ as a sum $\sum_v f_v(q)T_v$.
Multiplication uses the rule:
\[
 T_vT_{s_i}=
 \begin{cases}
 T_{vs_i},&v(i)<v(i+1),\\
 (q-1)T_v+qT_{vs_i},&v(i)>v(i+1).
 \end{cases}
\]
To multiply by $T_y$ for a general $y\in S_c$, sort the one-line
notation of $y$ to the identity by adjacent swaps that remove
inversions. If the recorded swaps give $ys_{i_1}\cdots s_{i_l}=1$,
then $y=s_{i_l}\cdots s_{i_1}$ is a reduced word, and we multiply by
$T_{s_{i_l}},\ldots,T_{s_{i_1}}$ in this order.

The proof of Lemma~\ref{lem:finiteA} uses the elements $E_1=1$ and
$E_j=T_{s_{c-1}}T_{s_{c-2}}\cdots T_{s_{c-j+1}}$, $2\le j\le c$, so
that $E_j=T_{(1,\ldots,1,j)}$ comes from the descending Coxeter word
of the last $j$ positions. Put $r=c-j+1$. Then $E_j$ commutes with
the subalgebra $H_{r-1}$ spanned by the $T_y$, $y\in S_{r-1}$, and
$E_jT_\mu=T_{(\mu,j)}$ for every composition $\mu$ of $r-1$.

A
\emph{state} of $H_c$ is a pair $(j,v)$ with $1\le j\le c$ and
$v\in S_r$, which stands for the element $E_jT_v$ of $E_jH_r$, and we
call $r$ its \emph{relative rank}. The routine computes the
coefficients $c_\eta(q)$ of every state it meets, beginning with the
state $(1,v)$, which stands for $T_v$. The space $E_jH_r$ is the left
side of the inclusion \eqref{eq:finiteA-step} in the proof of
Lemma~\ref{lem:finiteA}, with $n=c$ and $i=j-1$, and the following
two cases correspond to the two terms $E_jH_{r-1}$ and
$E_{j+1}H_{r-1}$ on its right side:
\begin{enumerate}[label=\textup{(\arabic*)}]
\item If $v(r)=r$, then $v\in S_{r-1}$. Apply the finite routine of
the algebra $H_{r-1}$ to $T_v$, and append the part $j$ to every
composition it returns. An expression of $T_v$ modulo
$[H_{r-1},H_{r-1}]$, multiplied by $E_j$, is an expression of
$E_jT_v$ modulo $[H_c,H_c]$, because $E_j[h_1,h_2]=[E_jh_1,h_2]$ for
$h_1,h_2\in H_{r-1}$ and $E_jT_\mu=T_{(\mu,j)}$. This is the
induction on the rank in the proof of Lemma~\ref{lem:finiteA}, see
\eqref{eq:finiteA-ind}. The routine of $H_0=\Z[q]$ returns the empty
composition with coefficient $1$, so the state $(c,1)$ gives the
composition $(c)$.
\item If $v(r)\ne r$, let $i=v^{-1}(r)$, so that $i<r$, and put:
\[
 d=s_{r-1}s_{r-2}\cdots s_i,\qquad a=vd^{-1},\qquad
 b=s_{r-2}\cdots s_i,
\]
with $b=1$ if $i=r-1$. By the proof of the Claim in
Lemma~\ref{lem:finiteA}, with $r$ in place of $j$, the permutation
$a$ fixes $r$, so that $a\in S_{r-1}$, and $T_v=T_aT_{s_{r-1}}T_b$.
Since $E_j$ commutes with $T_a$ and $E_jT_{s_{r-1}}=E_{j+1}$, we get:
\[
 E_jT_v=T_aE_{j+1}T_b=E_{j+1}T_bT_a+[T_a,E_{j+1}T_b].
\]
Expand $T_bT_a=\sum_zf_z(q)T_z$ by the multiplication rule, where
$z\in S_{r-1}$, and return the sum of the answers for the states
$(j+1,z)$, weighted by $f_z(q)$.
\end{enumerate}

Both cases lower the relative rank by one, so the recursion has at
most $c$ levels. Induction on that rank, using the two displayed
identities, proves that the returned expression is congruent to
$E_jT_v$ modulo $[H_c,H_c]$. All coefficients lie in $\Z[q]$.

The routine returns compositions, and these are sorted only after
the folding: the routine is applied to every term of the standard
basis expansion of $F$, each composition $\eta$ is replaced by the
partition $\bar\eta$ formed by its parts, and the coefficients of
equal partitions are added. By the folding
paragraph, which uses Lemma~\ref{lem:swap} for this sorting, the
result is the row $\mathbf a(w)$.

Inside $H_c$ the sorting is the
congruence $T_\eta\equiv T_{\eta'}$ for a rearrangement $\eta'$ of
$\eta$, which holds once $q$ is inverted, by the commutator identity
of Section~\ref{subsec:alg-finite-block}. It fails over $\Z[q]$:
at $q=0$ the one-dimensional representation of $H_3$ with
$T_{s_1}\mapsto-1$ and $T_{s_2}\mapsto0$ vanishes on $[H_3,H_3]$ and
separates $T_{(2,1)}=T_{s_1}$ from $T_{(1,2)}=T_{s_2}$.

\emph{Termination and correctness.}
The box of Stage~1, the decrease of the affine length in Stage~1
and the decrease of the relative rank in the finite routine prove
termination. The displayed Hecke identities prove correctness by
induction on the length, Proposition~\ref{prop:HN3} and
Lemma~\ref{lem:cert} justify the zero rows of Stage~2, and by
Theorem~\ref{thm:H} and Corollary~\ref{cor:basis} the completed row
does not depend on the choices made in the searches
(Section~\ref{subsec:alg-finite-block}). For efficiency, one may use all compatible
separators and reduce several smaller finite blocks, as explained
in Section~\ref{subsec:alg-discrete}.

\emph{Coefficient examples.}
The following examples apply the folding and the finite routine
directly to $w$. This is valid whenever $(w,p)$ passes the interval
test, even if $w$ is not terminal, in which case the coefficient
routine itself would begin with Stage~1. For $m=k=3$ and
$w=(6,5,4)=(s_1s_2s_1)\tau^3$ we have $N=3$ and $e=1$, the pair
$(w,0)$ passes the test, and $F=T_v$ with $v=\pi_0=(3,2,1)$.

The
finite routine starts at the state $(1,v)$ of $H_3$, where $r=3$ and
$v(3)=1$. Thus case~(2) applies with $i=1$, $d=s_2s_1$ and $a=b=s_1$,
and since $T_{s_1}^2=(q-1)T_{s_1}+q$, it leads to the states
$(2,s_1)$ and $(2,1)$ with the weights $q-1$ and $q$. At the state
$(2,s_1)$ case~(2) applies with $i=1$ and $a=b=1$ and leads to the
state $(3,1)$, which gives the composition $(3)$ by case~(1). At the
state $(2,1)$, case~(1) calls the routine of $H_1$, which returns the
composition $(1)$, and appends the part $2$. Thus:
\[
 T_{s_1s_2s_1}\equiv(q-1)T_{(3)}+qT_{(1,2)}
 =(q-1)T_{s_2s_1}+qT_{s_2}\pmod{[H_3,H_3]},
\]
which can also be seen directly from
$T_{s_1}\cdot T_{s_2}T_{s_1}\equiv T_{s_2}T_{s_1}\cdot T_{s_1}$. The
partition formed by the parts of $(1,2)$ is $(2,1)$, and therefore:
\[
 \mathbf a(6,5,4)=(q-1,q,0)
 \quad\text{in columns }(3),(2,1),(1,1,1).
\]

For $m=4$, $k=2$ and $w=(4,3,6,5)=(s_1s_3)\tau^2$ we have $N=e=2$
and $u=(2,1,4,3)$, the separator $p=0$ passes the test, and
$\pi_0=\pi_1=(2,1)$, so that $F=T_{s_1}^2=(q-1)T_{s_1}+q$. In $H_2$
the state $(1,s_1)$ leads by case~(2), with $i=1$ and $a=b=1$, to
the state $(2,1)$ and the composition $(2)$, while the state $(1,1)$
gives the composition $(1,1)$ by case~(1). Thus $T_{s_1}=T_{(2)}$
and $1=T_{(1,1)}$, and:
\[
 \mathbf a(4,3,6,5)=(q-1,q)
 \quad\text{in columns }(2),(1,1).
\]

For $e\ge3$ the cyclic order of the factors in $F$ matters. For
$m=9$, $k=3$ and $w=(5,4,6,8,9,7,10,12,11)$, the separator $p=0$
gives $\pi_0=s_2$, $\pi_1=s_1$ and $\pi_2=s_1s_2$, and since
$s_2s_1s_2s_1=s_1s_2$, the multiplication rule gives:
\[
 F=T_{s_2}T_{s_1s_2}T_{s_1}=(q-1)T_{s_1s_2s_1}+qT_{s_1s_2}.
\]
The first term was reduced in the case $m=k=3$. For
$v=s_1s_2=(2,3,1)$, case~(2) at the state $(1,v)$ applies with $i=2$,
$d=s_2$, $a=s_1$ and $b=1$, and it leads to the state $(2,s_1)$,
which gives $(3)$ as before. Hence $T_{s_1s_2}\equiv T_{(3)}$, so
that $F\equiv(q^2-q+1)T_{(3)}+q(q-1)T_{(1,2)}$ and:
\[
 \mathbf a(w)=(q^2-q+1,\ q(q-1),\ 0)
 \quad\text{in columns }(3),(2,1),(1,1,1).
\]
The product $T_{\pi_0}T_{\pi_1}T_{\pi_2}$ in the other cyclic order
would give $\bigl((q-1)^2,\ 2q(q-1),\ q^2\bigr)$ instead.

Some terminal classes $\mathcal C(w)$ contain no representative
$w_\rho$, so they require the block calculation. For $m=k=6$ the window
$w=(9,7,8,10,12,11)=(s_2s_1s_5)\tau^6$, the Coxeter element of the
blocks of sizes $(3,1,2)$, is terminal with central Newton point, and
$\mathcal C(w)$ contains no $w_\rho$. Indeed, $\tau^6$ is
central, a length-preserving conjugation by a simple reflection does
not change the support $\{s_1,s_2,s_5\}$, and conjugation by $\tau$
rotates it, whereas the support $\{s_1,s_2,s_4\}$ of $u_{(3,2,1)}$
is its mirror image and not a rotation of it. Here the separator
$p=0$ gives $F=T_{\pi_0}$ with $\pi_0=s_2s_1s_5=(3,1,2,4,6,5)$.

At the
state $(1,\pi_0)$ of $H_6$ case~(2) applies with $r=6$, $i=5$,
$d=s_5$, $a=(3,1,2,4,5)$ and $b=1$, and it leads to the state
$(2,a)$. There $a(5)=5$, so case~(1) calls the routine of $H_4$ on
$(3,1,2,4)$, which by case~(1) again calls the routine of $H_3$ on
$(3,1,2)$. The latter applies case~(2) with $i=1$, $a=1$ and $b=s_1$,
and it returns $(3)$ through the states $(2,s_1)$ and $(3,1)$.
Appending the parts $1$ and $2$ gives the composition $(3,1,2)$ with
coefficient $1$, agreeing with $T_{\pi_0}=T_{(3,1,2)}$, and
sorting gives $\mathbf a(w)=\mathbf e_{(3,2,1)}$.

The terminality requirement is
essential for the zero rule: when $m=3,k=1$, the window
$w=(6,2,1)$ has Newton point $(1/2,0,1/2)$, but
$a_{w,(1)}(q)=q^2(q-1)$. For example, writing $a(v)$ for this
single coefficient, length reduction gives
\[
\begin{aligned}
 a(6,2,1)&=(q-1)a(6,1,2)+q a(1,6,2),\\
 a(6,1,2)&=a(2,6,1)
          =(q-1)a(3,5,1)+q a(3,1,5),\\
 a(3,1,5)&=(q-1)a(3,2,4)+q a(2,3,4),\\
 a(1,6,2)&=(q-1)a(1,5,3)+q a(1,3,5).
\end{aligned}
\]
Here $a(2,3,4)=a(\tau)=1$; the other four terminal rows
$a(3,5,1)$, $a(3,2,4)$, $a(1,5,3)$ and $a(1,3,5)$ vanish.
The intermediate move $(6,1,2)\approx(2,6,1)$ also shows why
the search must include length-preserving conjugations.

\paragraph{\textbf{Constructing the transfer matrix.}}
With the rows $\mathbf a(w)$ now computed, put
$e=m/N$ and $v_\alpha(q)=\prod_i[\alpha_i]_q!$, where
$[r]_q!=\prod_{j=1}^r(1+q+\cdots+q^{j-1})$.
Let $\mathcal S_{\alpha,k}$ be the finite set of affine windows
carrying each lattice of the standard chain of type $\alpha$ into itself and having
determinant valuation $k$. The inequalities in \eqref{eq:alg-support} give a finite enumeration of this set.
Form
\begin{align}
 \mathsf A_{\alpha\rho}
   &=v_\alpha(q)^{-1}
        \sum_{w\in\mathcal S_{\alpha,k}}a_{w,\rho}(q),
       &&\alpha\vdash m,\quad\rho\vdash N,
       \label{eq:alg-A}\\
 \mathsf M_{\beta\rho}
   &=v_\beta(q)^{-e}
        \sum_{y\in W_{\beta^{[e]}}}a_{\tau^ky,\rho}(q),
       &&\beta,\rho\vdash N.
       \label{eq:alg-M}
\end{align}
Here $W_{\beta^{[e]}}$ is the product of symmetric groups on the
consecutive blocks of $\beta^{[e]}$. Corollary~\ref{cor:basis}
implies that $\mathsf M$ is invertible and that
\begin{equation}
 \mathsf C_{m,k}\mathsf M=\mathsf A,
 \qquad \mathsf C_{m,k}=\mathsf A\mathsf M^{-1}.
 \label{eq:alg-matrix}
\end{equation}
Thus solve $\mathsf M^{\mathsf T}\mathsf C_{m,k}^{\mathsf T}
=\mathsf A^{\mathsf T}$ over $\Q(q)$. Set
$\mathsf C_{1,k}=(1)$ directly.

The terminal flags are computed in
Section~\ref{subsec:alg-terminal-hall} using Gaussian polynomials
and the semisimple multiplication rule \eqref{eq:alg-hall-pieri}.
That section treats both stable flags $H_\alpha$ and flags $F_\eta$
lowered by the nilpotent operator.

\paragraph{\textbf{Computation by the Hall character.}}
When only the output $P_{\underline\mu,n}(q)=\#E_A(O^n)$ is needed,
the transfer from $A$ to the branch $T$ obtained by resolving its
first characteristic pair can be replaced by a linear system of size
$p(Dn)$ whose right side consists of torus counts.

Use the cabling
notation of Theorem~\ref{thm:cable}: $A=F[[x,y]]$ with $\val x=aD$,
$\val y=bD$ and $\gcd(a,b)=1$, the transformed framed branch
$(T;\varpi,\omega)$ of Corollary~\ref{cor:cabling-data} has
$\val\varpi=D$, and $N=Dn$. Here $D$ is the common source weight;
the source frame has valuation $aD$. For the frame $(x,y)$, of type $(a,b)$,
Theorem~\ref{thm:transfer} gives
$\mathbf V(aD,\ms(A);n)=\mathsf C_{aN,bN}\mathbf V(D,\ms(T);n)$, and
$\mathbf V(D,\ms(T);n)=\sM_NB_n^T$ by Lemma~\ref{lem:chainflag},
where $B_n^T(\lambda)$, $\lambda\vdash N$, is the number of
$\Lambda\in E_T(O^n)$ on which $\omega$ induces a nilpotent
endomorphism of $\Lambda/\varpi\Lambda$ of Jordan type $\lambda$.

Since $P_{\underline\mu,n}(q)$ is the entry of
$\mathbf V(aD,\ms(A);n)$ at $(aN)$, computing it requires only the row $(aN)$ of
$\mathsf C_{aN,bN}\sM_N$, which is the function of
Definition~\ref{def:Phi}:
\[
 \Phi^{a,b}_{q,N}(\lambda)
 =\sum_{\beta\vdash N}\bigl(\mathsf C_{aN,bN}\bigr)_{(aN),\beta}
  H_\beta(\lambda)\qquad(\lambda\vdash N),
\]
and $P_{\underline\mu,n}(q)=\sum_\lambda B_n^T(\lambda)
\Phi^{a,b}_{q,N}(\lambda)$, which is Theorem~\ref{thm:cable}.

Theorem~\ref{thm:hallchar} computes $\Phi^{a,b}_{q,N}$ from the
numbers $E_j=E_{a,b}(j)$, $j\le N$; these are lattice counts of the
torus branch, not the elements $E_j$ of the finite routine. Let
$\psi_\lambda=\varphi_{a,b}(u_\lambda)$ be the value of the character
$\varphi_{a,b}$ of that theorem on the basis element $u_\lambda$ of
the Hall algebra $\sH$ of Section~\ref{sub:hall}. By the definition
of $\varphi_{a,b}$ we have:
\[
 \psi_\lambda=\frac{q^{-ab\binom N2}\,\Phi^{a,b}_{q,N}(\lambda)}
                   {a_\lambda(q)},\qquad
 a_\lambda(q)=|\Aut(M_\lambda)|
 =q^{\sum_i(\lambda'_i)^2}
      \prod_{j\ge1}\prod_{r=1}^{m_j(\lambda)}(1-q^{-r}),
\]
where $m_j(\lambda)=\#\{i:\lambda_i=j\}$, so it suffices to compute
the numbers $\psi_\lambda$.

For a partition $\eta$ of $N$ let
$F_\eta(\lambda)$ be the coefficient of $u_\lambda$ in the product
$u_{(1^{\eta_1})}u_{(1^{\eta_2})}\cdots$, that is, the number of flags
$0=V_0\subseteq V_1\subseteq\cdots\subseteq F^N$ with
$\dim V_i/V_{i-1}=\eta_i$ and $XV_i\subseteq V_{i-1}$ for a fixed
nilpotent $X$ of Jordan type $\lambda$ (Section~\ref{sub:hall}). It is
computed by the Hall routine of
Section~\ref{subsec:alg-terminal-hall}. Since $\varphi_{a,b}$ is a
homomorphism with
$\varphi_{a,b}(u_{(1^j)})=q^{-ab\binom j2}E_j/|\GL_j(\F_q)|$,
applying it to this product gives:
\begin{equation}
 \sum_{\lambda\vdash N}F_\eta(\lambda)\psi_\lambda
 =\prod_i\frac{q^{-ab\binom{\eta_i}2}E_{\eta_i}}
                       {|\GL_{\eta_i}(\F_q)|}
 \quad(\eta\vdash N).
 \label{eq:alg-character-system}
\end{equation}
These are the identities of the last display of the proof of
Theorem~\ref{thm:hallchar}, $p(N)$ linear equations for the $p(N)$
unknowns $\psi_\lambda$. By Lemma~\ref{lem:hallbasics}(b) we have
$F_\eta(\lambda)\ne0$ only if $\lambda'\ctd\eta$, and
$F_\eta(\eta')=1$, so with the columns indexed by $\lambda'$ in an
order compatible with dominance the matrix is unitriangular, and the
system is solved by substitution. Then:
\begin{equation}
 \Phi^{a,b}_{q,N}(\lambda)
   =q^{ab\binom N2}a_\lambda(q)\psi_\lambda,
 \qquad
 P_{\underline\mu,n}(q)
   =\sum_{\lambda\vdash N}B_n^T(\lambda)
                       \Phi^{a,b}_{q,N}(\lambda).
 \label{eq:alg-character-output}
\end{equation}

For example, for $(a,b)=(2,3)$ and $N=2$ the matrix
$(F_\eta(\lambda))$ is the matrix $\mathsf F_2$ of
\eqref{eq:alg-small-hall}, and the two equations read
$\psi_{(1,1)}=q^{-6}E_2/|\GL_2(\F_q)|$ and
$\psi_{(2)}+(q+1)\psi_{(1,1)}=E_1^2/(q-1)^2$. They give
$\Phi^{2,3}_{q,2}((1,1))=E_2$ and
$\Phi^{2,3}_{q,2}((2))=(q^7E_1^2-E_2)/(q-1)$, the values used in
Section~\ref{subsec:alg-symbolic}.

If $B_n^T$ is not yet known, obtain its chain vector from
\eqref{eq:alg-recursion} and solve
$\sM_{Dn}B_n^T=\mathbf V(D,\ms(T);n)$.
The finite torus formula in Theorem~\ref{thm:torus}(i) provides
the $E_j$; an explicit enumeration is given in
Section~\ref{subsec:alg-torus}.
The general transfer recursion does not require this torus input.
Only finite torus counts enter the shortcut; the product formula for
$\cN_A$ concerns the subsequent rank limit.

\paragraph{\textbf{Cost.}}
The largest calculations are the support enumeration, matrix
inversion and rational arithmetic. Section~\ref{subsec:alg-checks}
gives bounds for the two supports and identities for checking the
calculation. The bounded searches in
Appendix~\ref{sec:algorithm-implementation} prove termination but
do not give a polynomial-time bound.

\paragraph{\textbf{A composition law.}}
The matrices $\mathsf C_{m,k}$ with $k>m$ are determined by those
with $k\le m$.

\begin{proposition}\label{prop:composition}
Let $e,d\ge1$ be coprime, $N\ge1$ and $j\ge0$. Then
\[
 \mathsf C_{eN,(d+je)N}=\mathsf C_{eN,eN}^{\,j}\;\mathsf C_{eN,dN}
\]
as matrices over $\Q(q)$ with rows indexed by the partitions of $eN$
and columns indexed by the partitions of $N$. In particular,
$\mathsf C_{N,hN}=\mathsf C_{N,N}^{\,h}$ for $h\ge1$.
\end{proposition}
\begin{proof}
Since $\gcd(e,d+ie)=1$ for all $i$, it suffices to treat $j=1$. Let
$(\varpi,\omega)$ be a framed pair in $\mathbf O=O^r$, and let
$\mathbf n=(n_1,\dots,n_r)$ be a tuple of positive integers with
$\sum_in_i\val_i(\varpi)=N$. Choose integers $A_0,B_0$ with
$A_0e+B_0(d+e)=1$, and put $\theta=1+\omega$ and
\[
 u=\varpi^e\theta^{-B_0},\qquad w=\varpi^{d+e}\theta^{A_0},\qquad
 w'=w/u=\varpi^{d}\theta^{A_0+B_0}.
\]
Recall that in Theorem~\ref{thm:transfer} the transform may be taken
with respect to any choice of the B\'ezout pair: another choice
multiplies $\varpi$ by a power of the unit $\theta$ of the
transformed ring (Proposition~\ref{prop:transform}(i)), which changes
neither this ring nor its periodic chains of lattices. We compute the
chain masses $\Ch_\alpha(u,w;\mathbf n)$, $\alpha\vdash eN$, in two
ways.

(a) By Lemma~\ref{lem:invtransform}, $(u,w)$ is a balanced framed
pair of type $(e,d+e)$ with transform $(\varpi,\omega)$, and
Theorem~\ref{thm:transfer} for $(m,k)=(eN,(d+e)N)$ gives
\[
 \Ch_\alpha(u,w;\mathbf n)=\sum_{\beta\vdash N}
 \bigl(\mathsf C_{eN,(d+e)N}\bigr)_{\alpha\beta}\,
 \Ch_\beta(\varpi,\omega;\mathbf n).
\]

(b) The pair $(u,w+u)$ is a framed pair with the same ring, because
$F((u))[w+u]=F((u))[w]$. Since $\val_i(w)=(d+e)\val_i(\varpi)>
\val_i(u)$ for all $i$, we have $\val_i(w+u)=\val_i(u)$ and
$(w+u)/u=1+w'\in1+\rad\mathbf O$. Hence $(u,w+u)$ is balanced of type
$(1,1)$, with $g_i=\val_i(u)=e\val_i(\varpi)$, and for the B\'ezout
pair $(1,0)$ its transform is $(u,w')$. The chain masses depend only
on the ring and on $u$ (see the note after
Definition~\ref{def:chainmass}), so Theorem~\ref{thm:transfer} for
$(m,k)=(eN,eN)$ gives
\[
 \Ch_\alpha(u,w;\mathbf n)=\sum_{\gamma\vdash eN}
 \bigl(\mathsf C_{eN,eN}\bigr)_{\alpha\gamma}\,
 \Ch_\gamma(u,w';\mathbf n).
\]
Next, $(u,w')$ is a framed pair by
Proposition~\ref{prop:transform}(i), we have
$\val_i(w')=d\val_i(\varpi)$, and $w'^e/u^d=w^e/u^{d+e}=\theta$. So
$(u,w')$ is balanced of type $(e,d)$, and since
$(A_0+B_0)e+B_0d=1$, its transform for the B\'ezout pair
$(A_0+B_0,B_0)$ is
$(u^{A_0+B_0}w'^{B_0},\theta-1)=(u^{A_0}w^{B_0},\omega)
=(\varpi,\omega)$. Theorem~\ref{thm:transfer} for $(m,k)=(eN,dN)$
gives
\[
 \Ch_\gamma(u,w';\mathbf n)=\sum_{\beta\vdash N}
 \bigl(\mathsf C_{eN,dN}\bigr)_{\gamma\beta}\,
 \Ch_\beta(\varpi,\omega;\mathbf n).
\]

Comparing (a) and (b), the matrix
$\mathsf X=\mathsf C_{eN,(d+e)N}-\mathsf C_{eN,eN}\mathsf C_{eN,dN}$
annihilates the vector
$\bigl(\Ch_\beta(\varpi,\omega;\mathbf n)\bigr)_{\beta\vdash N}$ for
every such $(\varpi,\omega;\mathbf n)$ and every prime power $q$. For
the split pairs $(\varpi,\omega_\beta)$ with ranks $\mathbf n=\beta$,
$\beta\vdash N$, of the proof of Proposition~\ref{prop:indep}, where
$\val_i(\varpi)=1$, these vectors are the columns of the matrix
$\sM_N\cdot\sN\cdot\diag(q^{c_\beta})$, which is invertible. Hence
$\mathsf X(q)=0$ for every prime power $q$, and since the entries of
$\mathsf X$ are rational functions of $q$
(Corollary~\ref{cor:basis}(ii)), we get $\mathsf X=0$.
\end{proof}
\noindent Consequently only the matrices with $k\le m$ have to be
computed, and the output of the recursion \eqref{eq:alg-recursion}
does not depend on whether a frame with $\val w>\val u$ is first
normalized by $w\mapsto w+u$, as in Step~2 of the proof of
Theorem~\ref{thm:E}, or used directly. For example, for the
multiplicity sequence $(4,2,2)$ the canonical recursion contracts
$\mathsf C_{4n,4n}\mathsf C_{4n,2n}$, while a frame with valuations
$(4,6)$ uses $\mathsf C_{4n,6n}$; the two agree by the case
$(e,d,N,j)=(2,1,2n,1)$.

\section{The universal coefficients and the Hall algebra}\label{sec:hall}

We determine the scalar transfer coefficients from torus counts by
evaluating the transfer identity on split collections of shifted
torus branches. Computing the same count by branch multiplicativity
and by residual type gives the Hall-algebra character of
Theorem~\ref{thm:hallchar}.

Throughout this section we fix coprime integers $a,b\ge1$ and integers $A_0,B_0$ with $A_0a+B_0b=1$. Let us denote the lattice counts of the branch $F[[t^a,t^b]]$ by:
\[
E_j=E_{a,b}(j)=\#E_{F[[t^a,t^b]]}(O^j)\quad(j\ge1),\qquad E_0=1 .
\]
Theorem~\ref{thm:hallchar} determines the functions $\Phi^{a,b}_{q,N}$
of Definition~\ref{def:Phi} from these numbers. To compute a general
pattern one also needs its residual counts $B_n^T$, or equivalently
its chain vector (Lemma~\ref{lem:chainflag}); see
Section~\ref{subsec:explicit-counting-algorithm}.

\begin{theorem}\label{thm:hallchar}
The linear form $\varphi=\varphi_{a,b}:\sH\to\Q$ defined by
\[
\varphi(u_\mu)=q^{-ab\binom{|\mu|}2}\ \frac{\Phi^{a,b}_{q,|\mu|}(\mu)}{a_\mu}\qquad(\mu\ne\emptyset),\qquad\varphi(u_\emptyset)=1,
\]
is a homomorphism of $\Q$-algebras, and $\varphi(u_{(1^j)})=q^{-ab\binom j2}E_j/|\GL_j(F)|$ for all $j\ge1$. Since the elements $u_{(1^j)}$ generate $\sH$, this determines $\Phi^{a,b}_{q,N}$ for all $N$ in terms of the numbers $E_j$, $j\le N$. Equivalently, for every composition $\mathbf n=(n_1,\dots,n_r)$ of $N$ we have:
\begin{equation}\label{eq:hallmean}
q^{ab\,e_2(\mathbf n)}\prod_{i=1}^rE_{n_i}=q^{-e_2(\mathbf n)}\sum_{\mu\vdash N}\#\big(\fn_{\mathbf n}(F)\cap\bO_\mu\big)\ \Phi^{a,b}_{q,N}(\mu).
\end{equation}
\end{theorem}
\noindent In the proof we evaluate the orbit-sum identity of Corollary~\ref{cor:orbitsum} on split pairs. One side of the identity is then given by Proposition~\ref{prop:split}, while the other side is a count of lattices over a reducible plane curve singularity, and this count factors over the branches.

\subsection{Shifted companions}
For $\eta\in tO$ we consider the elements $x_\eta,y_\eta\in O$ and the ring $A_\eta$ defined by:
\[
x_\eta=t^a(1+\eta)^{-B_0},\qquad y_\eta=t^b(1+\eta)^{A_0},\qquad A_\eta=F[[x_\eta,y_\eta]].
\]
By Lemma~\ref{lem:invtransform}, applied with $r=1$, $\varpi=t$, $\omega=\eta$ and $c=1$, the pair $(x_\eta,y_\eta)$ is a balanced framed pair of type $(a,b)$ and weight $1$ with transform $(t,\eta)$, whose ring is $O$. Proposition~\ref{prop:transform} then shows that $A_\eta$ is a plane branch with $\ms(A_\eta)=\eu(a,b)$ and $\delta_{A_\eta}=\delta_{a,b}$. For $\eta=0$ the branch $A_\eta$ is $F[[t^a,t^b]]$, and in general we call $A_\eta$ a \emph{shifted companion} of $F[[t^a,t^b]]$.

\begin{lemma}\label{lem:scalar}
For all $\eta\in tO$ and $j\ge1$ we have $\#E_{A_\eta}(O^j)=\Phi^{a,b}_{q,j}((1^j))$, and in particular $\#E_{A_\eta}(O^j)=E_j$.
\end{lemma}
\begin{proof}
Corollary~\ref{cor:orbitsum}, applied to the framed pair $(x_\eta,y_\eta)$ and its transform $(t,\eta)$, gives $\#E_{A_\eta}(O^j)=\sum_{\mu\vdash j}B_{(t,\eta)}(\mu;j)\Phi^{a,b}_{q,j}(\mu)$. The only lattice in $E_O(O^j)$ is $O^j$, and the element $\eta$ acts by zero on $O^j/tO^j$. Thus the number $B_{(t,\eta)}(\mu;j)$ is equal to $1$ for $\mu=(1^j)$ and vanishes otherwise. This proves the first statement, and the second one follows by comparison with the case $\eta=0$, in which $A_\eta=F[[t^a,t^b]]$.
\end{proof}

\subsection{Multiplicativity over branches}
Let $P=F[[X,Y]]$, let $\mathbf O=O^r$, and let $(U,W)$ be a framed pair in $\mathbf O$. Its ring $R=F[[U,W]]$ is the image of the homomorphism $P\to\mathbf O$, $X\mapsto U$, $Y\mapsto W$. The kernel of the $i$-th component $P\to O$ is a prime ideal of height one, and we denote by $f_i\in P$ a generator of this kernel and by $A_i=F[[U_i,W_i]]\cong P/(f_i)$ the corresponding plane branch.

Since $R$ has finite colength in $\mathbf O$ (Lemma~\ref{lem:gen}), it contains zero divisors supported on a single component, hence the $f_i$ are pairwise non-associate and $R\cong P/(f_1\cdots f_r)$. We denote by $I(f_i,f_j)=\dim_FP/(f_i,f_j)$ the intersection multiplicity, which is finite for $i\ne j$. Finally, for positive integers $\mathbf n=(n_1,\dots,n_r)$ let us recall the set $E_R(\mathbf O^{\mathbf n})$ of Section~\ref{sub:framedmod}.

\begin{proposition}\label{prop:branchmult}
With these notations we have:
\[
\#E_R(\mathbf O^{\mathbf n})=q^{\sum_{i<j}n_in_jI(f_i,f_j)}\ \prod_{i=1}^r\#E_{A_i}(O^{n_i}).
\]
\end{proposition}
\begin{proof}
By induction on $r$ it suffices to prove the following identity, in which $f'=f_2\cdots f_r$, $R_2=P/(f')\subseteq O^{r-1}$, $W_1=\cK^{n_1}$ and $W_2=\prod_{i\ge2}\cK^{n_i}$:
\begin{equation}\label{eq:twofactor}
\#E_R(\mathbf O^{\mathbf n})=q^{n_1\sum_{j\ge2}n_jI(f_1,f_j)}\ \#E_{A_1}(O^{n_1})\ \#E_{R_2}\Big(\prod_{i\ge2}O^{n_i}\Big).
\end{equation}

\emph{Step 1. The map and its fibres.} For an $R$-lattice $L\subseteq W_1\times W_2$ we set $L_1=\pr_1(L)\subseteq W_1$ and $L_2=L\cap(0\times W_2)$, regarded as a subset of $W_2$. Then $L_1$ is an $A_1$-lattice, and $L_2\supseteq f_1\pr_2(L)$ because $f_1$ acts as zero on $W_1$. As $f_1$ acts invertibly on $W_2$, it follows that $L_2$ is an $R_2$-lattice in $W_2$.

For fixed $(L_1,L_2)$ and $\overline W=W_2/L_2$, the lattices $L$ with $\pr_1(L)=L_1$ and $L\cap(0\times W_2)=L_2$ are in bijection with $\Hom_P(L_1,\overline W)$. Namely, to $L$ corresponds the map $l_1\mapsto l_2+L_2$ for $(l_1,l_2)\in L$, and to $\phi$ corresponds the lattice $\{(l_1,l_2):\ l_2+L_2=\phi(l_1)\}$, which is a finitely generated $P$-module annihilated by $f_1f'$ and spans $W_1\times W_2$.

We count the fibres of this map of sets of lattices and apply the
equivariant mass formula \eqref{eq:massfibre}.

\emph{Step 2. The fibres are finite and of constant size.} If $M$ is a $P$-module annihilated by $f_1$, then $\Ext^j_P(M,W_2)=0$ for all $j$, because $f_1$ acts on this group as zero through the first argument and bijectively through the second. We apply this to $M=L_1$ and to the exact sequence $0\to L_2\to W_2\to\overline W\to0$, and we obtain:
\[
\Hom_P(L_1,\overline W)\cong\Ext^1_P(L_1,L_2),\qquad\Ext^1_P(L_1,\overline W)\cong\Ext^2_P(L_1,L_2)=0 .
\]
The last group vanishes because $L_1$ has projective dimension one over the regular local ring $P$ of dimension two. Indeed, $L_1$ is torsion-free over the one-dimensional domain $A_1$, hence it has depth one, and the Auslander--Buchsbaum formula \cite[Thm.~1.3.3]{BH} applies. The module $\Ext^1_P(L_1,L_2)$ has finite length, because it is finitely generated and annihilated by $f_1$ and by $f'$, and we denote $\hbar(L_1,L_2)=\dim_F\Hom_P(L_1,\overline W)$.

Now let $L_1'\subseteq L_1$ be an $A_1$-lattice with $L_1/L_1'\cong F$. Applying the functor $\Hom_P(-,\overline W)$ to the exact sequence $0\to L_1'\to L_1\to F\to0$ and using that $\Ext^1_P(L_1,\overline W)=0$, we get:
\begin{align*}
\hbar(L_1,L_2)-\hbar(L_1',L_2)&=\dim\Hom_P(F,\overline W)-\dim\Ext^1_P(F,\overline W)\\
&=\dim\Ext^1_P(F,L_2)-\dim\Ext^2_P(F,L_2),
\end{align*}
where the second equality uses that $\Ext^j_P(F,W_2)=0$ for all $j$, as $f_1\in(X,Y)$. The framed parameter $U$ acts injectively on $L_2\subseteq W_2$,
since its component on each branch is nonzero. Hence $L_2$ has
depth at least one over $P$. Its support has dimension at most one because $f'$ annihilates it.
Thus its depth is exactly one.

The Auslander--Buchsbaum formula gives projective dimension
one and a resolution $0\to P^s\to P^{s'}\to L_2\to0$.
Here $s=s'$ because $L_2$ is a torsion $P$-module. Note that $\Ext^j_P(F,P)$ is $F$ for $j=2$ and zero otherwise, thus the resolution gives an exact sequence $0\to\Ext^1_P(F,L_2)\to F^s\to F^{s'}\to\Ext^2_P(F,L_2)\to0$ and the difference above is $s-s'=0$. Any two $A_1$-lattices in $W_1$ are joined by a chain of lattices with successive quotients $F$, thus $\hbar(L_1,L_2)=\hbar(A_1^{n_1},L_2)$.

Next, we observe that $\Hom_P(A_1,\overline W)=\{\bar w\in\overline W:f_1\bar w=0\}=f_1^{-1}L_2/L_2\cong L_2/f_1L_2$, and therefore:
\[
\hbar(L_1,L_2)=n_1\dim_FL_2/f_1L_2 .
\]
For a sublattice $L_2'\subseteq L_2$ of finite colength, the snake lemma for multiplication by $f_1$ on $0\to L_2'\to L_2\to L_2/L_2'\to0$ gives $\dim L_2/f_1L_2=\dim L_2'/f_1L_2'$, because the kernel and the cokernel of an endomorphism of a finite-dimensional space have the same dimension. Thus we may replace $L_2$ by the commensurable lattice $\bigoplus_{j\ge2}A_j^{n_j}$, for which $\dim L_2/f_1L_2=\sum_{j\ge2}n_j\dim A_j/f_1A_j=\sum_{j\ge2}n_jI(f_1,f_j)$. Consequently, every fibre has $q^{n_1\sum_{j\ge2}n_jI(f_1,f_j)}$ elements.

\emph{Step 3. Masses.} The map $L\mapsto(L_1,L_2)$ is equivariant for the group $Z=\GL_{n_1}(\cK)\times\prod_{i\ge2}\GL_{n_i}(\cK)$. By Lemma~\ref{lem:massE}, applied to $R$, $A_1$ and $R_2$, the $Z$-mass of the set of all $R$-lattices is $\#E_R(\mathbf O^{\mathbf n})$, while the mass of the set of pairs $(L_1,L_2)$ is $\#E_{A_1}(O^{n_1})\cdot\#E_{R_2}(\prod_{i\ge2}O^{n_i})$. By Step~2 all fibres have the same number of elements, and \eqref{eq:massfibre} implies \eqref{eq:twofactor}.
\end{proof}

\subsection{Intersection numbers of shifted companions}
\begin{lemma}\label{lem:intersection}
Let $\eta\ne\eta'$ be elements of $tO$, and let $f_\eta,f_{\eta'}\in P$ be generators of the kernels of the maps $P\to O$, $(X,Y)\mapsto(x_\eta,y_\eta)$ and $(X,Y)\mapsto(x_{\eta'},y_{\eta'})$, respectively, that is, equations of the branches $A_\eta$ and $A_{\eta'}$. Then
\[
I(f_\eta,f_{\eta'})=ab+\val(\eta-\eta').
\]
\end{lemma}
\begin{proof}
In $\mathbf O=O^2$ we put $U=(x_\eta,x_{\eta'})$, $W=(y_\eta,y_{\eta'})$, $\varpi=(t,t)$ and $\omega=(\eta,\eta')$. The pair $(\varpi,\omega)$ is a framed pair (Section~\ref{sub:split}), and Lemma~\ref{lem:invtransform} shows that $(U,W)$ is a balanced framed pair of type $(a,b)$ with $g_1=g_2=1$ and weight $2$, whose transform is $(\varpi,\omega)$. For $R=F[[U,W]]$ and $T_1=F[[\varpi,\omega]]$, Proposition~\ref{prop:transform}(ii) then gives $\dim_FT_1/R=\sum_{m\in\eu(a,b)}\binom{2m}2$.

The ring $T_1$ contains the diagonal copy of $O$ and the element $(0,\eta'-\eta)$, and for $k=\val(\eta-\eta')$ it is contained in the closed subring $\{(f,f'):f\equiv f'\bmod t^{k}\}$, which is generated by these elements, hence $\dim_F\mathbf O/T_1=k$. On the other hand, the exact sequence
\[
0\longrightarrow P/(f_\eta f_{\eta'})\longrightarrow P/(f_\eta)\times P/(f_{\eta'})\longrightarrow P/(f_\eta,f_{\eta'})\longrightarrow0
\]
shows that $I(f_\eta,f_{\eta'})=\dim(A_\eta\times A_{\eta'})/R=\dim\mathbf O/R-2\delta_{a,b}$, and together with \eqref{eq:eusums} we obtain:
\[
I(f_\eta,f_{\eta'})=k+\sum_{m\in\eu(a,b)}\Big(\binom{2m}2-2\binom m2\Big)=k+\sum_{m\in\eu(a,b)}m^2=k+ab .\qedhere
\]
\end{proof}

\subsection{Proof of Theorem~\ref{thm:hallchar}}
Let $\mathbf n=(n_1,\dots,n_r)$ be a composition of $N$, let us choose pairwise distinct $\eta_1,\dots,\eta_r\in tO$, for instance $\eta_i=t^i$, and set $k_{ij}=\val(\eta_i-\eta_j)$. In $\mathbf O=O^r$ we consider the split pair $(\varpi,\omega_\eta)$ of Section~\ref{sub:split} and the elements:
\[
U=(x_{\eta_1},\dots,x_{\eta_r}),\qquad W=(y_{\eta_1},\dots,y_{\eta_r}),\qquad R=F[[U,W]].
\]
By Lemma~\ref{lem:invtransform}, $(U,W)$ is a balanced framed pair of type $(a,b)$ with $g_i=1$ and transform $(\varpi,\omega_\eta)$, and here $n'=\sum_in_i=N$. Thus Corollary~\ref{cor:orbitsum} and Proposition~\ref{prop:split} give:
\[
\#E_R(\mathbf O^{\mathbf n})=q^{\sum_{i<j}n_in_j(k_{ij}-1)}\sum_{\mu\vdash N}\#\big(\fn_{\mathbf n}(F)\cap\bO_\mu\big)\,\Phi^{a,b}_{q,N}(\mu).
\]
On the other hand, the components of $(U,W)$ are the shifted companions $A_{\eta_i}$, and Proposition~\ref{prop:branchmult} together with Lemmas~\ref{lem:intersection} and~\ref{lem:scalar} gives:
\[
\#E_R(\mathbf O^{\mathbf n})=q^{\sum_{i<j}n_in_j(ab+k_{ij})}\prod_{i=1}^rE_{n_i}.
\]

Cancelling the common contact factor $q^{\sum_{i<j}n_in_jk_{ij}}$
leaves the exponents $-e_2(\mathbf n)$ and $ab\,e_2(\mathbf n)$,
respectively, and gives \eqref{eq:hallmean}. By Lemma~\ref{lem:hallbasics}(a), the right side of \eqref{eq:hallmean} is equal to $\prod_i|\GL_{n_i}(F)|\cdot\sum_\mu F_{\mathbf n}(\mu)\Phi^{a,b}_{q,N}(\mu)/a_\mu$, and if we multiply by $q^{-ab\binom N2}$ and use $\binom N2=e_2(\mathbf n)+\sum_i\binom{n_i}2$, we get:
\[
\varphi\Big(\prod_{i=1}^ru_{(1^{n_i})}\Big)=\sum_{\mu\vdash N}F_{\mathbf n}(\mu)\,\varphi(u_\mu)=\prod_{i=1}^r\frac{q^{-ab\binom{n_i}2}E_{n_i}}{|\GL_{n_i}(F)|}.
\]
For $r=1$ this is the stated value of $\varphi(u_{(1^j)})$, where we use $a_{(1^j)}=|\GL_j(F)|$ and Lemma~\ref{lem:scalar}. In general it says that $\varphi$ is multiplicative on monomials in the elements $u_{(1^j)}$, and since these monomials form a basis of $\sH$ by Lemma~\ref{lem:hallbasics}(c), $\varphi$ is a homomorphism of algebras. \qed

\begin{remark}\label{rem:Phi}
(1) For $(a,b)=(1,h)$ the branch $F[[t,t^h]]$ is $O$, hence $E_j=1$ for all $j$ and the theorem gives $\varphi_{1,h}(u_{(1^j)})=q^{-h\binom j2}/|\GL_j(F)|$.

(2) By \eqref{eq:hallmean} and Lemma~\ref{lem:hallbasics}, $\Phi^{a,b}_{q,N}$ is obtained from
the products $\prod_i E_{\lambda_i}$, $\lambda\vdash N$, by
inverting a triangular matrix; in particular,
$\Phi^{a,b}_{q,N}=\Phi^{b,a}_{q,N}$. Its values are rational functions
of $q$ by Corollary~\ref{cor:basis}.
In the setting of Corollary~\ref{cor:cabling-data}, Conjecture~S predicts that this is
the cardinality of each saturation fibre over a lattice of type $\mu$.
Corollary~\ref{cor:orbitsum} proves equality after summing over all transformed
lattices; it does not identify the sum over one prescribed residual
type or an individual fibre.

(3) The identity \eqref{eq:hallmean} does not involve the contacts $k_{ij}$, because they cancel between the intersection numbers of Lemma~\ref{lem:intersection} and the Jacobian of Proposition~\ref{prop:split}.
\end{remark}

\section{The cabling formula}\label{sec:limit}

The cabling formula $\cN_A=\cN_{a,b}^{\,D}\cN_S$ follows by
comparing the lattice counts of $A$ and $S$ at finite rank and then
passing to the limit. Here $A=F[[x,y]]$ has $\val x=Da$,
$\val y=Db$ and $\gcd(a,b)=1$, and $T$, $S$ and $h=a+b-1$ are as in
Corollary~\ref{cor:cabling-data}. The framed branches $(A;x,y)$ and
$(S;\varpi,\varpi^h(1+\omega))$ have the same transform
$(T;\varpi,\omega)$, so their orbit-sum identities involve the same
residual multiplicities. After normalization, their coefficient
functions solve Hall systems with the same matrix $A^{(n)}$,
which depends on $(a,b)$ only through $h$.

Near the rectangle $(n^D)$, the products of normalized torus counts
converge to $\cN_{a,b}^{\,D}$. To pass this comparison through
$(A^{(n)})^{-1}$, we need an estimate uniform in $n$, since the matrix
size grows with $n$. Lemma~\ref{lem:inverse} bounds the rows of the
inverse after multiplication by their Gaussian weights and shows
that the tails tend to zero uniformly. We split the sum at distance
$M$ from the rectangle and let first $n\to\infty$, then
$M\to\infty$, obtaining Theorem~\ref{thm:limit}. Iterating the
formula by induction on multiplicity gives Theorem~\ref{thm:product}.

\subsection{The finite-rank identities}\label{sub:finite-rank-id}
Let $A=F[[x,y]]$ be a plane branch with $\val x=Da$, $\val y=Db$ and $\gcd(a,b)=1$, and let $\varpi,\omega,T,S$ and $h=a+b-1$ be as in Corollary~\ref{cor:cabling-data}. For $\mu\vdash Dn$ we denote by $B^T_n(\mu)=B_{(\varpi,\omega)}(\mu;n)$ the number of $\Lambda\in E_T(O^n)$ such that $\omega$ has Jordan type $\mu$ on $\Lambda/\varpi\Lambda$.

\begin{theorem}\label{thm:cable}
For every $n\ge1$ we have
\begin{equation}\label{eq:twoidentities}
\#E_A(O^n)=\sum_{\mu\vdash Dn}B^T_n(\mu)\,\Phi^{a,b}_{q,Dn}(\mu),\qquad\#E_S(O^n)=\sum_{\mu\vdash Dn}B^T_n(\mu)\,\Phi^{1,h}_{q,Dn}(\mu),
\end{equation}
where $\Phi^{a,b}$ and $\Phi^{1,h}$ are given by Theorem~\ref{thm:hallchar} and Remark~\ref{rem:Phi}(1). Moreover, $\sum_\mu B^T_n(\mu)=\#E_T(O^n)$, and $B^T_n(\mu)=0$ unless all parts of $\mu$ are at most $D$. In the case $D=1$ we have $T=S=O$, $\delta_A=\delta_{a,b}$ and $\#E_A(O^n)=E_{a,b}(n)$.
\end{theorem}
\begin{proof}
By Corollary~\ref{cor:cabling-data}(i), $(T;\varpi,\omega)$ is the transform of the framed branch $(A;x,y)$ of type $(a,b)$, and it is also the transform of the framed branch $(S;\varpi,\varpi^h(1+\omega))$ of type $(1,h)$, both framed branches having weight $D$. Thus \eqref{eq:twoidentities} is Corollary~\ref{cor:orbitsum}, applied twice with $r=1$ and $n'=Dn$, and the statement on the parts of $\mu$ is Corollary~\ref{cor:cabling-data}(iv).

Finally, if $D=1$, then $\varpi$ is a uniformizer, $T=S=O$ and $E_O(O^n)=\{O^n\}$, and since $\omega$ acts by zero on $O^n/\varpi O^n$, the first identity reduces to $\#E_A(O^n)=\Phi^{a,b}_{q,n}((1^n))$. It remains to recall that $\Phi^{a,b}_{q,n}((1^n))=E_{a,b}(n)$ by Lemma~\ref{lem:scalar}, and that $\delta_A=\delta_{a,b}$ by Corollary~\ref{cor:cabling-data}(iii).
\end{proof}

\subsection{The normalized system}\label{sub:normalized}
Fix coprime $a,b\ge1$, and put
$\delta_0=(a-1)(b-1)/2$ and $h=a+b-1$, so that
$ab-2\delta_0=h$. The following normalization makes the Hall
matrix depend only on $h$:
\[
\widehat E_j=Q^{j^2\delta_0}E_j=e_{F[[t^a,t^b]]}(j),\qquad\widehat\Phi(\mu)=Q^{\frac{ab}2N^2-\frac h2\sum_i(\mu'_i)^2}\ \Phi^{a,b}_{q,N}(\mu)\quad(\mu\vdash N).
\]
For a composition $\mathbf n$ of $N$ we denote by $\pi_{\mathbf n}(\mu)=q^{-e_2(\mathbf n)}\#(\fn_{\mathbf n}(F)\cap\bO_\mu)$ the probability that a uniformly random element of $\fn_{\mathbf n}(F)$ has type $\mu$. In terms of these probabilities the identity \eqref{eq:hallmean} is equivalent to the following one, because $ab\,e_2(\mathbf n)+\delta_0\sum_in_i^2=\frac{ab}2N^2-\frac h2\sum_in_i^2$:
\begin{equation}\label{eq:normalized}
\prod_i\widehat E_{n_i}=\sum_{\mu\vdash N}\pi_{\mathbf n}(\mu)\ Q^{\frac h2\left(\sum_i(\mu'_i)^2-\sum_in_i^2\right)}\ \widehat\Phi(\mu).
\end{equation}

Now we fix $D\ge2$ and $n\ge1$ and denote by $\cP_n$ the set of partitions of $Dn$ with at most $D$ parts, written as $\lambda=(\lambda_1\ge\dots\ge\lambda_D\ge0)$. We write
\[
\lambda=(n+\alpha_1,\dots,n+\alpha_D),\qquad\textstyle\sum_i\alpha_i=0,\qquad\|\alpha\|^2=\sum_i\alpha_i^2,\qquad\|\lambda\|^2=\sum_i\lambda_i^2=Dn^2+\|\alpha\|^2,
\]
and for $\lambda,\nu\in\cP_n$ we set
\[
\widehat\Phi_\nu=\widehat\Phi(\nu'),\qquad\pi_{\lambda\nu}=\pi_\lambda(\nu'),\qquad A_{\lambda\nu}=\pi_{\lambda\nu}\,Q^{\frac h2(\|\nu\|^2-\|\lambda\|^2)},
\]
where the zero parts of $\lambda$ are omitted in $\fn_\lambda$. An element of $\fn_\lambda(F)$ has nilpotency index at most $D$, hence its type is $\nu'$ with $\nu\in\cP_n$ and $\sum_{\nu\in\cP_n}\pi_{\lambda\nu}=1$. Thus \eqref{eq:normalized} reads:
\begin{equation}\label{eq:system}
\prod_{i=1}^D\widehat E_{\lambda_i}=\sum_{\nu\in\cP_n}A_{\lambda\nu}\,\widehat\Phi_\nu\qquad(\lambda\in\cP_n),
\end{equation}
with $\widehat E_0=1$. The matrix $A=A^{(n)}$ is square, with both
indices in $\cP_n$; it differs from the generally rectangular
chain-transfer matrix $\mathsf C_{Dn,dn}$ of
Section~\ref{subsec:explicit-counting-algorithm}.
It depends on $(a,b)$ only through $h$, which lets us compare
the systems for $(a,b)$ and $(1,h)$.

All powers in these normalizations are integral. For a partition
$\nu$ of $N$,
\[
 \frac{abN^2-h\|\nu\|^2}{2}
   =\delta_0N^2+h\,e_2(\nu),\qquad
 \|\nu\|^2\equiv N\pmod2.
\]
Thus $\|\nu\|^2-\|\lambda\|^2$ is even for partitions of
the same weight. Likewise, $\|\alpha\|^2$ is even when
$\alpha\in\Z^D$ has coordinate sum zero. Thus the normalizations in
Section~\ref{sec:modular} use the ordinary power-series ring $\Q[[Q]]$.

\begin{lemma}\label{lem:triangular}
Let $\lambda=n+\alpha$ and $\nu=n+\beta$ be elements of $\cP_n$.
\begin{enumerate}[label=\textup{(\roman*)}]
\item If $\pi_{\lambda\nu}\ne0$, then $\nu\ctd\lambda$.
\item $\pi_{\lambda\lambda}=\prod_{i=1}^{D-1}\dfrac{(Q;Q)_{\lambda_i}}{(Q;Q)_{\lambda_i-\lambda_{i+1}}}\ge(Q;Q)_\infty^{D-1}$.
\item If $\nu\ctd\lambda$, then $\|\beta\|\ge\|\alpha\|$, and if moreover $\nu\ne\lambda$, then $\mathsf h(\beta)\ge\mathsf h(\alpha)+1$, where $\mathsf h(\beta)=\sum_{i=1}^D(D-i)\beta_i$.
\item $0\le\mathsf h(\beta)\le\varrho\,\|\beta\|$, where $\varrho^2=\sum_{i=1}^D(D-i)^2$.
\end{enumerate}
\end{lemma}
\begin{proof}
Part (i) follows from Lemma~\ref{lem:hallbasics}(a),(b).

For the proof of (ii) we consider a flag with $\dim V_i/V_{i-1}=\lambda_i$. By Lemma~\ref{lem:hallbasics}(a), $\pi_{\lambda\lambda}$ is also the probability that a uniformly random element $X$ of the space $\{X:XV_i\subseteq V_{i-1}\}$ has type $\lambda'$, and by Lemma~\ref{lem:hallbasics}(b) this happens if and only if the induced maps $V_i/V_{i-1}\to V_{i-1}/V_{i-2}$ are injective. For uniformly random $X$ these maps are independent and uniformly distributed, and it remains to note that a uniformly random $r\times s$ matrix with $r\ge s$ is injective with probability $\prod_{j=0}^{s-1}(1-Q^{r-j})=(Q;Q)_r/(Q;Q)_{r-s}$.

For the proof of (iii) we put $S_k(v)=\sum_{i\le k}v_i$, and summation by parts gives:
\[
\|\nu\|^2-\|\lambda\|^2=\sum_{k=1}^{D-1}\big(S_k(\nu)-S_k(\lambda)\big)\big(\nu_k+\lambda_k-\nu_{k+1}-\lambda_{k+1}\big)\ge0,
\]
where the left side is equal to $\|\beta\|^2-\|\alpha\|^2$. The second assertion of (iii) holds because $\mathsf h(\beta)-\mathsf h(\alpha)=\sum_{k=1}^{D-1}(S_k(\nu)-S_k(\lambda))$ is a nonnegative integer, positive if $\nu\ne\lambda$.

Finally, for (iv) we note that $S_k(\beta)\ge0$ for a decreasing vector of sum zero, hence $\mathsf h(\beta)=\sum_{k<D}S_k(\beta)\ge0$, and the upper bound is the Cauchy--Schwarz inequality.
\end{proof}

To control the inverse uniformly in $n$, we expand it along
increasing chains of partitions. The factors $\pi_{\lambda\nu}$
are transition probabilities, whose products sum to at most one
at each fixed chain length. Gaussian decay then bounds both long
chains and partitions far from the rectangle $(n^D)$.

\begin{lemma}[Uniform inverse estimate]\label{lem:inverse}
Let $\mathsf c=(Q;Q)_\infty^{-(D-1)}$, and for $r\ge0$ let us define
\[
\tau(r)=\sum_{k\ge0}\mathsf c^{\,k+1}\,Q^{\frac h2\max(r^2,\,k^2/\varrho^2)} .
\]
Then $\tau(r)<\infty$, $\tau$ is nonincreasing, $\tau(r)\to0$ as $r\to\infty$, and the function $\tau$ depends only on $(D,h,q)$. For every $n$ the matrix $A=A^{(n)}$ is invertible, and for all $\lambda=n+\alpha\in\cP_n$ and $M\ge0$ we have
\begin{equation}\label{eq:inverse}
Q^{\frac h2\|\alpha\|^2}\sum_{\nu=n+\beta\in\cP_n,\ \|\beta\|\ge M}\big|(A^{-1})_{\lambda\nu}\big|\ \le\ \tau\big(\max(M,\|\alpha\|)\big).
\end{equation}
\end{lemma}
\begin{proof}
The series converges because $Q^{\frac h2k^2/\varrho^2}$ decays faster than any geometric sequence, and the inequality $\max(x,y)\ge(x+y)/2$ gives $\tau(r)\le\mathsf K\,Q^{hr^2/4}$ with $\mathsf K=\sum_k\mathsf c^{k+1}Q^{hk^2/(4\varrho^2)}$.

With respect to a linear extension of the dominance order on $\cP_n$, Lemma~\ref{lem:triangular}(i),(ii) gives $A=\Delta+A^+$, where $\Delta$ is diagonal with $\Delta_{\lambda\lambda}=\pi_{\lambda\lambda}\ge\mathsf c^{-1}$ and $A^+$ is strictly triangular. Thus $A$ is invertible, and expanding the finite sum $A^{-1}=\sum_{k\ge0}(-\Delta^{-1}A^+)^k\Delta^{-1}$ we see that the entry $(A^{-1})_{\lambda\nu}$ is a sum over the chains $\lambda=\nu^0\triangleleft\nu^1\triangleleft\dots\triangleleft\nu^k=\nu$, $k\ge0$, and that the term of such a chain is
\[
\pm\ \frac1{\pi_{\nu\nu}}\ \prod_{j=1}^k\frac{\pi_{\nu^{j-1}\nu^j}}{\pi_{\nu^{j-1}\nu^{j-1}}}\ \ Q^{\frac h2(\|\nu\|^2-\|\lambda\|^2)},
\]
where the Gaussian factors have telescoped. The reciprocals of the diagonal entries contribute at most $\mathsf c^{k+1}$, and for fixed $k$ the sum of the products $\prod_j\pi_{\nu^{j-1}\nu^j}$ over all chains of length $k$ starting at $\lambda$ is at most $1$, because the matrix $(\pi_{\lambda\nu})$ is stochastic.

For the end point $\nu=n+\beta$ of a chain of length $k$, Lemma~\ref{lem:triangular}(iii),(iv) gives $\|\beta\|\ge\|\alpha\|$ and $\varrho\|\beta\|\ge\mathsf h(\beta)\ge\mathsf h(\alpha)+k\ge k$. Hence, after multiplication by $Q^{\frac h2\|\alpha\|^2}$, the Gaussian factor is $Q^{\frac h2\|\beta\|^2}\le Q^{\frac h2\max(M^2,\|\alpha\|^2,k^2/\varrho^2)}$ whenever $\|\beta\|\ge M$, and summing over $k$ we obtain \eqref{eq:inverse}.
\end{proof}

\subsection{The limit}\label{sub:limit}
Passing to the limit $n\to\infty$ in the finite-rank identities, we obtain the cabling formula for the numbers $\cN_A$.

\begin{theorem}\label{thm:limit}
Let $A=F[[x,y]]$, $(a,b,D)$, $T$ and $S$ be as in Corollary~\ref{cor:cabling-data}. Then the numbers $\cN_A$, $\cN_{a,b}$ and $\cN_S$ are finite, and
\[
\cN_A=\cN_{a,b}^{\,D}\cdot\cN_S .
\]
\end{theorem}
\begin{proof}
Let us denote $C=\cN_{a,b}=\lim_j\widehat E_j$, $B_E=\sup_{j\ge0}\widehat E_j\ge\widehat E_0=1$ and $B_T=\sup_ne_T(n)$. These numbers and $\cN_S$ are finite by Theorem~\ref{thm:finite}, because $F[[t^a,t^b]]$, $T$ and $S$ are plane branches (Corollary~\ref{cor:cabling-data}(i)), and we have $0\le C\le B_E$. If $D=1$, then $e_A(n)=\widehat E_n$ and $\cN_S=\cN_O=1$ by Theorem~\ref{thm:cable}, and there is nothing more to prove, so let $D\ge2$ and $N=Dn$.

\emph{Step 1. The two systems.} We apply \eqref{eq:system} to the pair $(a,b)$ and to the pair $(1,h)$, with the normalized coefficients:
\[
\widehat\Phi_\lambda=Q^{\frac{ab}2N^2-\frac h2\|\lambda\|^2}\,\Phi^{a,b}_{q,N}(\lambda'),\qquad\widehat\Phi^{\circ}_\lambda=Q^{\frac h2N^2-\frac h2\|\lambda\|^2}\,\Phi^{1,h}_{q,N}(\lambda').
\]
For the pair $(1,h)$ the torus branch is $O$ and $\delta_0=0$,
so its normalized counts are all $1$ (Remark~\ref{rem:Phi}(1)).
Since $1+h-1=h$, both pairs give the same matrix $A=A^{(n)}$ in
\eqref{eq:system}. Thus, for $\lambda\in\cP_n$, we have:
\begin{equation}\label{eq:two-systems}
\prod_{i=1}^D\widehat E_{\lambda_i}=\sum_{\nu\in\cP_n}A_{\lambda\nu}\widehat\Phi_\nu,\qquad1=\sum_{\nu\in\cP_n}A_{\lambda\nu}\widehat\Phi^{\circ}_\nu .
\end{equation}
By Lemma~\ref{lem:inverse} the matrix $A$ is invertible, so the two systems determine $\widehat\Phi$ and $\widehat\Phi^{\circ}$.

\emph{Step 2. The normalized counts.} By Theorem~\ref{thm:cable}, $B^T_n(\mu)=0$ unless all parts of $\mu$ are at most $D$, that is, unless $\mu=\lambda'$ with $\lambda\in\cP_n$. Thus the residual types are the partitions $\lambda'$ with $\lambda\in\cP_n$, and the numbers $b_n(\lambda)=Q^{n^2\delta_T}B^T_n(\lambda')\ge0$ satisfy:
\[
\sum_{\lambda\in\cP_n}b_n(\lambda)=Q^{n^2\delta_T}\,\#E_T(O^n)=e_T(n)\le B_T .
\]

Multiply the two identities in \eqref{eq:twoidentities} by
$Q^{n^2\delta_A}$ and $Q^{n^2\delta_S}$, respectively, and substitute
the normalized coefficients. The remaining powers of $Q$ have the
exponents below. For $\lambda=n+\alpha$, we have $N=Dn$ and
$\|\lambda\|^2=Dn^2+\|\alpha\|^2$, so
Corollary~\ref{cor:cabling-data}(iii) gives:
\begin{align*}
n^2(\delta_A-\delta_T)-\Big(\tfrac{ab}2N^2-\tfrac h2\|\lambda\|^2\Big)&=n^2\Big(D^2\delta_0+\tbinom D2h\Big)-\tfrac{ab}2D^2n^2+\tfrac h2\big(Dn^2+\|\alpha\|^2\big)=\tfrac h2\|\alpha\|^2,\\
n^2(\delta_S-\delta_T)-\Big(\tfrac h2N^2-\tfrac h2\|\lambda\|^2\Big)&=n^2\tbinom D2h-\tfrac h2D^2n^2+\tfrac h2\big(Dn^2+\|\alpha\|^2\big)=\tfrac h2\|\alpha\|^2,
\end{align*}
where the first line uses $2\delta_0-ab=-h$. Hence \eqref{eq:twoidentities} becomes:
\begin{equation}\label{eq:normalized-counts}
e_A(n)=\sum_{\lambda\in\cP_n}b_n(\lambda)\,Q^{\frac h2\|\alpha\|^2}\,\widehat\Phi_\lambda,\qquad e_S(n)=\sum_{\lambda\in\cP_n}b_n(\lambda)\,Q^{\frac h2\|\alpha\|^2}\,\widehat\Phi^{\circ}_\lambda .
\end{equation}
Thus $e_A(n)$ and $e_S(n)$ are sums of the normalized coefficients with the same nonnegative weights $b_n(\lambda)Q^{\frac h2\|\alpha\|^2}$, and the sum of these weights is at most $B_T$.

\emph{Step 3. The estimate.} Subtracting $C^D$ times the second system of \eqref{eq:two-systems} from the first one, we see that the differences $\epsilon_\lambda=\widehat\Phi_\lambda-C^D\,\widehat\Phi^{\circ}_\lambda$ solve the system with the matrix $A$ and the left sides $\prod_i\widehat E_{\lambda_i}-C^D$. Hence:
\[
\epsilon_\lambda=\sum_{\nu\in\cP_n}(A^{-1})_{\lambda\nu}\Big(\prod_{i=1}^D\widehat E_{\nu_i}-C^D\Big).
\]
Since $0\le\widehat E_j,C\le B_E$, the bracketed factor is at most
$2B_E^D$ in absolute value for every $\nu\in\cP_n$.

Next we fix an integer $M$ with $0\le M\le n$. If $\nu=n+\beta$ satisfies $\|\beta\|\le M$, then $|\beta_i|\le M$, and hence $\nu_i\ge n-M$ for all $i$. The telescoping identity
\[
\prod_{i=1}^D\widehat E_{\nu_i}-C^D=\sum_{k=1}^D\Big(\prod_{i<k}\widehat E_{\nu_i}\Big)\big(\widehat E_{\nu_k}-C\big)\,C^{D-k}
\]
then gives the bound:
\[
\Big|\prod_{i=1}^D\widehat E_{\nu_i}-C^D\Big|\le\vartheta_n(M):=D\,B_E^{D-1}\sup_{j\ge n-M}\big|\widehat E_j-C\big|\qquad(\|\beta\|\le M).
\]

Now let $\lambda=n+\alpha\in\cP_n$. We split the sum defining $\epsilon_\lambda$ into the terms with $\|\beta\|<M$ and the terms with $\|\beta\|\ge M$, and we multiply by $Q^{\frac h2\|\alpha\|^2}$. For the first part we use the bound $\vartheta_n(M)$ and \eqref{eq:inverse} with $M=0$, which bounds $Q^{\frac h2\|\alpha\|^2}\sum_\nu|(A^{-1})_{\lambda\nu}|$ by $\tau(\|\alpha\|)$. For the second part we use the bound $2B_E^D$ and \eqref{eq:inverse} with the given $M$. Since $\tau$ is nonincreasing, we obtain:
\[
Q^{\frac h2\|\alpha\|^2}|\epsilon_\lambda|\le\vartheta_n(M)\,\tau(\|\alpha\|)+2B_E^D\,\tau\big(\max(M,\|\alpha\|)\big)\le\vartheta_n(M)\,\tau(0)+2B_E^D\,\tau(M).
\]

\emph{Step 4. Conclusion.} Subtracting $C^D$ times the second identity of \eqref{eq:normalized-counts} from the first one, we get $e_A(n)-C^De_S(n)=\sum_{\lambda\in\cP_n}b_n(\lambda)\,Q^{\frac h2\|\alpha\|^2}\epsilon_\lambda$. The weights $b_n(\lambda)$ are nonnegative with sum at most $B_T$, so the estimate of Step~3 gives, for all $n\ge1$ and all integers $0\le M\le n$:
\[
\big|e_A(n)-C^De_S(n)\big|\le\sum_{\lambda\in\cP_n}b_n(\lambda)\,Q^{\frac h2\|\alpha\|^2}|\epsilon_\lambda|\le B_T\big(\vartheta_n(M)\,\tau(0)+2B_E^D\,\tau(M)\big).
\]
Let us fix $M$. Since $\widehat E_j\to C$ as $j\to\infty$, we have $\vartheta_n(M)\to0$ as $n\to\infty$, and therefore:
\[
\limsup_{n\to\infty}\big|e_A(n)-C^De_S(n)\big|\le2B_TB_E^D\,\tau(M).
\]
Letting $M\to\infty$ and using Lemma~\ref{lem:inverse} gives
$e_A(n)-C^De_S(n)\to0$. Since $e_S(n)\to\cN_S<\infty$, we obtain
$e_A(n)\to C^D\cN_S$, proving the cabling formula.
\end{proof}

\subsection{Iteration}\label{sub:iterate}
For a plane branch $A$ with $m=\mult(A)\ge2$, Lemma~\ref{lem:normalgen} gives generators with $\val x=m$ and $\val y=B$, $m\nmid B$, and we set:
\[
D=\gcd(m,B),\qquad a=m/D,\qquad b=B/D.
\]
Then $2\le a<b$, the branch $S=F[[\varpi,\varpi^{a+b-1}\omega]]$ of Corollary~\ref{cor:cabling-data} has multiplicity $D<m$, and we call $(a,b,D)$ a \emph{first cabling triple} of $A$, the ring $T$ of Corollary~\ref{cor:cabling-data} the \emph{pattern ring}, and $S$ a \emph{reduction} of $A$. Repeating the construction with $S$ until we reach the branch $O$, which happens when the multiplicity becomes $1$ (Lemma~\ref{lem:noether}), we obtain a finite list
\[
\mathcal L(A)=\big((a_1,b_1,D_1),\dots,(a_g,b_g,D_g)\big),
\]
with $2\le a_i<b_i$, $\gcd(a_i,b_i)=1$, $D_ia_i=D_{i-1}$ (where $D_0=m$) and $D_g=1$, which we call the \emph{cabling list} of $A$. For $A=O$ the cabling list is empty, and in general it is determined by the multiplicity sequence of $A$, as the following proposition shows.

\begin{proposition}\label{prop:list}
The first cabling triple of $A$ and the multiplicity sequence of its reduction are determined by $\ms(A)=(\mu_0,\mu_1,\dots)$ as follows: we have $m=\mu_0$, and if $k\ge1$ is the number of initial entries equal to $m$, then $B=km+\mu_k$ and $\ms(S)=\big(D^{\,a+b-1},\mu_{\ell(a,b)},\mu_{\ell(a,b)+1},\dots\big)$. Consequently, the cabling list $\mathcal L(A)$ depends only on $\ms(A)$.
\end{proposition}
\begin{proof}
By Corollary~\ref{cor:cabling-data}(ii) the multiplicity sequences of $A$ and $S$ are $\ms(A)=(D\cdot\eu(a,b),\ms(T))$ and $\ms(S)=(D^{a+b-1},\ms(T))$. Since $2\le a<b$, the sequence $\eu(a,b)$ begins with $\lfloor b/a\rfloor$ entries equal to $a$, followed by $b-\lfloor b/a\rfloor a<a$, hence $k=\lfloor b/a\rfloor$ and $\mu_k=D(b-ka)=B-km$.
\end{proof}

\begin{theorem}\label{thm:product}
For every plane branch $A$ over $\F_q$, in any characteristic, we have:
\[
\cN_A=\prod_{i=1}^g\cN_{a_i,b_i}^{\,D_i},\qquad\mathcal L(A)=\big((a_i,b_i,D_i)\big)_{1\le i\le g}.
\]
In particular, $\cN_A$ depends only on $q$ and on the multiplicity sequence $\ms(A)$.
\end{theorem}
\begin{proof}
We use induction on $m=\mult(A)$. If $m=1$, then $A=O$ by Lemma~\ref{lem:noether}, and $\cN_O=1$, so let $m\ge2$, and let $(a,b,D)$, $T$, $S$ be as above. By Theorem~\ref{thm:limit} we have $\cN_A=\cN_{a,b}^D\cN_S$, and the plane branch $S$ has multiplicity $D<m$, since its displayed generators have valuations $D$ and greater than $D$. Thus the product formula follows from the induction hypothesis for $S$, because $\mathcal L(A)=((a,b,D),\mathcal L(S))$. The last assertion follows from Proposition~\ref{prop:list}, and also from Theorem~\ref{thm:E}.
\end{proof}

Let us compare this result with the classical invariants of plane branches. Let $\beta_0<\beta_1<\dots<\beta_g$ be positive integers such that the numbers $e_i=\gcd(\beta_0,\dots,\beta_i)$ satisfy $\beta_0=e_0>e_1>\dots>e_g=1$, put $n_i=e_{i-1}/e_i$, and consider the sequence:
\[
\operatorname{Enr}(\beta_0;\beta_1,\dots,\beta_g)=\Big(e_1\cdot\eu\big(n_1,\tfrac{\beta_1}{e_1}\big),\ e_2\cdot\eu\big(n_2,\tfrac{\beta_2-\beta_1}{e_2}\big),\ \dots,\ e_g\cdot\eu\big(n_g,\tfrac{\beta_g-\beta_{g-1}}{e_g}\big)\Big).
\]
In characteristic zero, the Enriques formula identifies this sequence,
up to trailing $1$'s, with the multiplicity sequence of a plane
branch of Puiseux characteristic $(\beta_0;\beta_1,\dots,\beta_g)$;
see \cite[Thm.~3.2 and its proof, p.~7]{BDF}.
In positive characteristic, Hamburger--Noether expansions provide
the appropriate resolution data; see \cite[Ch.~III, \S\S3.2--3.3,
pp.~76--87]{Campillo}. Exponents visible in a chosen parametrization
need not recover these data when the characteristic divides the
multiplicity, as Section~\ref{sub:wild} illustrates.

The next corollary assumes only the numerical equality
$\ms(A)=\operatorname{Enr}(\beta_0;\beta_1,\dots,\beta_g)$.
Its proof in every characteristic uses this equality, which can
also be checked by successive quadratic transforms.

For $g=0$, we interpret $\operatorname{Enr}(1;)$ as the smooth multiplicity sequence.

\begin{corollary}\label{cor:puiseux}
If $\ms(A)=\operatorname{Enr}(\beta_0;\beta_1,\dots,\beta_g)$, then the cabling list of $A$ and the limit $\cN_A$ are given by
\[
\mathcal L(A)=\Big(\big(n_i,\ b_i,\ e_i\big)\Big)_{1\le i\le g},\qquad b_i=\frac{\beta_0+\beta_i-e_{i-1}}{e_i},\qquad\text{and}\qquad\cN_A=\prod_{i=1}^g\cN_{n_i,b_i}^{\,e_i}.
\]
\end{corollary}
\begin{proof}
We use induction on $g$. For $g=0$, the branch is smooth and both
lists are empty. For $g=1$, we have $e_1=1$, the first cabling triple
is $(\beta_0,\beta_1,1)$, and its reduction is $O$, so the assertion
holds.

Assume $g\ge2$. The first block $e_1\cdot\eu(n_1,\beta_1/e_1)$ begins with $k=\lfloor\beta_1/\beta_0\rfloor$ entries equal to $\beta_0$, followed by $\beta_1-k\beta_0<\beta_0$, thus by Proposition~\ref{prop:list} the first cabling triple is $(n_1,\beta_1/e_1,e_1)$, in agreement with the formula for $b_1$. With $h_1=n_1+\beta_1/e_1-1$, the reduction $S$ has:
\[
\ms(S)=\Big(e_1^{\,h_1},\ e_2\cdot\eu\big(n_2,\tfrac{\beta_2-\beta_1}{e_2}\big),\dots\Big)=\operatorname{Enr}\big(e_1;\ h_1e_1+\beta_2-\beta_1,\ \dots,\ h_1e_1+\beta_g-\beta_1\big).
\]
Indeed, $\eu(n_2,h_1n_2+c)=(n_2^{\,h_1},\eu(n_2,c))$, and the successive differences and the greatest common divisors of the new data are $\beta_{i+1}-\beta_i$ and $e_{i+1}$. By induction, the cabling list of $S$ consists of the triples $(n_{i+1},b'_i,e_{i+1})$ with $b'_i=(e_1+h_1e_1+\beta_{i+1}-\beta_1-e_i)/e_{i+1}$, and finally $h_1e_1=\beta_0+\beta_1-e_1$ gives $b'_i=b_{i+1}$.
\end{proof}

\begin{example}\label{ex:46}
For the Puiseux characteristic $(4;6,7)$ we have $e=(4,2,1)$, $\operatorname{Enr}=(4,2,2,1,1)$ and $\mathcal L=((2,3,2),(2,9,1))$, thus $\cN_A=\cN_{2,3}^{\,2}\,\cN_{2,9}$, and for $(12;16,22,25)$ we have $e=(12,4,2,1)$ and $\mathcal L=((3,4,4),(2,15,2),(2,35,1))$.

For comparison with the link of a complex branch of this
characteristic, let $\overline\beta_0,\ldots,\overline\beta_g$
be its minimal semigroup generators. Define $\delta^{\mathrm{tr}}_{i-1}$
using the primitive truncation with characteristic
\[
 \left(\frac{\beta_0}{e_{i-1}};
       \frac{\beta_1}{e_{i-1}},\ldots,
       \frac{\beta_{i-1}}{e_{i-1}}\right),
 \qquad \delta^{\mathrm{tr}}_0=0.
\]
Put $S_{i-1}=\sum_{1\le j<i}(n_j-1)\overline\beta_j$.
The semigroup recurrence and conductor formula give
\[
 \overline\beta_i=\beta_i+S_{i-1},\qquad
 2\delta^{\mathrm{tr}}_{i-1}
    =\frac{S_{i-1}-\beta_0+e_{i-1}}{e_{i-1}};
\]
see \cite[Thm.~4.1, Lemma~4.2 and Cor.~4.3, pp.~9--10]{BDF}.
Consequently,
\[
 \frac{\overline\beta_i}{e_i}-2n_i\delta^{\mathrm{tr}}_{i-1}
    =\frac{\beta_0+\beta_i-e_{i-1}}{e_i}=b_i.
\]

We use the cabling convention in which $(n,m)$ represents
$n$ times the preferred longitude plus $m$ times the meridian.
The topological pairs are then $(n_i,\overline\beta_i/e_i)$
\cite[Rem.~2.1, pp.~6--7]{Casals}, whereas our pairs are
$(n_i,b_i)$. In positive characteristic, this comparison is
with the complex characteristic data labeling the multiplicity
sequence, not with a topological link of the finite-field germ.
\end{example}


\section{Torus knots}\label{sec:torus}

The reduction of the preceding section expresses the stable invariant
in terms of torus branches. The finite-rank gap-poset formula and
infinite-product identity of Theorem~\ref{thm:torus} therefore give an
explicit product for every plane branch. We also derive
$\dinv(\mathbf n)\ge n_f$ from the module-mass estimate of
Theorem~\ref{thm:finite}; this bound controls the individual
power-series coefficients in Section~\ref{sec:modular}.

Let $2\le a<b$ be coprime and put $\delta=\delta_{a,b}$.
Write $G=\N\smallsetminus\langle a,b\rangle$ for the set of gaps,
partially ordered by $i\preceq j$ if $j-i\in\langle a,b\rangle$,
and let $f=ab-a-b$ be its largest element. The semigroup
$\langle a,b\rangle$ is symmetric: $f-i\in\langle a,b\rangle$
for every $i\in G$. Indeed, if $i=ua+vb$ with $0\le u<b$, then
$v<0$ and $f-i=(b-1-u)a+(-1-v)b$. Hence $i\preceq f$ for all
$i\in G$. Let $C$ be the set of
$\mathbf n=(n_i)_{i\in G}\in\Z_{\ge0}^G$ with $n_i\le n_j$ for
$i\preceq j$; in particular, $n_i\le n_f$ for $\mathbf n\in C$
and $i\in G$. For $\mathbf n\in C$ and $n\in\N\cup\{\infty\}$ put
\[
\qbinom n{\mathbf n}_{x;G}=\frac{(x;x)_n}{(x;x)_{n-n_f}}\prod_{i\in G}\frac{(x;x)_{n_i-n_{i-a-b}}}{(x;x)_{n_i-n_{i-a}}(x;x)_{n_i-n_{i-b}}},\qquad\dinv(\mathbf n)=\sum_{i,j\in G}K(j-i)\,n_in_j,
\]
where $n_j=0$ for $j<0$ and $K(d)=\mathbf 1_{d\ge0}-\mathbf 1_{d\ge a}-\mathbf 1_{d\ge b}+\mathbf 1_{d\ge a+b}$. For $n=\infty$ the first factor is omitted, and if $n_f>n$, the coefficient is $0$. The coefficient $\qbinom n{\mathbf n}_{x;G}$ is a product of $\delta$ Gaussian binomial coefficients \cite[Lemma~5.9]{HJO}, namely:
\[
\qbinom n{\mathbf n}_{x;G}=\prod_{d\in G}\qbinom{\hat n_{d+b}-\hat n_{d-a}}{n_d-\hat n_{d-a}}_x,
\]
where $\hat n_j=n_j$ for $j\in G$, $\hat n_j=n$ for
$j\in\langle a,b\rangle$ and $\hat n_j=0$ for $j<0$.
At $n=\infty$ we use the limiting convention
$\qbinom\infty r_Q=(Q;Q)_r^{-1}$.

We verify this product form for use in
Propositions~\ref{prop:dinvmodules} and~\ref{prop:torus-positive} and
Lemma~\ref{lem:torusformal}. For $i\in G$ the numbers
$i-a$, $i-b$ and $i-a-b$ are gaps or negative. Expanding the Gaussian
coefficients, the product form reads
\[
 \prod_{d\in G}\frac{(x;x)_{\hat n_{d+b}-\hat n_{d-a}}}
   {(x;x)_{n_d-\hat n_{d-a}}\,(x;x)_{\hat n_{d+b}-n_d}} .
\]
The factors $(x;x)_{n_d-\hat n_{d-a}}^{-1}$ are the factors
$(x;x)_{n_i-n_{i-a}}^{-1}$ of the definition. If $d\in G$ and
$i=d+b\in G$, the two remaining factors for $d$ are
$(x;x)_{n_i-n_{i-a-b}}/(x;x)_{n_i-n_{i-b}}$, the remaining factors of
the definition for the gap $i$. Every gap $i>b$ arises in this way
exactly once, and for a gap $i<b$ the corresponding factors of the
definition are $(x;x)_{n_i}/(x;x)_{n_i}=1$.

The gaps $d$ with
$d+b\in\langle a,b\rangle$ are the numbers $d_s=sa-b$ with
$b/a<s\le b-1$: indeed, a gap has a unique expression $sa+vb$ with
$0\le s<b$ and $v<0$, and $d+b$ lies in the semigroup if and only if
$v=-1$. They satisfy $d_s-a=d_{s-1}$, the smallest of them has
$d_s-a<0$, and the largest is $d_{b-1}=f$. For these gaps
$\hat n_{d+b}=n$, and their remaining factors
$(x;x)_{n-\hat n_{d_{s-1}}}/(x;x)_{n-n_{d_s}}$ telescope to
$(x;x)_n/(x;x)_{n-n_f}$, the first factor of the definition. All
indices are nonnegative, because $d-a\preceq d\preceq d+b$ whenever
these numbers are gaps, and $n_d\le n_f\le n$.

For $n=\infty$ the product form follows from the case of finite
$n\ge n_f$ by letting $n\to\infty$. Indeed, for $|x|<1$, and also
coefficientwise in $\Z[[x]]$, the first factor $(x;x)_n/(x;x)_{n-n_f}$
of the definition tends to $1$. Every factor
$\qbinom{n-\hat n_{d-a}}{n_d-\hat n_{d-a}}_x$ with $d+b\in\langle a,b\rangle$
tends to $(x;x)_{n_d-\hat n_{d-a}}^{-1}$, as stipulated by the
limiting convention. The remaining factors on both sides do not
depend on $n$.

The exponents in the product formula below are given by
\[
\rho_{a,b}(j)=\min\{a,\ b,\ \operatorname{dist}(aj,(a+b)\Z)\}-1+\mathbf 1_{(a+b)\mid j}\qquad(j\ge1).
\]

\Needspace{10\baselineskip}
\begin{theorem}[Huang--Jiang--Oblomkov, Huang, Huang--Lau--Ono]\label{thm:torus}
Let $2\le a<b$ be coprime and $F=\F_q$, $Q=q^{-1}$.
\begin{enumerate}[label=\textup{(\roman*)}]
\item \textup{\cite[Corollary~3.2]{HJO}} For every $n\ge1$ we have:
\[
E_{a,b}(n)=\#E_{F[[t^a,t^b]]}(O^n)=q^{n^2\delta}\sum_{\mathbf n\in C,\ n_f\le n}\qbinom n{\mathbf n}_{Q;G}\,Q^{\dinv(\mathbf n)} .
\]
\item \textup{\cite[Theorem~1.3]{Huang2}} For all $\mathbf n\in C$ we have $\dinv(\mathbf n)\ge n_f^2/\delta$.
\item \textup{\cite[Cor.~1.2]{HLO}} $\displaystyle\cN_{a,b}=\prod_{j\ge1}\big(1-Q^j\big)^{-\rho_{a,b}(j)}$.
\end{enumerate}
\end{theorem}
\noindent Part~(i) is \cite[Corollary~3.2]{HJO}, specialized from motives to point
counts and to $t=1$. For a fixed dimension vector $\mathbf n$, the
homogeneous initial modules form the poset flag variety
$\operatorname{Fl}_G(\mathbf n;n)$, with
$\left[\begin{smallmatrix}n\\\mathbf n\end{smallmatrix}\right]_{q;G}$
rational points. The map $L\mapsto\operatorname{in}(L)$ has an
affine-space fibre over each fixed initial flag; its dimension
$B(\mathbf n;n)$ depends only on $\mathbf n$
\cite[Theorem~6.9]{HJO}. Thus the stratum $X_{\mathbf n}$ with this
initial dimension vector has
\[
 \#X_{\mathbf n}(\mathbb F_q)=q^{B(\mathbf n;n)}
       \left[\begin{matrix}n\\\mathbf n\end{matrix}\right]_{q;G}.
\]
The affine-space assertion concerns a fibre over a fixed flag; the
whole stratum need not be affine space. The rectangular-grid
acyclicity theorem computes the fibre dimension, and rewriting the
sum in $Q=q^{-1}$ gives~(i).
Since $\left[\begin{smallmatrix}n\\\mathbf n\end{smallmatrix}\right]_{Q;G}$
is nondecreasing in $n$ and tends to
$\left[\begin{smallmatrix}\infty\\\mathbf n\end{smallmatrix}\right]_{Q;G}$,
part~(i) and monotone convergence give

\begin{equation}\label{eq:Nabsum}
\cN_{a,b}=\sum_{\mathbf n\in C}\qbinom\infty{\mathbf n}_{Q;G}\,Q^{\dinv(\mathbf n)},
\end{equation}
and the sum converges by Theorem~\ref{thm:finite}.
The positive-definiteness in (ii), conjectured in \cite{HJO} and
proved by Huang \cite{Huang2}, also implies this convergence;
this is the finiteness argument used in \cite{HJO}. A further
proof of boundedness follows from the interpretation of the sums
in (i) as generating functions of bounded cylindric partitions
\cite[Thm.~1.1]{HLO}. Part~(ii) is not needed below.

Part~(iii) is the power-series identity $Z_{a,b}=P_{a,b}$ of
\cite[Cor.~1.2]{HLO}, where $P_{a,b}$ denotes the product in (iii).
The series $Z_{a,b}$ of \cite[(1.6)--(1.10)]{HLO} is the right
side of \eqref{eq:Nabsum}: its quadratic form is $\dinv$, and its
factor $P_G(\mathbf n)$ is $\qbinom\infty{\mathbf n}_{Q;G}$.
The identity holds for $|Q|<1$ \cite[\S2.3]{HLO}, so it applies
at $Q=q^{-1}$. Conjectured in \cite{HJO}, it specializes to the
Andrews--Gordon identity for $a=2$ and the first Rogers--Ramanujan
identity for $(a,b)=(2,3)$. Lau and Ono \cite{LO} proved the
complete case $a=3$ before the general proof in \cite{HLO}.

Combining Theorem~\ref{thm:torus}(iii) with
Theorem~\ref{thm:product} gives the product formula of
Theorem~\ref{thm:introD}. Its assertion about the Puiseux
characteristic is Corollary~\ref{cor:puiseux}.

\begin{corollary}\label{cor:finalproduct}
For every plane branch $A$ over $\F_q$, in any characteristic, with cabling list $\mathcal L(A)=((a_i,b_i,D_i))_{1\le i\le g}$, we have:
\[
\cN_A=\prod_{i=1}^g\ \prod_{j\ge1}\big(1-q^{-j}\big)^{-D_i\,\rho_{a_i,b_i}(j)} .
\]
\end{corollary}

\subsection{The inequality \texorpdfstring{$\dinv(\mathbf n)\ge n_f$}{dinv(n) >= n\_f}}
Formal stabilization requires the bound $\dinv(\mathbf n)\ge n_f$
for every $n_f$, including its small values. We prove it from the
conductor-module estimate by interpreting the gap sum at fixed $n_f$.

\Needspace{16\baselineskip}
\begin{proposition}\label{prop:dinvmodules}
Let $A=F[[t^a,t^b]]$ and $\fc=t^{2\delta}O$, and for $s\ge0$ let us put:
\[
\mathsf Z_s(Q)=\sum_{\mathbf n\in C,\ n_f=s}\qbinom\infty{\mathbf n}_{Q;G}\,Q^{\dinv(\mathbf n)} .
\]
Then for every $s\ge0$ we have, in the notation of Proposition~\ref{prop:conductor}:
\[
\mathsf Z_s(q^{-1})=\sum_{[M]:\ \fc M=0,\ s(M)=s}\frac1{|\Aut_A(M)|} .
\]
Consequently $\dinv(\mathbf n)\ge n_f$ for all $\mathbf n\in C$.
For every real $0<Q<1$ the infinite gap sum satisfies
\begin{equation}\label{eq:gap-sum-bound}
 \sum_{\mathbf n\in C}
   \qbinom\infty{\mathbf n}_{Q;G}Q^{\dinv(\mathbf n)}
 \le (Q;Q)_\infty^{-\delta}
        \sum_{s\ge0}(s+1)^{\delta-1}Q^s<\infty.
\end{equation}
\end{proposition}
\begin{proof}
Denote the right-hand side by $\mathsf m_s$. This is a finite sum:
the proof of Proposition~\ref{prop:conductor}(ii) shows that every
module $M$ with $\fc M=0$ and $s(M)=s$ is isomorphic to a submodule
of $(O/A)^s$. Also, $\mathsf m_0=1=\mathsf Z_0$. Since
$\prod_{i=0}^{s-1}(1-q^{i-n})$ equals $(Q;Q)_n/(Q;Q)_{n-s}$ for
$s\le n$ and $0$ for $s>n$, Proposition~\ref{prop:conductor}(ii)
and Theorem~\ref{thm:torus}(i) give, for all $n\ge1$ and $Q=q^{-1}$,
\[
\sum_{s=0}^n\frac{(Q;Q)_n}{(Q;Q)_{n-s}}\,\mathsf m_s=e_A(n)=\sum_{s=0}^n\frac{(Q;Q)_n}{(Q;Q)_{n-s}}\,\mathsf Z_s(Q),
\]
and the identity $\mathsf m_s=\mathsf Z_s(q^{-1})$ follows by induction on $s$.

The inequality $\dinv(\mathbf n)\ge n_f$ is immediate when $n_f=0$.
Now let $\mathbf n\in C$ with $n_f=s\ge1$. All terms of
$\mathsf Z_s$ are positive for $0<Q<1$, and by
Remark~\ref{rem:finite}(2), the bound \eqref{eq:soclebound} applies
to $A$ with $k=a-1$. Hence for every prime power $q$ and $Q=q^{-1}$,
\[
\qbinom\infty{\mathbf n}_{Q;G}\,Q^{\dinv(\mathbf n)}\ \le\ \mathsf Z_s(Q)=\mathsf m_s\ \le\ Q^s\,(Q;Q)_\infty^{-(a-1)}\sum_{j\ge0}\binom{j+s+a-2}{a-2}Q^j .
\]
As $q\to\infty$, the coefficient $\qbinom\infty{\mathbf n}_{Q;G}$ tends to $1$, and the factor of $Q^s$ on the right tends to $\binom{s+a-2}{a-2}$. If we had $\dinv(\mathbf n)<s$, the left side divided by $Q^s$ would tend to infinity, which is impossible.

For \eqref{eq:gap-sum-bound}, each generalized Gaussian coefficient
is a product of $\delta$ Gaussian binomial coefficients, by the
product form proved before
Theorem~\ref{thm:torus}. Each of these factors is at most
$(Q;Q)_\infty^{-1}$ for $0<Q<1$, including the limiting factors.
There are at most $(s+1)^{\delta-1}$ gap vectors with $n_f=s$,
because every coordinate lies between $0$ and $s$.
Summing the resulting bound gives \eqref{eq:gap-sum-bound}.
\end{proof}

\begin{remark}\label{rem:dinv}
(1) The proof shows that the inequality $\dinv(\mathbf n)\ge n_f$ is equivalent to a bound for the number of modules: by Proposition~\ref{prop:dinvmodules}, the modules of finite length over $A/\fc$ with a socle of dimension $s$, counted with the weights $|\Aut_A(M)|^{-1}$, have total mass $\mathsf Z_s(q^{-1})$, and this mass is $O(q^{-s})$ as $q\to\infty$ by Theorem~\ref{thm:finite}. The stronger inequality of Theorem~\ref{thm:torus}(ii) says that the mass is $O(q^{-s^2/\delta})$.

(2) By Proposition~\ref{prop:conductor}(iii), the sum \eqref{eq:Nabsum} is the sum of the numbers $|\Aut_A(M)|^{-1}$ over the modules of finite length over the finite ring $A/\fc$. For $(a,b)=(2,2k+1)$ this ring is $F[x]/(x^k)$, and by the formula for the order of the automorphism group of a module over this ring \cite[Ch.~II]{Mac}, the sum over its modules is the sum side of the Andrews--Gordon identity. For $(a,b)=(3,4)$ it is $F[u,v]/(u,v)^2$, and Theorem~\ref{thm:torus}(iii) becomes the identity:
\[
\sum_{d\ge0}\frac1{|\GL_d(F)|}\sum_{r=0}^d\qbinom dr_q\prod_{i=0}^{r-1}\big(q^{2(d-r)}-q^i\big)=\prod_{j\ge1}\big(1-Q^j\big)^{-\rho_{3,4}(j)},
\]
in which the inner sum is the number of pairs $(U,V)$ of matrices in $\Mat_d(F)$ with $U^2=UV=VU=V^2=0$.
\end{remark}

\subsection{Coefficient positivity for torus branches}
Theorem~\ref{thm:torus}(i) implies Conjecture~P for the torus
sequences.

\Needspace{12\baselineskip}
\begin{proposition}\label{prop:torus-positive}
Let $2\le a<b$ be coprime and $n\ge1$. For $\mathbf n\in C$ with
$n_f\le n$ and $d\in G$ put $N_d=\hat n_{d+b}-\hat n_{d-a}$ and
$k_d=n_d-\hat n_{d-a}$, and let
$B'(\mathbf n;n)=n^2\delta-\dinv(\mathbf n)-\sum_{d\in G}k_d(N_d-k_d)$.
Then
\[
 P_{\eu(a,b),n}(q)=\sum_{\mathbf n\in C,\ n_f\le n}
   q^{B'(\mathbf n;n)}\prod_{d\in G}\qbinom{N_d}{k_d}_q
\]
as Laurent polynomials in $q$. In particular,
$P_{\eu(a,b),n}\in\Z_{\ge0}[X]$.
\end{proposition}
\begin{proof}
Since $\qbinom Nk_{q^{-1}}=q^{-k(N-k)}\qbinom Nk_q$, the product
form displayed before Theorem~\ref{thm:torus} and
Theorem~\ref{thm:torus}(i) show that the right-hand side of the
asserted identity equals $E_{a,b}(n)$ at every prime power $q$. By
Example~\ref{ex:torus} and Corollary~\ref{cor:poly} we have
$E_{a,b}(n)=P_{\eu(a,b),n}(q)$. Two Laurent polynomials which agree
at all prime powers are equal. Every summand on the right is a
Laurent polynomial with nonnegative coefficients, hence the
polynomial $P_{\eu(a,b),n}$ has nonnegative coefficients.
\end{proof}

\section{Power series and modularity}\label{sec:modular}

The product of Corollary~\ref{cor:finalproduct} is the
coefficientwise limit of the normalized counting polynomials,
with agreement modulo $Q^{n+1}$ in rank $n$
(Theorem~\ref{thm:formal}). After multiplication by an explicit
rational power of $Q$, it is also a modular unit
(Theorem~\ref{thm:modular}). We first use polynomiality of the
refined counts to study the cabling systems in the $Q$-adic
topology. We then set $Q=e^{2\pi i\tau}$ and use the periodic
product exponents and transformation formulas for generalized
eta functions.

\subsection{The cabling formula for power series}
For coprime $2\le a<b$ we use the notation of Section~\ref{sec:torus}, and we put:
\[
P_{a,b}(Q)=\prod_{j\ge1}(1-Q^j)^{-\rho_{a,b}(j)},\qquad F_{a,b;j}(Q)=\sum_{\mathbf n\in C,\ n_f\le j}\qbinom j{\mathbf n}_{Q;G}Q^{\dinv(\mathbf n)}\quad(j\ge0).
\]
By Theorem~\ref{thm:torus}(i),
$e_{F[[t^a,t^b]]}(j)=F_{a,b;j}(q^{-1})$; also $F_{a,b;0}=1$.

\begin{lemma}\label{lem:torusformal}
\begin{enumerate}[label=\textup{(\roman*)}]
\item We have $\dinv(\mathbf n)\ge n_f$ for all $\mathbf n\in C$.
\item $F_{a,b;j}\in\Z[Q]$, and $\sum_{\mathbf n\in C}\qbinom\infty{\mathbf n}_{Q;G}Q^{\dinv(\mathbf n)}=P_{a,b}(Q)$ in $\Z[[Q]]$.
\item $F_{a,b;j}\equiv P_{a,b}$ modulo $Q^{\,j+1}$ for all $j\ge0$.
\end{enumerate}
\end{lemma}
\begin{proof}
(i) This is Proposition~\ref{prop:dinvmodules}.

(ii) Since $\dinv\ge0$ on $C$ and the coefficients are products
of Gaussian binomial coefficients by the product form proved before
Theorem~\ref{thm:torus}, we have $F_{a,b;j}\in\Z[Q]$. By (i),
only finitely many $\mathbf n\in C$ contribute to each coefficient
of the left side. Both sides are therefore power series with
nonnegative coefficients, and both converge for $|Q|<1$, the
left side by \eqref{eq:gap-sum-bound}. They agree at $Q=q^{-1}$
for all prime powers $q$ by \eqref{eq:Nabsum} and
Theorem~\ref{thm:torus}(iii), so they are equal.

(iii) For $\mathbf n\in C$ with $n_f\le j$ we have:
\[
\qbinom j{\mathbf n}_{Q;G}-\qbinom\infty{\mathbf n}_{Q;G}=\qbinom\infty{\mathbf n}_{Q;G}\Big(\frac{(Q;Q)_j}{(Q;Q)_{j-n_f}}-1\Big)\in Q^{\,j-n_f+1}\Z[[Q]].
\]
Thus, by (ii), $F_{a,b;j}-P_{a,b}$ is a sum of terms divisible by $Q^{\dinv(\mathbf n)+j-n_f+1}$, for $1\le n_f\le j$, or by $Q^{\dinv(\mathbf n)}$, for $n_f>j$, and by (i) the exponent is at least $j+1$ in both cases.
\end{proof}

Let $\underline\mu$ be the multiplicity sequence of a plane branch
over a finite field. Write $\delta(\underline\mu)=\sum_j\binom{\mu_j}2$
for its $\delta$-invariant (Lemma~\ref{lem:noether}),
$P_{\underline\mu,n}$ for the polynomial of Corollary~\ref{cor:poly},
and $\mathcal L(\underline\mu)=((a_i,b_i,D_i))_{1\le i\le g}$ for
its cabling list (Proposition~\ref{prop:list}). Put
\[
p_{\underline\mu,n}(Q)=Q^{n^2\delta(\underline\mu)}P_{\underline\mu,n}(Q^{-1}),\qquad\cN_{\underline\mu}(Q)=\prod_{i=1}^gP_{a_i,b_i}(Q)^{D_i}\in1+Q\,\Z_{\ge0}[[Q]],
\]
so that, for every plane branch $A$ over $\F_q$ with $\ms(A)=\underline\mu$,
Corollaries~\ref{cor:poly} and~\ref{cor:finalproduct} give
$e_A(n)=p_{\underline\mu,n}(q^{-1})$ and
$\cN_A=\cN_{\underline\mu}(q^{-1})$, respectively.

\begin{theorem}\label{thm:formal}
For every $n\ge1$, $p_{\underline\mu,n}$ is a polynomial in $Q$ with constant term $1$, and
\[
p_{\underline\mu,n}(Q)\equiv\cN_{\underline\mu}(Q)\pmod{Q^{\,n+1}} .
\]
In particular, the normalized counting polynomials $p_{\underline\mu,n}$ converge coefficientwise to $\cN_{\underline\mu}$ as $n\to\infty$.
\end{theorem}
\begin{proof}
Fix a plane branch $A_0$ over $\F_{q_0}$ with
$\ms(A_0)=\underline\mu$, and for $r\ge1$ let
$A_r=A_0\otimes\F_{q_0^r}$. Generators of $A_0$ as in
Lemma~\ref{lem:normalgen} remain such generators after scalar
extension. The branches $T$ and $S$ of Section~\ref{sub:iterate}
for $A_r$ are likewise obtained from those of $A_0$ by extension
of scalars. We use repeatedly the fact that two rational functions of $Q$
which agree at $Q=q_0^{-r}$ for every $r\ge1$ are equal.
The integrality observations in Section~\ref{sec:limit} show
that all the normalizations belong to $\Q(Q)$.

\emph{Step 1. Polynomials and their leading terms.}
For any plane branch $A=F[[x,y]]$, the summand $L=A^n$ in
Proposition~\ref{prop:finitesum} gives $e_A(n)\ge1$.
Choose $k\ge1$ with $y^k$ in the conductor. Theorem~\ref{thm:finite}
then gives the bounds
\begin{equation}\label{eq:sandwich}
1\ \le\ e_A(n)\ \le\ \cN_A\ \le\ (Q;Q)_\infty^{-k},
\qquad Q=q^{-1}.
\end{equation}
For the extensions $A_r$ of the fixed branch $A_0$, the same $k$
works for every $r$, since the generators and conductor extend by
scalars. As $Q=q_0^{-r}\to0$, the upper bound tends to $1$.
Consequently $p_{\underline\mu,n}(q_0^{-r})\to1$.
Corollary~\ref{cor:poly} makes $p_{\underline\mu,n}$ a Laurent
polynomial in $Q$. A Laurent polynomial with this limit has no
negative powers and has constant term $1$. Equivalently,
$P_{\underline\mu,n}$ is monic of degree $n^2\delta(\underline\mu)$.
Only the conductor bound and polynomiality are used here.

The same bound handles residual types. Fix a framed branch
$(T;\varpi,\omega)$ over $\F_{q_0}$ and a partition
$\mu\vdash n\val\varpi$. By Corollary~\ref{cor:poly}, its
normalized refined counts over $\F_{q_0^r}$ are values of one Laurent
polynomial in $Q$. These values lie between $0$ and $e_T(n)$.
Apply \eqref{eq:sandwich} to $T$ with a fixed $k_T$ valid under
scalar extension. Boundedness as $r\to\infty$ excludes negative powers, so
every such normalized refined count is also a polynomial in $Q$.

\emph{Step 2. Reduction.} We induct on $m=\mu_0$. When $m=1$,
Lemma~\ref{lem:noether} gives $A_0=O$ and
$p_{\underline\mu,n}=1=\cN_{\underline\mu}$. Let $m\ge2$, let
$(a,b,D)$ be the first cabling triple, and let $T$ and $S$ be the
pattern ring and the reduction of $A_0$. Put $h=a+b-1\ge4$,
$C=P_{a,b}(Q)$, $p^A_n=p_{\underline\mu,n}$ and
$p^S_n=p_{\ms(S),n}$. We claim that
\begin{equation}\label{eq:formalstep}
p^A_n\equiv C^D\,p^S_n\pmod{Q^{\,n+1}} .
\end{equation}
Since $\cN_{\underline\mu}=C^D\cN_{\ms(S)}$, the theorem follows
from \eqref{eq:formalstep} and the induction hypothesis for $S$:
write $p^A_n-C^D\cN_{\ms(S)}=(p^A_n-C^Dp^S_n)+C^D(p^S_n-\cN_{\ms(S)})$.
If $D=1$, then $S=O$ and $p^S_n=1$, while
$p^A_n=F_{a,b;n}$ by Theorem~\ref{thm:cable}. In this case,
\eqref{eq:formalstep} is Lemma~\ref{lem:torusformal}(iii).

\emph{Step 3. The case $D\ge2$.} We use the notation of Section~\ref{sec:limit}, with $N=Dn$, and consider the quantities $\widehat E_j$, $\pi_{\lambda\nu}$, $\widehat\Phi_\lambda$, $\widehat\Phi^\circ_\lambda$, $b_n(\lambda)$, $e_A(n)$ and $e_S(n)$, computed over $\F_{q_0^r}$. They are the values at $Q=q_0^{-r}$ of elements of $\Q(Q)$ which do not depend on $r$: for $\widehat E_j=F_{a,b;j}$ this is Theorem~\ref{thm:torus}(i), and for $\pi_{\lambda\nu}$ it follows from Lemma~\ref{lem:hallbasics}(a) and the polynomiality of the Hall numbers and of the numbers $a_\mu$.

The numbers $\Phi^{a,b}_{q,N}(\mu)$ and $\Phi^{1,h}_{q,N}(\mu)$ are also the values at $Q=q_0^{-r}$ of elements of $\Q(Q)$ which do not depend on $r$. Indeed, by Definition~\ref{def:Phi} each of them is a sum of products $c_\beta H_\beta(\mu)$. The coefficients $c_\beta$ are the values at $q$ of rational functions which do not depend on the field by Corollary~\ref{cor:basis}(ii), and $H_\beta(\mu)$, the coefficient of $u_\mu$ in $v_{\beta_1}\cdots v_{\beta_l}$ (see the proof of Lemma~\ref{lem:hallbasics}(d)), is a sum of products of Hall numbers. The remaining quantities have the same property by Theorem~\ref{thm:E}. Therefore the identities \eqref{eq:two-systems} and \eqref{eq:normalized-counts} hold in $\Q(Q)$. We embed this field in $\Q((Q))$ and denote by $\val_Q$ the $Q$-adic valuation.

(a) Lemma~\ref{lem:triangular}(ii) gives
$\pi_{\lambda\lambda}\in1+Q\Z[Q]$. Also,
$\pi_{\lambda\nu}\in\Q[Q]$: the number
$q^{e_2(\lambda)}\pi_{\lambda\nu}=\#(\fn_\lambda(F)\cap\bO_{\nu'})$
is the value at $q$ of a rational function, hence of a polynomial
by Lemma~\ref{lem:ratpoly}, and is at most $q^{e_2(\lambda)}$.
By Lemma~\ref{lem:triangular}(i), the rational function
$\pi_{\lambda\nu}$ vanishes unless $\nu\ctd\lambda$. The
expansion of $A^{-1}$ in the proof of Lemma~\ref{lem:inverse}
is therefore valid in $\Q(Q)$. Its Gaussian factors telescope,
giving, for $\lambda=n+\alpha$ and $\nu=n+\beta$ in $\cP_n$,
\[
(A^{-1})_{\lambda\nu}\in Q^{\frac h2(\|\beta\|^2-\|\alpha\|^2)}\,\Q[[Q]] .
\]

(b) By \eqref{eq:two-systems}, the elements $\epsilon_\lambda=\widehat\Phi_\lambda-C^D\widehat\Phi^\circ_\lambda$ satisfy $\epsilon_\lambda=\sum_{\nu\in\cP_n}(A^{-1})_{\lambda\nu}\big(\prod_{i=1}^D\widehat E_{\nu_i}-C^D\big)$, and Lemma~\ref{lem:torusformal}(iii) together with the identity $\prod_ix_i-C^D=\sum_i(\prod_{l<i}x_l)(x_i-C)C^{D-i}$ in $\Z[[Q]]$ gives:
\[
\val_Q\Big(\prod_{i=1}^D\widehat E_{\nu_i}-C^D\Big)\ \ge\ \min_i\,(\nu_i+1)=n+\beta_D+1,
\]
which, together with (a) and $\|\beta\|^2\ge\beta_D^2$, gives
\[
\val_Q\big(Q^{\frac h2\|\alpha\|^2}\epsilon_\lambda\big)\ \ge\ n+1+\min_{\nu}\Big(\tfrac h2\beta_D^2+\beta_D\Big)\ \ge\ n+1,
\]
because $\beta_D$ is an integer and $h\ge2$.

(c) By \eqref{eq:normalized-counts},
$p^A_n-C^Dp^S_n=\sum_{\lambda\in\cP_n}b_n(\lambda)\,Q^{\frac h2\|\alpha\|^2}\epsilon_\lambda$,
where $b_n(\lambda)\in\Q[Q]$ by Step~1. Part~(b) now proves
\eqref{eq:formalstep}.
\end{proof}

\begin{remark}\label{rem:formal}
(1) This is the $Q$-adic counterpart of the proof of
Theorem~\ref{thm:limit}. Positivity enters only through the bounds
of Steps~1 and~3(a). The inverse-matrix estimate is immediate
here, whereas the real estimate required Lemma~\ref{lem:inverse}.
Theorem~\ref{thm:formal} does, however, require all terms to be
values of rational functions of $q$ independent of the field
(Theorem~\ref{thm:E} and Corollary~\ref{cor:poly});
Theorem~\ref{thm:limit} concerns a single finite field.

(2) By Theorem~\ref{thm:formal}, the power series $\cN_{\underline\mu}(Q)$ is intrinsically attached to the equisingularity class, and for a plane branch $A$ we write $\cN_A(Q)=\cN_{\ms(A)}(Q)$. This series converges for $|Q|<1$, has nonnegative integer coefficients, and satisfies $\cN_A=\cN_A(q^{-1})$.
\end{remark}

\subsection{The exponents}
\begin{lemma}\label{lem:rho}
Let $2\le a<b$ be coprime and $M=a+b$, let $\rho=\rho_{a,b}$, extended to all integers by the formula of Section~\ref{sec:torus}, and let $a^*$ be an inverse of $a$ modulo $M$.
\begin{enumerate}[label=\textup{(\roman*)}]
\item For all $g\in\Z$ we have:
\[
\rho(g)=(a-1)\,\mathbf 1_{M\nmid g}-\sum_{j=1}^{a-1}(a-j)\,\mathbf 1_{ag\equiv\pm j\ (\mathrm{mod}\ M)} .
\]
In particular, $\rho$ is $M$-periodic, even and nonnegative, and $\rho(0)=0$. Moreover, $\sum_{g=1}^{M-1}\rho(g)=(a-1)(b-1)$, and if $M$ is even then $\rho(M/2)=a-1$ is even.
\item With $g_j\in\{1,\dots,M-1\}$ such that $g_j\equiv a^*j$ modulo $M$ we have:
\[
P_{a,b}(Q)=\frac{(Q^M;Q^M)_\infty^{\,a-1}}{(Q;Q)_\infty^{\,a-1}}\ \prod_{j=1}^{a-1}\Big[(Q^{g_j};Q^M)_\infty\,(Q^{M-g_j};Q^M)_\infty\Big]^{a-j} .
\]
\item $\displaystyle\sum_{g=1}^{M-1}\rho(g)\,B_2\Big(\frac gM\Big)=\frac{(a-1)(b-1)}{6M^2}$, where $B_2(x)=x^2-x+\frac16$.
\end{enumerate}
\end{lemma}
\begin{proof}
Since the distance $d(x)$ from $x$ to $M\Z$ satisfies $d(x)\le M/2<b$, we have $\rho(g)=\min\{a,d(ag)\}-1+\mathbf 1_{M\mid g}$.

(i) Both sides vanish if $M\mid g$, so let $M\nmid g$ and $x=ag$. Then $d(x)\ge1$, and $\min\{a,d(x)\}-1$ is the number of $k\in\{1,\dots,a-1\}$ with $k<d(x)$. Next, $d(x)\le k$ if and only if $x\equiv j$ for some $j$ with $1\le|j|\le k$, and the residues of the integers $j$ with $|j|\le a-1$ are pairwise distinct because $2a-2<M$. For a fixed $j$ there are $a-|j|$ values of $k$ with $|j|\le k\le a-1$, and this gives the formula. The formula implies the evenness of $\rho$ and the identity $\sum_{g=1}^{M-1}\rho(g)=(a-1)(M-1)-a(a-1)$. Finally, if $M$ is even, then $a$ is odd, $a\cdot M/2\equiv M/2$ and $d(M/2)=M/2>a$.

(ii) By (i) we have:
\[
\prod_{n\ge1}(1-Q^n)^{-\rho(n)}=\prod_{M\nmid n}(1-Q^n)^{-(a-1)}\prod_{j=1}^{a-1}\prod_{n\equiv\pm g_j}(1-Q^n)^{a-j}.
\]

(iii) Let us put $\zeta=e^{2\pi i/M}$, $\mathsf t(x)=\min\{a,d(x)\}$ for $x\in\Z/M$, and $I=\{0,1,\dots,a-1\}\subseteq\Z/M$. Since $2a-1\le M$, the number of pairs $(i,i')\in I^2$ with $i-i'=x$ is $\max\{0,a-d(x)\}=a-\mathsf t(x)$, and taking Fourier transforms we obtain for $m\not\equiv0$:
\[
\sum_{x\in\Z/M}\mathsf t(x)\zeta^{mx}=-\Big|\sum_{i\in I}\zeta^{mi}\Big|^2=-\frac{\sin^2(\pi am/M)}{\sin^2(\pi m/M)},\qquad\text{and}\qquad\sum_{x}\mathsf t(x)=aM-a^2=ab .
\]
Since $\rho(g)=\mathsf t(ag)-1+\mathbf 1_{M\mid g}$, the function $\widehat\rho(n)=\sum_{g\in\Z/M}\rho(g)\zeta^{ng}$ is given by:
\[
\widehat\rho(n)=1-\frac{\sin^2(\pi n/M)}{\sin^2(\pi a^*n/M)}\quad(M\nmid n),\qquad\widehat\rho(n)=(a-1)(b-1)\quad(M\mid n).
\]
Next, let us recall that $B_2(\{y\})=\frac1{2\pi^2}\sum_{n\ne0}e^{2\pi iny}/n^2$, that $\sum_{n\equiv r}n^{-2}=\pi^2/(M^2\sin^2(\pi r/M))$ for $M\nmid r$, and that $\sum_{n\ne0,\,M\mid n}n^{-2}=\pi^2/(3M^2)$. Hence we get:
\[
\begin{aligned}
\sum_{g\in\Z/M}\rho(g)B_2\Big(\Big\{\frac gM\Big\}\Big)&=\frac1{2\pi^2}\sum_{n\ne0}\frac{\widehat\rho(n)}{n^2}\\
&=\frac{(a-1)(b-1)}{6M^2}+\frac1{2M^2}\sum_{r=1}^{M-1}\Big(\frac1{\sin^2(\pi r/M)}-\frac1{\sin^2(\pi a^*r/M)}\Big),
\end{aligned}
\]
where the last sum vanishes because $r\mapsto a^*r$ permutes the nonzero residues.
\end{proof}

\subsection{Modularity}
We use the generalized Dedekind eta functions in Yang's
normalization \cite{Yang1,Yang2}. Let $N\ge5$ be an integer,
$\tau\in\mathbb H$ and $Q=e^{2\pi i\tau}$; for rational $x$,
$Q^x$ denotes $e^{2\pi ix\tau}$. With
$B_2(x)=x^2-x+\frac16$, put, for integers $1\le g\le N-1$,
\[
E^{(N)}_g(\tau)=Q^{NB_2(g/N)/2}\prod_{m\ge1}\big(1-Q^{(m-1)N+g}\big)\big(1-Q^{mN-g}\big).
\]
In \cite{Yang2} the variable $Q$ is denoted by $q$, whereas in this paper $q$ is the order of the field. For $0<g<N$ we have $E^{(N)}_{N-g}=E^{(N)}_g$ because $B_2(1-x)=B_2(x)$, and this allows us to extend the definition uniquely to all integers $g\not\equiv0$ modulo $N$ in such a way that $E^{(N)}_{g+N}=E^{(N)}_{-g}=-E^{(N)}_g$. In terms of the Jacobi theta function, $E^{(N)}_g(\tau)=-ie^{\pi ig^2\tau/N}\vartheta_1(g\tau|N\tau)/\eta(N\tau)$ (Proposition~\ref{prop:Etransf}(i)). We use the following two results.
\begin{itemize}[leftmargin=2.2em]
\item For $g\not\equiv0$ modulo $N$ and $\gamma=\left(\begin{smallmatrix}\alpha&\beta\\ \gamma_0N&\delta_0\end{smallmatrix}\right)\in\Gamma_0(N)$ we have:
\begin{equation}\label{eq:yang}
E^{(N)}_g(\gamma\tau)=\varepsilon(\alpha,\beta N,\gamma_0,\delta_0)\ e^{\pi i(g^2\alpha\beta/N-g\beta)}\ E^{(N)}_{\alpha g}(\tau),
\end{equation}
where $\varepsilon$ is given by:
\[
\varepsilon(a,b,c,d)=\begin{cases}e^{\pi i(bd(1-c^2)+c(a+d-3))/6}&\text{if $c$ is odd,}\\[3pt]-i\,e^{\pi i(ac(1-d^2)+d(b-c+3))/6}&\text{if $d$ is odd.}\end{cases}
\]
\item Let $(e_g)$ be a finite family of integers, indexed by integers $g\not\equiv0$ modulo $N$, such that $\sum_ge_g\equiv0$ modulo $12$, $\sum_gge_g\equiv0$ modulo $2$ and $\sum_gg^2e_g\equiv0$ modulo $2N$. Then $\prod_g(E^{(N)}_g)^{e_g}$ is a modular unit on $\Gamma_1(N)$.
\end{itemize}
Yang states these results in \cite[Propositions~1 and~2]{Yang2}
for $N>4$, with the nontranslation transformation law given for
$\gamma_0>0$, and attributes them to
\cite[Corollaries~2 and~3]{Yang1}.
The proof of Proposition~\ref{prop:Etransf} applies in the present
conventions to every integer $N\ge2$ and every matrix in
$\Gamma_0(N)$. The function $\varepsilon(a,b,c,d)$ is the value at $\left(\begin{smallmatrix}a&b\\c&d\end{smallmatrix}\right)$ of a character $\varepsilon$ of $\mathrm{SL}_2(\Z)$ with $\varepsilon^{12}=1$ and $\varepsilon|_{\Gamma(12)}=1$, and $\eta(\gamma\tau)^2=\varepsilon(\gamma)(c\tau+d)\eta(\tau)^2$ (Lemma~\ref{lem:etachar}).

For coprime $2\le a<b$ we put:
\[
c_{a,b}=\frac{(a-1)(b-1)}{24(a+b)},\qquad f_{a,b}(\tau)=Q^{-c_{a,b}}\,P_{a,b}(Q),\qquad Q=e^{2\pi i\tau},
\]
where $Q^{-c}$ stands for $e^{-2\pi ic\tau}$, and for a plane branch $A\ne O$ with cabling list $((a_i,b_i,D_i))_{i\le g}$ we put:
\[
c_A=\sum_{i=1}^gD_i\,c_{a_i,b_i}=\frac1{24}\sum_{i=1}^gD_i\,\frac{(a_i-1)(b_i-1)}{a_i+b_i},\qquad\Theta_A(\tau)=Q^{-c_A}\,\cN_A(Q)=\prod_{i=1}^gf_{a_i,b_i}(\tau)^{D_i}.
\]
Here $\cN_A(Q)=\cN_{\ms(A)}(Q)$ is the power series of Remark~\ref{rem:formal}(2), and we denote by $L_A$ the least common multiple of the numbers $a_i+b_i$. A \emph{modular function} for a congruence subgroup $\Gamma$ is a meromorphic function on $\mathbb H$ which is invariant under $\Gamma$ and meromorphic at every cusp, and a \emph{modular unit} for $\Gamma$ is a modular function for $\Gamma$ without zeros and poles in $\mathbb H$.

\begin{proposition}\label{prop:etaform}
Let $2\le a<b$ be coprime and $M=a+b$. Then
\[
f_{a,b}=\prod_{1\le g<M/2}\big(E^{(M)}_g\big)^{-\rho_{a,b}(g)}\cdot\big(E^{(M)}_{M/2}\big)^{-(a-1)/2},
\]
where the last factor is present only if $M$ is even, and all exponents are integers.
\end{proposition}
\begin{proof}
For $0<g<M$ we consider the product of $(1-Q^n)$ over the integers $n\ge1$ with $n\equiv g$ or $n\equiv-g$ modulo $M$, where each $n$ is taken once. If $g\ne M/2$, this product is $Q^{-MB_2(g/M)/2}E^{(M)}_g$, while at the midpoint we have
\[
 E^{(M)}_{M/2}
 =Q^{MB_2(1/2)/2}
   \left(\prod_{\substack{n\ge1\\n\equiv M/2\; (\mathrm{mod}\,M)}}
      (1-Q^n)\right)^2.
\]
Since $\rho=\rho_{a,b}$ is even with $\rho(0)=0$, and $\rho(M/2)=a-1$ is even if $M$ is even (Lemma~\ref{lem:rho}(i)), we get:
\[
P_{a,b}(Q)=Q^{\frac M4\sum_{g=1}^{M-1}\rho(g)B_2(g/M)}\prod_{1\le g<M/2}\big(E^{(M)}_g\big)^{-\rho(g)}\cdot\big(E^{(M)}_{M/2}\big)^{-\rho(M/2)/2},
\]
where the exponent of $Q$ is $c_{a,b}$ by Lemma~\ref{lem:rho}(iii).
\end{proof}

\begin{theorem}\label{thm:modular}
Let $A\ne O$ be a plane branch over $\F_q$, in any characteristic, and $c=c_A$, $L=L_A$, $\Theta=\Theta_A$.
\begin{enumerate}[label=\textup{(\alph*)}]
\item $\Theta$ is holomorphic and nowhere zero on $\mathbb H$, and $\Theta(\gamma\tau)=\Theta(\tau)$ for all $\gamma\in\Gamma(12L)$.
\item $\Theta^{12L}$ is a modular unit on $\Gamma_1(L)$, and $\Theta$ is a modular unit on $\Gamma(12L)$.
\item $\Theta(\tau)=Q^{-c}(1+O(Q))$ and $\Theta(-1/\tau)=C_A\,Q^{-c}(1+O(Q^{1/L}))$ as $\operatorname{Im}\tau\to\infty$, for an explicit constant $C_A>0$. Thus $\Theta$ has poles of the same order $c$, measured in $Q$, at the cusps $\infty$ and $0$.
\end{enumerate}
By Corollary~\ref{cor:finalproduct} we have $\cN_A=q^{-c}\,\Theta\big(\tfrac{i\log q}{2\pi}\big)$, thus, up to the factor $q^{-c}$, the invariant $\cN_A$ is a value of a modular unit.
\end{theorem}
\begin{proof}
(a) The product $\cN_A(Q)$ converges locally uniformly on $|Q|<1$ and has no zeros, so it remains to prove the invariance. By Proposition~\ref{prop:etaform}, and since $\Gamma(12L)\subseteq\Gamma(12M_i)$ for $M_i=a_i+b_i$, it suffices to show that $E_g=E^{(M)}_g$ is invariant under $\Gamma(12M)$ for every $M\ge5$ and $g\not\equiv0$ modulo $M$.

Let $\gamma=\left(\begin{smallmatrix}\alpha&\beta\\ \gamma_0M&\delta_0\end{smallmatrix}\right)\in\Gamma(12M)$, that is, $\alpha\equiv\delta_0\equiv1$ and $\beta\equiv0$ modulo $12M$, and $\gamma_0\equiv0$ modulo $12$. We apply \eqref{eq:yang} and observe that the matrix $\left(\begin{smallmatrix}\alpha&\beta M\\ \gamma_0&\delta_0\end{smallmatrix}\right)$ lies in $\Gamma(12)$, hence $\varepsilon(\alpha,\beta M,\gamma_0,\delta_0)=1$ by Lemma~\ref{lem:etachar}, and that the integer $g^2\alpha\beta/M-g\beta$ is even, because $12M$ divides $\beta$. Finally, $\alpha g=g+12M\alpha'g$ with $\alpha'\in\Z$, thus $E_{\alpha g}=(-1)^{12\alpha'g}E_g=E_g$ and $E_g(\gamma\tau)=E_g(\tau)$.

(b) By Proposition~\ref{prop:etaform}, $f_{a,b}^{12M}$ is a product of functions $E^{(M)}_g$ with exponents divisible by $12M$, hence it is a modular function on $\Gamma_1(M)$ by the second result quoted above. Thus $\Theta^{12L}=\prod_i\big(f_{a_i,b_i}^{12M_i}\big)^{D_iL/M_i}$ is a modular function on $\Gamma_1(L)\subseteq\Gamma_1(M_i)$, and by (a) it has no zeros or poles in $\mathbb H$.

Finally, for $\sigma\in\mathrm{SL}_2(\Z)$ the function $\tau\mapsto\Theta(\sigma\tau)$ is invariant under $\tau\mapsto\tau+12L$ by (a), because the group $\Gamma(12L)$ is normal in $\mathrm{SL}_2(\Z)$ and contains this translation. In other words, $\Theta(\sigma\tau)$ is a holomorphic function of $Q^{1/(12L)}$ on the punctured unit disc, and it has at most a pole at the origin, because its $12L$-th power does. Thus $\Theta$ is a modular function for $\Gamma(12L)$.

(c) The expansion at $\infty$ follows from the defining product.
For the cusp $0$, fix a torus factor, put $M=a+b$, and set
$\widetilde\rho(g)=\rho_{a,b}(g)$ for $1\le g<M/2$ and
$\widetilde\rho(M/2)=(a-1)/2$ when $M$ is even.
Proposition~\ref{prop:etaform} expresses $f_{a,b}$ as the product of
$(E^{(M)}_g)^{-\widetilde\rho(g)}$, and
Lemma~\ref{lem:rho}(i) gives
$\sum_g\widetilde\rho(g)=\delta_{a,b}$.
Applying Proposition~\ref{prop:Etransf}(iv) gives
\[
\begin{aligned}
f_{a,b}(-1/\tau)
 &=C_{a,b}Q^{-c_{a,b}}
   \prod_{1\le g\le M/2}\prod_{m\ge1}
   \left(1-2\cos\frac{2\pi g}{M}\,Q^{m/M}
                 +Q^{2m/M}\right)^{-\widetilde\rho(g)},\\
C_{a,b}
 &=\prod_{1\le g\le M/2}
       \left(2\sin\frac{\pi g}{M}\right)^{-\widetilde\rho(g)}>0.
\end{aligned}
\]
The product on the right is $1+O(Q^{1/M})$ as
$\operatorname{Im}\tau\to\infty$. Multiplication of the finitely
many torus factors yields
\[
 \Theta_A(-1/\tau)
   =C_AQ^{-c_A}\bigl(1+O(Q^{1/L_A})\bigr),\qquad
 C_A=\prod_iC_{a_i,b_i}^{D_i}>0.
\]
Both cusps have width $12L_A$ for $\Gamma(12L_A)$. Thus their
multiplicities in the pole divisor on this modular curve are
\[
 12L_Ac_A
   =\sum_iD_i\frac{L_A}{a_i+b_i}\,\delta_{a_i,b_i}
   \in\Z_{>0}.
\]
\end{proof}

\begin{remark}\label{rem:modular}

(1) The product $P_{a,b}$ is the normalized character of the
balanced sector of $\mathcal W_a(a,a+b)$
\cite[Proposition~9.3]{HJO}. Proposition~\ref{prop:balanced-character-normalization}
proves that its ordinary character is $f_{a,b}=Q^{-c_{a,b}}P_{a,b}$
and that this sector has minimal conformal weight.
Corollary~\ref{cor:cabling-character-normalization} then identifies
$\Theta_A$ as the ordinary character of the tensor product of the
balanced sectors, with effective central charge $24c_A$.
The character conventions and weight calculation are given in
Section~\ref{supp:character}. This interpretation is not used in
the counting or modularity proofs.

(2) The proof of (c) shows that $\log\cN_A(e^{-t})\sim4\pi^2c_A/t$ as $t\to0^+$.

(3) Theorem~\ref{thm:modular} concerns the limit only. At finite
rank, Theorem~\ref{thm:formal} says that the polynomials
$p_{\underline\mu,n}$ approximate the normalized product
$\cN_A(Q)=Q^{c_A}\Theta_A(\tau)$ modulo $Q^{n+1}$.
The modular unit is $\Theta_A$; multiplication by $Q^{c_A}$ gives
the power series with constant term~$1$.

(4) The proof of Theorem~\ref{thm:modular}(c) gives the
positive constant $C_A$ explicitly. For the cusp singularity,
$C_{2,3}=(2\sin(\pi/5))^{-1}=0.850651\ldots$.

(5) For the torus branches, see also \cite{LO} on the modularity of
point counts for the curves $X^a=Y^b$.
\end{remark}

\section{Examples}\label{sec:examples}

We begin with rank-one calculations combining Hall coefficients
with the residual counts of the transformed branch. The examples
in wild characteristic show how the multiplicity sequence can change
with the characteristic. A Hecke-algebra calculation illustrates
cancellation among signed transfer coefficients.

\subsection{Branches with $D=2$ at rank one}\label{sub:rankone}
Let $A$ be a plane branch with first cabling triple $(a,b,2)$ and pattern $(T;\varpi,\omega)$, where $\val\varpi=2$, and let us write $E_j=E_{a,b}(j)$. Since $u_{(1)}^2=u_{(2)}+(q+1)u_{(1,1)}$ in $\sH$, Theorem~\ref{thm:hallchar} gives the following two values:
\begin{equation}\label{eq:Phi2}
\Phi^{a,b}_{q,2}((1,1))=E_2,\qquad\Phi^{a,b}_{q,2}((2))=\frac{q^{ab+1}E_1^2-E_2}{q-1}.
\end{equation}
If $T=O$, equivalently $\val\omega=1$, then the only lattice is $O$, and the element $\omega$ has type $(2)$ on $O/\varpi O$, thus $\#E_A(O)=\Phi^{a,b}_{q,2}((2))$.

If $T$ is a cusp, that is, $\ms(T)=(2)$, then $T=F+t^2O$ and $\val\omega\in\{2,3\}$; both values occur, the value $2$ in the example of Section~\ref{sub:wild}. Recall that the refined counts depend only on $T$ and $\varpi$ (see the note after Definition~\ref{def:chainmass}). In this case $E_T(O)$ consists of $O$, on which $\omega$ acts by zero modulo $\varpi$ because $\omega\in t^2O=\varpi O$, and of the $q$ free $T$-modules $(1+\xi t)T$, $\xi\in F$, of type $(2)$, hence $\#E_A(O)=E_2+q\,\Phi^{a,b}_{q,2}((2))$.

For $(a,b)=(2,3)$ Theorem~\ref{thm:torus}(i) gives $E_1=q+1$ and $E_2=q^4+q^3+q^2+1$, hence we have:
\begin{align*}
\ms(A)=(4,2,2):&\quad\#E_A(O)=q^8+3q^7+4q^6+4q^5+4q^4+3q^3+2q^2+q+1,\\
\ms(A)=(4,2,2,2):&\quad\#E_A(O)=q^9+3q^8+4q^7+4q^6+4q^5+4q^4+3q^3+2q^2+q+1.
\end{align*}
At $q=2$ these polynomials take the values $1123$ and $2275$, the point counts of the compactified Jacobians of the singularities with semigroups $\langle4,6,13\rangle$ and $\langle4,6,15\rangle$. In these two cases Corollary~\ref{cor:finalproduct} gives $\cN_A=\cN_{2,3}^{\,2}\cN_{2,9}$ and $\cN_A=\cN_{2,3}^{\,2}\cN_{2,11}$, respectively.

\subsection{Wild characteristic}\label{sub:wild}
The same written parametrization can have different multiplicity
sequences in different characteristics. In particular, one cannot
read its characteristic-zero Puiseux data unchanged when the
characteristic divides the multiplicity. Consider the branch $R_B=F[[t^4,t^{2B}(1+t)]]$, where $B\ge3$ is odd.

If $p\ne2$, the Puiseux characteristic of $R_B$ is $(4;2B,2B+1)$, thus $\ms(R_B)=\operatorname{Enr}(4;2B,2B+1)=(4^{(B-1)/2},2,2)$ and $\cN_{R_B}=\cN_{2,B}^{\,2}\,\cN_{2,2B+3}$ by Corollary~\ref{cor:puiseux}.

If $p=2$, we apply Corollary~\ref{cor:cabling-data} to the given generators, with $(a,b,D)=(2,B,2)$ and $(A_0,B_0)=((1-B)/2,1)$. Then $\varpi=x^{(1-B)/2}y=t^2(1+t)$ and $y^2/x^B=(1+t)^2=1+t^2$, so $\omega=t^2$ and $T=F[[t^2+t^3,t^2]]=F[[t^2,t^3]]$ is a cusp. The ring $S=F[[\varpi,\varpi^{B+1}t^2]]=F[[\varpi,\varpi^{B+1}(t^2+\varpi)]]$ has generators of valuations $2$ and $2B+5$, hence $\ms(R_B)=(4^{(B-1)/2},2,2,2)$, the cabling list is $((2,B,2),(2,2B+5,1))$, and we obtain:
\[
\cN_{R_B}=\cN_{2,B}^{\,2}\ \cN_{2,2B+5}\qquad(p=2).
\]


\subsection{A reduction in the Hecke algebra}\label{sub:hecke46}
The symbolic calculation in
Section~\ref{subsec:alg-symbolic} gives, at $q=2$,
\[
 \bigl(c_{(2)}(F_{(4),6}),c_{(1,1)}(F_{(4),6})\bigr)
       =(1670,-547).
\]
The independent affine-Hecke computation gives the same row.
By Proposition~\ref{prop:composition} it is the row of index $(4)$ of
$\mathsf C_{4,4}\mathsf C_{4,2}$.
For the cusp in rank two, the terminal flag vector is $(1,3)^{\mathsf T}$,
so the contraction gives $29$. For multiplicity sequence $(4,2,2)$
in rank one, the terminal frame has valuation two and the flag
vector is $(1,1)^{\mathsf T}$, giving $1123$.
The two contractions depend on the terminal framing and involve
cancellation between signed transfer coefficients.


\clearpage
\appendix
\section{Bounded models for lattice and flag counts}\label{app:integrality}

A conductor bound places normalized lattices in a finite
Grassmannian. On this locus we construct the quotient bundle
carrying the residual endomorphism and its stable flags, and
check compatibility with arbitrary base change. These are the
finite-type schemes used in Corollary~\ref{cor:poly}. Their
construction does not use polynomiality; the conclusion about
integral coefficients follows from Lemma~\ref{lem:completion-integral}.

\begin{proposition}[Bounded counting models]\label{cor:completion-integral-main}
Let $A=k[[u,w]]\subset O=k[[t]]$ be a plane branch over a finite
field $k$, with $u,w\in tO$ and $u\ne0$, and let $n\ge1$.
Each of the following
counts is the number of points of a fixed separated finite-type
$k$-scheme, uniformly over every finite extension $K/k$:
the total count $\#E_{A_K}(O_K^n)$, each residual-type count
$B_{(u,w)}(\mu;n)$, and each chain count
$\Ch_\alpha(A_K,u;n)$, with the generators and types extended
from $k$ to $K$. The constructions of these schemes commute
with extension of the base field.
\end{proposition}
\begin{proof}
Choose $c\ge1$ with $t^cO\subset A$ and put
$V_c=(O/t^cO)^n$. For $0\le r\le cn$, let $\cE_{A,n,r}^{(c)}$
be the locally closed subscheme of $\Gr(r,V_c)$ on which the
universal subbundle $W$ is $A$-stable and
$W\to V_c/tV_c$ is surjective. Stability may be imposed using
a finite $k$-basis of $A/t^cO$, so these are finite-type schemes.
Their finite disjoint union is separated.

For every finite extension $K/k$, its $K$-points are precisely
the normalized $A_K$-lattices $L\subset O_K^n$, where
$O_K=K[[t]]$, by $L\mapsto L/t^cO_K^n$. Indeed, normalization
implies $t^cO_K^n\subset L$: for $v\in O_KL$, multiplying an
expression $v=\sum a_i l_i$ by $t^c$ gives coefficients
$t^ca_i\in A_K$. Conversely, a stable subspace has a stable
preimage, and surjectivity modulo $t$ is equivalent to
$O_KL=O_K^n$ by Nakayama. The conductor bound and the defining
equations extend by scalars.

For the refined conditions, put $D=\val_tu$ and $H=c+D$.
On each of the preceding Grassmannian loci, regard the fixed
spaces $V_c$ and $V_H$ as trivial vector bundles. Lift $W$ to
its inverse image $\widetilde W\subset V_H$. It is a subbundle
containing the fixed subbundle
\[
 \mathcal K_u=\ker(u:V_H\to V_H)=t^cO^n/t^HO^n,
 \qquad \operatorname{rank}\mathcal K_u=Dn.
\]
Multiplication by $u$ induces an isomorphism of fixed vector bundles
\[
 V_H/\mathcal K_u\xrightarrow{\ \sim\ }\im(u)
       =t^DO^n/t^HO^n.
\]
Since $V_H/\widetilde W$ is locally free,
$\widetilde W/\mathcal K_u$ is a subbundle of $V_H/\mathcal K_u$.
Its image $u\widetilde W$ is therefore a subbundle of the fixed
subbundle $\im(u)$, and hence of $V_H$. The $A$-stability
condition gives $u\widetilde W\subset\widetilde W$, so there
is an exact sequence
\[
 0\longrightarrow\widetilde W/u\widetilde W
 \longrightarrow V_H/u\widetilde W
 \longrightarrow V_H/\widetilde W\longrightarrow0.
\]
The two right-hand terms are vector bundles, so the surjection
splits locally. Its kernel
$\mathcal Q=\widetilde W/u\widetilde W$ is consequently locally
free of rank $Dn$, the rank of $\mathcal K_u$. These locally split
sequences also show that the construction commutes with arbitrary
base change. On a field-valued point it gives $L/uL$:
$t^HO^n=u(t^cO^n)\subseteq uL$.

The other generator $w$ acts on $\mathcal Q$. If $w=0$, its
action is already zero. Otherwise, choose $M$ with
$M\val_t(w)-D\ge c$. Then $w^M/u\in t^cO\subset A$, so
$w^ML\subset uL$ in every family. Thus the endomorphism is
uniformly nilpotent. A specified Jordan type is imposed by
exact-rank determinantal conditions on its powers, giving a
locally closed subscheme. Stable flags of a fixed composition
are closed incidence conditions in the relative flag scheme
of $\mathcal Q$. These schemes are separated and of finite
type. Their field-valued points give the residual-type counts
and, by Lemma~\ref{lem:chainflag}, the chain counts. This proves
the proposition.
\end{proof}

\section{Finite subroutines for the counting algorithm}\label{sec:algorithm-implementation}

This appendix gives the finite subroutines for the framed recursion
of Section~\ref{subsec:explicit-counting-algorithm}, with matrix convention
$\mathsf C\mathsf M=\mathsf A$. We describe the support enumeration,
the finite search for affine cocenter reduction, the commutator
identities for finite blocks, and the Hall calculation of the
terminal flag vector. Torus counts give an alternative evaluation,
illustrated by a symbolic example. All calculations take place
in $\Q(q)$, and the bounds below prove termination. We retain
the notation and partition orders of
Section~\ref{subsec:explicit-counting-algorithm}.

\subsection{Finite supports and normalized rows}\label{subsec:alg-support}
Fix $m,k\ge1$, put $N=\gcd(m,k)$ and $e=m/N$, and use the
partition orders of Section~\ref{subsec:explicit-counting-algorithm}.
Both test functions in \eqref{eq:testfunctions} have the following
finite index sets.
For a composition $\alpha$ put
\[
 v_\alpha(q)=\prod_i[\alpha_i]_q!,\qquad
 [r]_q!=\prod_{j=1}^r(1+q+\cdots+q^{j-1}).
\]
Let
$\mathbf a(w)=(a_{w,\rho}(q))_{\rho\vdash N}$ be a row obtained
by the reduction of Section~\ref{sec:hecke}, so that
\[
 D(T_w)=\sum_{\rho\vdash N}a_{w,\rho}(q)D(T_{w_\rho})
 \quad\text{for every }D\in J_k.
\]
This is an identity on $J_k$, not necessarily an identity in the
entire cocenter after terms have been discarded.

Represent an affine permutation by its window
\[
 w(i)=\sigma(i)+mz_i\quad(1\leq i\leq m),\qquad
 \sigma\in S_m,\quad z_i\in\Z,
\]
and extend it by $w(i+m)=w(i)+m$. Let $b_\alpha(i)$ be the block index of $i$. Define
\begin{equation}
 \mathcal S_{\alpha,k}=
 \left\{w:\ \sum_i z_i=k,\quad
 z_i\geq\mathbf1_{b_\alpha(\sigma(i))<b_\alpha(i)}
 \text{ for every }i\right\}.
 \label{eq:alg-support}
\end{equation}
All lower bounds in this formula are nonnegative. These are precisely
the windows satisfying $w\Lambda_j^\alpha\subseteq\Lambda_j^\alpha$
for every $j$. Indeed, at a separator $s$, the lattice $\cL_s$ has
basis $\pi^{\mathbf1_{i\leq s}}b_i$, and preservation is equivalent to
\[
 z_i+\mathbf1_{i\leq s}\geq\mathbf1_{\sigma(i)\leq s}.
\]
Taking all separators gives \eqref{eq:alg-support}.
To enumerate this support, take each $\sigma\in S_m$, put
$\epsilon_i=\mathbf1_{b_\alpha(\sigma(i))<b_\alpha(i)}$,
and enumerate the weak compositions $(t_1,\dots,t_m)$ of
$k-\sum_i\epsilon_i$, setting $z_i=\epsilon_i+t_i$.
If $k-\sum_i\epsilon_i<0$, this permutation contributes no window.

The finite sums in \eqref{eq:alg-A} and \eqref{eq:alg-M} give
the normalized rows $\mathsf A$ and $\mathsf M$.
The rigid support is the finite group $W_{\beta^{[e]}}$, translated
by $\tau^k$. Both normalizing factors enter the equation
$\mathsf C\mathsf M=\mathsf A$.

\subsection{Finite search and block refinements}\label{subsec:alg-discrete}
For the reduction of $\mathbf a(w)$ in
Section~\ref{subsec:explicit-counting-algorithm}, we prove that the
affine search is finite and show how to refine the blocks. We retain
the notation
$N=\gcd(m,k)$, $e=m/N$ and $\delta=\operatorname{Ad}(\tau^k)$.

\emph{Finiteness of the search.}
Write a window of valuation $k$ as $w(i)=\sigma(i)+mz_i$, with
$\sigma\in S_m$ and $\sum_i z_i=k$. Then
\[
 \left\lfloor\frac{w(j)-w(i)}m\right\rfloor
 =z_j-z_i-\mathbf1_{\sigma(j)<\sigma(i)}\qquad(i<j).
\]
If $\ell(w)\le L$, the length formula in
Section~\ref{subsec:explicit-counting-algorithm}
therefore gives $|z_i-z_j|\le L+1$ for all $i,j$. Since
$mz_i-k=\sum_j(z_i-z_j)$, every such window lies in the finite box
\[
 \left\lceil\frac{k-(m-1)(L+1)}m\right\rceil
 \le z_i\le
 \left\lfloor\frac{k+(m-1)(L+1)}m\right\rfloor.
\]
There are only finitely many choices of $\sigma$ and of the integers
$z_i$ in this box. Thus the search through length-preserving
conjugations terminates if previously visited windows are not
revisited. Conjugation by $\tau$ also has period dividing $m$,
since $\tau^m$ is central.
At an input of length $L$, the recursive calls in \eqref{eq:GP}
have lengths $L-1$ and $L-2$.
For a central terminal element, Proposition~\ref{prop:HN3} supplies
a compatible block decomposition in the same finite search graph.
Together with the decrease in length, this proves termination of
the affine reduction. Proposition~\ref{prop:effective-standard-triple} gives a
second finite construction that applies to every affine permutation
of rank $m\ge2$:
it constructs the full standard triple $(K,x,u)$ using the rational
regularity test of Proposition~\ref{prop:regularity-test}.

\emph{The finest compatible separators.}
For a candidate $w'=u\tau^k$, let $S(u)$ be the set of
$p\in\{0,\ldots,N-1\}$ for which $u$ preserves every interval with
endpoints in $p+N\Z$. The singleton test in
Section~\ref{subsec:explicit-counting-algorithm} determines $S(u)$
using the lifted integer window. A nonempty set of separators
$S\subseteq\{0,\ldots,N-1\}$ is compatible with $u$ if and only if
$S\subseteq S(u)$.

Indeed, the decomposition with endpoints in $p+N\Z$ is a coarsening
of the decomposition with endpoints in $S+N\Z$ whenever $p\in S$.
Preservation of the latter therefore implies preservation of every
singleton decomposition. Conversely, if $u$ preserves the singleton
decompositions for all $p\in S$, it preserves their common refinement:
each refined block is an intersection of blocks that $u$ preserves.
That common refinement has endpoints in $S+N\Z$. Consequently a
compatible proper $\delta$-stable parabolic exists exactly when
$S(u)\ne\varnothing$, and $S(u)$ gives the finest compatible block
decomposition. Only $N$ singleton tests are needed, instead of a
search through $2^N-1$ nonempty separator sets. Testing residues
modulo $m$ alone would lose the interval-preservation condition.

\emph{Reduction with smaller blocks.}
Suppose that $w$ is central terminal and $w'\approx w$ has the
compatible form $w'=u\tau^k$. Choose the finest decomposition, or
any compatible coarsening, and let $R$ contain one block from each
$\delta$-orbit. Write a representative as
$B=\{p+1,\ldots,p+c_B\}$. Its orbit consists of the $e$ distinct
blocks $B+jk$ modulo $m$, and
$\sum_{B\in R}c_B=N$. For $0\le j<e$, identify the factor of $u$
on $B+jk$ with the permutation
\[
 \pi_{B,j}(r)=u(p+jk+r)-(p+jk),\qquad 1\le r\le c_B.
\]
Thus the full shift $jk$, including any multiple of $m$, is retained.
The product to reduce in $H_{c_B}$ is
\[
 F_B=T_{\pi_{B,0}}T_{\pi_{B,e-1}}\cdots T_{\pi_{B,1}},
\]
with $F_B=T_{\pi_{B,0}}$ if $e=1$. To see the order, write
$\phi(h)=[hT_{\tau^k}]$ and use
$\phi(hh')=\phi(h'\delta(h))$. First move the factor at $B+(e-1)k$
into $B$, where it is multiplied on the right; next move the factor
at $B+(e-2)k$ through the now empty blocks into $B$, and continue.
Since $\delta^e=1$, this gives the displayed product, without an
additional power of $q$.

Apply the finite routine to each $F_B$, initially keeping composition
indices:
\[
 F_B\equiv\sum_\eta c_{B,\eta}(q)T_\eta
 \pmod{[H_{c_B},H_{c_B}]},\qquad c_{B,\eta}(q)\in\Z[q].
\]
Combining the block reductions gives
\[
 a_{w,\rho}(q)=
 \sum_{\substack{(\eta_B)_{B\in R}\\
       \operatorname{sort}(\bigsqcup_{B\in R}\eta_B)=\rho}}
       \prod_{B\in R}c_{B,\eta_B}(q).
\]
Here $\bigsqcup$ means concatenation of the parts, followed by sorting
into a partition. Each orbit contributes once, so these parts sum
to $N$; they are not repeated $e$ times. Factors in distinct blocks
commute, and Lemma~\ref{lem:twist}(iii) allows the finite commutator
relations in one representative block while the other representative
factors are held fixed. Lemma~\ref{lem:swap} then identifies the
resulting refined block structure with $T_{w_\rho}$. This proves the
formula; using smaller blocks changes only the size of the finite
calculations. Sorting affine cycle types in place of
these Hecke calculations is not justified; see
Remark~\ref{rem:terminal}(1).

\subsection{Finite-block identities and uniqueness}\label{subsec:alg-finite-block}
The finite routine of Section~\ref{subsec:explicit-counting-algorithm}
computes the reduction of Lemma~\ref{lem:finiteA}. We give the
commutator identities proving this reduction and keep track of
the coefficient rings needed for reduction and sorting.

Let $H_c$ be the type-$A$ Hecke algebra over $\Q(q)$, with
$T_s^2=(q-1)T_s+q$, and use the descending block Coxeter words
$T_\eta$ of Lemma~\ref{lem:finiteA}. For a standard basis element,
the reduction can be verified by an identity in the full Hecke algebra:
\[
 T_v-\sum_{\rho\vdash c}c_{v,\rho}(q)T_\rho
 =\sum_t d_t[A_t,B_t].
\]
The normal-form coefficients $c_{v,\rho}$ lie in $\Z[q]$.
The right side can be constructed over $\Z[q]$ until composition
indices are sorted, and over $\Z[q,q^{-1}]$ after sorting. The
identity extends linearly to every finite Hecke element; polynomial
input coefficients still give polynomial normal-form coefficients.

\emph{Identities for the relative recursion.}
We use the elements $E_j$, the relative rank $r=c-j+1$ and, for
$i=v^{-1}(r)<r$, the elements $d$, $a$, $b$ of the finite routine in
Section~\ref{subsec:explicit-counting-algorithm}.
The recursive state is $E_jT_v$ with $v\in S_r$, and $E_j$ commutes
with $H_{r-1}$. If $v(r)=r$, a commutator $[A,B]$ from the
smaller algebra $H_{r-1}$ lifts by the identity
\[
 E_j[A,B]=[E_jA,B].
\]
The accompanying normal-form identity is
$E_jT_\eta=T_{(\eta,j)}$.

If $i<r$, the factorization $v=as_{r-1}b$ is
length-additive: the inversions involving the position $i$ of $r$
contribute $r-i$, and the remaining inversions are exactly those
of $a$. Since $E_j$ commutes with $T_a$ and
$E_{j+1}=E_jT_{s_{r-1}}$, we obtain the identity
\[
 E_jT_v-E_{j+1}T_bT_a=[T_a,E_{j+1}T_b]
\]
used there.
After expanding $T_bT_a$ in the standard basis, multiply the
identities for its recursive terms by the corresponding
coefficients and add this commutator. This proves the identity
for $E_jT_v$ and accounts for the order $T_bT_a$ in the algorithm.
Each step reduces the relative rank $r$, or calls the
ordinary routine at smaller rank, and all multiplication
coefficients lie in $\Z[q]$.

\emph{Identities for sorting.}
Suppose that $\eta'$ is obtained from $\eta$ by interchanging
adjacent blocks. Let $X$ move the first block past the second,
preserving the order within both and fixing all other positions.
For the corresponding block Coxeter elements,
\[
 Xw_\eta=w_{\eta'}X,\qquad
 \ell(Xw_\eta)=\ell(X)+\ell(w_\eta),\qquad
 \ell(w_{\eta'}X)=\ell(w_{\eta'})+\ell(X).
\]
The length equalities hold because $X$ is increasing on each old
block and $X^{-1}$ is increasing on each new block. Thus
$T_XT_\eta=T_{\eta'}T_X$, and the sorting identity is
\[
 T_\eta-T_{\eta'}=[T_\eta T_X^{-1},T_X].
\]
The inverse is computed using
\[
 T_s^{-1}=q^{-1}T_s+(q^{-1}-1),
\]
with the factors taken in reverse reduced-word order. This is an
identity after inverting $q$; it does not assert invertibility in
the algebra over $\Z[q]$. Adjacent interchanges sort every
composition and leave its polynomial coefficient unchanged.
The resulting commutator identity can be checked by expanding
both sides in the full standard basis.

\emph{Independence of choices.}
Theorem~\ref{thm:H} and Corollary~\ref{cor:basis} show that the
evaluations $D\mapsto D(T_{w_\rho})$ form a basis of $J_k^*$ for
every prime power $q$. Two valid symbolic reductions of the same
$T_w$ therefore agree at every prime power where their rational
expressions are defined. There are infinitely many such values,
so the resulting rows agree in $\Q(q)$. This applies to choices
of search paths, compatible separators, block representatives,
and finite reductions.

Moreover, \eqref{eq:GP}, finite Hecke multiplication, and the
relative recursion all have coefficients in $\Z[q]$, and sorting
only collects coefficients. Hence the resulting row satisfies
$\mathbf a(w)\in\Z[q]^{p(N)}$. The row indexed by partitions is
unique; unsorted composition coordinates and commutator
expressions need not be unique. The integrality assertion does
not extend to the normalized matrices $\mathsf A$, $\mathsf M$
and $\mathsf C$. Factorial support and growth of intermediate
polynomials limit practical running time; the argument gives no
polynomial-time complexity bound.

\subsection{An elementary Hall routine for terminal flags}\label{subsec:alg-terminal-hall}
The terminal flag counts can be computed from semisimple Hall
multiplication. Write
$\mathsf e_r=u_{(1^r)}$ and
$\mathsf v_j=\sum_{\lambda\vdash j}u_\lambda$, with
$\mathsf e_0=\mathsf v_0=u_\varnothing=1$. For partitions $\lambda,\nu$ with
$|\lambda|-|\nu|=r$, put
\[
 d_i=\lambda'_i-\nu'_i,\qquad
 m_i=\lambda'_i-\lambda'_{i+1}.
\]
The coefficient of $u_\lambda$ in $u_\nu\mathsf e_r$ is
\begin{equation}
 g^{\lambda}_{\nu,(1^r)}(q)=
 q^{\sum_{i<j}d_i(m_j-d_j)}
       \prod_i\qbinom{m_i}{d_i}_q,
 \label{eq:alg-hall-pieri}
\end{equation}
provided $0\leq d_i\leq m_i$ for every $i$, and is zero otherwise.
Here missing conjugate parts are zero. To see the formula, a
semisimple submodule $U$ lies in the socle of $M_\lambda$.
Filter this socle by
$F_i=\ker z\cap z^{i-1}M_\lambda$. The quotient $M_\lambda/U$
has type $\nu$ precisely when
\[
 \dim(U\cap F_i)-\dim(U\cap F_{i+1})=d_i.
\]
The Gaussian factors choose the subspaces in the graded pieces,
of dimensions $m_i$. The graph maps between an earlier chosen
piece and later unchosen pieces contribute
$\sum_{i<j}d_i(m_j-d_j)$ free parameters. This proves
\eqref{eq:alg-hall-pieri}.

The identity in the proof of Lemma~\ref{lem:hallbasics}(d) gives
\begin{equation}
 \mathsf v_j=\sum_{r=1}^j(-1)^{r+1}q^{\binom r2}
                     \mathsf e_r\mathsf v_{j-r}.
 \label{eq:alg-V}
\end{equation}
Expand each $\mathsf v_j$ as a polynomial in the commuting
$\mathsf e_r$, and
evaluate these monomials by repeated use of
\eqref{eq:alg-hall-pieri}. Then
\[
 H_\alpha(\lambda)=[u_\lambda]\prod_i\mathsf v_{\alpha_i},
 \qquad
 F_\eta(\lambda)=[u_\lambda]\prod_i\mathsf e_{\eta_i}.
\]
Gaussian binomials can themselves be computed without division, using
$ \qbinom ab_q=\qbinom{a-1}{b-1}_q+q^b\qbinom{a-1}b_q$,
with $\qbinom a0_q=\qbinom aa_q=1$ for $a\ge0$ and value zero
when $b<0$ or $b>a$. For example, with rows and columns in the order
$(2),(1,1)$,
\begin{equation}
 \sM_2=\begin{pmatrix}1&1\\1&q+1\end{pmatrix},\qquad
 \mathsf F_2=\begin{pmatrix}0&1\\1&q+1\end{pmatrix}.
 \label{eq:alg-small-hall}
\end{equation}
The two matrices are distinct: $H_\alpha$ counts stable flags,
whereas $F_\eta$ counts flags lowered by the nilpotent operator.

\subsection{Finite torus inputs for the Hall-character shortcut}\label{subsec:alg-torus}
The torus inputs $E_j$ in \eqref{eq:alg-character-system} can be
computed by the finite formula of Theorem~\ref{thm:torus}(i). If $a=1$ or $b=1$,
set $E_j=1$ directly. Otherwise interchange $a,b$ if necessary so
that $2\le a<b$, and form
the gaps $G=\N\setminus\langle a,b\rangle$, let
$f=ab-a-b$, $\delta=|G|$, and enumerate the vectors
$(r_g)_{g\in G}$ with $0\leq r_g\leq j$ which are increasing
for the gap order $g\preceq h$ if $h-g\in\langle a,b\rangle$.
For $h<0$ set $r_h=0$ and put
\begin{align*}
 \mathcal G_{j,\mathbf r}(Q)
  &=\frac{(Q;Q)_j}{(Q;Q)_{j-r_f}}
    \prod_{g\in G}
    \frac{(Q;Q)_{r_g-r_{g-a-b}}}
    {(Q;Q)_{r_g-r_{g-a}}(Q;Q)_{r_g-r_{g-b}}},\\
 \dinv(\mathbf r)
  &=\sum_{g,h\in G}
    \bigl(\mathbf1_{h-g\geq0}-\mathbf1_{h-g\geq a}
          -\mathbf1_{h-g\geq b}+\mathbf1_{h-g\geq a+b}\bigr)r_gr_h.
\end{align*}
Whenever a shifted index in the product is nonnegative it is again
a gap. Then
\[
 E_j=q^{j^2\delta}
      \sum_{\mathbf r}\mathcal G_{j,\mathbf r}(q^{-1})
                          q^{-\dinv(\mathbf r)}.
\]
This torus calculation uses the formula quoted in
Theorem~\ref{thm:torus}; the general algorithm
\eqref{eq:alg-recursion}--\eqref{eq:alg-matrix} is independent
of that formula.

\subsection{A symbolic two-stage example}\label{subsec:alg-symbolic}
For the cusp, the conductor is $t^2O$. An element of
$E_{\F_q[[t^2,t^3]]}(O^j)$ corresponds modulo
$t^2$ to a subspace surjecting onto the constant coefficients.
Choose its intersection with $t\F_q^j$ to have dimension
$r$, and then choose a graph over $\F_q^j$. This gives
\[
 E_{2,3}(j)=\sum_{r=0}^j\qbinom jr_q q^{j(j-r)},\qquad
 E_1=q+1,\quad E_2=q^4+q^3+q^2+1.
\]
With these values, \eqref{eq:Phi2} gives $\Phi^{2,3}_{q,2}((1,1))=E_2$ and
$\Phi^{2,3}_{q,2}((2))=(q^7E_1^2-E_2)/(q-1)$. The second one is the
polynomial $P_{(4,2,2),1}$ displayed in Section~\ref{sub:rankone}: for
a branch with multiplicity sequence $(4,2,2)$, choose its
first frame with valuations $(4,6)$; its transform is smooth with
framing valuation $2$, so its only residual type in rank one is
$(2)$.
One can also recover the symbolic row in Section~\ref{sub:hecke46}. Write
$c_2=c_{(2)}(F_{(4),6})$ and
$c_{11}=c_{(1,1)}(F_{(4),6})$. Since
$\Phi(\lambda)=c_2H_{(2)}(\lambda)+c_{11}H_{(1,1)}(\lambda)$,
\eqref{eq:alg-small-hall} gives
\[
 c_{11}=\frac{E_2-\Phi((2))}{q},\qquad c_2=\Phi((2))-c_{11}.
\]
Thus
\begin{align*}
 c_{11}&=-q^7-3q^6-4q^5-4q^4-3q^3-2q^2-q-1,\\
 c_2&=q^8+4q^7+7q^6+8q^5+8q^4+6q^3+4q^2+2q+2.
\end{align*}
Their values at $q=2$ and the contractions with the two terminal
vectors are discussed in Section~\ref{sub:hecke46}.

\subsection{Cost and exact checks}\label{subsec:alg-checks}
A direct enumeration of \eqref{eq:alg-support} has at most
$m!\binom{k+m-1}{m-1}$ terms; the rigid support has
$\prod_i(\beta_i!)^e$ terms. Thus the universal algorithm is
finite but is not a polynomial-time algorithm. In the canonical
recursion $k\leq m$, while a useful alternative frame may have
$k>m$; the two computations agree by
Proposition~\ref{prop:composition}, which expresses every
$\mathsf C_{m,k}$ with $k>m$ through matrices with $k\le m$.
Reusing reductions of affine windows, finite-block tables, Hall
monomials, and matrices for each pair $(m,k)$ avoids repeated
calculations. To compute a single output, one can propagate
one row and calculate only the source rows it requires.

Corollary~\ref{cor:poly} guarantees that the final rational functions for
actual chain counts, refined counts, and $P_{\underline\mu,n}$
cancel to polynomials in $\Z[q]$. Intermediate transfer
coefficients may be negative, as the preceding example shows,
and their exact values must be retained throughout the calculation.
Polynomiality and integrality concern the final counting expressions.
The matrix identity $\mathsf C\mathsf M=\mathsf A$, the terminal
Hall values, and the two contractions of Section~\ref{sub:hecke46}
give exact checks. Theorem~\ref{thm:formal} gives
\[
 \deg P_{\underline\mu,n}=n^2\delta(\underline\mu),\qquad
 [q^{n^2\delta}]P_{\underline\mu,n}=1,\qquad
 \delta(\underline\mu)=\sum_i\binom{\mu_i}2
\]
as further checks. These degree statements are not needed for the rational-function
algorithm or its proof of correctness.

\section{Affine-Weyl, Hecke, and modular arguments}\label{app:conventions}

We prove the affine-Weyl and Hecke results used in
Sections~\ref{sec:hecke}--\ref{sec:chains}, followed by the modular
transformation formulas of Section~\ref{sec:modular}. The final
subsection computes the character normalization of the cabling product.

\subsection{Reduction to standard triples}\label{app:HN3}

The reduction for $\widetilde W=\mathbb Z^m\rtimes S_m$ follows
the geometric argument of He and Nie
\cite[Propositions~2.5 and~2.8]{HN}. The key step is to find an
alcove whose closure contains a regular point of the affine space
on which $w$ acts by translation. Its stabilizer in the affine Weyl
group $W_a$ gives the finite parabolic factor, and the remaining
factor is straight.
We conclude with a finite procedure for finding these factors.

We work modulo the diagonal line, keeping the full Newton vector
and the valuation to recover the central direction.

\subsubsection{Setting}
Fix $m\ge2$. An affine permutation is a bijection $w:\Z\to\Z$
with $w(i+m)=w(i)+m$. We identify $t^\lambda\sigma\in
\tW=\Z^m\rtimes S_m$ with the affine permutation
$w(i)=\sigma(i)+m\lambda_{\sigma(i)}$ for $1\le i\le m$.
Put $\tau(i)=i+1$ and
$\val(w)=m^{-1}\sum_{i=1}^m(w(i)-i)=\sum_i\lambda_i$.
For $0\le i<m$, the simple reflection $s_i$ exchanges
$i+mj$ and $i+1+mj$ for all $j\in\mathbb Z$.
If $N$ is the order of $\sigma$, then
\[
 w^N=t^{\lambda+\sigma\lambda+\cdots+\sigma^{N-1}\lambda}.
\]
Define the \emph{full Newton vector}
\[
 \nu_w^{\mathrm{full}}
 =\frac1N\sum_{j=0}^{N-1}\sigma^j\lambda\in\mathbb Q^m.
\]
This definition is independent of the power used: replacing $N$ by
any multiple multiplies both the translation vector and the denominator
by the same integer, and any two admissible powers have a common
multiple. If $c$ is a cycle of $\sigma$, then
\begin{equation}\label{eq:full-newton-cycles}
 (\nu_w^{\mathrm{full}})_i
 =\frac{\sum_{j\in c}\lambda_j}{|c|}\quad(i\in c),
 \qquad
 \sum_{i=1}^m(\nu_w^{\mathrm{full}})_i=\val(w).
\end{equation}
The Newton point is \emph{central} when all these coordinates are equal.
For $a=t^\mu\pi\in\tW$, conjugating a translation power gives
\begin{equation}\label{eq:newton-conjugation}
 \nu_{a^{-1}wa}^{\mathrm{full}}=\pi^{-1}\nu_w^{\mathrm{full}},
 \qquad \val(a^{-1}wa)=\val(w).
\end{equation}

The group $\tW$ acts on $\R^m$ by $t^\lambda\sigma\cdot v=\sigma v+\lambda$, where $(\sigma v)_{\sigma(i)}=v_i$, and hence on $\mathbb V=\R^m/\R(1,\dots,1)$, which we identify with the Euclidean space $\{v:\sum_iv_i=0\}$.
When describing points or transformations of $\mathbb V$, coordinate
tuples in $\mathbb R^m$ are understood modulo $\mathbb R(1,\ldots,1)$;
their sum-zero representatives are obtained by subtracting the mean.
The full Newton vectors remain in $\mathbb Q^m$; the solution space
in Proposition~\ref{prop:regularity-test} is taken in $\mathbb R^m$.

The action on $\mathbb V$ is by affine isometries, and its kernel is generated by the central element $t^{(1,\dots,1)}=\tau^m$. Let us denote by $\fH$ the set of hyperplanes $\{v_i-v_j=c\}$ of $\mathbb V$, where $i<j$ and $c\in\Z$. This set is locally finite and is permuted by $\tW$, and the connected components of $\mathbb V\smallsetminus\bigcup\fH$ are called \emph{alcoves}. For $H\in\fH$ we denote by $s_H$ the orthogonal reflection in $H$, which for $H=\{v_i-v_j=c\}$ is the element $t^{c(e_i-e_j)}(i\,j)$ of $\tW$, of valuation $0$. For alcoves $\mathsf A,\mathsf A'$ we denote by $\fH(\mathsf A,\mathsf A')$ the set of hyperplanes separating them. An alcove is the intersection of the open half-spaces which are bounded by hyperplanes of $\fH$ and contain it, because this intersection is convex and meets no hyperplane, and consequently we have:
\begin{equation}\label{eq:alc-a}
\begin{aligned}
\fH(\mathsf A,\mathsf A')=\emptyset\ &\Longrightarrow\ \mathsf A=\mathsf A',\\
v\in\overline{\mathsf A}\cap\overline{\mathsf A'}\ \text{ and }\ H\in\fH(\mathsf A,\mathsf A')\ &\Longrightarrow\ v\in H .
\end{aligned}
\end{equation}
For the second statement, note that if $v\notin H$, then a neighbourhood of $v$ lies on one side of $H$ and meets both alcoves.

The fundamental alcove and its walls are
\begin{gather*}
\Delta=\{v\in\mathbb V:\ v_1>v_2>\dots>v_m>v_1-1\},\\
H_i=\{v_i=v_{i+1}\}\ \ (1\le i<m),\qquad H_0=\{v_1-v_m=1\}.
\end{gather*}
On $\Delta$ we have $0<v_i-v_j<1$ for all $i<j$, so $\Delta$ is convex and meets no hyperplane. It is also the intersection of the half-spaces bounded by $H_0,\dots,H_{m-1}$, and thus an alcove. The reflection in $H_i$ is the simple reflection $s_i$, and we put $S=\{s_0,\dots,s_{m-1}\}$. The closure $\overline\Delta$ is the simplex with vertices $\omega_k=e_1+\dots+e_k$, $0\le k<m$, and the hyperplane $H_i$ contains all of them except $\omega_i$. The barycenter of the vertices other than $\omega_i$ lies on $H_i$ and on no other hyperplane of $\fH$. Thus a small ball around it meets $\Delta$ and $s_i\Delta$ but no hyperplane other than $H_i$, and we obtain:
\begin{equation}\label{eq:alc-b}
\fH(\Delta,s_i\Delta)=\{H_i\}\qquad(0\le i<m).
\end{equation}
The element $\tau$ acts by $(v_1,\dots,v_m)\mapsto(v_m+1,v_1,\dots,v_{m-1})$ and maps $\Delta$ onto itself. Define $\ell(w)$ to be the number of inversions of the affine permutation $w$, that is, of the pairs $(a,a')$ with $1\le a\le m$, $a<a'$ and $w(a)>w(a')$. We prove the following standard facts from this definition of $\ell$;
see also \cite[Ch.~V, \S3; Ch.~VI, \S2]{Bou},
\cite[Ch.~1, 4, 5]{Hum}, and \cite[\S1.5]{HN}.

\begin{lemma}\label{lem:standard}
\begin{enumerate}[label=\textup{(S\arabic*)},leftmargin=3em]
\item The group $W_a=\langle S\rangle$ acts simply transitively on the alcoves, $W_a=\ker\val$, and the stabilizer of $\Delta$ in $\tW$ is $\langle\tau\rangle$. For an alcove $\mathsf A$ we write $\mathsf A=x_{\mathsf A}\Delta$ with $x_{\mathsf A}\in W_a$.
\item $\ell(w)=\#\fH(\Delta,w\Delta)$ for $w\in\tW$.
\item For $w\in\tW$ and $H\in\fH$ we have $\ell(s_Hw)<\ell(w)$ if and only if $H\in\fH(\Delta,w\Delta)$.
\item For $v\in\overline\Delta$ we put $K_v=\{s\in S:s(v)=v\}$. Then $K_v\subsetneq S$, and every proper subset of $S$ is of this form. The stabilizer of $v$ in $W_a$ is the group $W_{K_v}$ generated by $K_v$. This group is finite and contains $s_H$ for every $H\in\fH$ with $v\in H$.
\item Let $K\subsetneq S$. Every $w\in\tW$ has a unique decomposition $w=ux$ with $u\in W_K$ and $x$ of minimal length in $W_Kw$, and $\ell(yx)=\ell(y)+\ell(x)$ for all $y\in W_K$. An element $x'\in\tW$ is of minimal length in $W_Kx'$ if and only if $\ell(sx')>\ell(x')$ for all $s\in K$.
\end{enumerate}
\end{lemma}
\begin{proof}
\emph{The length formula (S2).} Let $b\in\Delta$ and $w=t^\lambda\sigma$. Let us recall that $w(k)=\sigma(k)+m\lambda_{\sigma(k)}$ and $(wb)_{\sigma(k)}=b_k+\lambda_{\sigma(k)}$, and for $1\le k<l\le m$ we put $i=\sigma(k)$, $j=\sigma(l)$ and $n=\lambda_i-\lambda_j$. Since $0<b_k-b_l<1$, the number $(wb)_i-(wb)_j$ lies in the open interval $(n,n+1)$. On the other hand, $b_i-b_j$ lies in $(0,1)$ if $i<j$ and in $(-1,0)$ if $i>j$, thus the number of hyperplanes $\{v_i-v_j=c\}$, $c\in\Z$, separating $b$ from $wb$ is $|n|$ if $i<j$ and $|n+1|$ if $i>j$.

Next we count the inversions, noting first that no inversion $(a,a')$ of $w$ has $a'\equiv a$ modulo $m$. The number $x=(w(k)-w(l))/m=n+(i-j)/m$ is not an integer, and the inversions with $\{a,a'\}\equiv\{k,l\}$ are the pairs $(k,l+mc)$ with $0\le c<x$ and the pairs $(l,k+mc)$ with $1\le c<-x$, whose number is $\lceil x\rceil$ if $x>0$ and $\lfloor-x\rfloor$ if $x<0$. If $i<j$, then $n-1<x<n$ and this number is $|n|$, while if $i>j$, then $n<x<n+1$ and this number is $|n+1|$. In both cases the two counts agree. In particular,
\begin{equation}\label{eq:affine-window-length}
 \ell(t^\lambda\sigma)
 =\sum_{1\le k<l\le m}
   \left|\lambda_{\sigma(k)}-\lambda_{\sigma(l)}
           +\mathbf1_{\sigma(k)>\sigma(l)}\right|.
\end{equation}
Summing the separating-hyperplane counts proves (S2).

\emph{The alcove action and its stabilizer (S1).} First, we show by induction on $N=\#\fH(\Delta,\mathsf B)$ that every alcove $\mathsf B$ is $g\Delta$ for some $g\in\langle S\rangle$. If $N=0$, then $\mathsf B=\Delta$ by \eqref{eq:alc-a}, while if $N>0$, then $\mathsf B\not\subseteq\Delta$, so some $H_i$ separates $\mathsf B$ from $\Delta$. In this case \eqref{eq:alc-b} gives $\fH(s_i\Delta,\mathsf B)=\fH(\Delta,\mathsf B)\smallsetminus\{H_i\}$, and applying $s_i$ we get $\#\fH(\Delta,s_i\mathsf B)=N-1$, thus $s_i\mathsf B=g'\Delta$ by induction and $\mathsf B=s_ig'\Delta$.

Next, if $w\in\tW$ satisfies $w\Delta=\Delta$, then $w$ has no inversions by (S2), so $w(1)<w(2)<\dots<w(m)<w(1)+m$, and consequently $w(a)=w(1)+a-1$ for all $a$ and $w=\tau^{w(1)-1}$. As $\val(s_i)=0$ and $\val(\tau)=1$, we have $\langle S\rangle\cap\langle\tau\rangle=1$, thus $\langle S\rangle$ acts simply transitively on the alcoves. Finally, every $w\in\tW$ is $g\tau^c$ with $g\in\langle S\rangle$ and $c=\val(w)$, and it follows that $\langle S\rangle=\ker\val$.

\emph{Reflection descent (S3).} For alcoves $\mathsf A,\mathsf B$ and $H\in\fH(\mathsf A,\mathsf B)$ we claim that $\#\fH(s_H\mathsf A,\mathsf B)<\#\fH(\mathsf A,\mathsf B)$. Let us choose $a\in\mathsf A$ and $b\in\mathsf B$ and put $a'=s_H(a)$. Then the hyperplane $H$ separates $a$ from $b$ but does not separate $a'$ from $b$. Let $H'\ne H$ be a hyperplane of $\fH$ with an affine equation $\phi'=0$, and put $H''=s_HH'$ and $\phi''=\phi'\circ s_H$. If $H''=H'$, then $\phi'$ and $\phi''$ are proportional and agree on $H\not\subseteq H'$, so $\phi''=\phi'$ and $H'$ separates $a'$ from $b$ if and only if it separates $a$ from $b$.

Suppose now that $H''\ne H'$, and let $\epsilon_1,\epsilon_2,\epsilon_3,\epsilon_4$ be the signs of $\phi'(a),\phi''(a),\phi'(b),\phi''(b)$. Since $\phi'(a')=\phi''(a)$ and $\phi''(a')=\phi'(a)$, the pair $\{H',H''\}$ contributes $[\epsilon_1\ne\epsilon_3]+[\epsilon_2\ne\epsilon_4]$ to $\#\fH(\mathsf A,\mathsf B)$ and $[\epsilon_2\ne\epsilon_3]+[\epsilon_1\ne\epsilon_4]$ to $\#\fH(s_H\mathsf A,\mathsf B)$.

These two numbers are equal if $\epsilon_1=\epsilon_2$ or $\epsilon_3=\epsilon_4$, and otherwise $a$ and $b$ lie in the union of the convex sets $C_1=\{\phi'>0>\phi''\}$ and $C_2=\{\phi'<0<\phi''\}=s_HC_1$. These sets do not meet $H$, because $\phi'=\phi''$ on $H$, so they lie on opposite sides of $H$. As $a$ and $b$ lie on opposite sides of $H$ as well, one of them is in $C_1$ and the other one is in $C_2$. Then $\epsilon_1\ne\epsilon_3$, $\epsilon_2\ne\epsilon_4$, $\epsilon_2=\epsilon_3$ and $\epsilon_1=\epsilon_4$, and the contributions are $2$ and $0$, which proves the claim. Now by (S2) we have:
\[
\ell(s_Hw)=\#\fH(\Delta,s_Hw\Delta)=\#\fH(s_H\Delta,w\Delta).
\]
If $H\in\fH(\Delta,w\Delta)$, then the claim for $(\Delta,w\Delta)$ gives $\ell(s_Hw)<\ell(w)$, and otherwise $H$ separates $s_H\Delta$ from $w\Delta$ and the claim for $(s_H\Delta,w\Delta)$ gives $\ell(w)<\ell(s_Hw)$.

\emph{Point stabilizers (S4).} The equations of $H_1,\dots,H_{m-1}$ give $v_1=v_m$, which contradicts the equation of $H_0$. Thus the hyperplanes $H_0,\dots,H_{m-1}$ have no common point, and $K_v\ne S$. Let $K\subsetneq S$, and let $v$ be the barycenter of the vertices $\omega_k$ with $s_k\notin K$. Then $v$ lies on $H_i$ if and only if $s_i\in K$, so $K=K_v$.

Now let $v\in\overline\Delta$, and let us show by induction on $N=\#\fH(\Delta,\mathsf B)$ that $x_{\mathsf B}\in W_{K_v}$ for every alcove $\mathsf B$ with $v\in\overline{\mathsf B}$. For $N=0$ we have $\mathsf B=\Delta$ and $x_{\mathsf B}=1$, while for $N>0$ some $H_i$ separates $\mathsf B$ from $\Delta$, as in the proof of (S1), and $v\in H_i$ by \eqref{eq:alc-a}. In this case $s_i\in K_v$ and $v=s_i(v)\in\overline{s_i\mathsf B}$, and since $\#\fH(\Delta,s_i\mathsf B)=N-1$, the induction hypothesis gives $x_{\mathsf B}=s_i\,x_{s_i\mathsf B}\in W_{K_v}$. Now if $g\in W_a$ fixes $v$, then $v\in\overline{g\Delta}$ and $g=x_{g\Delta}\in W_{K_v}$. Conversely, $W_{K_v}$ fixes $v$, thus $W_{K_v}$ is the stabilizer of $v$ in $W_a$, and in particular $s_H\in W_{K_v}$ if $v\in H$. Finally, the stabilizer is finite: an element $t^\lambda\sigma$ of $W_a$ fixing $v$ is determined by $\sigma$, because $\lambda$ is congruent to $v-\sigma v$ modulo $\R(1,\dots,1)$ and has coordinate sum $0$.

\emph{Parabolic coset representatives (S5).} By (S4), $K=K_v$ for some $v\in\overline\Delta$, and $W_K$ is the stabilizer of $v$ in $W_a$. Let $x$ be an element of minimal length in $W_Kw$, and note that no hyperplane $H\ni v$ separates $\Delta$ from $x\Delta$: otherwise $s_H\in W_K$ by (S4) and $\ell(s_Hx)<\ell(x)$ by (S3). For $y\in W_K$ we have $y(v)=v$, so $v$ lies in the closures of $\Delta$, $y\Delta$ and $y^{-1}\Delta$. Let us count the hyperplanes $H$ separating $\Delta$ from $yx\Delta$, or equivalently the hyperplanes $H$ such that $y^{-1}H$ separates $y^{-1}\Delta$ from $x\Delta$. First, if $v\notin H$, then $v\notin y^{-1}H$, and $\Delta$ and $y^{-1}\Delta$ lie on the same side of $y^{-1}H$ by \eqref{eq:alc-a}, so the condition is that $y^{-1}H$ separates $\Delta$ from $x\Delta$. The number of such $H$ is the number of hyperplanes in $\fH(\Delta,x\Delta)$ not containing $v$, which is $\ell(x)$.

Next, if $v\in H$, then $v\in y^{-1}H$, and $\Delta$ and $x\Delta$ lie on the same side of $y^{-1}H$, so the condition is that $y^{-1}H$ separates $y^{-1}\Delta$ from $\Delta$, or equivalently that $H\in\fH(\Delta,y\Delta)$. By \eqref{eq:alc-a} every element of $\fH(\Delta,y\Delta)$ contains $v$, hence the number of such $H$ is $\ell(y)$.

Thus $\ell(yx)=\ell(y)+\ell(x)$. An element of $W_a$ of length $0$ fixes $\Delta$ and is trivial by (S1), so $\ell(yx)>\ell(x)$ for $y\ne1$, and in particular $x$ is unique.

Finally, let $x'\in\tW$ with $\ell(sx')>\ell(x')$ for all $s\in K$, and write $x'=yx$ with $y\in W_K$ and $x$ of minimal length. Suppose that $y\ne1$. Then $y\Delta\ne\Delta$ and $v\in\overline{y\Delta}$, so by the proof of (S4) some $H_i$ with $s_i\in K$ separates $\Delta$ from $y\Delta$, hence $\ell(s_iy)<\ell(y)$ by (S3) and $\ell(s_ix')=\ell(s_iy)+\ell(x)<\ell(x')$, a contradiction. Thus $y=1$.
\end{proof}

\begin{lemma}[Expressions of prescribed length]\label{lem:alcove-reduced-expression}
If $w\in\tW$, $k=\val(w)$ and $L=\ell(w)$, then
\[
 w=s_{i_1}\cdots s_{i_L}\tau^k
 \qquad(s_{i_j}\in S).
\]
Consequently, for fixed integers $k$ and $L\ge0$, the set
\[
 \mathcal W_{k,\le L}
 =\{w\in\tW:\val(w)=k,\ \ell(w)\le L\}
\]
is finite, of cardinality at most $\sum_{j=0}^{L}m^j$.
\end{lemma}
\begin{proof}
If $L=0$, (S1) and (S2) give $w=\tau^k$. If $L>0$, then
$w\Delta\ne\Delta$, and one of the defining walls $H_i$ of $\Delta$
separates $\Delta$ from $w\Delta$, as in the proof of (S1).
Since $\Delta$ and $s_i\Delta$ are separated only by $H_i$,
\[
 \ell(s_iw)
 =\#\fH(s_i\Delta,w\Delta)=\ell(w)-1.
\]
Also $\val(s_iw)=k$. Induction gives the stated expression.
There are at most $m^j$ words of length $j$ in the alphabet $S$,
which proves the bound. In particular, an element of $W_a$ of
length one lies in $S$.
\end{proof}

For $w\in\tW$ and an alcove $\mathsf A$ we put $w_{\mathsf A}=x_{\mathsf A}^{-1}wx_{\mathsf A}$ and note that $\ell(w_{\mathsf A})=\#\fH(\mathsf A,w\mathsf A)$ by (S2). We say that alcoves $\mathsf A\ne\mathsf A'$ are \emph{adjacent along} $H\in\fH$ if $\mathsf A'=s_H\mathsf A$ and $\overline{\mathsf A}\cap\overline{\mathsf A'}$ contains a point $h$ lying on $H$ and on no other hyperplane, and we call $\overline{\mathsf A}\cap\overline{\mathsf A'}$ their \emph{common face}. In this case $x_{\mathsf A}^{-1}h$ is a point of the boundary of $\Delta$ lying on the single hyperplane $x_{\mathsf A}^{-1}H$, hence this hyperplane is one of the $H_i$ and $s:=x_{\mathsf A}^{-1}s_Hx_{\mathsf A}=s_i\in S$. Consequently $x_{\mathsf A'}=x_{\mathsf A}s$, $w_{\mathsf A'}=sw_{\mathsf A}s$, and $\fH(\mathsf A,\mathsf A')=\{H\}$ by \eqref{eq:alc-b}.

For $s\in S$ and $w'\in\tW$ the lengths of $w'$ and $sw'$ differ by $1$, by (S2) and \eqref{eq:alc-b}, and the same holds for $w'$ and $w's$. Thus the lengths of $w'$ and $sw's$ differ by $0$ or $\pm2$. We write $w\to w'$ if $w'$ is obtained from $w$ by a sequence of conjugations by simple reflections none of which increases the length.
We write $w\approx w'$ if a sequence of length-preserving simple
conjugations and conjugations by powers of $\tau$ connects them.
An element $w$ is \emph{terminal} if no $w'\approx w$ admits
$s\in S$ with $\ell(sw's)<\ell(w')$. An element $x$ is
\emph{straight} if $\ell(x^n)=n\ell(x)$ for every $n\ge1$.
A \emph{standard triple} $(K,x,u)$ consists of $K\subsetneq S$,
$u\in W_K$, and $x\in\tW$ with $xKx^{-1}=K$ and
$\ell(sx)>\ell(x)$ for every $s\in K$.

For $w\in\tW$, write $w(v)=\bar wv+w(0)$ for its action on
$\mathbb V$, with $\bar w$ orthogonal. Define its \emph{reduced Newton
vector} by
\begin{equation}\label{eq:reduced-newton}
 \nu_w^{\mathrm{red}}
 =\nu_w^{\mathrm{full}}
   -\frac{\val(w)}m(1,\ldots,1)\in\mathbb V.
\end{equation}
Thus $w^N$ acts on $\mathbb V$ as translation by
$N\nu_w^{\mathrm{red}}$, and
\[
 \nu_w^{\mathrm{red}}=0
 \quad\Longleftrightarrow\quad
 \nu_w^{\mathrm{full}}\text{ is central}.
\]
Although $\langle\tau\rangle$ is infinite, its action on $\mathbb V$
factors through the finite cyclic group
$\langle\tau\rangle/\langle\tau^m\rangle$ of order $m$.
Length and all alcove calculations factor through
$\tW/\langle\tau^m\rangle$; the valuation records the central
information lost in this quotient. In the gradient-flow and straightness arguments we
abbreviate $\nu_w^{\mathrm{red}}$ by $\nu$ when $w$ is fixed.

\subsubsection{The function \texorpdfstring{$f_w$}{f-w} and its gradient flow}
\begin{lemma}[{\cite[Lemmas~2.2, 2.3, Cor.~2.4]{HN}}]\label{lem:Vw}
\begin{enumerate}[label=\textup{(\roman*)}]
\item We have $\bar w\nu=\nu$, the set $\mathbb V_w=\{v\in\mathbb V:w(v)=v+\nu\}$ is a nonempty affine subspace with direction $\mathsf F=\ker(\bar w-1)$, and $w(\mathbb V_w)=\mathbb V_w$.
\item Let $p\in\mathbb V_w$. The map $M=(\bar w-1)|_{\mathsf F^\perp}$ is an automorphism of $\mathsf F^\perp$, and for $x\in\mathsf F$ and $y\in\mathsf F^\perp$ we have:
\[
f_w(p+x+y):=\|w(p+x+y)-(p+x+y)\|^2=\|My\|^2+\|\nu\|^2 .
\]
\item The maps $\phi_s(p+x+y)=p+x+e^{-2sM^{T}M}y$, $s\ge0$, satisfy $\frac d{ds}\phi_s(z)=-\nabla f_w(\phi_s(z))$, and $\phi_s(z)\to\varpi_w(z):=p+x$ as $s\to+\infty$.
\end{enumerate}
\end{lemma}
\begin{proof}
(i) Since $w$ commutes with $w^N$, the translation by $N\nu$, we have $w(v+N\nu)=w(v)+N\nu$, while the left side equals $w(v)+N\bar w\nu$. For any $q\in\mathbb V$ the point $p_0=\frac1N\sum_{i=0}^{N-1}w^i(q)$ satisfies the equation defining $\mathbb V_w$:
\[
w(p_0)-p_0=\frac1N(w^N(q)-q)=\nu.
\]
The equation $w(v)=v+\nu$ is equivalent to $(\bar w-1)v=\nu-w(0)$, hence $\mathbb V_w$ is an affine subspace with direction $\mathsf F$, and $w(v)=v+\nu\in\mathbb V_w$ for $v\in\mathbb V_w$, because $\nu\in\mathsf F$.

(ii) We have:
\[
w(p+x+y)-(p+x+y)=(w(p)-p)+(\bar w-1)(x+y)=\nu+My,
\]
and $My\in\mathsf F^\perp$ is orthogonal to $\nu$, because the space $\mathsf F^\perp$ is $\bar w$-stable, $\bar w$ being orthogonal. Finally, $M$ is injective on $\mathsf F^\perp$.

(iii) We have $\nabla f_w(p+x+y)=2M^{T}My$, and $M^{T}M$ is positive definite on $\mathsf F^\perp$.
\end{proof}

\begin{lemma}[{\cite[Lemma~2.1]{HN}}]\label{lem:deriv}
Let $\mathsf A,\mathsf A'$ be alcoves adjacent along $H\in\fH$, let $h$ be a point of their common face which lies on no other hyperplane of $\fH$, and let $\mathsf n$ be the unit normal of $H$ pointing out of $\mathsf A$. If $\ell(w_{\mathsf A'})=\ell(w_{\mathsf A})+2$, then $\langle\nabla f_w(h),\mathsf n\rangle>0$.
\end{lemma}
\begin{proof}
The alcoves $\mathsf A,\mathsf A'$ are separated only by $H$, and the alcoves $w\mathsf A,w\mathsf A'$ are separated only by $wH$, thus a hyperplane other than $H,wH$ separates $\mathsf A'$ from $w\mathsf A'$ if and only if it separates $\mathsf A$ from $w\mathsf A$. Since $\#\fH(\mathsf A',w\mathsf A')=\#\fH(\mathsf A,w\mathsf A)+2$, we get $H\ne wH$, and neither $H$ nor $wH$ separates $\mathsf A$ from $w\mathsf A$. The point $w(h)$ lies on $wH$ and on no other hyperplane, hence $w(h)\notin H$ and $h\notin wH$.

The alcove $\mathsf A$ lies in $\{z:\langle z-h,\mathsf n\rangle<0\}$, and so does $w\mathsf A$, thus $\langle w(h)-h,\mathsf n\rangle<0$, because $w(h)\in\overline{w\mathsf A}\smallsetminus H$. On the other hand, the vector $\bar w\mathsf n$ is the unit normal of $wH$ pointing out of $w\mathsf A$. Since the alcove $\mathsf A$ lies on the same side of $wH$ as $w\mathsf A$ and $h\in\overline{\mathsf A}\smallsetminus wH$, we get $\langle h-w(h),\bar w\mathsf n\rangle<0$. Finally, $w(h+s\mathsf n)-(h+s\mathsf n)=(w(h)-h)+s(\bar w\mathsf n-\mathsf n)$, and we obtain the following expression for the derivative of $f_w$ at $h$ in the direction $\mathsf n$:
\[
\langle\nabla f_w(h),\mathsf n\rangle=2\langle w(h)-h,\bar w\mathsf n-\mathsf n\rangle=2\langle w(h)-h,\bar w\mathsf n\rangle-2\langle w(h)-h,\mathsf n\rangle>0 .\qedhere
\]
\end{proof}

For an affine subspace $\mathsf Y\subseteq\mathbb V$ we denote by $\fH_{\mathsf Y}$ the set of hyperplanes containing $\mathsf Y$, and we call a point of $\mathsf Y$ \emph{regular} if it lies on no hyperplane of $\fH$ which does not contain $\mathsf Y$.

\begin{lemma}\label{lem:null}
Let $d\ge2$, let $C\subseteq\R^{d-1}$ be a compact cube and let $\Psi:C\to\R^d$ be a Lipschitz map. Then $\Psi(C)$ has Lebesgue measure zero.
\end{lemma}
\begin{proof}
Let $R$ be the side length of $C$, let $L$ be a Lipschitz constant of $\Psi$, and for $N\ge1$ let us divide $C$ into $N^{d-1}$ cubes of side $R/N$. The image of each of them lies in a ball of radius $L\sqrt{d-1}\,R/N$, thus $\Psi(C)$ is covered by sets of total measure at most $N^{d-1}\,\omega_d\,(L\sqrt{d-1}\,R/N)^d$, where $\omega_d$ is the volume of the unit ball, and this bound tends to zero as $N\to\infty$.
\end{proof}

\begin{proposition}[{\cite[Prop.~2.5]{HN}}]\label{prop:flow}
For every $w\in\tW$ and every alcove $\mathsf A$ there is an alcove $\mathsf A'$ with $w_{\mathsf A}\to w_{\mathsf A'}$ whose closure contains a regular point of $\mathbb V_w$.
\end{proposition}
\begin{proof}
We follow \cite[proof of Prop.~2.5]{HN}, and we put $d=\dim\mathbb V=m-1$ and $\gamma_y(s)=\phi_s(y)$ for $y\in\mathbb V$ and $s\ge0$.

\emph{Step 1. Choice of a generic starting point.} Let $Z_1$ be the set of $y\in\mathbb V$ such that $\varpi_w(y)$ is not a regular point of $\mathbb V_w$. The non-regular points of $\mathbb V_w$ are the points of the sets $H\cap\mathbb V_w$, $H\in\fH\smallsetminus\fH_{\mathbb V_w}$, each of which is empty or an affine hyperplane of $\mathbb V_w$. As $\varpi_w:\mathbb V\to\mathbb V_w$ is affine and surjective, $Z_1$ is a countable union of affine hyperplanes of $\mathbb V$, hence a null set.

Let us denote by $\mathbb V^{(2)}$ the set of points lying on at least two hyperplanes of $\fH$. It is the union of the countably many affine subspaces $P=H\cap H'$ of dimension $d-2$, where $H\ne H'$ are not parallel, and in particular it is empty if $d=1$. Let $Z_2$ be the set of $y$ such that $\gamma_y(s)\in\mathbb V^{(2)}$ for some $s\ge0$. The map $\phi_s$ is bijective, and its inverse $\phi_{-s}$ is given by the formula of Lemma~\ref{lem:Vw}(iii) with $-s$ in place of $s$. Hence $Z_2$ is the union of the sets $\Psi(P\times[0,\infty))$, where $\Psi(z,s)=\phi_{-s}(z)$ is a smooth map $\mathbb V\times\R\to\mathbb V$. Each $P\times[0,\infty)$ is a countable union of compact cubes in an affine space of dimension $d-1$, and $\Psi$ is Lipschitz on each of them. If $d=1$, then $Z_2$ is empty. If $d\ge2$, Lemma~\ref{lem:null} applied to each bounded parameter cube shows that $Z_2$ is a null set.

\emph{Step 2. The limit of the path.} As the alcove $\mathsf A$ is open and nonempty, there is a point $y_0\in\mathsf A\smallsetminus(Z_1\cup Z_2)$, and we put $\gamma=\gamma_{y_0}$ and $v=\varpi_w(y_0)$. Then $v$ is a regular point of $\mathbb V_w$, every point $\gamma(s)$, $s\ge0$, lies on at most one hyperplane of $\fH$, and by Lemma~\ref{lem:Vw}(iii) we have $\gamma(s)\to v$ as $s\to\infty$ and $\gamma'(s)=-\nabla f_w(\gamma(s))$.

Since $\fH$ is locally finite, there is $\varepsilon>0$ such that every hyperplane of $\fH$ meeting the ball $B(v,\varepsilon)$ contains $v$. We claim that $v\in\overline{\mathsf B}$ for every alcove $\mathsf B$ whose closure meets $B(v,\varepsilon)$. Indeed, such an alcove $\mathsf B$ contains a point $z\in B(v,\varepsilon)$, and no hyperplane $H\in\fH$ meets the segment from $z$ to $v$ in a point other than $v$. Otherwise $H$ contains $v$ and a second point of the line through $z$ and $v$, and therefore also $z$, which is impossible. Thus this segment without its end point $v$ is a connected subset of $\mathbb V\smallsetminus\bigcup\fH$ containing $z$, hence it lies in $\mathsf B$, and $v\in\overline{\mathsf B}$. Finally, since $\gamma(s)\to v$ as $s\to\infty$, we can fix $t_\ast>0$ with $\gamma(t_\ast)\in B(v,\varepsilon)$.

\emph{Step 3. The alcoves along the path.}
Taking the last exit from each alcove will give a finite gallery
with no repetitions. Define alcoves $\mathsf A_0,\mathsf A_1,\dots$ and numbers $0<t_0<t_1<\cdots\le t_\ast$, starting with $\mathsf A_0=\mathsf A$. If $\mathsf A_i$ is defined and $\gamma(t)\in\overline{\mathsf A_i}$ for some $t\in[0,t_\ast]$, we put:
\[
t_i=\max\{t\in[0,t_\ast]:\ \gamma(t)\in\overline{\mathsf A_i}\}.
\]
If $t_i=t_\ast$, we stop, and otherwise we put $h_i=\gamma(t_i)$. The point $h_i$ is not in the open set $\mathsf A_i$, because $\gamma(t)\notin\overline{\mathsf A_i}$ for $t>t_i$. Thus $h_i$ lies on a hyperplane of $\fH$, and by the choice of $y_0$ it lies on exactly one of them, which we denote by $H^{(i)}$. We choose $r>0$ such that the ball $B(h_i,r)$ meets no other hyperplane. Then $B(h_i,r)\smallsetminus H^{(i)}$ has two connected components, each of which is contained in an alcove, and these two alcoves are different, being separated by $H^{(i)}$. One of the components, $B^-$, lies in $\mathsf A_i$, because $h_i\in\overline{\mathsf A_i}$. The other one, $B^+$, lies in the alcove $\mathsf A_{i+1}:=s_{H^{(i)}}\mathsf A_i$, because $s_{H^{(i)}}$ exchanges $B^-$ and $B^+$. Thus $\mathsf A_i$ and $\mathsf A_{i+1}$ are adjacent along $H^{(i)}$, and the point $h_i$ lies in their common face and on no other hyperplane. Moreover, the closure of $\mathsf A_i$ meets the ball $B(h_i,r)$ in the following set:
\[
\overline{\mathsf A_i}\cap B(h_i,r)=B^-\cup(H^{(i)}\cap B(h_i,r)).
\]
For $t>t_i$ close to $t_i$ the point $\gamma(t)$ lies in $B(h_i,r)$ but not in $\overline{\mathsf A_i}$, so it lies in $B^+\subseteq\mathsf A_{i+1}$, and in particular $t_{i+1}$ is defined and $t_{i+1}>t_i$.

For $j>i$ the alcove $\mathsf A_j$ contains a point $\gamma(t)$ with $t>t_{j-1}\ge t_i$, and this point is not in $\overline{\mathsf A_i}$, so the alcoves $\mathsf A_0,\mathsf A_1,\dots$ are pairwise distinct. Since the compact set $\gamma([0,t_\ast])$ meets the closures of only finitely many alcoves, the construction stops, and there is $N\ge0$ with $t_N=t_\ast$. We put $\mathsf A'=\mathsf A_N$, and then $\gamma(t_\ast)\in\overline{\mathsf A'}\cap B(v,\varepsilon)$, hence $v\in\overline{\mathsf A'}$ by Step~2.

\emph{Step 4. Length does not increase.} Let $0\le i<N$, and let $\mathsf n$ be the unit normal of $H^{(i)}$ pointing out of $\mathsf A_i$. Then $B^+=\{z\in B(h_i,r):\langle z-h_i,\mathsf n\rangle>0\}$. The function $\psi(t)=\langle\gamma(t)-h_i,\mathsf n\rangle$ vanishes at $t_i$ and is positive for $t>t_i$ close to $t_i$, thus we get:
\[
0\le\psi'(t_i)=\langle\gamma'(t_i),\mathsf n\rangle=-\langle\nabla f_w(h_i),\mathsf n\rangle ,
\]
and Lemma~\ref{lem:deriv} shows that $\ell(w_{\mathsf A_{i+1}})\ne\ell(w_{\mathsf A_i})+2$. On the other hand, the alcoves $\mathsf A_i$ and $\mathsf A_{i+1}$ are adjacent, hence $w_{\mathsf A_{i+1}}$ is the conjugate of $w_{\mathsf A_i}$ by a simple reflection, and the two lengths differ by $0$ or $\pm2$. Therefore $\ell(w_{\mathsf A_{i+1}})\le\ell(w_{\mathsf A_i})$ for every $0\le i<N$, which means that $w_{\mathsf A}=w_{\mathsf A_0}\to w_{\mathsf A_N}=w_{\mathsf A'}$.
\end{proof}

\begin{remark}\label{rem:flowremarks}
The generic choice of the starting point ensures that each crossing
lies in the relative interior of a common face, as required by
Lemma~\ref{lem:deriv}; the lemma makes no assertion about
crossings at an intersection of walls. Only the segment up to $t_\ast$ is used;
the last-exit construction does not require a transversality assumption
on the trajectory.

\end{remark}

\subsubsection{Straightness and the standard triple}
\begin{lemma}\label{lem:straight}
Let $x\in\tW$, let $\mathsf Y\subseteq\mathbb V$ be an affine subspace, and let $\nu_1\in\mathbb V$ be such that $\mathsf Y+\nu_1=\mathsf Y$ and $x(z)=z+\nu_1$ for all $z\in\mathsf Y$. Let $v_1\in\mathsf Y\cap\overline\Delta$ be a regular point of $\mathsf Y$, and assume that no hyperplane of $\fH_{\mathsf Y}$ separates $\Delta$ from $x\Delta$. Then $\ell(x^n)=n\,\ell(x)$ for all $n\ge1$, and $\ell(x)=0$ if $\nu_1=0$.
\end{lemma}
\begin{proof}
Since $\mathsf Y+\nu_1=\mathsf Y$, the element $x$ maps $\mathsf Y$ onto itself, and therefore it permutes $\fH_{\mathsf Y}$ and $\fH\smallsetminus\fH_{\mathsf Y}$. For $i\ge0$ we put $\mathsf A_i=x^i\Delta$ and $v^{(i)}=x^i(v_1)=v_1+i\nu_1\in\overline{\mathsf A_i}$, and we note that every $v^{(i)}$ is a regular point of $\mathsf Y$. Applying $x^{i-1}$ to the hypothesis, we see that no hyperplane of $\fH_{\mathsf Y}$ separates $\mathsf A_{i-1}$ from $\mathsf A_i$, and consequently none of them separates $\mathsf A_0$ from $\mathsf A_n$.

For $H\in\fH\smallsetminus\fH_{\mathsf Y}$ no $v^{(i)}$ lies on $H$, hence $\mathsf A_i$ lies on the side of $H$ containing $v^{(i)}$, and $H$ separates $\mathsf A_i$ from $\mathsf A_j$ if and only if it separates the points $v^{(i)}$ and $v^{(j)}$.

These points are equal (if $\nu_1=0$) or equally spaced on a line. Hence $H$ separates $v^{(0)}$ from $v^{(n)}$ if and only if it separates $v^{(i-1)}$ from $v^{(i)}$ for exactly one $i\in\{1,\dots,n\}$, and it never does so for two values of $i$. Therefore we have:
\[
\#\fH(\mathsf A_0,\mathsf A_n)=\sum_{i=1}^n\#\fH(\mathsf A_{i-1},\mathsf A_i)=n\,\#\fH(\Delta,x\Delta),
\]
which is the first claim by (S2). If $\nu_1=0$, then no hyperplane separates $\Delta$ from $x\Delta$, and $\ell(x)=0$ by (S2).
\end{proof}

\begin{proposition}[After He--Nie, Propositions~2.5 and~2.8]\label{prop:HN3app}
For every $w\in\tW$ there are $K\subsetneq S$, $u\in W_K$ and $x\in\tW$ such that $w\to ux$, $\ell(sx)>\ell(x)$ for all $s\in K$, and $xKx^{-1}=K$. Moreover, $x$ is straight, $\val(x)=\val(w)$, and $\nu_x^{\mathrm{full}}$ is the image of $\nu_w^{\mathrm{full}}$ under an element of $S_m$. Finally, if the Newton point of $w$ is central, then $x\in\langle\tau\rangle$.
\end{proposition}
\begin{proof}
By Proposition~\ref{prop:flow} with $\mathsf A=\Delta$, there is an alcove $\mathsf A'$ with $w\to w_1:=w_{\mathsf A'}$ such that $\overline{\mathsf A'}$ contains a regular point $v$ of $\mathbb V_w$. Write $\pi_{\mathsf A'}\in S_m$ for the permutation part of $x_{\mathsf A'}$. We put $\mathsf Y=x_{\mathsf A'}^{-1}\mathbb V_w$ and $v_1=x_{\mathsf A'}^{-1}v\in\overline\Delta$, and then $\mathsf Y=\mathbb V_{w_1}$, with $\nu_{w_1}^{\mathrm{red}}=\pi_{\mathsf A'}^{-1}\nu_w^{\mathrm{red}}=:\nu_1$, and $v_1$ is a regular point of $\mathsf Y$. Let $K=K_{v_1}$. Every $H\in\fH_{\mathsf Y}$ contains $v_1$,
so $s_H\in W_K$ by (S4). Conversely, for $s_i\in K$ the wall
$H_i$ contains $v_1$ and hence contains $\mathsf Y$ by regularity.
Therefore
\[
 W_K=\langle s_H:H\in\fH_{\mathsf Y}\rangle,
\]
and this group fixes $\mathsf Y$ pointwise. Since
$w_1(\mathsf Y)=\mathsf Y$, the element $w_1$ permutes
$\fH_{\mathsf Y}$, hence $w_1W_Kw_1^{-1}=W_K$.
Write $w_1=ux$ as in (S5).

\emph{$xKx^{-1}=K$.} The element $x=u^{-1}w_1$ normalizes $W_K$, and for $s\in K$ the element $y=xsx^{-1}\in W_K$ satisfies:
\[
\ell(y)+\ell(x)=\ell(yx)=\ell(xs)\le\ell(x)+1,
\]
hence $\ell(y)\le1$. Since $y\ne1$ and $y\in W_a$, (S1)--(S2)
give $\ell(y)=1$, and Lemma~\ref{lem:alcove-reduced-expression}
gives $y\in S$. The element $y\in W_K$ fixes $v_1$, so
$y\in K_{v_1}=K$. Thus $xKx^{-1}\subseteq K$; conjugation is
injective and $K$ is finite, so equality holds.

\emph{No hyperplane of $\fH_{\mathsf Y}$ separates $\Delta$ from $x\Delta$.} Otherwise we would have $\ell(s_Hx)<\ell(x)$ by (S3), with $s_H\in W_K$, contradicting (S5).

\emph{$x(z)=z+\nu_1$ on $\mathsf Y$.} For $z\in\mathsf Y$ we have:
\[
x(z)=u^{-1}(w_1(z))=u^{-1}(z+\nu_1)=z+\nu_1,
\]
because $z+\nu_1=w_1(z)\in\mathsf Y$ is fixed by $u^{-1}\in W_K$.

In particular $\mathsf Y+\nu_1=w_1(\mathsf Y)=\mathsf Y$, and by Lemma~\ref{lem:straight}, $x$ is straight, and if $\nu_1=0$, then $\ell(x)=0$ and hence $x\in\langle\tau\rangle$ by (S1) and (S2). A power $x^N$ of $x$ acts on $\mathbb V$ as a translation, and on $\mathsf Y$ it is $z\mapsto z+N\nu_1$, so $\nu_x^{\mathrm{red}}=\nu_1$. Also
$\val(x)=\val(w_1)=\val(w)$, because $u,x_{\mathsf A'}\in W_a$.
Together with \eqref{eq:reduced-newton} and
\eqref{eq:newton-conjugation}, this proves
$\nu_x^{\mathrm{full}}=
 \pi_{\mathsf A'}^{-1}\nu_w^{\mathrm{full}}$.
\end{proof}

\begin{corollary}[Terminal form of the reduction]\label{cor:terminal-standard}
If $w\in\tW$ is terminal, there is a standard triple $(K,x,u)$
with $x$ straight and $w\approx ux$. The full Newton point of
$x$ is central if and only if that of $w$ is central. In this
case $x=\tau^{\val(w)}$.
\end{corollary}
\begin{proof}
Take $K,u,x$ from Proposition~\ref{prop:HN3app}. Along the chain
$w\to ux$ no step increases length. If a first strictly decreasing
step existed, all its preceding steps would preserve length and
would produce $w'\approx w$ admitting a length-decreasing simple
conjugation, contrary to terminality. Thus $w\approx ux$.
The full Newton vectors of $w$ and $x$ differ by a coordinate
permutation, so centrality is equivalent. In the central case,
Proposition~\ref{prop:HN3app} gives $x\in\langle\tau\rangle$ and
$\val(x)=\val(w)$; since $\val(\tau)=1$, this gives
$x=\tau^{\val(w)}$.
\end{proof}

\subsubsection{Regularity and effective reduction}

The flow argument proves that a suitable alcove exists. To find one
by a finite search, we first express regularity as a system of rational
difference inequalities.

\begin{proposition}[A regularity criterion]\label{prop:regularity-test}
Let $w'=t^\lambda\sigma\in\tW$, and let $c$ run over the $r$ cycles
of $\sigma$. Choose $v^0\in\mathbb Q^m$ satisfying
\[
 w'(v^0)=v^0+\nu_{w'}^{\mathrm{full}},
\]
and put $d_{ij}=v_i^0-v_j^0$. Such a choice is effective: if $N$
is the order of $\sigma$, take
\[
 v^0=\frac1N\sum_{a=0}^{N-1}(w')^a(0).
\]
Choose a positive integer $D$ such that every $Dd_{ij}$ is integral,
and set $\varepsilon=1/[D(r+1)]$.

Form a directed multigraph with the cycles of $\sigma$ as vertices.
For every $i<j$ in different cycles, insert the two edges
\[
 \begin{array}{c|c}
 \text{edge}&\text{weight}\\ \hline
 c(j)\longrightarrow c(i)&1-d_{ij}\\
 c(i)\longrightarrow c(j)&d_{ij}.
 \end{array}
\]
\Needspace{8\baselineskip}
Let the \emph{modified weight} of each edge be its weight minus
$\varepsilon$. Then $\overline\Delta$ contains a regular point of
$\mathbb V_{w'}$ if and only if both of the following conditions hold:
\begin{enumerate}[label=\textup{(\roman*)}]
\item $0\le d_{ij}\le1$ whenever $i<j$ lie in the same cycle;
\item the graph with modified weights has no negative-weight
      directed cycle.
\end{enumerate}
When these conditions hold, adjoining a source with zero-weight
edges to every vertex and taking shortest-path distances $x_c$
gives a rational regular point: take the image in $\mathbb V$ of
\[
 v^0+\sum_c x_c\chi_c,
\]
where $\chi_c$ is the indicator vector of $c$.
\end{proposition}
\begin{proof}
The averaging identity used in Lemma~\ref{lem:Vw}(i), now in
$\mathbb R^m$, gives the displayed formula for $v^0$.
The full solution space of
$w'(v)=v+\nu_{w'}^{\mathrm{full}}$ is
\begin{equation}\label{eq:regularity-cycle-parameters}
 \left\{v^0+\sum_c x_c\chi_c:x_c\in\mathbb R\right\},
\end{equation}
because $\ker(\sigma-1)$ consists of vectors constant on each cycle.
Its image in $\mathbb V$ is exactly $\mathbb V_{w'}$.
Indeed, equality after quotienting by the diagonal line means that
$w'(v)-v-\nu_{w'}^{\mathrm{full}}$ is diagonal; its coordinate sum
is zero by \eqref{eq:full-newton-cycles}, so that diagonal vector
is zero.

An affine root hyperplane $\{v_i-v_j=a\}$ contains this affine
space if and only if $i,j$ are in the same cycle and $d_{ij}=a$.
Also $v\in\overline\Delta$ is equivalent to
$0\le v_i-v_j\le1$ for every $i<j$. Thus a regular point in the
closure exists precisely when condition~(i) holds and the following
strict system is feasible:
\begin{equation}\label{eq:regularity-strict-system}
 0<d_{ij}+x_{c(i)}-x_{c(j)}<1
 \qquad(i<j,\ c(i)\ne c(j)).
\end{equation}
The two inequalities for a pair encode the two edges above, with
an edge $a\to b$ of weight $w_e$ meaning $x_b-x_a<w_e$.
Summing around a directed cycle shows that feasibility implies
strict positivity of every cycle weight.

Conversely, suppose every nonempty directed cycle has positive
weight. Each simple directed cycle has at most $r$ edges, and
its weight is a positive multiple of $1/D$, hence at least $1/D$.
After subtracting $\varepsilon$ from each edge, its weight is at
least
\[
 \frac1D-r\varepsilon=\frac1{D(r+1)}>0.
\]
Every closed directed walk decomposes into simple directed cycles,
so the modified graph has no negative-weight cycle. If an original
cycle has nonpositive weight, subtracting $\varepsilon$ from each
of its edges makes it negative. This proves that positivity of all
original cycle weights is equivalent to condition~(ii).

When (ii) holds, add the source as stated. A shortest path to any
vertex exists: removing a repeated-vertex segment cannot increase
the weight, because there is no negative cycle, so it is enough to
consider the finitely many simple paths. The shortest-path distances
satisfy
\[
 x_b\le x_a+w_e-\varepsilon
 \quad\text{for every edge }a\to b.
\]
They are rational and satisfy all the strict inequalities of
\eqref{eq:regularity-strict-system}. The resulting point is in
$\overline\Delta$; every root hyperplane through it contains
$\mathbb V_{w'}$, by the same-cycle/different-cycle distinction above.
It is therefore regular. If $r=1$, there are no edges, condition~(ii)
is vacuous, and the source gives $x_c=0$; the argument remains valid.
\end{proof}

\begin{example}\label{ex:regularity-graph}
For $m=3$ and $w'=(1\,3)$, take $v^0=0$. The two cycles are
$A=\{1,3\}$ and $B=\{2\}$. The pairs $(1,2)$ and $(2,3)$ require
both $0<x_A-x_B<1$ and $0<x_B-x_A<1$, which is impossible.
Correspondingly, the graph has a zero-weight directed two-cycle,
which becomes negative after modifying the weights. For
$w'=(1\,2)$, the point represented by $(0,0,-1/3)$ is regular in
$\mathbb V_{w'}\cap\overline\Delta$: the sole wall through it is
$v_1=v_2$, and that wall contains $\mathbb V_{w'}$.
\end{example}

\begin{proposition}[Finite construction of a standard triple]\label{prop:effective-standard-triple}
Given an affine-permutation window for $w\in\tW$, one can compute,
using integer and rational arithmetic, a finite simple-conjugation
chain $w\to w_1$, a rational regular point of
$\mathbb V_{w_1}\cap\overline\Delta$, and a standard triple
$(K,x,u)$ with $w_1=ux$ satisfying Proposition~\ref{prop:HN3app}.
If $w$ is terminal, the returned chain preserves length and
satisfies Corollary~\ref{cor:terminal-standard}.
\end{proposition}
\begin{proof}
Put $k=\val(w)$ and $L=\ell(w)$. Starting with $w$, enumerate
reachable vertices by the following finite queue procedure: for
each vertex $z$ already reached and each $s\in S$, form $szs$;
retain it if $\ell(szs)\le\ell(z)$, and add it to the queue only
if it has not previously been reached. Record a predecessor and
the simple reflection for each newly added vertex. Every reached
vertex belongs to $\mathcal W_{k,\le L}$, which has at most
$\sum_{j=0}^{L}m^j$ elements by
Lemma~\ref{lem:alcove-reduced-expression}. Thus the queue procedure
terminates. All group operations and length comparisons use the
finite window and formula~\eqref{eq:affine-window-length}.

Apply Proposition~\ref{prop:regularity-test} to the reached vertices.
At least one passes the test: Proposition~\ref{prop:flow}, with
$\mathsf A=\Delta$, supplies a reachable
$w_1=x_{\mathsf A'}^{-1}wx_{\mathsf A'}$ such that
$\overline\Delta$ contains a regular point of $\mathbb V_{w_1}$.
Tracing the predecessors gives a chain $w\to w_1$.
The test in Proposition~\ref{prop:regularity-test} requires only
the finitely many simple directed cycles. After adjoining the source,
the distances are the minima of the weights of the finitely many
simple paths. Both computations use rational arithmetic.

Let $v_1$ be the rational regular point obtained in this way, and
set
\[
 K=\{s_i\in S:v_1\in H_i\}.
\]
Enumerate $W_K$ by multiplication by its generators; this procedure
terminates by (S4). For each $a\in W_K$ compute $\ell(aw_1)$,
and let $x$ be the unique minimum-length element of $W_Kw_1$,
whose uniqueness is (S5). Put $u=w_1x^{-1}\in W_K$.
The proof of Proposition~\ref{prop:HN3app}, with
$\mathsf Y=\mathbb V_{w_1}$ and this $v_1$, applies:
$W_K$ fixes $\mathsf Y$ pointwise, $xKx^{-1}=K$, $x$ is straight,
and the required valuation and Newton-vector identities hold.
If $w$ is terminal, the argument of
Corollary~\ref{cor:terminal-standard} shows that every step of the chain
preserves length.
\end{proof}

\subsection{Bruhat decomposition and Hecke normalization}\label{app:IM}

We derive the Bruhat decomposition and Iwahori--Matsumoto relations
\cite{IM} from the lattice-chain description of $\GL_m$.
The affine permutations and length function are those of
Section~\ref{app:HN3}. A rank-one calculation gives the convolution
coefficients and double-coset volumes, from which the finite
parabolic algebras and their normalizations follow.

Fix a non-archimedean local field $\K$ with ring of integers
$\fo$, uniformizer $\pi$, and residue field $\F_q$.
For $m=1$, the cosets of $I=\fo^\times$ in $\K^\times$ are
$\pi^jI$, $j\in\Z$, and the normalized characteristic functions
multiply according to addition of their exponents. Thus the Hecke
algebra is $\C[\Z]$, and the finite parabolic subgroup is trivial.
In the proofs below we assume $m\ge2$.

Put $G=\GL_m(\K)$ and choose a Haar measure on $G$. Extend the
standard basis $b_1,\ldots,b_m$ of $\K^m$ by
$b_{i+m}=\pi b_i$, and put
\[
 \cL_j=\sum_{i>j}\fo b_i,\qquad
 I=\{g\in G:g\cL_j=\cL_j\text{ for every }j\}.
\]
Here $\cL_0=\fo^m$, and for $0\le j\le m$ the image of
$\cL_j$ in $\cL_0/\pi\cL_0$ is spanned by
$\bar b_{j+1},\ldots,\bar b_m$. Consequently $I$ is the
inverse image of the \emph{lower triangular} Borel subgroup under
$\GL_m(\fo)\to\GL_m(\F_q)$. In particular, $I$ is compact and open.
An affine permutation acts by $w b_i=b_{w(i)}$. The algebra
$\cH=C_c(I\backslash G/I)$ consists of complex-valued functions
and has convolution defined by the chosen Haar measure; we put
$T_w=\vol(I)^{-1}\mathbf1_{IwI}$. Thus $T_1$ is its unit.

A \emph{complete chain} is a family $M=(M_j)_{j\in\Z}$ of
lattices in $\K^m$ with $M_j\supsetneq M_{j+1}$ and
$M_{j+m}=\pi M_j$ for all $j$. Each successive quotient is
killed by $\pi$, and the $m$ positive dimensions in one period
sum to $\dim_{\F_q}M_j/\pi M_j=m$. Hence
$\dim_{\F_q}M_j/M_{j+1}=1$ for every $j$.
The group $G$ acts on complete chains, and $I$ is the stabilizer
of $\cL=(\cL_j)$. For an affine permutation $w$ we have:
\[
w\cL_j=\sum_{a:\,w^{-1}(a)>j}\fo\,b_a,
\]
in particular $\tau\cL_j=\cL_{j+1}$, thus $\tau$ normalizes $I$, and $\tau s_i\tau^{-1}=s_{i+1}$ because $\tau s_i\tau^{-1}(a)=s_i(a-1)+1$ for $a\in\Z$.

The affine-permutation valuation agrees with the determinant
valuation in this realization. Indeed, write $w(i)=r_i+mk_i$
with $1\le r_i\le m$. Then $wb_i=\pi^{k_i}b_{r_i}$ and
$\det w=\pm\pi^{\sum_i k_i}$. Since the $r_i$ are a
permutation of $1,\ldots,m$, we obtain
\[
 \val(\det w)=\sum_i k_i
 =\frac1m\sum_{i=1}^m(w(i)-i)=\val(w).
\]

\subsubsection{Lengths}
\begin{lemma}\label{lem:IMlength}
Let $w\in\tW$ and $0\le i<m$.
\begin{enumerate}[label=\textup{(\roman*)}]
\item The number $\ell(w)$ is finite, and we have $\ell(w^{-1})=\ell(w)$ and $\ell(\tau w)=\ell(w\tau)=\ell(w)$.
\item We have $\ell(s_iw)=\ell(w)+1$ if $w^{-1}(i)<w^{-1}(i+1)$ and $\ell(s_iw)=\ell(w)-1$ otherwise, while $\ell(ws_i)=\ell(w)+1$ if $w(i)<w(i+1)$ and $\ell(ws_i)=\ell(w)-1$ otherwise.
\item We have $\ell(w)=0$ if and only if $w\in\langle\tau\rangle$, and if $\ell(w)>0$, then $\ell(s_jw)<\ell(w)$ for some $j$. Hence $w=s_{i_1}\cdots s_{i_r}\tau^k$ with $r=\ell(w)$ and $k=\val(w)$. Moreover, $\ell(w)$ is the smallest number $r$ for which such an expression exists, and in particular $W_a=\ker\val$ and $\tW=W_a\rtimes\langle\tau\rangle$.
\item If $\ell(s_iw)>\ell(w)$, $\ell(ws_i)>\ell(w)$ and $\ell(s_iws_i)=\ell(w)$, then $s_iws_i=w$.
\end{enumerate}
\end{lemma}
\begin{proof}
Let us denote by $\mathrm{Inv}(w)$ the set of inversions of $w$:
\[
\mathrm{Inv}(w)=\{(a,a')\in\Z^2:\ a<a',\ w(a)>w(a')\}.
\]
This set is stable under the translation $(a,a')\mapsto(a+m,a'+m)$, and $\ell(w)$ is the number of its orbits under this translation.

(i) By Lemma~\ref{lem:standard}(S2), $\ell(w)$ counts the
hyperplanes separating $\Delta$ and $w\Delta$. This number is
finite by local finiteness of the arrangement. Applying $w$ to
the pair $(\Delta,w^{-1}\Delta)$ gives
$\ell(w^{-1})=\ell(w)$. Since $\tau\Delta=\Delta$ and
$\tau$ permutes the reflecting hyperplanes, the same description
gives $\ell(\tau w)=\ell(w\tau)=\ell(w)$.

(ii) Let $s=s_i$, and note that for integers $x<y$ we have $s(x)>s(y)$ if and only if $(x,y)=(i+mk,i+1+mk)$ for some $k$. Hence a pair $(a,a')$ with $a<a'$ lies in exactly one of the sets $\mathrm{Inv}(w)$ and $\mathrm{Inv}(sw)$ if and only if $\{w(a),w(a')\}=\{i+mk,i+1+mk\}$ for some $k$, and these pairs form one orbit under the translation. If $A=w^{-1}(i)$ is smaller than $B=w^{-1}(i+1)$, then the pair $(A,B)$ lies in $\mathrm{Inv}(sw)$ and not in $\mathrm{Inv}(w)$, and if $A>B$, then the pair $(B,A)$ lies in $\mathrm{Inv}(w)$ and not in $\mathrm{Inv}(sw)$. This proves the first statement, and the second one follows from the first one, applied to $w^{-1}$, and from (i), because $\ell(ws)=\ell(sw^{-1})$.

(iii) By Lemma~\ref{lem:standard}(S1),(S2), the elements of
length zero are exactly the powers of $\tau$, and
$W_a=\ker\val$ with $\tW=W_a\rtimes\langle\tau\rangle$.
If $w$ has no left descent, (ii) gives
$w^{-1}(j)<w^{-1}(j+1)$ for every $j\in\Z$.
An increasing bijection of $\Z$ is a translation, so in this
case $w\in\langle\tau\rangle$ and $\ell(w)=0$.
Thus every element of positive length has a left descent.
The asserted expression follows from
Lemma~\ref{lem:alcove-reduced-expression}. Conversely, if
$w=s_{j_1}\cdots s_{j_t}\tau^k$, repeated application of
(i) and (ii) gives $\ell(w)\le t$, which proves minimality.

(iv) Let $s=s_i$, $A=w^{-1}(i)$ and $B=w^{-1}(i+1)$. By (ii) the inequality $\ell(sw)>\ell(w)$ means that $A<B$, and the assumptions give $\ell(ws)=\ell(w)+1$ and $\ell(s\cdot ws)=\ell(ws)-1$, so (ii), applied to $ws$, gives $(ws)^{-1}(i)>(ws)^{-1}(i+1)$, that is, $s(A)>s(B)$. As in the proof of (ii), this is possible only if $(A,B)=(i+mk,i+1+mk)$ for some $k$, and then $w(a)=a-mk$ for all $a$ which are congruent to $i$ or to $i+1$ modulo $m$. For these $a$ we get $sw(a)=ws(a)$, while the other residue classes are permuted by $w$ and fixed pointwise by $s$. Thus $sw=ws$.
\end{proof}

\subsubsection{The Bruhat decomposition}
\begin{lemma}\label{lem:IMbruhat}
\begin{enumerate}[label=\textup{(\roman*)}]
\item If $v_a\in\cL_{a-1}\smallsetminus\cL_a$ for $1\le a\le m$, then the linear map $h$ with $hb_a=v_a$ for $1\le a\le m$ lies in $I$.
\item For every complete chain $M$ there are $h\in I$ and $w\in\tW$ with $M=hw\cL$.
\item For $g\in G$ and $i,j\in\Z$ let $\mathsf d_g(i,j)$ be the length of the $\fo$-module $\cL_i/(\cL_i\cap g\cL_j)$. Then $\mathsf d_{hgh'}=\mathsf d_g$ for $h,h'\in I$, and for $w\in\tW$ we have $\mathsf d_w(i-1,j)-\mathsf d_w(i,j)=1$ if $w^{-1}(i)\le j$ and $\mathsf d_w(i-1,j)-\mathsf d_w(i,j)=0$ otherwise.
\item The group $G$ is the disjoint union of the double cosets of the elements of $\tW$: $G=\bigsqcup_{w\in\tW}IwI$.
\end{enumerate}
\end{lemma}
\begin{proof}
(i) If we put $v_{a+mk}=\pi^kv_a$, then $v_a\in\cL_{a-1}\smallsetminus\cL_a$ and $hb_a=v_a$ for all $a\in\Z$, thus $h\cL_i\subseteq\cL_i$ for all $i$. The map induced by $h$ on $\cL_i/\cL_{i+m}$ preserves the filtration by the subspaces $\cL_a/\cL_{i+m}$, $i\le a\le i+m$, and it is nonzero on each of the one-dimensional quotients $\cL_{a-1}/\cL_a$, because it maps the class of $b_a$ to the class of $v_a$, so it is bijective. Hence $h\cL_i+\pi\cL_i=\cL_i$, and Nakayama's lemma gives $h\cL_i=\cL_i$ for all $i$, thus $h\in G$ and $h\in I$.

(ii) For $a\in\Z$ we have $M_j\supseteq\cL_{a-1}$ for $j\ll0$ and $M_j\subseteq\cL_a$ for $j\gg0$, thus there is a largest integer $\beta(a)$ with $\cL_{a-1}\cap M_{\beta(a)}\not\subseteq\cL_a$. We have $\cL_{a-1}\cap M_j\not\subseteq\cL_a$ if and only if $j\le\beta(a)$, and moreover $\beta(a+m)=\beta(a)+m$. For $1\le a\le m$ we choose $v_a\in\cL_{a-1}\cap M_{\beta(a)}$ with $v_a\notin\cL_a$, we put $v_{a+mk}=\pi^kv_a$, and we denote by $h\in I$ the map of (i). For $j\in\Z$ let us consider the module spanned by the vectors $v_a$ with $\beta(a)\ge j$:
\[
N_j=\sum_{a:\,\beta(a)\ge j}\fo\,v_a .
\]
In every residue class $r$ modulo $m$ the integers $a$ with $\beta(a)\ge j$ are the integers $a_r(j)+mk$, $k\ge0$, for some $a_r(j)$. Hence $N_j=\bigoplus_r\fo\,v_{a_r(j)}$, and this is a lattice because $h$ is invertible, while $N_j\subseteq M_j$ by the choice of the vectors $v_a$.

We claim that $N_j=M_j$. Suppose to the contrary that $\xi\in M_j\smallsetminus N_j$, and for a nonzero vector $\zeta$ let $\operatorname{ht}(\zeta)$ be the integer with $\zeta\in\cL_{\operatorname{ht}(\zeta)-1}\smallsetminus\cL_{\operatorname{ht}(\zeta)}$. The lattice $N_j$ contains $\cL_{c_0}$ for some integer $c_0$, so the coset $\xi+N_j$ does not meet $\cL_{c_0}$. Thus $\operatorname{ht}\le c_0$ on this coset, and we can choose $\zeta\in\xi+N_j$ such that $a_0=\operatorname{ht}(\zeta)$ is maximal. We have $\cL_{a_0-1}=\fo v_{a_0}+\cL_{a_0}$, so there is $\epsilon\in\fo$ with $\zeta-\epsilon v_{a_0}\in\cL_{a_0}$. If $\beta(a_0)\ge j$, then $v_{a_0}\in N_j$, so $\zeta-\epsilon v_{a_0}$ lies in $\xi+N_j$, it is nonzero, and $\operatorname{ht}(\zeta-\epsilon v_{a_0})>a_0$, contrary to the choice of $\zeta$. Hence $\beta(a_0)<j$. On the other hand, $\zeta$ lies in $\cL_{a_0-1}\cap M_j$ and not in $\cL_{a_0}$, so $j\le\beta(a_0)$, and this contradiction proves the claim.

If $\beta(a_r(j))>j$, then $a_r(j+1)=a_r(j)$, and if
$\beta(a_r(j))=j$, then $a_r(j+1)=a_r(j)+m$.
Hence $\dim N_j/N_{j+1}$ is the number of integers $a$ with
$\beta(a)=j$. Since $N_j/N_{j+1}=M_j/M_{j+1}$ is
one-dimensional, $\beta$ is a bijection of $\Z$; together with
$\beta(a+m)=\beta(a)+m$, this says that $\beta$ is an affine
permutation. To identify the required $w$, observe that in the
chain $w\cL$ the last layer containing the vector $b_a$ has
index $w^{-1}(a)-1$. Thus we must set
$w^{-1}(a)=\beta(a)+1$, or equivalently
$w^{-1}=\tau\beta$ and $w=\beta^{-1}\tau^{-1}$.
Then $w^{-1}(a)>j$ if and only if $\beta(a)\ge j$, so
\[
 w\cL_j=\sum_{a:\,\beta(a)\ge j}\fo\,b_a,
 \qquad hw\cL_j=N_j=M_j.
\]

(iii) The first statement holds because $h'\cL_j=\cL_j$ and $h^{-1}\cL_i=\cL_i$. For the second one, let us observe that the lattices $\cL_i$ and $w\cL_j$ are direct sums of their intersections with the lines $\K b_r$, $1\le r\le m$. These intersections are $\fo b_{a_r}$ and $\fo b_{a'_r}$, where $a_r$ is the smallest integer $a\equiv r$ with $a>i$, and $a'_r$ is the smallest integer $a\equiv r$ with $w^{-1}(a)>j$. Here we use that $w^{-1}(a+m)=w^{-1}(a)+m$, and we obtain:
\[
\cL_i/(\cL_i\cap w\cL_j)\cong\bigoplus_{r=1}^m\fo\,b_{a_r}\big/\fo\,b_{\max(a_r,a'_r)}\cong\bigoplus_{r=1}^m\fo/\pi^{(\max(a_r,a'_r)-a_r)/m}\fo .
\]
The length of this module is the number of integers $a>i$ with $w^{-1}(a)\le j$, and the formula follows.

(iv) For $g\in G$ the family $g\cL$ is a complete chain, and by (ii) we have $g\cL=hw\cL$ with $h\in I$ and $w\in\tW$, thus $w^{-1}h^{-1}g\in I$ and $g\in IwI$. If $IwI=Iw'I$, then $\mathsf d_w=\mathsf d_{w'}$ by (iii), and since, again by (iii), $w^{-1}(i)$ is the smallest integer $j$ with $\mathsf d_w(i-1,j)-\mathsf d_w(i,j)=1$, we get $w=w'$.
\end{proof}

\subsubsection{The rank-one computation}
Fix $0\le i<m$ and put $s=s_i$. For $\mu\in\fo$ we denote by $\mathsf u_\mu$ and $\bar{\mathsf u}_\mu$ the linear maps with
\[
\mathsf u_\mu b_a=b_a+\mu\,b_{a+1}\ \ (a\equiv i),\qquad\bar{\mathsf u}_\mu b_a=b_a+\mu\,b_{a-1}\ \ (a\equiv i+1),
\]
where the congruences are modulo $m$, and which fix all other $b_a$. We also fix a set $\mathcal R\subset\fo$ of representatives of $\fo/\pi\fo$ with $0\in\mathcal R$.

\begin{lemma}\label{lem:IMrank1}
\begin{enumerate}[label=\textup{(\roman*)}]
\item We have $\mathsf u_\mu\in I$ and $\mathsf u_\mu\mathsf u_{\mu'}=\mathsf u_{\mu+\mu'}$ for $\mu,\mu'\in\fo$, and $IsI=\bigsqcup_{\mu\in\mathcal R}\mathsf u_\mu sI$, in particular $|IsI/I|=q$.
\item Let $z\in\tW$ and $\mu\in\fo$. If $z^{-1}(i)<z^{-1}(i+1)$, then $s\,\mathsf u_\mu z\in szI$, and if $z^{-1}(i)>z^{-1}(i+1)$, then $s\,\mathsf u_\mu z\in IszI$ for $\mu\in\pi\fo$ and $s\,\mathsf u_\mu z\in IzI$ for $\mu\in\fo^\times$.
\end{enumerate}
\end{lemma}
\begin{proof}
(i) We have $\mathsf u_\mu\in I$ by Lemma~\ref{lem:IMbruhat}(i), and since $m\ge2$, the map $\mathsf u_\mu$ fixes $b_{a+1}$ for $a\equiv i$, so $\mathsf u_\mu\mathsf u_{\mu'}=\mathsf u_{\mu+\mu'}$. Let us call a complete chain $M$ an $i$-neighbour of $\cL$ if $M_j=\cL_j$ for all $j\not\equiv i$ modulo $m$ and $M_i\ne\cL_i$. For $j\not\equiv i$ the set of integers $a>j$ is stable under $s$, and we get $s\cL_j=\cL_j$, while $s\cL_i=\fo b_i+\cL_{i+1}\ne\cL_i$. Thus $s\cL$ is an $i$-neighbour of $\cL$, and so is $hs\cL$ for every $h\in I$. Conversely, if $M$ is an $i$-neighbour of $\cL$, then $\cL_{i-1}\supsetneq M_i\supsetneq\cL_{i+1}$. Since $m\ge2$, we have $\pi\cL_{i-1}=\cL_{i-1+m}\subseteq\cL_{i+1}$, so $\cL_{i-1}/\cL_{i+1}$ is a plane over $\F_q$, and $M_i/\cL_{i+1}$ is a line in it. This line is different from the line $\cL_i/\cL_{i+1}$ spanned by the class of $b_{i+1}$, thus for a unique $\mu\in\mathcal R$ we have:
\[
M_i=\fo(b_i+\mu b_{i+1})+\cL_{i+1}=\mathsf u_\mu s\cL_i .
\]
For $j\not\equiv i$ we have $\mathsf u_\mu s\cL_j=\cL_j=M_j$, and we get $M=\mathsf u_\mu s\cL$. Therefore the chains $hs\cL$, $h\in I$, are the $q$ chains $\mathsf u_\mu s\cL$, $\mu\in\mathcal R$, and (i) follows, because $I$ is the stabilizer of $\cL$.

(ii) Let $A=z^{-1}(i)$ and $B=z^{-1}(i+1)$. The map $z^{-1}\mathsf u_\mu z$ sends $b_a$ to $b_a+\mu b_{a+B-A}$ for $a\equiv A$ and fixes the other $b_a$, so if $A<B$, it lies in $I$ by Lemma~\ref{lem:IMbruhat}(i), and we get:
\[
s\,\mathsf u_\mu z=sz\,(z^{-1}\mathsf u_\mu z)\in szI.
\]

Now let $A>B$. Then the map $z^{-1}\bar{\mathsf u}_\mu z$ sends $b_a$ to $b_a+\mu b_{a+A-B}$ for $a\equiv B$ and fixes the other $b_a$, so it lies in $I$ for every $\mu\in\fo$. For $\mu\in\pi\fo$ we have $s\,\mathsf u_\mu s=\bar{\mathsf u}_\mu$ and $\bar{\mathsf u}_\mu b_a=b_a+(\mu/\pi)\,b_{a-1+m}$ for $a\equiv i+1$, and since $m\ge2$, this shows that $\bar{\mathsf u}_\mu\in I$ by Lemma~\ref{lem:IMbruhat}(i), thus $s\,\mathsf u_\mu z=\bar{\mathsf u}_\mu sz\in IszI$. Next, for $\mu\in\fo^\times$ we have:
\[
s\,\mathsf u_\mu=\mathsf u_{\mu^{-1}}\,\mathsf d\,\bar{\mathsf u}_{\mu^{-1}} ,
\]
where $\mathsf d$ is the linear map with $\mathsf d\,b_a=\mu b_a$ for $a\equiv i$ and $\mathsf d\,b_a=-\mu^{-1}b_a$ for $a\equiv i+1$ which fixes the other $b_a$.
Indeed, for $a\equiv i$ both sides map $b_a$ to $\mu b_a+b_{a+1}$. For $a\equiv i+1$ the right side maps $b_a$ first to $b_a+\mu^{-1}b_{a-1}$, then to $-\mu^{-1}b_a+b_{a-1}$, and then to $b_{a-1}$, which is $s\,\mathsf u_\mu b_a$. The map $z^{-1}\mathsf dz$ multiplies every $b_a$ by a unit, so it lies in $I$, and we get:
\[
s\,\mathsf u_\mu z=\mathsf u_{\mu^{-1}}\,z\,(z^{-1}\mathsf dz)(z^{-1}\bar{\mathsf u}_{\mu^{-1}}z)\in IzI .\qedhere
\]
\end{proof}

\subsubsection{The relations}
Every double coset $IwI$ is compact and open, and a compact set meets only finitely many of them, hence by Lemma~\ref{lem:IMbruhat}(iv) the functions $T_w=\vol(I)^{-1}\mathbf 1_{IwI}$, $w\in\tW$, form a basis of $\cH$. Since $I\subseteq\GL_m(\fo)$, the function $\val\det$ takes the value $\val(w)$ on $IwI$, and $\cH$ is graded by $\val$ because $\val\det$ is a homomorphism. For $x,y\in\tW$ the function $T_x*T_y$ is invariant under $I$ on both sides and supported on the compact set $IxI\cdot IyI$. So it lies in $\cH$, and we have $T_x*T_y=\sum_zc_zT_z$ with $c_z=\vol(I)\,(T_x*T_y)(z)$.

\Needspace{6\baselineskip}
\begin{proposition}\label{prop:IM}
Let $s$ be a simple reflection and $x,y\in\tW$.
\begin{enumerate}[label=\textup{(\roman*)}]
\item We have $T_sT_y=T_{sy}$ if $\ell(sy)>\ell(y)$ and $T_sT_y=(q-1)T_y+q\,T_{sy}$ if $\ell(sy)<\ell(y)$, and the same holds for $T_yT_s$, with $ys$ in place of $sy$.
\item For $k\in\Z$ we have $T_{\tau^k}T_y=T_{\tau^ky}$ and $T_yT_{\tau^k}=T_{y\tau^k}$, and if $\ell(xy)=\ell(x)+\ell(y)$, then $T_xT_y=T_{xy}$. In particular, $T_x=T_{s_{i_1}}\cdots T_{s_{i_r}}T_{\tau^k}$ whenever $x=s_{i_1}\cdots s_{i_r}\tau^k$ with $r=\ell(x)$, and every $T_x$ is invertible.
\item We have $IsI\cdot IyI=IsyI$ if $\ell(sy)>\ell(y)$ and $IsI\cdot IyI=IyI\sqcup IsyI$ if $\ell(sy)<\ell(y)$, and the same holds for $IyI\cdot IsI$, with $ys$ in place of $sy$. Moreover, if $\ell(xy)=\ell(x)+\ell(y)$, then we have $IxI\cdot IyI=IxyI$.
\item For all $w\in\tW$ the double coset $IwI$ is the union of $q^{\ell(w)}$ cosets of $I$, that is, $|IwI/I|=q^{\ell(w)}$.
\end{enumerate}
\end{proposition}
\begin{proof}
(i) Let $s=s_i$. Lemma~\ref{lem:IMrank1}(i) shows that the set $IsI$ is the disjoint union of the cosets $\mathsf u_\mu sI$, $\mu\in\mathcal R$, and since $h^{-1}g\in I\,s\,\mathsf u_{-\mu}g$ for $h\in\mathsf u_\mu sI$, we get:
\[
(T_s*T_y)(g)=\vol(I)^{-2}\int_{IsI}\mathbf 1_{IyI}(h^{-1}g)\,dh=\vol(I)^{-1}\,\#\{\mu\in\mathcal R:\ s\,\mathsf u_{-\mu}g\in IyI\}.
\]
Thus $c_z$ is the number of $\mu\in\mathcal R$ with $s\,\mathsf u_{-\mu}z\in IyI$, and to compute it we use Lemma~\ref{lem:IMrank1}(ii), Lemma~\ref{lem:IMlength}(ii) and the fact that distinct elements of $\tW$ have disjoint double cosets. If $\ell(sz)>\ell(z)$, then $s\,\mathsf u_{-\mu}z\in IszI$ for all $\mu$, so $c_z=q$ if $y=sz$, and $c_z=0$ otherwise. If $\ell(sz)<\ell(z)$, then $s\,\mathsf u_{-\mu}z$ lies in $IszI$ for $\mu=0$ and in $IzI$ for the $q-1$ other elements $\mu$ of $\mathcal R$, so $c_z=1$ if $y=sz$, $c_z=q-1$ if $y=z$, and $c_z=0$ otherwise. In particular, $c_z=0$ unless $z=y$ or $z=sy$, and we get $c_y=0$ and $c_{sy}=1$ if $\ell(sy)>\ell(y)$, and $c_y=q-1$ and $c_{sy}=q$ if $\ell(sy)<\ell(y)$. This proves the first statement.

The map $f\mapsto f^\vee$, $f^\vee(g)=f(g^{-1})$, satisfies $(f_1*f_2)^\vee=f_2^\vee*f_1^\vee$, as follows from the substitution $h\mapsto gh$ in the integral, and it maps $T_w$ to $T_{w^{-1}}$. Applying it to the formula for $T_sT_{y^{-1}}$ and using Lemma~\ref{lem:IMlength}(i), we obtain the formula for $T_yT_s$.

(ii) Since $\tau^k$ normalizes $I$, we have $I\tau^kI=\tau^kI$ and $(T_{\tau^k}*T_y)(g)=\vol(I)^{-1}\mathbf 1_{IyI}(\tau^{-k}g)$, and the right side is $T_{\tau^ky}(g)$ because $\tau^kIyI=I\tau^kyI$. The formula for $T_yT_{\tau^k}$ follows in the same way. Now let $\ell(xy)=\ell(x)+\ell(y)$, let $x=s_{i_1}\cdots s_{i_r}\tau^k$ with $r=\ell(x)$ be an expression given by Lemma~\ref{lem:IMlength}(iii), and put $x_t=s_{i_t}\cdots s_{i_r}\tau^k$ for $1\le t\le r+1$. By Lemma~\ref{lem:IMlength}(i),(ii) we have $\ell(x_{r+1}y)=\ell(y)$ and $\ell(x_ty)\le\ell(x_{t+1}y)+1$ for $1\le t\le r$, and together with the first equality these $r$ inequalities give $\ell(x_1y)\le r+\ell(y)$. As $\ell(x_1y)=r+\ell(y)$, all of them are equalities, that is, $\ell(x_ty)=\ell(x_{t+1}y)+1$ for $1\le t\le r$. In the same way $\ell(x_t)=\ell(x_{t+1})+1$, and by (i) and the first part of (ii) we get $T_x=T_{s_{i_1}}\cdots T_{s_{i_r}}T_{\tau^k}$ and:
\[
T_xT_y=T_{s_{i_1}}\cdots T_{s_{i_r}}T_{\tau^ky}=T_{xy} .
\]
The formula for $T_x$ holds for every such expression of $x$. By (i) we have the quadratic relation:
\[
T_s\big(T_s-(q-1)\big)=q,
\]
so $T_s$ is invertible, and since $T_{\tau^k}T_{\tau^{-k}}=T_1$, every $T_x$ is invertible.

(iii) Let us observe that:
\[
\vol(I)^2\,(T_x*T_y)(g)=\vol\{h\in IxI:\ h^{-1}g\in IyI\},
\]
where the set on the right is open, and it is nonempty if and only if $g\in IxI\cdot IyI$. Since a Haar measure is positive on nonempty open sets, $IxI\cdot IyI$ is the union of the double cosets $IzI$ with $c_z\ne0$, so (iii) follows from (i) and (ii).

(iv) We have $\int_GT_w=|IwI/I|$ and $\int_G(f_1*f_2)=\int_Gf_1\int_Gf_2$, and $|IsI/I|=q$ by Lemma~\ref{lem:IMrank1}(i), while $|I\tau^kI/I|=1$. Thus (iv) follows from the formula for $T_x$ in (ii).
\end{proof}

\subsubsection{Parabolic subgroups}
Let $\Sigma\subseteq\Z/m$ be nonempty, and put
$K=K_\Sigma=\{s_i:i\notin\Sigma\}$ and $W_K=\langle K\rangle$.
Every proper subset of $\{s_0,\ldots,s_{m-1}\}$ arises in this way.
Write $\widehat\Sigma\subseteq\Z$ for the set of integers
whose residues lie in $\Sigma$. A \emph{window} is an interval
$\{p+1,\ldots,p+c\}$ for consecutive elements $p<p+c$ of
$\widehat\Sigma$. These windows partition $\Z$, have length
at most $m$, and are permuted by translation by $m$.
Consecutive integers $a,a+1$ lie in different windows precisely
when $a\in\widehat\Sigma$.

The residue set of a window is called a \emph{block} of $\Sigma$;
the window is said to lie over that block. If
$p_1\in\widehat\Sigma$, the windows contained in
$\{p_1+1,\ldots,p_1+m\}$ represent all windows modulo
translation by $m$, and their residue sets are pairwise disjoint.
Thus two windows over the same block differ by a multiple of $m$.
Choose a set $\mathcal B$ containing one window over each block.
For a window $B$, let $\mathrm{Sym}(B)$ be its permutation group
and let $\mathrm{inv}(\sigma)$ count the inversions of
$\sigma\in\mathrm{Sym}(B)$ in the inherited order.
For $B=\{p+1,\ldots,p+c\}$ define
\[
 S_B=\langle s_{p+1},\ldots,s_{p+c-1}\rangle,
\]
where indices are read modulo $m$. This subgroup depends only
on the block below $B$.

\begin{lemma}\label{lem:parabolic}
\begin{enumerate}[label=\textup{(\roman*)}]
\item An element $y\in\tW$ lies in $W_K$ if and only if $y$ maps every window onto itself, and the restriction to the windows in $\mathcal B$ is an isomorphism $W_K\to\prod_{B\in\mathcal B}\mathrm{Sym}(B)$. This isomorphism maps $S_B$ onto $\mathrm{Sym}(B)$, and it maps $s_a$, $p<a<p+c$, to the transposition of $a$ and $a+1$ in $\mathrm{Sym}(B)$. In particular, $W_K=\prod_{B\in\mathcal B}S_B$, and $W_K$ is finite.
\item For $y\in W_K$ we have $\ell(y)=\sum_{B\in\mathcal B}\mathrm{inv}(y|_B)$, and $y$ is a product of $\ell(y)$ elements of $K$. If $y\in S_B$, then $y$ is a product of $\ell(y)$ of the generators of $S_B$.
\item The set $P_K=\bigcup_{y\in W_K}IyI$ is a compact open subgroup of $G$, and its index is given by:
\[
[P_K:I]=\sum_{y\in W_K}q^{\ell(y)}=\prod_{B\in\mathcal B}\,[\,|B|\,]_q!\,,\qquad\text{where}\quad[n]_q!=\prod_{i=1}^n\frac{q^i-1}{q-1}.
\]
\item The span $\cH_K$ of the $T_y$, $y\in W_K$, is a subalgebra of $\cH$, and so is the span $\cH_B$ of the $T_y$, $y\in S_B$, for every $B\in\mathcal B$. The subalgebras $\cH_B$ commute with each other, and $T_y=\prod_{B\in\mathcal B}T_{y_B}$ for $y=\prod_By_B$ with $y_B\in S_B$.
\end{enumerate}
\end{lemma}
\begin{proof}
(i) If $s_i\in K$, then the integers $i+mj$ do not lie in $\widehat\Sigma$, so $i+mj$ and $i+1+mj$ lie in the same window for every $j$, and $s_i$ maps every window onto itself. Thus every element of $W_K$ maps every window onto itself. Let $W_\Sigma$ be the group of affine permutations which map every window onto itself. Since every window is a translate of a unique window in $\mathcal B$ and $y(a+m)=y(a)+m$ for $y\in\tW$, an element of $W_\Sigma$ is determined by its restrictions to the windows in $\mathcal B$. Conversely, every family $(\sigma_B)\in\prod_B\mathrm{Sym}(B)$ is the family of restrictions of the affine permutation $y$ with $y(a+mk)=\sigma_B(a)+mk$ for $a\in B\in\mathcal B$ and $k\in\Z$, so the restriction is an isomorphism $W_\Sigma\to\prod_B\mathrm{Sym}(B)$.

Now let $B=\{p+1,\dots,p+c\}\in\mathcal B$ and $p<a<p+c$. Then $a\notin\widehat\Sigma$, so $s_a\in K$, and since the residues of the windows in $\mathcal B$ are pairwise disjoint, $s_a$ restricts to the transposition of $a$ and $a+1$ on $B$ and to the identity on the other windows of $\mathcal B$. These transpositions generate $\prod_B\mathrm{Sym}(B)$. Indeed, a permutation $\sigma\ne1$ of $B$ has a descent $\sigma(a)>\sigma(a+1)$, and composing $\sigma$ with the transposition of $a$ and $a+1$ lowers $\mathrm{inv}(\sigma)$ by one. Thus $W_K=W_\Sigma$, and the remaining statements follow.

(ii) Let $y\in W_K$ and $(a,a')\in\mathrm{Inv}(y)$, with the notation of the proof of Lemma~\ref{lem:IMlength}, and suppose that $a$ and $a'$ lie in different windows. Since the windows are intervals and $a<a'$, every element of the window of $a$ is smaller than every element of the window of $a'$. This applies to $y(a)$ and $y(a')$, which lie in these windows, and it contradicts $y(a)>y(a')$, so $a$ and $a'$ lie in the same window. Every translation orbit of such pairs contains exactly one pair in a window of $\mathcal B$, and this proves the formula for $\ell(y)$.

If $\ell(y)=0$, then no $y|_B$ has inversions, so $y=1$ by (i). If $\ell(y)>0$, then by Lemma~\ref{lem:IMlength}(ii),(iii) there is $0\le j<m$ with $y^{-1}(j)>y^{-1}(j+1)$. If $j\in\widehat\Sigma$, then the window of $j$ lies to the left of the window of $j+1$, and $y^{-1}$ maps both windows onto themselves, thus $y^{-1}(j)<y^{-1}(j+1)$, which is impossible. Hence $s_j\in K$. Then $s_jy\in W_K$ and $\ell(s_jy)=\ell(y)-1$, and by induction on $\ell(y)$ the element $y$ is a product of $\ell(y)$ elements of $K$. If $y\in S_B$, then $y$ is the identity on the windows over the other blocks, so $j$ and $j+1$ lie in a window over $B$, and $s_j$ is one of the generators of $S_B$. Thus in this case $y$ is a product of $\ell(y)$ generators of $S_B$.

(iii) The set $P_K$ contains $I$, it is open, and it is compact because $W_K$ is finite. For $y\in W_K$ and $s\in K$, Proposition~\ref{prop:IM}(iii) gives $IsI\cdot IyI\subseteq IyI\cup IsyI\subseteq P_K$, thus $IsI\cdot P_K\subseteq P_K$. If $x\in W_K$ and $x=s_{i_1}\cdots s_{i_r}$ with $s_{i_t}\in K$ and $r=\ell(x)$ is an expression given by (ii), then $\ell(s_{i_t}\cdots s_{i_r})=r-t+1$ for all $t$, and $IxI=Is_{i_1}I\cdots Is_{i_r}I$ by Proposition~\ref{prop:IM}(iii), so $IxI\cdot P_K\subseteq P_K$. It follows that $P_K\cdot P_K\subseteq P_K$, and since $(IyI)^{-1}=Iy^{-1}I$ and $y^{-1}\in W_K$, the set $P_K$ is a group. Next, $P_K/I$ is the disjoint union of the sets $IyI/I$, and $|IyI/I|=q^{\ell(y)}$ by Proposition~\ref{prop:IM}(iv), so by (i) and (ii) we get:
\[
\sum_{y\in W_K}q^{\ell(y)}=\prod_{B\in\mathcal B}\ \sum_{\sigma\in\mathrm{Sym}(B)}q^{\mathrm{inv}(\sigma)} .
\]
Finally, we identify $\mathrm{Sym}(B)$ with the symmetric group of $\{1,\dots,|B|\}$ by the increasing bijection, which preserves inversions. Let $\sigma$ be a permutation of $\{1,\dots,n\}$, and for $1\le i\le n$ let $\kappa_i(\sigma)$ be the number of $j<i$ with $\sigma(j)>\sigma(i)$. Then $0\le\kappa_i(\sigma)\le i-1$, $\mathrm{inv}(\sigma)=\sum_i\kappa_i(\sigma)$ and $\sigma(n)=n-\kappa_n(\sigma)$, while the numbers $\kappa_i(\sigma)$ with $i<n$ are those of the permutation of $\{1,\dots,n-1\}$ obtained from $\sigma|_{\{1,\dots,n-1\}}$ by the increasing bijection of $\{1,\dots,n\}\smallsetminus\{\sigma(n)\}$ onto $\{1,\dots,n-1\}$. By induction on $n$, the map $\sigma\mapsto(\kappa_1(\sigma),\dots,\kappa_n(\sigma))$ is a bijection onto $\prod_{i=1}^n\{0,\dots,i-1\}$, and we obtain:
\[
\sum_\sigma q^{\mathrm{inv}(\sigma)}=\prod_{i=1}^n(1+q+\dots+q^{i-1})=[n]_q!.
\]

(iv) Let $x,y\in W_K$. If $y=s_{i_1}\cdots s_{i_r}$ with $s_{i_t}\in K$ and $r=\ell(y)$ is an expression given by (ii), then $T_y=T_{s_{i_1}}\cdots T_{s_{i_r}}$ by Proposition~\ref{prop:IM}(ii). For $z\in W_K$ and $s\in K$ the element $T_zT_s$ is a combination of $T_z$ and $T_{zs}$ by Proposition~\ref{prop:IM}(i), and $zs\in W_K$, thus $T_xT_y\in\cH_K$. By the last statement of (ii), the same argument applies to $\cH_B$. Next, if $y=\prod_By_B$ with $y_B\in S_B$, then the factors $y_B$ commute by (i), and by (ii) the length of a product of some of them is the sum of their lengths. Hence $T_y=\prod_BT_{y_B}$ by Proposition~\ref{prop:IM}(ii), and in particular, for $x\in S_B$ and $x'\in S_{B'}$ with $B\ne B'$ we have:
\[
T_xT_{x'}=T_{xx'}=T_{x'x}=T_{x'}T_x.\qedhere
\]
\end{proof}

The normalization formulas are
\begin{equation}\label{eq:IM-normalizations}
\begin{aligned}
 T_1&=\vol(I)^{-1}\mathbf1_I,
 & T_s^2&=(q-1)T_s+qT_1,\\
 \vol(IwI)&=q^{\ell(w)}\vol(I),
 & \vol(P_K)&=\vol(I)\prod_{B\in\mathcal B}[|B|]_q!.
\end{aligned}
\end{equation}
In particular, the characteristic function of $P_K$ normalized
to have integral one is
\[
 \frac{\mathbf1_{P_K}}{\vol(P_K)}
 =\frac{1}{\prod_{B\in\mathcal B}[|B|]_q!}
   \sum_{y\in W_K}T_y.
\]
These formulas hold for any choice of Haar measure.

\subsubsection{The generic finite Hecke algebra}
We now construct the Hecke algebra of $S_n$ over $\Z[X]$ from
the convolution algebras for $\GL_n(\F_q((\pi)))$.
The parameter $X$ is an indeterminate; specialization $X=q$
recovers the finite parabolic algebra above.
Write $s_1,\ldots,s_{n-1}$ for the simple transpositions of
$S_n$ and $\mathrm{inv}(w)$ for the inversion number, and put
$H_1=\Z[X]$.

\Needspace{9\baselineskip}
\begin{proposition}\label{prop:genericHecke}
Let $n\ge1$.
\begin{enumerate}[label=\textup{(\roman*)}]
\item There is a unique associative $\Z[X]$-algebra $H_n$ with unit $T_1$ which is a free $\Z[X]$-module with basis $(T_w)_{w\in S_n}$ and satisfies
\[
 T_wT_{s_i}=
 \begin{cases}
 T_{ws_i},&\mathrm{inv}(ws_i)>\mathrm{inv}(w),\\
 (X-1)T_w+X\,T_{ws_i},&\mathrm{inv}(ws_i)<\mathrm{inv}(w),
 \end{cases}
\]
for $w\in S_n$ and $1\le i<n$. The same rules hold for $T_{s_i}T_w$, with $s_iw$ in place of $ws_i$, and $T_xT_y=T_{xy}$ if $\mathrm{inv}(xy)=\mathrm{inv}(x)+\mathrm{inv}(y)$.
\item Let $\K$ be a local field with residue field $\F_q$, let $m\ge2$, let $\Sigma\subseteq\Z/m$ be nonempty, and let $B$ be a block of $\Sigma$ with $n$ elements. Let $\{p+1,\dots,p+n\}$ be a window over $B$, and let $\iota_B:S_n\to S_B$ be the isomorphism with $\iota_B(s_i)=s_{p+i}$. Then the linear map $H_n\otimes_{\Z[X]}\C\to\cH_B$ with $T_w\mapsto T_{\iota_B(w)}$ is an isomorphism of algebras, where $X$ acts on $\C$ by $q$.
\end{enumerate}
\end{proposition}
\begin{proof}
For $n=1$ we have $S_1=\{1\}$, $H_1=\Z[X]$ and $\cH_B=\C T_1$, and there is nothing to prove, so let $n\ge2$. First let $\K$, $m$, $\Sigma$, $B$ and $\iota_B$ be as in (ii). By Lemma~\ref{lem:parabolic}(i),(ii), the map $\iota_B$ is an isomorphism with $\ell(\iota_B(w))=\mathrm{inv}(w)$, and by Lemma~\ref{lem:IMbruhat}(iv), the elements $T_{\iota_B(w)}$, $w\in S_n$, form a basis of $\cH_B$. By Proposition~\ref{prop:IM}(i) we have:
\begin{equation}\label{eq:heckerule}
T_{\iota_B(w)}T_{\iota_B(s_i)}=\begin{cases}T_{\iota_B(ws_i)}&\text{if }\mathrm{inv}(ws_i)>\mathrm{inv}(w),\\(q-1)\,T_{\iota_B(w)}+q\,T_{\iota_B(ws_i)}&\text{if }\mathrm{inv}(ws_i)<\mathrm{inv}(w),\end{cases}
\end{equation}
and the same holds for left multiplication.

For every $y\in S_n$ we fix an expression $y=s_{i_1}\cdots s_{i_r}$ with $r=\mathrm{inv}(y)$, as provided by Lemma~\ref{lem:parabolic}(ii). Let $\mathsf M$ be the free $\Z[X]$-module with basis $(T_w)_{w\in S_n}$, and note that $\mathrm{inv}(ws_i)=\mathrm{inv}(w)\pm1$, because $ws_i$ is obtained from $w$ by exchanging its values at $i$ and $i+1$. For $1\le i<n$ let $\mathsf r_i$ be the $\Z[X]$-linear endomorphism of $\mathsf M$ with $\mathsf r_i(T_w)=T_{ws_i}$ if $\mathrm{inv}(ws_i)>\mathrm{inv}(w)$ and $\mathsf r_i(T_w)=(X-1)T_w+XT_{ws_i}$ if $\mathrm{inv}(ws_i)<\mathrm{inv}(w)$. We define a $\Z[X]$-bilinear product on $\mathsf M$ by the formula
\[
T_x\cdot T_y=\mathsf r_{i_r}\cdots\mathsf r_{i_1}(T_x),
\]
where we use the fixed expression of $y$.

For a prime power $q$ we take $\K=\F_q((\pi))$, $m=n$, $\Sigma=\{0\}$ and the block $B=\Z/n$ with the window $\{1,\dots,n\}$, and we denote by $\mathrm{ev}_q:\mathsf M\to\cH_B$ the map $\sum_wf_w(X)T_w\mapsto\sum_wf_w(q)T_{\iota_B(w)}$. By \eqref{eq:heckerule} we have $\mathrm{ev}_q(\mathsf r_i(f))=\mathrm{ev}_q(f)\,T_{\iota_B(s_i)}$ for $f\in\mathsf M$, and since $T_{\iota_B(y)}=T_{\iota_B(s_{i_1})}\cdots T_{\iota_B(s_{i_r})}$ by Proposition~\ref{prop:IM}(ii), we get:
\[
\mathrm{ev}_q(f\cdot g)=\mathrm{ev}_q(f)\,\mathrm{ev}_q(g)\qquad(f,g\in\mathsf M).
\]

For $x,y,z,u\in S_n$ the coefficient of $T_u$ in $(T_x\cdot T_y)\cdot T_z-T_x\cdot(T_y\cdot T_z)$ is a polynomial in $X$, and by the last display, its value at every prime power $q$ is a coefficient of $(T_{\iota_B(x)}T_{\iota_B(y)})T_{\iota_B(z)}-T_{\iota_B(x)}(T_{\iota_B(y)}T_{\iota_B(z)})=0$, thus the polynomial is zero, and the product on $\mathsf M$ is associative. In the same way we get $T_1\cdot T_y=T_y$, because $T_1$ is the unit of $\cH$, while $T_x\cdot T_1=T_x$ and the rules in (i) for right multiplication hold by definition. The rules for left multiplication and the formula $T_xT_y=T_{xy}$ hold after applying $\mathrm{ev}_q$ for every $q$, by Proposition~\ref{prop:IM}(i),(ii). Hence they hold in $\mathsf M$, and the existence in (i) follows. The rules determine $T_xT_y$ by induction on $\mathrm{inv}(y)$, because $T_y=T_{y'}T_{s_i}$ for $y=y's_i$ with $\mathrm{inv}(y')<\mathrm{inv}(y)$, and this proves the uniqueness.

(ii) The linear map $\psi$ of (ii) is bijective because the $T_{\iota_B(w)}$ form a basis of $\cH_B$, and by \eqref{eq:heckerule} we have $\psi(f\,T_{s_i})=\psi(f)\,\psi(T_{s_i})$ for every $f$ and every $i$. If $y=s_{i_1}\cdots s_{i_r}$ with $r=\mathrm{inv}(y)$, then $T_y=T_{s_{i_1}}\cdots T_{s_{i_r}}$ in $H_n$ by the rules of (i), and we get $\psi(fT_y)=\psi(f)\,\psi(T_{s_{i_1}})\cdots\psi(T_{s_{i_r}})$ for every $f$. For $f=1$ this gives:
\[
\psi(T_y)=\psi(T_{s_{i_1}})\cdots\psi(T_{s_{i_r}}),
\]
thus $\psi(fT_y)=\psi(f)\,\psi(T_y)$, and $\psi$ is multiplicative.
\end{proof}

The invertibility assertion in Proposition~\ref{prop:IM}(ii)
concerns the convolution algebra over $\C$, where $q$ is
invertible. For the generic algebra the corresponding formula
is
\[
 T_s^{-1}=X^{-1}\bigl(T_s-(X-1)T_1\bigr)
 \quad\text{in}\quad
 H_n\otimes_{\Z[X]}\Z[X,X^{-1}].
\]
The construction and polynomial structure constants of
Proposition~\ref{prop:genericHecke} are over $\Z[X]$;
inverse generators are used in this localization or after
specialization to a nonzero parameter.

\subsection{Modular transformations and cusp expansions}\label{app:siegel}

The Jacobi triple product and Poisson summation give the
transformation law and modularity criterion for $E_g^{(N)}$.
Yang states these formulas for $N>4$ in
\cite[Propositions~1 and~2]{Yang2}, attributing them to
\cite[Corollaries~2 and~3]{Yang1}. The nontranslation formula
there assumes a positive lower-left entry. The argument below applies
to every integer $N\ge2$ and every matrix in $\Gamma_0(N)$, with the
Siegel-function normalization of \cite[\S2.1]{Yang2}.
Here $a,b,c,d$ denote matrix entries and $S,T$ the standard generators
of $\mathrm{SL}_2(\Z)$.

\subsubsection{The theta function}
Throughout this section $\tau\in\mathbb H$ and $Q=e^{2\pi i\tau}$, and for a real number $x$ we write $Q^x$ for $e^{2\pi ix\tau}$. We denote by $(-i\tau)^{1/2}$ the branch of the square root which is positive for $\tau\in i\R_{>0}$. Let us recall that $\eta(\tau)=Q^{1/24}\prod_{n\ge1}(1-Q^n)$ and $B_2(x)=x^2-x+\frac16$, and that $\eta$ has no zeros on $\mathbb H$, because the product converges absolutely and all its factors are nonzero. For $z\in\C$ we consider the theta function:
\[
\vartheta(z,\tau)=\sum_{k\in\Z}(-1)^k\,e^{\pi i(k+\frac12)^2\tau+2\pi i(k+\frac12)z}.
\]
To compare with the usual odd theta function, our convention is
\[
 \vartheta_1(z\mid\tau)=-i\,\vartheta(z,\tau).
\]
Here $B_2$ is the Bernoulli polynomial, with periodic version
$\widetilde B_2(x)=B_2(\{x\})$.

The modulus of the term of index $k$ is $e^{-\pi(k+\frac12)^2\operatorname{Im}\tau-2\pi(k+\frac12)\operatorname{Im}z}$, thus the series converges absolutely and locally uniformly, and $\vartheta$ is holomorphic on $\C\times\mathbb H$. We shall use the following identities.

\begin{lemma}\label{lem:theta}
\begin{enumerate}[label=\textup{(\roman*)}]
\item The function $\vartheta$ satisfies
\[
 \vartheta(z+1,\tau)=-\vartheta(z,\tau),\qquad
 \vartheta(z+\tau,\tau)=-e^{-\pi i\tau-2\pi iz}\vartheta(z,\tau),
 \qquad \vartheta(-z,\tau)=-\vartheta(z,\tau).
\]
\item The function $\vartheta$ has the product expansion $\vartheta(z,\tau)=Q^{1/8}\big(e^{\pi iz}-e^{-\pi iz}\big)\prod_{n\ge1}(1-Q^n)(1-Q^ne^{2\pi iz})(1-Q^ne^{-2\pi iz})$, and in particular $\frac{\partial\vartheta}{\partial z}(0,\tau)=2\pi i\,\eta(\tau)^3$.
\item The transformations under the generators are
\[
 \vartheta(z,\tau+1)=e^{\pi i/4}\vartheta(z,\tau),\qquad
 \vartheta(z/\tau,-1/\tau)
   =-i\,(-i\tau)^{1/2}e^{\pi iz^2/\tau}\vartheta(z,\tau).
\]
\item The Dedekind eta function satisfies
\[
 \eta(\tau+1)=e^{\pi i/12}\eta(\tau),\qquad
 \eta(-1/\tau)=(-i\tau)^{1/2}\eta(\tau).
\]
\end{enumerate}
\end{lemma}
\begin{proof}
(i) The first identity holds because $e^{2\pi i(k+\frac12)}=-1$, and the third one follows if we replace $k$ by $-k-1$. For the second one, the identity $(k+\frac12)^2+2(k+\frac12)=(k+\frac32)^2-1$ shows that the term of index $k$ of $\vartheta(z+\tau,\tau)$ is equal to $e^{-\pi i\tau-2\pi iz}(-1)^ke^{\pi i(k+\frac32)^2\tau+2\pi i(k+\frac32)z}$, and it remains to replace $k$ by $k-1$.

(ii) Let us recall the Jacobi triple product identity \cite[Theorem~14.6]{Ap}:
\[
\prod_{n\ge1}(1-X^{2n})(1+X^{2n-1}Y^2)(1+X^{2n-1}Y^{-2})=\sum_{k\in\Z}X^{k^2}Y^{2k}\qquad(X,Y\in\C,\ |X|<1,\ Y\ne0).
\]
If we put $X=e^{\pi i\tau}$ and choose $Y$ with $Y^2=-e^{\pi i\tau+2\pi iz}$, it becomes:
\[
\prod_{n\ge1}(1-Q^n)(1-Q^ne^{2\pi iz})(1-Q^{n-1}e^{-2\pi iz})=\sum_{k\in\Z}(-1)^kQ^{k(k+1)/2}e^{2\pi ikz}.
\]
The right side is $Q^{-1/8}e^{-\pi iz}\vartheta(z,\tau)$ because $(k+\frac12)^2=k(k+1)+\frac14$, and on the left side the factor $1-Q^{n-1}e^{-2\pi iz}$ with $n=1$ is $1-e^{-2\pi iz}$, so the product formula follows.

The factor $e^{\pi iz}-e^{-\pi iz}$ vanishes at $z=0$ and its derivative at $z=0$ is $2\pi i$, thus $\frac{\partial\vartheta}{\partial z}(0,\tau)=2\pi i\,Q^{1/8}\prod_{n\ge1}(1-Q^n)^3=2\pi i\,\eta(\tau)^3$.

(iii) The first identity holds because $(k+\frac12)^2=k(k+1)+\frac14$ and $k(k+1)$ is even. For the second one let us fix $z$ and $\tau$ and write $\vartheta(z,\tau)=\sum_{k\in\Z}f(k)$, where $f(v)=e^{\pi iv}e^{\pi i(v+\frac12)^2\tau+2\pi i(v+\frac12)z}$ is a Schwartz function on $\R$ because $\operatorname{Im}\tau>0$. By the Poisson summation formula \cite[Ch.~VII, Prop.~15]{Serre} we have $\sum_{k\in\Z}f(k)=\sum_{m\in\Z}\widehat f(m)$, where $\widehat f(m)=\int_\R f(v)e^{-2\pi imv}dv$, and the substitution $v=\xi-\frac12$ gives:
\[
\widehat f(m)=-i\,(-1)^m\int_\R e^{\pi i\tau\xi^2+2\pi i\xi(z+\frac12-m)}\,d\xi=-i\,(-1)^m\,(-i\tau)^{-1/2}\,e^{-\pi i(z+\frac12-m)^2/\tau}.
\]
Here we used the formula $\int_\R e^{\pi i\tau\xi^2+2\pi i\xi\lambda}d\xi=(-i\tau)^{-1/2}e^{-\pi i\lambda^2/\tau}$ for $\lambda\in\C$. The integral converges locally uniformly in $(\tau,\lambda)\in\mathbb H\times\C$, hence both sides of the formula are holomorphic in $(\tau,\lambda)$. For $\tau\in i\R_{>0}$ and $\lambda\in\R$ it is the formula for the Fourier transform of a Gaussian, and it extends to all $(\tau,\lambda)\in\mathbb H\times\C$ by analytic continuation, first in $\lambda$ and then in $\tau$. Next, the identity $(z+\frac12-m)^2=z^2-2(m-\frac12)z+(m-\frac12)^2$ gives:
\[
e^{-\pi i(z+\frac12-m)^2/\tau}=e^{-\pi iz^2/\tau}\ e^{\pi i(m-\frac12)^2(-1/\tau)+2\pi i(m-\frac12)(z/\tau)} ,
\]
and putting $m=k+1$ and summing over $k$, we obtain the second identity in the following form:
\[
\vartheta(z,\tau)=i\,(-i\tau)^{-1/2}e^{-\pi iz^2/\tau}\vartheta(z/\tau,-1/\tau).
\]

(iv) The first identity holds because $Q$ does not change and $Q^{1/24}$ is multiplied by $e^{2\pi i/24}$. For the second one we differentiate the second identity of (iii) with respect to $z$ at $z=0$, and since $\vartheta(0,\tau)=0$, we get:
\[
\frac1\tau\,\frac{\partial\vartheta}{\partial z}(0,-1/\tau)=-i\,(-i\tau)^{1/2}\,\frac{\partial\vartheta}{\partial z}(0,\tau).
\]
By (ii) this says that $\eta(-1/\tau)^3=\big((-i\tau)^{1/2}\eta(\tau)\big)^3$. Thus $\eta(-1/\tau)/\big((-i\tau)^{1/2}\eta(\tau)\big)$ is a holomorphic function on $\mathbb H$ whose cube is $1$, so it is constant, and its value at $\tau=i$ is $1$.
\end{proof}

\subsubsection{Siegel functions}
Put
\[
S=\begin{pmatrix}0&-1\\1&0\end{pmatrix},\qquad T=\begin{pmatrix}1&1\\0&1\end{pmatrix}.
\]
The group $\mathrm{SL}_2(\Z)$ acts on $\mathbb H$ by $\gamma\tau=(a\tau+b)/(c\tau+d)$ for $\gamma=\left(\begin{smallmatrix}a&b\\c&d\end{smallmatrix}\right)$, and it acts from the right on row vectors $u=(u_1,u_2)\in\R^2$ by $u\gamma=(u_1a+u_2c,\,u_1b+u_2d)$. For $u\in\R^2$ the Siegel function $\mathfrak s_u$ is defined as follows:
\[
\mathfrak s_u(\tau)=e^{\pi iu_1(u_1\tau+u_2)}\ \frac{\vartheta(u_1\tau+u_2,\tau)}{\eta(\tau)} .
\]

\begin{lemma}\label{lem:siegel}
\begin{enumerate}[label=\textup{(\roman*)}]
\item The function $\mathfrak s_u$ satisfies
\[
 \mathfrak s_{-u}=-\mathfrak s_u,\qquad
 \mathfrak s_{(u_1,u_2+1)}=-e^{\pi iu_1}\mathfrak s_u,\qquad
 \mathfrak s_{(u_1+1,u_2)}=-e^{-\pi iu_2}\mathfrak s_u.
\]
\item With the notation $w=e^{2\pi i(u_1\tau+u_2)}$ we have:
\[
\mathfrak s_u(\tau)=-e^{\pi iu_2(u_1-1)}\,Q^{B_2(u_1)/2}\,(1-w)\prod_{n\ge1}(1-Q^nw)(1-Q^nw^{-1}),
\]
and if $u\notin\Z^2$, then $\mathfrak s_u$ has no zeros on $\mathbb H$.
\item There is a character $\varepsilon$ of $\mathrm{SL}_2(\Z)$ with $\varepsilon(T)=e^{\pi i/6}$ and $\varepsilon(S)=-i$ such that:
\[
\mathfrak s_u(\gamma\tau)=\varepsilon(\gamma)\,\mathfrak s_{u\gamma}(\tau)\qquad\text{for all }\gamma\in\mathrm{SL}_2(\Z),\ u\in\R^2\text{ and }\tau\in\mathbb H .
\]
\end{enumerate}
\end{lemma}
\begin{proof}
Throughout the proof we put $z=u_1\tau+u_2$, thus $w=e^{2\pi iz}$ in the notation of (ii).

(i) The first identity follows from $\vartheta(-z,\tau)=-\vartheta(z,\tau)$, because replacing $u$ by $-u$ replaces $z$ by $-z$ and does not change $u_1z$. If we replace $u_2$ by $u_2+1$, then $z$ is replaced by $z+1$, and by Lemma~\ref{lem:theta}(i) the function $e^{\pi iu_1z}\vartheta(z,\tau)$ is multiplied by $-e^{\pi iu_1}$. If we replace $u_1$ by $u_1+1$, then $z$ is replaced by $z+\tau$, and Lemma~\ref{lem:theta}(i) gives:
\[
e^{\pi i(u_1+1)(z+\tau)}\,\vartheta(z+\tau,\tau)=-e^{\pi i(u_1+1)(z+\tau)-\pi i\tau-2\pi iz}\,\vartheta(z,\tau)=-e^{\pi i(u_1\tau-z)}\,e^{\pi iu_1z}\,\vartheta(z,\tau),
\]
and it remains to note that $u_1\tau-z=-u_2$.

(ii) By Lemma~\ref{lem:theta}(ii) and the definition of $\eta$ we have:
\[
\mathfrak s_u(\tau)=Q^{1/12}\,e^{\pi iu_1z}\big(e^{\pi iz}-e^{-\pi iz}\big)\prod_{n\ge1}(1-Q^nw)(1-Q^nw^{-1}) .
\]
Moreover, $e^{\pi iu_1z}(e^{\pi iz}-e^{-\pi iz})=-e^{\pi i(u_1-1)z}(1-w)$ and $e^{\pi i(u_1-1)z}=Q^{u_1(u_1-1)/2}e^{\pi iu_2(u_1-1)}$, so the formula follows from the identity $\frac1{12}+\frac12u_1(u_1-1)=\frac12B_2(u_1)$.

Now let $u\notin\Z^2$. By (i) we may assume that $0\le u_1<1$, and in this case $|Q|<|w|\le1$, so $|Q^nw|<1$ and $|Q^nw^{-1}|<1$ for all $n\ge1$. On the other hand, $w=1$ would mean that $u_1\tau+u_2\in\Z$, and since $\tau$ is not real, this happens only for $u_1=0$ and $u_2\in\Z$. Thus no factor vanishes, and the product converges to a nonzero number.

(iii) Let us first treat $T$ and $S$, for which $uT=(u_1,u_1+u_2)$ and $uS=(u_2,-u_1)$, and Lemma~\ref{lem:theta}(iii),(iv) gives:
\[
\mathfrak s_u(\tau+1)=e^{\pi iu_1(u_1\tau+u_1+u_2)}\ \frac{e^{\pi i/4}\,\vartheta(u_1\tau+u_1+u_2,\tau)}{e^{\pi i/12}\,\eta(\tau)}=e^{\pi i/6}\,\mathfrak s_{uT}(\tau).
\]
Next, for $\zeta=u_2\tau-u_1$ we have $u_1\cdot(-1/\tau)+u_2=\zeta/\tau$, and we apply Lemma~\ref{lem:theta}(iii),(iv) once more:
\[
\mathfrak s_u(-1/\tau)=e^{\pi iu_1\zeta/\tau}\ \frac{-i\,(-i\tau)^{1/2}\,e^{\pi i\zeta^2/\tau}\,\vartheta(\zeta,\tau)}{(-i\tau)^{1/2}\,\eta(\tau)}=-i\,e^{\pi i\zeta(u_1+\zeta)/\tau}\ \frac{\vartheta(\zeta,\tau)}{\eta(\tau)}=-i\,\mathfrak s_{uS}(\tau),
\]
where the last equality holds because $u_1+\zeta=u_2\tau$.

Now let $\Gamma'$ be the set of all $\gamma\in\mathrm{SL}_2(\Z)$ for which there is a constant $\varepsilon(\gamma)$ with $\mathfrak s_u(\gamma\tau)=\varepsilon(\gamma)\mathfrak s_{u\gamma}(\tau)$ for all $u$ and $\tau$. For $u\notin\Z^2$ we have $u\gamma\notin\Z^2$, and the functions $\mathfrak s_u$ and $\mathfrak s_{u\gamma}$ have no zeros by (ii), hence the constant $\varepsilon(\gamma)$ is unique and nonzero. We have shown that $S,T\in\Gamma'$, and for $\gamma,\tilde\gamma\in\Gamma'$ we have:
\[
\mathfrak s_u(\gamma\tilde\gamma\tau)=\varepsilon(\gamma)\,\mathfrak s_{u\gamma}(\tilde\gamma\tau)=\varepsilon(\gamma)\varepsilon(\tilde\gamma)\,\mathfrak s_{u\gamma\tilde\gamma}(\tau),
\]
so $\gamma\tilde\gamma\in\Gamma'$ and $\varepsilon(\gamma\tilde\gamma)=\varepsilon(\gamma)\varepsilon(\tilde\gamma)$. If $\gamma\in\Gamma'$, then replacing $u$ by $u\gamma^{-1}$ and $\tau$ by $\gamma^{-1}\tau$ we get $\gamma^{-1}\in\Gamma'$. Thus $\Gamma'$ is a subgroup and $\varepsilon$ is a character of $\Gamma'$, and it remains to recall that $S$ and $T$ generate $\mathrm{SL}_2(\Z)$.

Indeed, for $\gamma=\left(\begin{smallmatrix}a&b\\c&d\end{smallmatrix}\right)$ with $c=0$ we have $\gamma=T^b$ or $\gamma=S^2T^{-b}$, and if $c\ne0$, we choose an integer $n$ with $|a+nc|<|c|$ and observe that the lower left entry of $ST^n\gamma$ is $a+nc$, so we conclude by induction on $|c|$.
\end{proof}

For $\gamma=\left(\begin{smallmatrix}a&b\\c&d\end{smallmatrix}\right)$
we also write $\varepsilon(a,b,c,d)=\varepsilon(\gamma)$.
The next lemma computes this character and identifies it with the
multiplier of $\eta^2$.

\begin{lemma}\label{lem:etachar}
For $\gamma=\left(\begin{smallmatrix}a&b\\c&d\end{smallmatrix}\right)
\in\mathrm{SL}_2(\Z)$ the following statements hold.
\begin{enumerate}[label=\textup{(\roman*)}]
\item We have $\varepsilon^{12}=1$ and
$\varepsilon(-\gamma)=-\varepsilon(\gamma)$.
\item The character is given by
\[
\varepsilon(\gamma)=\begin{cases}
 e^{\pi i(bd(1-c^2)+c(a+d-3))/6}&\text{if $c$ is odd,}\\[3pt]
 -i\,e^{\pi i(ac(1-d^2)+d(b-c+3))/6}&\text{if $d$ is odd.}
\end{cases}
\]
In particular, $\varepsilon|_{\Gamma(12)}=1$.
\item The transformation formula for $\eta^2$ is
\begin{equation}\label{eq:eta-square-character}
 \eta(\gamma\tau)^2
   =\varepsilon(\gamma)(c\tau+d)\eta(\tau)^2.
\end{equation}
\end{enumerate}
\end{lemma}
\begin{proof}
Since $\varepsilon(T)$ and $\varepsilon(S)$ are $12$th roots of
unity and $S,T$ generate $\mathrm{SL}_2(\Z)$, we have
$\varepsilon^{12}=1$. Also $-I=S^2$ and
$\varepsilon(-I)=(-i)^2=-1$, which proves the second assertion of (i).

For an integer $n$ we put $\mathbf e(n)=e^{\pi in/6}$, a number depending only on $n$ modulo $12$, and we consider the following two functions on $\mathrm{SL}_2(\Z)$:
\[
e_1(\gamma)=bd(1-c^2)+c(a+d-3),\qquad e_2(\gamma)=ac(1-d^2)+d(b-c+3).
\]
We have to show that $\varepsilon(\gamma)=\mathbf e(e_1(\gamma))$ if $c$ is odd and that $\varepsilon(\gamma)=\mathbf e(e_2(\gamma)-3)$ if $d$ is odd. Note that at least one of these two cases occurs, because $c$ and $d$ are coprime, and that the relation $ad-bc=1$ gives:
\[
e_1(\gamma)-e_2(\gamma)+3=cd\,(ad-bc)+2cd-3c-3d+3=3(c-1)(d-1),
\]
and if $c$ and $d$ are both odd, then this number is divisible by $12$. Hence the two expressions agree in this case, and they define a function $\varepsilon'$ on $\mathrm{SL}_2(\Z)$. Since $e_2(1)=3$, we have $\varepsilon'(1)=1$, and we claim that:
\begin{equation}\label{eq:etachar}
\varepsilon'(\gamma S)=\mathbf e(-3)\,\varepsilon'(\gamma)\qquad\text{and}\qquad\varepsilon'(\gamma T)=\mathbf e(1)\,\varepsilon'(\gamma)\qquad\text{for all }\gamma\in\mathrm{SL}_2(\Z).
\end{equation}
Since $\varepsilon(S)=\mathbf e(-3)$ and $\varepsilon(T)=\mathbf e(1)$, the claim shows that $\varepsilon'(\gamma)/\varepsilon(\gamma)$ does not change if $\gamma$ is multiplied from the right by $S^{\pm1}$ or by $T^{\pm1}$. As $S$ and $T$ generate $\mathrm{SL}_2(\Z)$ and $\varepsilon'(1)=\varepsilon(1)$, this gives $\varepsilon'=\varepsilon$.

\emph{The claim for $S$.} Here $\gamma S=\left(\begin{smallmatrix}b&-a\\d&-c\end{smallmatrix}\right)$, and we first let $d$ be even. Then $c$ is odd and the definitions give $e_2(\gamma S)=bd(1-c^2)+c(a+d-3)=e_1(\gamma)$, and since the lower right entry $-c$ of $\gamma S$ is odd, we get $\varepsilon'(\gamma S)=\mathbf e(e_1(\gamma)-3)=\mathbf e(-3)\,\varepsilon'(\gamma)$. Now let $d$ be odd. Then the lower left entry $d$ of $\gamma S$ is odd and the definitions give $e_1(\gamma S)=ac(1-d^2)+d(b-c-3)=e_2(\gamma)-6d$, hence $\varepsilon'(\gamma S)=\mathbf e(e_2(\gamma)-6)=\mathbf e(-3)\,\varepsilon'(\gamma)$, again because $d$ is odd.

\emph{The claim for $T$.} Here $\gamma T=\left(\begin{smallmatrix}a&a+b\\c&c+d\end{smallmatrix}\right)$, and in both cases below we use the relation $ad=1+bc$. If $c$ is odd, then we have:
\[
e_1(\gamma T)-e_1(\gamma)=(1-c^2)(ac+ad+bc)+c^2=1-(c-1)\,c\,(c+1)\,(a+2b).
\]
Since $c$ is odd, $8$ divides $(c-1)(c+1)$, and $3$ divides $(c-1)c(c+1)$, so $\varepsilon'(\gamma T)=\mathbf e(1)\,\varepsilon'(\gamma)$. Now let $c$ be even. Then $a$, $d$ and $c+d$ are odd, and for $\phi=a(c+1)(c+2d)+c-2$ we get:
\[
e_2(\gamma T)-e_2(\gamma)=-ac^2(c+2d)+c(a+b-c+3)+ad=1-(c-1)\,\phi .
\]
Let us show that $12$ divides $(c-1)\phi$: if $4$ divides $c$, then $\phi\equiv2ad-2\equiv0$ modulo $4$, because $ad$ is odd, and if $c\equiv2$ modulo $4$, then $4$ divides $c+2d$ and $c-2$. Thus $4$ divides $\phi$ in both cases.

Modulo $3$ there are three cases: if $c\equiv1$ modulo $3$, then $3$ divides $c-1$; if $c\equiv0$ modulo $3$, then $\phi\equiv2ad-2=2bc\equiv0$ modulo $3$; and if $c\equiv-1$ modulo $3$, then $\phi\equiv c-2\equiv0$ modulo $3$. Hence $12$ divides $(c-1)\phi$, so $\varepsilon'(\gamma T)=\mathbf e(1)\,\varepsilon'(\gamma)$ also for even $c$, and \eqref{eq:etachar} is proved.

For $\gamma\in\Gamma(12)$ the entry $d$ is odd and
$e_2(\gamma)\equiv3$ modulo $12$, so
$\varepsilon(\gamma)=\mathbf e(0)=1$. This proves (ii).

For (iii), Lemma~\ref{lem:theta}(iv) gives
\[
 \eta(T\tau)^2=e^{\pi i/6}\eta(\tau)^2,\qquad
 \eta(S\tau)^2=(-i\tau)\eta(\tau)^2.
\]
These are \eqref{eq:eta-square-character} for $T$ and $S$.
Writing $j(\gamma,\tau)=c\tau+d$, direct matrix multiplication gives
\[
 j(\gamma\gamma',\tau)
  =j(\gamma,\gamma'\tau)j(\gamma',\tau).
\]
Together with the character property of $\varepsilon$, this shows that
\eqref{eq:eta-square-character} is preserved under multiplication
and inversion. The generator formulas therefore prove it for every
$\gamma\in\mathrm{SL}_2(\Z)$.
\end{proof}

\subsubsection{The functions \texorpdfstring{$E^{(N)}_g$}{E}}
For integers $N\ge2$ and $1\le g\le N-1$, define
\[
 E_g^{(N)}(\tau)=Q^{N B_2(g/N)/2}
 \prod_{n\ge1}(1-Q^{(n-1)N+g})(1-Q^{nN-g}),
\]
and extend to $g\in\Z\setminus N\Z$ by
$E_{g+N}^{(N)}=E_{-g}^{(N)}=-E_g^{(N)}$.
Write $E_g=E_g^{(N)}$. A modular function on $\Gamma_1(N)$ is a
meromorphic function on $\mathbb H$ that is invariant under
$\Gamma_1(N)$ and is meromorphic at every cusp. A \emph{modular unit}
is a modular function with no zero or pole on $\mathbb H$.
At a cusp represented by $\sigma\in\mathrm{SL}_2(\Z)$, meromorphy is
measured in a local parameter $e^{2\pi i\tau/w}$, where $w$ is the
cusp width. An exponent of $Q=e^{2\pi i\tau}$ must therefore be
multiplied by $w$ to give the order in that local parameter.

\Needspace{6\baselineskip}
\begin{proposition}\label{prop:Etransf}
Let $g$ be an integer which is not divisible by $N$.
\begin{enumerate}[label=\textup{(\roman*)}]
\item The Siegel-function and theta-function expressions are
\[
 \begin{aligned}
 E_g(\tau)&=-\mathfrak s_{(g/N,0)}(N\tau)\\
 &=-i\,e^{\pi ig^2\tau/N}
      \frac{\vartheta_1(g\tau\mid N\tau)}{\eta(N\tau)}.
 \end{aligned}
\]
\item For $\gamma=\left(\begin{smallmatrix}\alpha&\beta\\ \gamma_0N&\delta_0\end{smallmatrix}\right)\in\Gamma_0(N)$ we have:
\[
E_g(\gamma\tau)=\varepsilon(\alpha,\beta N,\gamma_0,\delta_0)\ e^{\pi i(g^2\alpha\beta/N-g\beta)}\ E_{\alpha g}(\tau).
\]
\item Let $(e_g)$ be a finite family of integers indexed by
$g\in\Z\setminus N\Z$. If
\[
 \sum_g e_g\equiv0\pmod{12},\qquad
 \sum_g ge_g\equiv0\pmod2,\qquad
 \sum_g g^2e_g\equiv0\pmod{2N},
\]
then $\prod_gE_g^{e_g}$ is a modular unit on $\Gamma_1(N)$.
\item For $0<g<N$ the function $E_g(-1/\tau)$ is given by the following product:
\[
E_g(-1/\tau)=2\sin\frac{\pi g}N\ Q^{\frac1{12N}}\ \prod_{n\ge1}\Big(1-2\cos\frac{2\pi g}N\,Q^{n/N}+Q^{2n/N}\Big).
\]
\end{enumerate}
\end{proposition}
\begin{proof}
(i) For $0<g<N$ we apply Lemma~\ref{lem:siegel}(ii) to $u=(g/N,0)$ with $N\tau$ in place of $\tau$, in which case $Q$ is replaced by $Q^N$ and $w=Q^g$, and we get:
\[
\mathfrak s_{(g/N,0)}(N\tau)=-Q^{NB_2(g/N)/2}\,(1-Q^g)\prod_{n\ge1}(1-Q^{nN+g})(1-Q^{nN-g})=-E_g(\tau).
\]
If $g$ is replaced by $g+N$ or by $-g$, then $E_g$ changes sign
by definition and the Siegel expression changes sign by
Lemma~\ref{lem:siegel}(i). Thus the first equality holds for all
$g\in\Z\setminus N\Z$. Substituting the definition of
$\mathfrak s_{(g/N,0)}$ and $\vartheta=i\vartheta_1$ gives the second
equality. This proves (i).

(ii) The matrix $\gamma'=\left(\begin{smallmatrix}\alpha&\beta N\\ \gamma_0&\delta_0\end{smallmatrix}\right)\in\mathrm{SL}_2(\Z)$ satisfies $N\cdot\gamma\tau=\gamma'(N\tau)$ and $(g/N,0)\gamma'=(\alpha g/N,\beta g)$. By Lemma~\ref{lem:siegel}(i) we have $\mathfrak s_{(u_1,u_2+n)}=(-1)^ne^{\pi iu_1n}\mathfrak s_u$ for all integers $n$, and together with (i) and Lemma~\ref{lem:siegel}(iii) this gives:
\[
\begin{aligned}
E_g(\gamma\tau)=-\varepsilon(\gamma')\,\mathfrak s_{(\alpha g/N,\beta g)}(N\tau)&=-\varepsilon(\gamma')\,(-1)^{\beta g}e^{\pi i\alpha\beta g^2/N}\,\mathfrak s_{(\alpha g/N,0)}(N\tau)\\
&=\varepsilon(\gamma')\,e^{\pi i(g^2\alpha\beta/N-g\beta)}\,E_{\alpha g}(\tau).
\end{aligned}
\]
Note that $\alpha g$ is not divisible by $N$, because $\alpha$ is prime to $N$.

(iii) By (i) and Lemma~\ref{lem:siegel}(ii), the function $\Pi=\prod_gE_g^{e_g}$ is holomorphic on $\mathbb H$ and has no zeros. For $\gamma\in\Gamma_1(N)$, with the notation of (ii), we have $\alpha=1+N\alpha'$ with $\alpha'\in\Z$, and (ii) gives:
\[
\Pi(\gamma\tau)=\varepsilon(\gamma')^{\sum_ge_g}\ e^{\pi i(\alpha\beta/N)\sum_gg^2e_g}\ e^{-\pi i\beta\sum_gge_g}\ \prod_gE_{\alpha g}(\tau)^{e_g}.
\]
The first factor is $1$ because $\varepsilon^{12}=1$ by Lemma~\ref{lem:etachar}, the second and third factors are $1$ by the assumptions on $\sum_gg^2e_g$ and on $\sum_gge_g$, and the last factor is $(-1)^{\alpha'\sum_gge_g}\Pi(\tau)=\Pi(\tau)$, because $E_{\alpha g}=E_{g+N\alpha'g}=(-1)^{\alpha'g}E_g$. Thus $\Pi$ is invariant under $\Gamma_1(N)$.

It remains to show that $\Pi$ is meromorphic at the cusps. Let $\sigma\in\mathrm{SL}_2(\Z)$, and let $h>0$ be the greatest common divisor of the two entries of the first column of the integer matrix $\left(\begin{smallmatrix}N&0\\0&1\end{smallmatrix}\right)\sigma$. If we choose $\sigma_1\in\mathrm{SL}_2(\Z)$ whose first column is this column divided by $h$, then $\sigma_1^{-1}\left(\begin{smallmatrix}N&0\\0&1\end{smallmatrix}\right)\sigma$ is an integer matrix of determinant $N$ with first column $(h,0)$, that is:
\[
\begin{pmatrix}N&0\\0&1\end{pmatrix}\sigma=\sigma_1\begin{pmatrix}h&j\\0&N/h\end{pmatrix}\qquad\text{with an integer }j .
\]
Let us put $\tau_1=h(h\tau+j)/N$. Then $N\cdot\sigma\tau=\sigma_1\tau_1$, and by (i) and Lemma~\ref{lem:siegel}(iii) we get $E_g(\sigma\tau)=-\varepsilon(\sigma_1)\,\mathfrak s_u(\tau_1)$ with $u=(g/N,0)\sigma_1\notin\Z^2$. By Lemma~\ref{lem:siegel}(i) we may replace $u$ by a vector $u'\in u+\Z^2$ with $0\le u'_1<1$, at the cost of a nonzero constant factor. Let us put $Q_1=e^{2\pi i\tau_1}$ and $w=e^{2\pi i(u'_1\tau_1+u'_2)}$, and observe that $|Q_1|=|Q|^{h^2/N}\to0$ as $\operatorname{Im}\tau\to\infty$.

If $u'_1>0$, then $|w|=|Q_1|^{u'_1}\to0$ and $\sum_{n\ge1}\big(|Q_1^nw|+|Q_1^nw^{-1}|\big)\le2\sum_{n\ge1}|Q_1|^{n-u'_1}\to0$, and if $u'_1=0$, then $w=e^{2\pi iu'_2}$ is a constant different from $1$, because $u'\notin\Z^2$.

In both cases $(1-w)\prod_{n\ge1}(1-Q_1^nw)(1-Q_1^nw^{-1})$ tends to a nonzero constant, and the convergence is uniform in $\operatorname{Re}\tau$ because these bounds depend only on $|Q_1|$. Together with Lemma~\ref{lem:siegel}(ii) this gives $\mathfrak s_{u'}(\tau_1)=\kappa\,e^{\pi iB_2(u'_1)\tau_1}(1+o(1))$ with a constant $\kappa\ne0$, and thus:
\[
E_g(\sigma\tau)=\kappa_\sigma\,Q^{r_\sigma}\,(1+o(1))\qquad\text{with }\kappa_\sigma\ne0\text{ and }r_\sigma=h^2B_2(u'_1)/(2N).
\]
Hence $\Pi(\sigma\tau)=\kappa'Q^{r'}(1+o(1))$ with a constant $\kappa'\ne0$ and a rational number $r'$, uniformly in $\operatorname{Re}\tau$. The function $\tau\mapsto\Pi(\sigma\tau)$ is invariant under $\tau\mapsto\tau+N$, because $\sigma T^N\sigma^{-1}\in\Gamma(N)\subseteq\Gamma_1(N)$, so it is a holomorphic function of $Q^{1/N}$ on the punctured unit disc. By the estimate it is bounded by a power of $|Q^{1/N}|$ near the origin, hence it has at most a pole there.

(iv) We have $-N/\tau=S(\tau/N)$ and $(g/N,0)S=(0,-g/N)$, and by (i) and Lemma~\ref{lem:siegel}(iii) we get:
\[
E_g(-1/\tau)=-\mathfrak s_{(g/N,0)}\big(S(\tau/N)\big)=i\,\mathfrak s_{(0,-g/N)}(\tau/N).
\]
We apply Lemma~\ref{lem:siegel}(ii) to $u=(0,-g/N)$ with $\tau/N$ in place of $\tau$. Then $w=e^{-2\pi ig/N}$, the variable $Q$ is replaced by $Q^{1/N}$, and the constant factor is $-e^{\pi ig/N}(1-e^{-2\pi ig/N})=-2i\sin(\pi g/N)$. The factors of the product are $(1-Q^{n/N}w)(1-Q^{n/N}w^{-1})=1-2\cos(2\pi g/N)\,Q^{n/N}+Q^{2n/N}$, and $B_2(0)=\frac16$, so together with the factor $i$ we obtain (iv).
\end{proof}

\subsection{Character normalization of the cabling product}\label{supp:character}

The torus product $P_{a,b}$ is the normalized character of a
balanced minimal-model sector. We compute its conformal exponent and then
apply the cabling product formula to obtain a tensor-product character.

Fix coprime integers $2\le a<b$, and put $M=a+b$. Recall the product
\[
 \rho_{a,b}(j)=\min\{a,b,\operatorname{dist}(aj,M\Z)\}-1
                  +\mathbf1_{M\mid j},\qquad
 P_{a,b}(Q)=\prod_{j\ge1}(1-Q^j)^{-\rho_{a,b}(j)}.
\]
Here $\operatorname{dist}(aj,M\Z)$ is the distance to the nearest
multiple of $M$. For $\tau\in\mathbb H$ set
\[
 Q=e^{2\pi i\tau},\qquad
 c_{a,b}=\frac{(a-1)(b-1)}{24M},\qquad
 f_{a,b}(\tau)=Q^{-c_{a,b}}P_{a,b}(Q),
\]
where $Q^x=e^{2\pi ix\tau}$ for rational $x$.

We use the minimal model $\mathcal W_a(a,M)$. Its irreducible sectors
are indexed, up to cyclic rotation, by positive integer tuples
$j=(j_0,\ldots,j_{a-1})$ of sum $M$; see
\cite[\S9.1]{HJO} and \cite[\S\S3.2--3.3]{FW}.
Let $V_j$ denote the corresponding module, $h_j$ its lowest
conformal weight, and $c_{\mathrm{VOA}}$ the central charge of this
minimal model. Our character conventions are
\[
 \overline\chi_j(Q)=\operatorname{Tr}_{V_j}Q^{L_0-h_j},\qquad
 \chi_j(\tau)=\operatorname{Tr}_{V_j}Q^{L_0-c_{\mathrm{VOA}}/24}
       =Q^{h_j-c_{\mathrm{VOA}}/24}\overline\chi_j(Q).
\]
Thus $\overline\chi_j$ is the normalized character, with constant
term $1$, and $\chi_j$ is the ordinary character. The number
$c_{a,b}$ is a normalizing exponent, not the VOA central charge.

The balanced, or Euclidean-rhythm, tuple is
\[
 j_i=\left\lceil\frac{(i+1)M}{a}\right\rceil
          -\left\lceil\frac{iM}{a}\right\rceil,
 \qquad 0\le i<a.
\]
This agrees with the convention of \cite{HJO}: writing
$M=\ell a+r$ with $0<r<a$, the entries equal to $\ell+1$ occur
at $i=\lfloor ka/r\rfloor$ modulo $a$.
The character identity of \cite[Proposition~9.3]{HJO} is
\begin{equation}\label{eq:balanced-normalized-character}
 P_{a,b}(Q)=\overline\chi_j(Q).
\end{equation}

\begin{proposition}[Balanced-sector normalization]
\label{prop:balanced-character-normalization}
For the balanced tuple $j$, one has
\[
 h_j-\frac{c_{\mathrm{VOA}}}{24}=-c_{a,b},\qquad
 f_{a,b}(\tau)=\chi_j(\tau).
\]
Moreover, $h_j$ is the lowest conformal weight among all irreducible
sectors of $\mathcal W_a(a,M)$. Consequently the effective central
charge of this minimal model is
\[
 c_{\mathrm{eff}}
    :=c_{\mathrm{VOA}}-24\min_{j'}h_{j'}
     =24c_{a,b}=\frac{(a-1)(b-1)}{a+b}.
\]
\end{proposition}
\begin{proof}
Let $\varrho$ and $\omega_i$ be the Weyl vector and fundamental
weights of $A_{a-1}$, with roots of squared length $2$, and set
$\lambda_j=\sum_{i=1}^{a-1}j_i\omega_i$.
In the notation of Foda--Welsh, take $r=p=a$, $p'=M$ and $\xi=0$.
Their affine weight $\zeta$ has Dynkin labels $j_i-1$, so its finite
part satisfies $\zeta_{\mathrm{fin}}+\varrho=\lambda_j$.
The conformal exponent in \cite[\S3.2, (24), and \S3.3]{FW} is
\begin{equation}\label{eq:balanced-conformal-exponent}
 \Delta_j:=h_j-\frac{c_{\mathrm{VOA}}}{24}
    =\frac{\|M\varrho-a\lambda_j\|^2}{2aM}
        -\frac{a-1}{24}.
\end{equation}

For the balanced tuple, set
$s_i=a\lceil iM/a\rceil-iM$, $1\le i\le a$.
Each $s_i$ lies in $\{0,\ldots,a-1\}$ and satisfies
$s_i\equiv-iM\pmod a$. Since $\gcd(a,M)=1$, these residues are
distinct, and the $s_i$ permute $0,\ldots,a-1$.
In the sum-zero realization of the weight space, the vector
$v=M\varrho-a\lambda_j$ has successive differences
$v_i-v_{i+1}=M-aj_i$. The coordinates
\[
 v_i=s_i-\frac{a-1}{2}\qquad(1\le i\le a)
\]
have sum zero and precisely these successive differences, because
$s_i-s_{i+1}=M-aj_i$. They therefore give the coordinates of $v$.
It follows that $v$ is a Weyl permutation of $\varrho$, and hence
\[
 \|v\|^2=\|\varrho\|^2=\frac{a(a^2-1)}{12}.
\]
Substitution into \eqref{eq:balanced-conformal-exponent} gives
\[
 \Delta_j
   =\frac{a^2-1}{24M}-\frac{a-1}{24}
   =-\frac{(a-1)(b-1)}{24M}=-c_{a,b}.
\]
Together with \eqref{eq:balanced-normalized-character}, this proves
the ordinary-character identity.

Now let $j'$ be any positive cyclic tuple of sum $M$, and set
$v'=M\varrho-a\lambda_{j'}$. For $i<k$,
\[
 v'_i-v'_k=M(k-i)-a\sum_{t=i}^{k-1}j'_t
       \equiv M(k-i)\pmod a.
\]
These differences are nonzero integers because $0<k-i<a$ and
$\gcd(a,M)=1$. If the coordinates are sorted as
$x_1<\cdots<x_a$, their successive gaps are at least $1$, and
$x_k-x_i\ge k-i$. Since their sum is zero,
\[
 \|v'\|^2
   =\frac1a\sum_{i<k}(x_k-x_i)^2
   \ge\frac1a\sum_{i<k}(k-i)^2
   =\|\varrho\|^2.
\]
The central charge is independent of the sector, so
\eqref{eq:balanced-conformal-exponent} proves that $h_j$ is minimal.
Finally $c_{\mathrm{eff}}=-24\Delta_j=24c_{a,b}$.
\end{proof}

\begin{corollary}[Tensor-product normalization]
\label{cor:cabling-character-normalization}
Let $((a_i,b_i,D_i))_{1\le i\le g}$ be a finite list with
$2\le a_i<b_i$, $\gcd(a_i,b_i)=1$, and
$D_i\in\Z_{\ge0}$. Suppose that a series $\cN_A(Q)$ satisfies
the product identity
\[
 \cN_A(Q)=\prod_{i=1}^gP_{a_i,b_i}(Q)^{D_i}.
\]
Put $c_A=\sum_iD_i c_{a_i,b_i}$ and
$\Theta_A(\tau)=Q^{-c_A}\cN_A(Q)$.
For each $i$, let $V_i$ be the balanced sector of
$\mathcal W_{a_i}(a_i,a_i+b_i)$. Then $\Theta_A$ is the ordinary
character of
\[
 V_A=\bigotimes_{i=1}^g V_i^{\otimes D_i}
 \quad\text{over}\quad
 \mathcal W_A=
   \bigotimes_{i=1}^g\mathcal W_{a_i}(a_i,a_i+b_i)^{\otimes D_i},
\]
and the effective central charge is
\[
 c_{\mathrm{eff}}(\mathcal W_A)=24c_A
     =\sum_iD_i\frac{(a_i-1)(b_i-1)}{a_i+b_i}.
\]
For the empty list the tensor product is $\C$, the character is $1$,
and $c_A=0$. Factors with $D_i=0$ are omitted.
\end{corollary}
\begin{proof}
On a tensor product, $L_0$ is the sum of the grading operators and
the central charges add. The graded pieces are finite-dimensional,
so multiplication of their character series gives
\[
 \operatorname{Tr}_{V_A}Q^{L_0-c_{\mathrm{VOA}}(\mathcal W_A)/24}
   =\prod_i f_{a_i,b_i}(\tau)^{D_i}
   =Q^{-c_A}\cN_A(Q).
\]
The lowest conformal weights also add, and each chosen factor
attains the minimum for its minimal model by
Proposition~\ref{prop:balanced-character-normalization}. The effective
central charges therefore add, giving the asserted formula.
\end{proof}

\begin{example}[The $(2,3)$ normalization]
For $(a,b)=(2,3)$, the balanced tuple is $(3,2)$ and
$\mathcal W_2(2,5)$ is the corresponding Virasoro minimal model.
For its vacuum label $(4,1)$, one has $\lambda_j=\omega_1$ and
$\|5\varrho-2\lambda_j\|^2=9/2$, so
\eqref{eq:balanced-conformal-exponent} gives
$-c_{\mathrm{VOA}}/24=9/40-1/24=11/60$.
Hence $c_{\mathrm{VOA}}=-22/5$, and the balanced exponent gives
$h_{\min}=-1/5$ and
\[
 h_{\min}-\frac{c_{\mathrm{VOA}}}{24}=-\frac1{60},\qquad
 f_{2,3}(\tau)=Q^{-1/60}P_{2,3}(Q),\qquad
 c_{\mathrm{eff}}=\frac25.
\]
\end{example}


\begin{thebibliography}{HJO}
\bibitem[Ap]{Ap} T. M. Apostol, \emph{Introduction to analytic number theory}, Undergrad. Texts Math., Springer, 1976.

\bibitem[BDF]{BDF} V. Barucci, M. D'Anna and R. Fr\"oberg,
On plane algebroid curves, arXiv:math/0302224v1 (2003).

\bibitem[BeDe]{BeDe} J. N. Bernstein, r\'edig\'e par P. Deligne, Le ``centre'' de Bernstein, in: \emph{Repr\'esentations des groupes r\'eductifs sur un corps local}, Travaux en Cours, Hermann, Paris, 1984, 1--32.

\bibitem[Bo]{Borel} A. Borel, Admissible representations of a semi-simple group over a local field with vectors fixed under an Iwahori subgroup, \emph{Invent. Math.} 35 (1976), 233--259.

\bibitem[Bou]{Bou} N. Bourbaki, \emph{Groupes et alg\`ebres de Lie, Chapitres 4, 5 et 6}, Hermann, 1968.

\bibitem[BH]{BH} W. Bruns and J. Herzog, \emph{Cohen--Macaulay rings}, Cambridge Stud. Adv. Math. 39, Cambridge Univ. Press, 1993.

\bibitem[Ca]{Campillo} A. Campillo, \emph{Algebroid curves in positive characteristic}, Lecture Notes in Math. 813, Springer, 1980.

\bibitem[Cas]{Casals} R. Casals,
Lagrangian skeleta and plane curve singularities, arXiv:2009.06737.
\url{https://www.math.ucdavis.edu/~casals/Viterbo60_RCasals.pdf}.

\bibitem[Ch]{Ch25} I. Cherednik, Superpolynomials of algebraic links,
arXiv:2410.14703v2 (2025).

\bibitem[CD]{CD} I. Cherednik and I. Danilenko, DAHA and iterated torus knots,
\emph{Algebr. Geom. Topol.} 16 (2016), 843--898.

\bibitem[CP]{CP} I. Cherednik and I. Philipp,
Modules over plane curve singularities in any ranks and DAHA,
arXiv:1711.07605v2 (2018).

\bibitem[CL]{CL} H. Cohen and H. W. Lenstra, Jr., Heuristics on class groups of number fields, in: \emph{Number theory, Noordwijkerhout 1983}, Lecture Notes in Math. 1068, Springer, 1984, 33--62.

\bibitem[DHS]{DHS} D.-E. Diaconescu, Z. Hua and Y. Soibelman,
HOMFLY polynomials, stable pairs and motivic Donaldson--Thomas invariants,
\emph{Commun. Number Theory Phys.} 6 (2012), no.~3, 517--600.

\bibitem[Dw]{Dwork} B. Dwork, On the rationality of the zeta function of an algebraic variety, \emph{Amer. J. Math.} 82 (1960), 631--648.

\bibitem[FW]{FW} O. Foda and T. A. Welsh,
Cylindric partitions, $\mathcal W_r$ characters and the
Andrews--Gordon--Bressoud identities,
\emph{J. Phys. A: Math. Theor.} 49 (2016), no.~16, 164004;
arXiv:1510.02213v1.

\bibitem[GP]{GP} M. Geck and G. Pfeiffer, \emph{Characters of finite Coxeter groups and Iwahori--Hecke algebras}, London Math. Soc. Monogr. (N.S.) 21, Oxford Univ. Press, 2000.

\bibitem[GM]{GM} E. Gorsky and M. Mazin, Compactified Jacobians and $q,t$-Catalan numbers, I, \emph{J. Combin. Theory Ser. A} 120 (2013), 49--63.

\bibitem[GN]{GN} E. Gorsky and A. Negu\c{t},
Refined knot invariants and Hilbert schemes,
\emph{J. Math. Pures Appl.} 104 (2015), 403--435;
arXiv:1304.3328v3.

\bibitem[He]{He} X. He, Cocenters of $p$-adic groups, I: Newton decomposition, \emph{Forum Math. Pi} 6 (2018), e2; arXiv:1610.04791.

\bibitem[HN]{HN} X. He and S. Nie, Minimal length elements of extended affine Weyl groups, \emph{Compos. Math.} 150 (2014), 1903--1927; arXiv:1112.0824v1.

\bibitem[Hu]{Huang} Y. Huang, Mutually annihilating matrices, and a Cohen--Lenstra series for the nodal singularity, \emph{J. Algebra} 619 (2023), 26--50.

\bibitem[Hu2]{Huang2} Y. Huang, A quadratic form generalization of rational dinv, with an appendix by K. Lau, \emph{Res. Math. Sci.} 13 (2026), 44; arXiv:2604.13238.

\bibitem[HJ]{HJ} Y. Huang and R. Jiang, Motivic Coh and Quot zeta functions of singular curves, arXiv:2312.12528.

\bibitem[HJO]{HJO} Y. Huang, R. Jiang and A. Oblomkov, Quot scheme of points on torus knot singularities, arXiv:2608.16086.

\bibitem[HLO]{HLO} Y. Huang, K. Lau and K. Ono, Rogers--Ramanujan identities from the geometry of $X^a=Y^b$, arXiv:2609.20567v2 (2026).

\bibitem[Hum]{Hum} J. E. Humphreys, \emph{Reflection groups and Coxeter groups}, Cambridge Stud. Adv. Math. 29, Cambridge Univ. Press, 1990.

\bibitem[IM]{IM} N. Iwahori and H. Matsumoto, On some Bruhat decomposition and the structure of the Hecke rings of $\mathfrak p$-adic Chevalley groups, \emph{Publ. Math. Inst. Hautes \'Etudes Sci.} 25 (1965), 5--48.

\bibitem[KT]{KT} O. Kivinen and C.-C. Tsai, Shalika germs for tamely ramified elements in $\mathrm{GL}_n$, arXiv:2209.02509v3 (2023).

\bibitem[Ku]{Kunz} E. Kunz, The value-semigroup of a one-dimensional Gorenstein ring, \emph{Proc. Amer. Math. Soc.} 25 (1970), 748--751.

\bibitem[LO]{LO} K. Lau and K. Ono, Modularity of point counts for the curves $X^a=Y^b$: New Rogers--Ramanujan identities, arXiv:2608.05480.

\bibitem[Mac]{Mac} I. G. Macdonald, \emph{Symmetric functions and Hall polynomials}, 2nd ed., Oxford Univ. Press, 1995.

\bibitem[Ma]{Maulik} D. Maulik, Stable pairs and the HOMFLY polynomial, \emph{Invent. Math.} 204 (2016), 787--831.

\bibitem[Ne1]{NegShuffle} A. Negu\c{t}, The shuffle algebra revisited,
\emph{Int. Math. Res. Not. IMRN} 2014, no.~22, 6242--6275;
arXiv:1209.3349v3.

\bibitem[Ne2]{Neg} A. Negu\c{t}, Operators on symmetric polynomials,
arXiv:1310.3515 (2013).

\bibitem[Ne3]{NegPieri} A. Negu\c{t}, The $m/n$ Pieri rule,
\emph{Int. Math. Res. Not. IMRN} 2016, no.~1, 219--257;
arXiv:1407.5303v5.

\bibitem[Ne4]{NegIntegral} A. Negu\c{t},
An integral form of quantum toroidal $\mathfrak{gl}_1$,
\emph{Math. Rep. (Bucur.)} 26(76) (2024), no.~3--4, 183--205.

\bibitem[ORS]{ORS} A. Oblomkov, J. Rasmussen and V. Shende, The Hilbert scheme of a plane curve singularity and the HOMFLY homology of its link, \emph{Geom. Topol.} 22 (2018), 645--691.

\bibitem[OS]{OS} A. Oblomkov and V. Shende, The Hilbert scheme of a plane curve singularity and the HOMFLY polynomial of its link, \emph{Duke Math. J.} 161 (2012), 1277--1303.

\bibitem[Ra]{Ram} A. Ram, A Frobenius formula for the characters of the Hecke algebras, \emph{Invent. Math.} 106 (1991), 461--488.

\bibitem[Sa]{Sa} P. Samuelson, Iterated torus knots and double affine Hecke algebras,
arXiv:1408.0483v4.

\bibitem[Se]{Serre} J.-P. Serre, \emph{A course in arithmetic}, Grad. Texts in Math. 7, Springer, 1973.

\bibitem[Ts]{Tsai} C.-C. Tsai, Inductive structure of Shalika germs and affine Springer fibers, arXiv:1512.00445v1 (2015).

\bibitem[Wa1]{Wal89} J.-L. Waldspurger, Sur les germes de Shalika pour les groupes lin\'eaires, \emph{Math. Ann.} 284 (1989), 199--221.

\bibitem[Wa2]{Wal91} J.-L. Waldspurger, Sur les int\'egrales orbitales tordues pour les groupes lin\'eaires: un lemme fondamental, \emph{Canad. J. Math.} 43 (1991), 852--896.

\bibitem[Ya1]{Yang1} Y. Yang, Transformation formulas for generalized Dedekind eta functions, \emph{Bull. London Math. Soc.} 36 (2004), 671--682.

\bibitem[Ya2]{Yang2} Y. Yang, Modular units and cuspidal divisor class groups of $X_1(N)$, \emph{J. Algebra} 322 (2009), 514--553; arXiv:0712.0629.

\bibitem[Yun]{Yun} Z. Yun, Orbital integrals and Dedekind zeta functions, in: \emph{The legacy of Srinivasa Ramanujan}, Ramanujan Math. Soc. Lect. Notes Ser. 20, 2013, 399--420.
\end{thebibliography}
\end{document}